\documentclass{article}

\usepackage{amsfonts}
\usepackage{dsfont}
\usepackage{geometry}
\usepackage{soul}
\usepackage{booktabs}
\usepackage[round]{natbib}
\usepackage[utf8]{inputenc}
\usepackage{enumitem}
\usepackage{amsmath}
\usepackage{amsthm}
\usepackage{amssymb}
\usepackage{bbm}
\usepackage{comment}
\usepackage{float}
\usepackage{needspace}

\usepackage{environ}
\usepackage{etoolbox}

\usepackage{graphicx}
\graphicspath{{fig/}}
\usepackage[svgnames,dvipsnames]{xcolor}

\usepackage[bookmarksopen=true,bookmarks=true,pdfencoding=auto,psdextra,
colorlinks,
hypertexnames=false,
linkcolor=RoyalBlue,
citecolor=ForestGreen,
urlcolor=RubineRed!90!black,
]{hyperref}

\hypersetup{
  unicode=true,    
  pdftoolbar=true,    
  pdfmenubar=true      
}

\usepackage[capitalize,noabbrev]{cleveref}

\usepackage{crossreftools}
\newtheorem{theorem}{Theorem}
\newtheorem{definition}{Definition}
\newtheorem{corollary}{Corollary}
\newtheorem{proposition}{Proposition}
\newtheorem{lemma}{Lemma}

\newtheorem{conjecture}{Conjecture}
\crefname{conjecture}{Conjecture}{Conjectures}
\Crefname{conjecture}{Conjecture}{Conjectures}

\theoremstyle{definition}
\newtheorem{remark}{Remark}

\newcommand{\bbE}{\mathbb E}

\newcommand{\bbQ}{\mathbb Q}
\newcommand{\bbP}{\mathbb P}

\newcommand{\bbR}{\mathbb R}

\newcommand{\bfA}{\mathbf A}

\newcommand{\bfI}{\mathbf I}
\newcommand{\bfi}{\mathbf i}
\newcommand{\bfj}{\mathbf j}
\newcommand{\bfW}{\mathbf W}   
\newcommand{\bfX}{\mathbf X}
\newcommand{\bfY}{\mathbf Y}

\newcommand{\bfy}{\mathbf y}

\newcommand{\calE}{\mathcal E}

\newcommand{\rarr}{\rightarrow}

\newcommand{\ee}{\mathbb E}
\newcommand{\vv}{\mathrm {Var}}

\usepackage[font=small,labelfont=bf]{caption}
\usepackage{subcaption}

\usepackage[normalem]{ulem}

\allowdisplaybreaks

\usepackage{microtype}

\usepackage[xspace]{ellipsis}
\usepackage{titletoc}

\newcommand{\footremember}[2]{
    \footnote{#2}
    \newcounter{#1}
    \setcounter{#1}{\value{footnote}}
}
\newcommand{\footrecall}[1]{
    \footnotemark[\value{#1}]
}

\title{On Statistical-Computational Gaps in\\Graph Combinatorial Testing Problems}

\author{
Ikjun Choi\footremember{utsds}{Department of Statistics and Data Sciences, University of Texas at Austin, TX 78712, USA.} \\ {\footnotesize \texttt{\url{ikjunchoi@utexas.edu}}}
\and
Pratik Patil\footrecall{utsds} \\ {\footnotesize \texttt{\url{pratikpatil@utexas.edu}}}
\and
Alessandro Rinaldo\footrecall{utsds} \\ {\footnotesize \texttt{\url{alessandro.rinaldo@austin.utexas.edu}}}
}

\date{\vspace{-15pt}}

\begin{document}

\maketitle

\begin{abstract}
We study high-dimensional graph combinatorial hypothesis testing problems in which a signal is supported on the edge set of an unknown subgraph drawn from a class invariant under vertex relabeling.
For Gaussian observations, we characterize sharp statistical and computational thresholds (at polynomial order) for the signal magnitude required for strong detection across various combinatorial graph alternatives.
Our analysis is based on bounding the second moments of the full and low-degree likelihood ratios in terms of the overlap of two independent latent structures and comparing the resulting lower bounds with scan, average, and coordinatewise maximum tests.
Under an appropriate low-degree conjecture for computational hardness, we obtain planted-clique thresholds matching those known for Gaussian submatrix detection and identify new statistical-computational gaps for multiple planted cliques and hypercliques, thick cycles and thick paths.
We also establish a novel transfer principle across natural exponential family observation models: matching the one-coordinate chi-square divergence preserves full likelihood-ratio second moments exactly and makes low-degree second moments equivalent up to constant factors under a degree condition.
Under local regularity, we obtain corresponding testing upper bounds in terms of an effective signal strength.
When these upper and lower bounds match, they yield universal threshold scales in the local alternative regime.
As a consequence, statistical-computational gaps for planted-clique detection in  Erd\H{o}s--R\'enyi settings can be deduced from the Gaussian model.
\end{abstract}

\section{Introduction}

Signal detection problems are central to modern statistics and arise in applications ranging from network analysis and neuroscience to genomics and spatial statistics.
A recurring feature of these problems is that a weak signal is embedded in a high-dimensional ambient space, but is present only on an unknown subset of coordinates with prescribed combinatorial structure.
In graph-valued data, for example, the signal may be supported on the edge set of an unknown clique, collection of cliques, path, cycle, or other structured subgraph.
A scan over all admissible structures, or an equivalent generalized likelihood-ratio test, often achieves rate-optimal detection but can be computationally infeasible when the class of possible supports is large.
This raises a basic question: when can statistically optimal detection be achieved by computationally efficient tests?

A growing body of work shows that the answer can depend sharply on the geometry of the latent structure.
For many high-dimensional testing problems, there is a range of signal strengths in which detection is information-theoretically possible but no known polynomial-time method succeeds.
Such a \emph{statistical-computational gap} is most familiar from planted clique and related subgraph detection problems, but analogous phenomena appear across sparse principal component analysis, submatrix detection, tensor problems, and structured normal means models.
Despite this broad picture, the mechanisms that determine the size of the gap remain only partially understood for many graph combinatorial testing problems, especially under observation models beyond the classical Erd\H{o}s--R\'enyi setting.

This paper studies these questions through a common testing framework for signals supported on unknown graph-indexed combinatorial structures.
Our starting point is the Gaussian observation model, in which the full and low-degree likelihood-ratio second moments can both be expressed in terms of the overlap distribution of two independent latent structures.
The former yields information-theoretic lower bounds, while the latter provides evidence for computational lower bounds under the low-degree heuristic.
We combine these lower bounds with standard scan, average, and coordinatewise maximum tests to identify statistical and conjectural computational thresholds for several graph testing problems.
We also show that the same overlap calculations extend beyond Gaussian noise to natural exponential family models by matching the one-coordinate chi-square divergence.

\begin{figure}[t]
\centering
\begin{subfigure}[t]{0.49\textwidth}
\centering
\includegraphics[width=\linewidth]{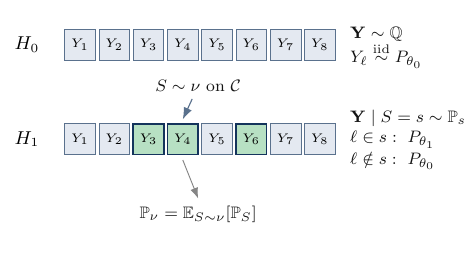}
\caption{General average-case testing model.}
\label{fig:testing-framework-general}
\end{subfigure}\hfill
\begin{subfigure}[t]{0.49\textwidth}
\centering
\includegraphics[width=\linewidth]{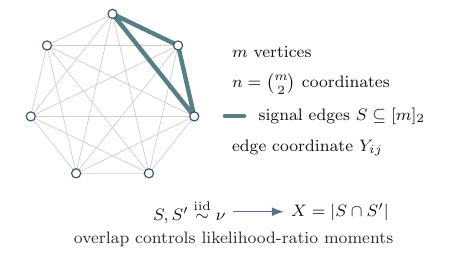}
\caption{Graph specialization and overlap.}
\label{fig:testing-framework-graph}
\end{subfigure}
\caption{Testing framework. Panel (a) shows the general mixture model: under $H_0$, the coordinates follow $P_{\theta_0}$, whereas under $H_1$, a latent structure $S\sim\nu$ determines which coordinates follow $P_{\theta_1}$ instead of $P_{\theta_0}$, and mixing over $S$ yields $\bbP_\nu$. Panel (b) specializes the framework to edge-indexed observations on $K_m$; the overlap $X=|S\cap S'|$, for $S,S'\stackrel{\mathrm{iid}}{\sim}\nu$, controls the likelihood-ratio moments used in our lower bounds.}
\label{fig:testing-framework}
\end{figure}

\subsection{Summary of results}

We study sequences of combinatorial hypothesis-testing problems indexed by the number $n=\binom{m}{2}$ of edges in $K_m$, the complete graph on $m$ vertices.
The null hypothesis represents unstructured noise, whereas the alternative is a mixture of structured distributions indexed by a class of edge-induced subgraphs of $K_m$; \Cref{fig:testing-framework} illustrates this distributional hierarchy and its graph specialization.
Let $\{P_\theta:\theta\in\Theta\}$ be a parametric family of distributions on a sample space $\Omega\subseteq\mathbb R$, where $\Theta\subseteq\mathbb R$ is open.
We observe a random vector $\bfY=(Y_1,\ldots,Y_n)\in\Omega^n$ whose coordinates are indexed by the edges of $K_m$.
Under the null hypothesis, the coordinates of $\bfY$ are independent draws from a reference distribution $P_{\theta_0}$, for $\theta_0\in\Theta$.
Under the alternative, a latent edge set is sampled from a combinatorial class $\mathcal C\subseteq 2^{[n]}$, and the coordinates indexed by that set instead follow $P_{\theta_1}$ for some $\theta_1\neq\theta_0$.
More precisely, for each fixed structure $s\in\mathcal C$, let $\bbP_s$ denote the distribution under which the coordinates of $\bfY$ are independent and
\begin{align}\label{equation: alternative distribution}
Y_\ell\sim
\begin{cases}
P_{\theta_1}, & \ell\in s,\\
P_{\theta_0}, & \ell\notin s,
\end{cases}
\qquad \ell\in[n].
\end{align}
Let $\nu$ be a mixing measure on $\mathcal C$, taken to be uniform in the basic fixed-size models.
The alternative distribution is the mixture $\bbP_\nu$ of $\{\bbP_s:s\in\mathcal C\}$ with respect to $\nu$, as illustrated in \Cref{fig:testing-framework-general}; that is,
\begin{align}\label{eq:mixture-alternative}
\bbP_\nu(A):=\bbE_{S\sim\nu}\big[\bbP_S(A)\big],
\qquad A\subseteq\Omega^n\ \text{measurable}.
\end{align}
The class $\mathcal C$ controls the combinatorial complexity of the alternative and can be viewed as a collection of subgraphs of $K_m$ indexing the possible anomalous coordinates of $\bfY$.
The corresponding edge-indexed graph model and the overlap variable used in our lower-bound arguments appear in \Cref{fig:testing-framework-graph}.

We study the resulting graph combinatorial testing problem
\begin{align}\label{eq:general-testing-problem}
H_0:\ \bfY\sim\bbQ
\qquad\text{versus}\qquad
H_1:\ \bfY\sim\bbP_\nu,
\end{align}
where $\bbQ=P_{\theta_0}^{\otimes n}$ is the $n$-fold product of the reference distribution.
We consider asymptotic regimes in which the dimension $n$ diverges and the combinatorial class $\mathcal C$ and signal parameter $\theta_1$ may vary with $n$.
The reference parameter $\theta_0$ is fixed in the Gaussian model and in the local alternative universality results of \Cref{sec: examples}.
Our general bounds for natural exponential family models also allow $\theta_0$ to vary with $n$, as in the sparse Erd\H{o}s--R\'enyi applications in \Cref{section: applications to Erdos--Renyi models}.
Our primary interest is in regimes where the separation between the null and alternative vanishes as $n$ grows.
We suppress dependence on $n$ when it is clear from context.

A test is a measurable function $\Delta:\Omega^n\to\{0,1\}$, with $\Delta(\bfY)=1$ indicating rejection of the null hypothesis.
We measure its performance by the average-case risk
\begin{align}\label{eq:average-case-risk}
R_\nu(\Delta)
:=\bbQ\{\Delta(\bfY)=1\}+\bbP_\nu\{\Delta(\bfY)=0\},
\end{align}
the sum of its type I error and its type II error averaged over the latent structure.

\begin{definition}[Strong detection, statistical and computational distinguishability]\label{definition: strong detection}
Consider a sequence of testing problems $H_0:\bfY\sim\bbQ_n$ versus $H_1:\bfY\sim\bbP_{\nu_n}$.
A sequence of tests $\Delta_n:\Omega^n\to\{0,1\}$ achieves \emph{strong detection} if
\[
R_{\nu_n}(\Delta_n)\mathrel{\to} 0.
\]
The sequence of testing problems is \emph{statistically distinguishable} if there exists a sequence of tests achieving strong detection, and \emph{computationally distinguishable} if such a sequence can be evaluated in time polynomial in $n$.
If no sequence of tests achieves strong detection, the testing problems are statistically indistinguishable; if no polynomial-time sequence does so, they are computationally indistinguishable.
\end{definition}

\begin{remark}[Strong vs weak detection] For comparison, a sequence of tests $\Delta_n$ achieves \emph{weak detection} if $\limsup_{n\to\infty}R_{\nu_n}(\Delta_n)<1$, meaning that it achieves a nonvanishing improvement over the trivial risk of one. Our distinguishability statements and threshold results concern strong detection.
\end{remark}

We review the second-moment and low-degree likelihood-ratio methods used to establish statistical and computational lower bounds in \Cref{sec:preliminaries}.

Our goal is to characterize, as a function of $n$ and the complexity of $\mathcal C$, the signal scale at which statistical and computational distinguishability emerge.
Most of the paper studies the Gaussian sequence model of \citet{addario2010combinatorial}, in which $P_{\theta_0}=N(0,1)$ and $P_{\theta_1}=N(\mu,1)$ for some nonzero $\mu$.
In this model, the separation between the hypotheses is measured by the signal strength $|\mu|$, and we denote the statistical and computational thresholds by $\mu_{\mathrm{stat}}^\ast$ and $\mu_{\mathrm{comp}}^\ast$, respectively.
These thresholds, defined formally in \Cref{definition: statistical and computational thresholds}, identify the signal scales for strong detection up to subpolynomial factors.
Our low-degree lower bounds rule out strong separation by polynomial statistics of polylogarithmic degree (\Cref{definition:strong-separation}); under \Cref{conjecture: low-degree polynomial conjecture}, we interpret low-degree indistinguishability as evidence that no polynomial-time test achieves strong detection.

Our main results for the Gaussian sequence model are collected in \Cref{tab:summary-of-problems}.
The threshold scales summarize the polynomially separated parameter regimes (see \Cref{sec:organization-and-notation} for the definition of polynomial separation) specified in the corresponding results.
As in \citet{addario2010combinatorial}, we focus on combinatorial classes $\mathcal C$ that are invariant under permutations of the vertices of $K_m$.
For the nonclique classes studied there, the full and low-degree likelihood-ratio lower bounds match polynomial-time testing upper bounds at polynomial order, and hence no statistical-computational gap arises; the corresponding results appear in the appendix sections listed in the first block of \Cref{tab:summary-of-problems}.
We obtain analogous no-gap results for $k$-stars and for paths and cycles of length $k$.
By contrast, \Cref{sec:results} identifies statistical-computational gaps for planted cliques, multiple planted cliques and hypercliques, and thick cycles and paths.

For $k$-vertex planted cliques, \Cref{theorem: fixed-size clique problem} identifies the statistical threshold $\mu_{\mathrm{stat}}^\ast\asymp k^{-1/2}\wedge m/k^2$ and the conjectural computational threshold $\mu_{\mathrm{comp}}^\ast\asymp 1\wedge m/k^2$; see \Cref{fig:planted-clique-phase-diagram}.
These thresholds agree with those established for Gaussian submatrix detection by \citet{butucea2013detection,ma2015computational}, with computational lower bounds obtained through a polynomial-time reduction from planted clique under the planted-clique conjecture \citep{alon2007testing}.
In contrast, our derivation uses the low-degree likelihood-ratio method.
Across the Gaussian families summarized in \Cref{tab:summary-of-problems}, the low-degree lower bounds match the best polynomial-time upper bounds at polynomial order, covering both the no-gap examples and the models that exhibit statistical-computational gaps.
\Cref{corollary: random-size clique problem} shows that the same thresholds hold when the planted clique has a binomially distributed vertex set with mean size $k$.
In the intermediate regime, the exhaustive scan test succeeds below the conjectural computational threshold, while the low-degree obstruction predicts, under \Cref{conjecture: low-degree polynomial conjecture}, that no polynomial-time test achieves strong detection.

\Cref{theorem: multiple cliques problem} extends the analysis to a union of $\gamma$ disjoint planted $k$-vertex cliques and gives $\mu_{\mathrm{stat}}^\ast\asymp k^{-1/2}\wedge m/(\gamma k^2)$ and $\mu_{\mathrm{comp}}^\ast\asymp 1\wedge m/(\gamma k^2)$ for both fixed-size and random-size constructions.
\Cref{corollary: multiple planted hypercliques problem} treats the corresponding $t$-uniform hypergraph model, with statistical threshold of order $k^{-(t-1)/2}\wedge m^{t/2}/(\gamma k^t)$ and computational threshold of order $1\wedge m^{t/2}/(\gamma k^t)$ in the regimes stated there.
For the random-size multiple-clique and hyperclique models, we derive the computational lower bounds by comparing the overlap moments with those of auxiliary single-clique and single-hyperclique models, respectively.

For length-$\ell$, $\Delta$-thick cycles, \Cref{theorem: thick cycle problem} identifies $\mu_{\mathrm{stat}}^\ast\asymp \Delta^{-1/2}\wedge m/(\ell\Delta)$ and $\mu_{\mathrm{comp}}^\ast\asymp 1\wedge m/(\ell\Delta)$.
\Cref{corollary: thick path problem} shows that the same thresholds hold for the corresponding thick path model.
The lower-bound argument reduces the relevant overlap calculation to the multiple-clique analysis by decomposing each thick cycle or path into local patches.

Finally, \Cref{sec:nef-equivalence} extends the framework beyond Gaussian observations to natural exponential families.
\Cref{theorem: equivalence between NEF models} shows that two natural exponential family models with matching one-coordinate chi-square divergence have identical full likelihood-ratio second moments and equivalent low-degree second moments under the stated degree condition.
\Cref{proposition: upper bound for NEF models} shows that the average and scan upper bounds depend on the observation model through its effective standardized mean shift.
Together, these results transfer Gaussian thresholds to models such as Erd\H{o}s--R\'enyi observations after the appropriate chi-square calibration and under the stated regularity conditions.

\begin{table}[t]
\centering
\caption{Summary of problems, with thresholds and test optimality stated at polynomial order.
The threshold scales summarize the polynomially separated parameter regimes (see \Cref{sec:organization-and-notation} for the definition of polynomial separation) specified in the corresponding results.
$^\dagger$ denotes the problems whose statistical thresholds were studied by \citet{addario2010combinatorial}.
The statistical and computational thresholds for $k$-vertex cliques match those for Gaussian submatrix detection, where \citet{ma2015computational} obtained computational lower bounds through a planted-clique reduction.
A checkmark indicates a statistical-computational gap (conditional on \Cref{conjecture: low-degree polynomial conjecture}) in some polynomial parameter regime.}
\label{tab:summary-of-problems}

\resizebox{\linewidth}{!}{
\begin{tabular}{lcccccc}
\toprule
Problem & Section & $\mu_{\mathrm{stat}}^\ast$ & $\mu_{\mathrm{comp}}^\ast$ & Rate-optimal test & Polytime test & Comp-stat gap \\
\midrule

$K$-Disjoint sets$^\dagger$
& \Cref{appendix: disjoint sets}
& $\sqrt{\frac{1}{K}}$
& $\sqrt{\frac{1}{K}}$
& scan
& scan
& $\times$ \\[5pt]

$K$-sets$^\dagger$
& \Cref{appendix: K-sets}
& $\sqrt{\frac{n}{K^2}}\wedge 1$
& $\sqrt{\frac{n}{K^2}}\wedge 1$
& scan, average
& scan, average
& $\times$ \\[5pt]

Perfect matching$^\dagger$
& \Cref{appendix: perfect matchings}
& 1
& 1
& average
& average
& $\times$ \\[5pt]

Star$^\dagger$
& \Cref{appendix: star}
& $\sqrt{\frac{1}{m}}$
& $\sqrt{\frac{1}{m}}$
& scan
& scan
& $\times$ \\[5pt]

Spanning tree$^\dagger$
& \Cref{appendix: spanning tree}
& 1
& 1
& average
& average
& $\times$ \\[5pt]

$k$-vertex star
& \Cref{appendix: k-size star}
& $\sqrt{\frac{m}{k^2}}\wedge 1$
& $\sqrt{\frac{m}{k^2}}\wedge 1$
& degree, outlier 
& degree, outlier
& $\times$ \\[5pt]

$k$-length path, cycle
& \Cref{appendix: k-length path}
& 1
& 1
& outlier
& outlier
& $\times$ \\

\midrule

$k$-vertex clique$^\dagger$
& \Cref{section: planted clique}
& $\frac{m}{k^2}\wedge \sqrt{\frac{1}{k}}$
& $\frac{m}{k^2}\wedge 1$
& scan, average
& outlier, average
& $\checkmark$ \\

$\gamma$ multiple $k$-vertex cliques
& \Cref{section: multiple planted cliques}
& $\frac{m}{\gamma k^2}\wedge \sqrt{\frac{1}{k}}$
& $\frac{m}{\gamma k^2}\wedge 1$
& scan, average
& outlier, average
& $\checkmark$ \\

$\gamma$ multiple $k$-vertex hypercliques
& \Cref{section: multiple planted cliques}
& $\frac{m^{t/2}}{k^t\gamma} \wedge \sqrt{\frac{1}{k^{t-1}}}$
& $\frac{m^{t/2}}{k^t\gamma} \wedge 1$
& scan, average
& outlier, average
& $\checkmark$ \\

length-$\ell$, $\Delta$-thick cycle, path
& \Cref{section: thick cycle and path}
& $\frac{m}{\ell\Delta} \wedge \sqrt{\frac{1}{\Delta}}$
& $\frac{m}{\ell\Delta} \wedge 1$
& scan, average
& outlier, average
& $\checkmark$ \\

\bottomrule
\end{tabular}
}

\end{table}

\subsection{Related work}\label{sec:related-work}

\paragraph{Gaussian combinatorial testing and graph scans.}
The Gaussian framework for testing structured means in this paper is closest to \citet{addario2010combinatorial}, who analyze combinatorial classes through scan and averaging tests and through lower bounds based on the intersections of latent supports.
Related work on anomalous paths and connected clusters includes \citet{arias2008searching,arias2011detection}, while \citet{sharpnack2013near,sharpnack2016detecting} develop computationally tractable relaxations of graph scan statistics.
Those graph-scan models generally place the signal on vertices or vertex-indexed observations, whereas our principal examples place the signal on the edge set of an unknown subgraph.
Distribution-free permutation and rank-based scans were studied by \citet{arias2018distribution}, including performance guarantees for one-parameter natural exponential families.

\paragraph{Planted cliques and submatrices.}
The planted-clique problem is a canonical source of statistical-computational gaps, with early evidence for algorithmic limitations due to \citet{jerrum1992large}.
Classical and subsequent algorithmic upper bounds include the spectral method of \citet{alon1998finding} and the nearly linear-time method of \citet{deshpande2015finding}.
Lower bounds for restricted computational models include spectral, sum-of-squares, and statistical-query results \citep{deshpande2015improved,montanari2017limitation,barak2019nearly,feldman2017statistical}, and \citet{hirahara2024planted} establish equivalences among many common formulations of the planted-clique conjecture.
For weighted observations, \citet{chatterjee2025detecting} study a general hidden-clique model and derive information-theoretic bounds and spectral guarantees, while our clique results focus on vanishing signals and on matching overlap calculations for both the full and low-degree likelihood ratios.
Sparse submatrix detection has sharp statistical thresholds \citep{butucea2013detection,goswami2026sharp}, planted-clique-based computational lower bounds \citep{ma2015computational}, and distributional universality results based on average-case reductions \citep{brennan2019universality}.
The closest submatrix comparison for our multiple-clique results is \citet{dadon2024detection}, who study detection and recovery for one or multiple hidden Gaussian submatrices using information-theoretic and low-degree methods.
Related work addresses localization and heterogeneous multiple submatrices \citep{chen2016statistical,cai2017computational,orenloberman2026inhomogeneous}.
These models allow distinct row and column supports, whereas our clique models impose symmetric principal edge sets.

\paragraph{Multiple structures and hypergraphs.}
Multiple planted structures have also been studied through growing-cluster recovery guarantees \citep{chen2016statistical}, low-degree lower bounds for planted colorings or disjoint cliques covering most vertices \citep{kothari2023planted}, and recent sum-of-squares and statistical-query lower bounds for multiple planted cliques and bicliques \citep{mosievskiy2026multiple}.
Our multiple-clique problem instead concerns detecting a prescribed number of vertex-disjoint planted cliques and tracks how the number of planted components changes both the statistical and low-degree thresholds.
For higher-order observations, related statistical-computational phenomena have been studied in planted tensor clustering \citep{luo2022tensor}, planted dense subhypergraph recovery \citep{corinzia2022statistical}, and low-degree detection of dense subhypergraphs \citep{dhawan2025detection}.
The matrix-only observation of a planted hypergraph clique considered by \citet{alaluusua2026planted} is complementary to our setting, where the full collection of independent hyperedge weights is observed.

\paragraph{General planted subgraphs, cycles, and paths.}
For dense and sparse planted communities in Erd\H{o}s--R\'enyi graphs, \citet{arias2014community,verzelen2015community} establish information-theoretic thresholds, while \citet{hajek2015computational} derive planted-clique-based computational lower bounds.
More general planted-subgraph detection has been studied in \citet{huleihel2022inferring,elimelech2025detecting}, with specialized results for bipartite graphs \citep{rotenberg2024planted}, constant-degree polynomial tests \citep{yu2024counting}, and planted matchings and spanning trees \citep{addario2026statistical}.
The corresponding recovery literature includes sharp thresholds for planted subgraphs \citep{mossel2025sharp,huleihel2026recovery} and for hidden weighted sparse graphs \citep{hou2026recovery}.
For path- and cycle-like structures, \citet{arias2008searching} analyze anomalous trails, while \citet{mao2023detection,mao2025information} establish computational and information-theoretic thresholds for planted dense cycles.
The random dense-cycle prior considered in these works coincides, after rescaling the bandwidth, with our random-size thick-cycle prior when $\ell=m$.
They consider Bernoulli observations and an Erd\H{o}s--R\'enyi null whose edge density matches the marginal edge density under the alternative.
In our testing framework, the null coordinate distribution instead agrees with the distribution outside the planted structure.
Mean-matched formulations also arise in the planted matching and spanning-tree problems studied by \citet{addario2026statistical}; we discuss this distinction further in \Cref{sec:discussion}.

\paragraph{Frameworks for statistical and computational lower bounds.}
The second moment and conditional second moment methods provide information-theoretic lower bounds in many planted models \citep{banks2017information,perry2018optimality,perry2020statistical}.
Evidence for computational hardness has been developed through reductions from planted clique \citep{brennan2018reducibility,brennan2020reducibility}, sum-of-squares lower bounds \citep{hopkins2017power,barak2019nearly}, statistical-query lower bounds \citep{feldman2017statistical}, and the low-degree polynomial framework \citep{hopkins2018statistical,kunisky2019notes,schramm2022computational,carpentier2025low,wein2025computational}.
The Franz--Parisi criterion gives a related overlap-based perspective \citep{bandeira2022franz}, and \citet{chen2025optimized} connect an optimized version of this criterion to statistical-query lower bounds.
Recent work also develops sharp low-degree thresholds beyond the planted-versus-null setting \citep{skeja2026sharp}.
At the same time, counterexamples show that a general implication from low-degree indistinguishability to polynomial-time hardness requires additional assumptions \citep{holmgren2021counterexamples,buhai2025quasi,mao2026polynomial}.
Accordingly, our low-degree lower bounds rule out strong separation by low-degree polynomial statistics and are interpreted as conjectural evidence, rather than unconditional lower bounds, for arbitrary polynomial-time algorithms.

\paragraph{Distributional universality and natural exponential families.}
Distributional extensions of planted detection problems have been obtained through Gaussianization and related reductions \citep{brennan2018reducibility,brennan2019universality,brennan2020reducibility}, through distribution-free scan statistics \citep{arias2018distribution}, and through general weighted planted-clique models \citep{chatterjee2025detecting}.
More generally, \citet{kunisky2025low} develops low-coordinate-degree lower bounds for broad latent-variable channels and proves channel-universality results governed by the scalar Fisher information of the observation channel.
Closely related to our natural exponential family results, \citet{kunisky2021hypothesis} derives full and low-degree likelihood-ratio formulas and channel-monotonicity inequalities for the Morris class of natural exponential families with quadratic variance functions.
For the models studied here, whose alternatives mix over latent supports, our comparison extends beyond the quadratic-variance setting.
With a common latent prior, matching the one-coordinate chi-square divergence gives identical full likelihood-ratio second moments and low-degree second moments equivalent up to constant factors, subject to the degree condition in \Cref{theorem: equivalence between NEF models}.
We also analyze the average, scan, and degree tests under local regularity conditions in \Cref{proposition: upper bound for NEF models}.
When these upper bounds match the second-moment lower bounds, our local alternative result transfers the resulting statistical and conjectural computational thresholds across observation families under the conditions of \Cref{theorem: universality of tightness in the local alternative regime via chi-square matching}.

\subsection{Organization and notation}\label{sec:organization-and-notation}

\paragraph{Organization.}
The rest of the paper is organized as follows.
\Cref{sec:gaussian-testing} develops the lower-bound framework based on the full and low-degree likelihood ratios, specializes it to Gaussian graph combinatorial testing, and records the overlap identities and the outlier, average, and scan upper bounds used throughout the paper.
\Cref{sec:results} applies these tools to planted cliques, multiple planted cliques and hypercliques, and thick cycles and paths, for both fixed-size and random-size constructions.
\Cref{sec:nef-equivalence} extends the framework to natural exponential family noise models through chi-square matching, develops corresponding upper bounds, and studies the universality of the resulting thresholds in the local alternative regime.
\Cref{sec:discussion} summarizes the main conclusions, limitations, and directions for future work.
The appendices collect technical lemmas and complete proofs (\Cref{sec:app-technical-lemmas,sec:app-proofs}) and analyze additional latent structures that complement the main examples (\Cref{sec:app-additional-examples}).

\paragraph{Notation.}
For a positive integer $n$, we write $[n] := \{1,\ldots,n\}$.
For $r\in[n]$, we write $[n]_r:=\{(i_1,\ldots,i_r)\in[n]^r: i_1<\cdots<i_r\}$, which we identify with the collection of $r$-element subsets of $[n]$.
For a finite set $A$, we write $|A|$ for the cardinality of $A$.
For real numbers $a$ and $b$, we write $a\wedge b:=\min\{a,b\}$ and $a\vee b:=\max\{a,b\}$.
The indicator of an event $\mathcal E$ is denoted by $\mathds{1}\{\mathcal E\}$.
Expectations and variances are denoted by $\bbE$ and $\operatorname{Var}$, respectively.
Probabilities are denoted by $\Pr$, with the source of randomness displayed in a subscript when useful; for example, $\Pr_{S\sim\nu}\{\ell\in S\}$.
For a named probability measure, we instead write its probability of an event directly, as in $\bbQ(A)$, $\bbP_s(A)$, and $\bbP_\nu(A)$.
We use plain capital letters such as $P_{\theta_0}$ and $P_{\theta_1}$ for one-coordinate probability measures and blackboard-bold letters such as $\bbQ$ and $\bbP$ for joint probability measures on observation vectors.
Calligraphic letters are used for classes, families, events, and function spaces, such as the latent-structure class $\mathcal C$.
The notation $S,S'\sim\nu$ means that $S$ and $S'$ are independent draws from $\nu$; a fixed realization of $S$ is denoted by $s\in\mathcal C$.
We write $H\sim\operatorname{Hypergeometric}(N,M,r)$ when $H$ is the number of marked items in a sample of size $r$ drawn uniformly without replacement from a population of $N$ items, of which $M$ are marked.
We use boldface symbols for vectors and multi-indices and nonbold symbols for scalars and sets; in particular, $\mathbf{Y}$ denotes a random observation vector and $\mathbf{y}$ denotes one of its realizations.

Unless stated otherwise, all asymptotic statements are taken as the ambient dimension tends to infinity.
In the general testing framework, the ambient dimension is $n$, whereas in the graph models, $m$ denotes the number of vertices and $n=\binom{m}{2}$ denotes the number of observed edges.
Model parameters, latent-structure classes, priors, and probability distributions may depend on $n$, or on $m$ in the graph setting, even when this dependence is suppressed.
Unless stated otherwise, constants such as $c,C>0$ do not depend on the ambient dimension and may change from line to line.

For two positive sequences $\{a_n\}$ and $\{b_n\}$, we write $a_n \lesssim b_n$ if $a_n \leq Cb_n$ eventually for a constant $C>0$ independent of $n$.
We write $a_n \gtrsim b_n$ if $b_n \lesssim a_n$.
We write $a_n \asymp b_n$ if both $a_n \lesssim b_n$ and $a_n \gtrsim b_n$ hold.
We write $a_n \ll b_n$ if there exists a constant $c>0$ such that $a_n\lesssim n^{-c}b_n$, and write $a_n\gg b_n$ if $b_n\ll a_n$.
Thus, $\ll$ and $\gg$ denote separation by a polynomial factor, whereas $o(\cdot)$ and $\omega(\cdot)$ retain their usual asymptotic meanings.

\Needspace{6\baselineskip}
\section{Gaussian graph combinatorial testing}\label{sec:gaussian-testing}

Throughout this section, we study the isotropic Gaussian specialization of the graph combinatorial testing problem in \Cref{eq:general-testing-problem}, a model that has been studied extensively in the signal detection literature; see, e.g., \citet{arias2008searching,arias2011detection,addario2010combinatorial,ma2015computational,montanari2017limitation,dadon2024detection}. 
 In our principal graph models, the data consist of the vector $\bfY=(Y_{ij}:(i,j) \in [m]_2)\in \mathbb R^{\binom m2}$, whose coordinates are indexed by the edges of the complete graph $K_m$. For a latent structure $S \subseteq [m]_2$ randomly selected from a class $\mathcal{C}$ of edge sets of $K_m$, the null and alternative hypotheses are
\begin{align}\label{equation: gaussian graph testing problem}
H_0: Y_{ij} \stackrel{iid}{\sim} N(0,1)\quad \text{versus} \quad H_1: Y_{ij}\,|\,S \stackrel{ind}{\sim}
\begin{cases}
\mathrm{N}(\mu,1), & \text{if } (i,j) \in S, \\
\mathrm{N}(0,1), & \text{if } (i,j) \notin S,
\end{cases},
\end{align}
where $\mu > 0$ is the signal strength. Recall that $\bbQ$ and $\bbP_\nu$ denote the distributions of $\bfY$ under the null and alternative hypotheses, respectively, where $\nu$ is an appropriate distribution on $\mathcal{C}$.

We consider two related types of latent structures: a structure drawn uniformly from a class $\mathcal{C}$ of fixed-size subgraphs of $K_m$, and a randomized counterpart whose size is not fixed. The former includes the class of all $k$-vertex cliques, while the latter includes cliques with random vertex sets. Inspired by the permutation and independent sampling schemes of \citet{carpentier2025low}, we construct these classes and their associated priors using auxiliary random vectors $\dot \bfA=(\dot A_1,\ldots,\dot A_m)$ and $\tilde \bfA=(\tilde A_1,\ldots,\tilde A_m)$ whose coordinates determine vertex membership in the latent structure. Specifically, $\dot \bfA$ is a uniformly random permutation of $[m]$, while $\tilde \bfA$ has i.i.d. coordinates drawn uniformly from $[0,m]$.
The latent structure is then defined as a problem-specific function of $\dot \bfA$ or $\tilde \bfA$. We denote the corresponding distributions of $S$ by $\dot \nu$ and $\tilde \nu$ and their supports by $\dot {\mathcal C}$ and $\tilde {\mathcal C}$, respectively. Both constructions are symmetric: the distribution of $S$ is invariant under permutations of the vertex labels.

\subsection{Lower bounds for testing}\label{sec:preliminaries}

Below, we briefly review the second-moment and low-degree likelihood-ratio methods used to establish lower bounds for the Gaussian testing problem above; see, e.g., \citet{banks2017information,perry2018optimality,hopkins2018statistical,kunisky2019notes,schramm2022computational,wein2025computational,carpentier2025low}.

For each $n$, we assume throughout that the mixture alternative $\bbP_{\nu_n}$ is absolutely continuous with respect to the null distribution $\bbQ_n$. This condition holds when the one-coordinate signal distribution and the null distribution belong to the same regular exponential family with common support, as they do in every model considered in this paper. Define the likelihood ratio
\begin{align*}
\bfy \in \Omega^n \mapsto L_{\nu_n}(\bfy):=\frac{d\bbP_{\nu_n}}{d\bbQ_n}(\bfy).
\end{align*}

For a fixed structure $s\in\mathcal C_n$, similarly write $L_{s,n}:=d\bbP_{s,n}/d\bbQ_n$.
When the dependence on $n$ is clear, we suppress it and write $\bbP_\nu$, $\bbQ$, $L_\nu$, and $L_s$.

\subsubsection{Statistical lower bounds}\label{sec:statistical-lower-bounds}
For simple-versus-simple testing, it is well known that the minimum risk equals one minus the total variation distance between the two distributions.
Directly analyzing this distance, or the optimal likelihood-ratio test, is often difficult in high-dimensional models.
The second moment method instead controls the risk through the null second moment of $L_{\nu_n}$ \citep{arias2008searching,addario2010combinatorial,montanari2017limitation,banks2017information,perry2018optimality,perry2020statistical}.
The proof of the following lemma is given in \Cref{proof: second moment method}.

\begin{lemma}[Second moment method]\label{lemma: second moment method}
Suppose that, for some constant $C > 1$ and all sufficiently large $n$,
\begin{align*}
\bbE_{\bfY\sim\bbQ_n}\big[L_{\nu_n}(\bfY)^2\big]\leq C.
\end{align*}
Then every test $\Delta_n$ satisfies
\begin{align*}
R_{\nu_n}(\Delta_n)
\geq 1-\sqrt{1-C^{-1}},
\end{align*}
and hence the sequence of testing problems is statistically indistinguishable.
\end{lemma}

\subsubsection{Low-degree lower bounds}
We next describe the low-degree likelihood-ratio framework,  which characterizes the ability of polynomial statistics to distinguish the null and alternative in a sense made precise below. This framework has been widely used to predict computational hardness in testing problems; see, e.g., \citep{hopkins2017power,hopkins2018statistical,barak2019nearly,kunisky2019notes,wein2025computational}.
For this subsection, suppose that $\Omega\subseteq\bbR$ and that the polynomials under consideration are square-integrable under the null distribution. We recall the notion of strong separation from \citet[Definition~1.8]{bandeira2022franz}.

\begin{definition}[Strong separation]\label{definition:strong-separation}
A sequence of functions $f_n:\Omega^n\to\bbR$ strongly separates $\bbP_{\nu_n}$ and $\bbQ_n$ if
\begin{align*}
\frac{
\left|\bbE_{\bfY\sim\bbP_{\nu_n}}[f_n(\bfY)]
-\bbE_{\bfY\sim\bbQ_n}[f_n(\bfY)]\right|
}{
\sqrt{\max\left\{
\operatorname{Var}_{\bfY\sim\bbP_{\nu_n}}[f_n(\bfY)],
\operatorname{Var}_{\bfY\sim\bbQ_n}[f_n(\bfY)]
\right\}}
}
\mathrel{\to}\infty.
\end{align*}
\end{definition}

If $f_n$ strongly separates the hypotheses, thresholding $f_n$ at the midpoint of its two expectations achieves strong detection by Chebyshev's inequality.

Let $\mathcal V_n^{\leq D}\subseteq L^2(\bbQ_n)$ be the linear space of polynomials in the coordinates of $\bfY$ that have total degree at most $D$ and are square-integrable under the null distribution, and let $\operatorname{Proj}_n^{\leq D}$ denote the orthogonal projection onto this space.
Define the degree-$D$ likelihood ratio by
\begin{align*}
L_{\nu_n}^{\leq D}:=\operatorname{Proj}_n^{\leq D}L_{\nu_n}.
\end{align*}
The following lemma shows that boundedness of its second moment rules out strong separation by degree-$D$ polynomial statistics; see \citet[Proposition~6.2]{bandeira2022franz} and its proof.

\begin{lemma}[Low-degree second moment method]\label{lemma: low-degree second moment method}
For every $D\geq1$,
\begin{align}\label{eq:low-degree-variational-principle}
\bbE_{\bfY\sim\bbQ_n}\big[L_{\nu_n}^{\leq D}(\bfY)^2\big]-1
=\sup_{\substack{f\in\mathcal V_n^{\leq D}\\
\operatorname{Var}_{\bfY\sim\bbQ_n}[f(\bfY)]>0}}
\frac{
\left(\bbE_{\bfY\sim\bbP_{\nu_n}}[f(\bfY)]
-\bbE_{\bfY\sim\bbQ_n}[f(\bfY)]\right)^2
}{
\operatorname{Var}_{\bfY\sim\bbQ_n}[f(\bfY)]
}.
\end{align}
Consequently, if $D=D_n$ and
\begin{align*}
\bbE_{\bfY\sim\bbQ_n}\big[L_{\nu_n}^{\leq D_n}(\bfY)^2\big]=O(1),
\end{align*}
then no sequence of degree-$D_n$ polynomials strongly separates $\bbP_{\nu_n}$ and $\bbQ_n$.
\end{lemma}

The lemma is a statement about degree-$D$ polynomial statistics.
Interpreting it as evidence against all polynomial-time tests requires an additional hypothesis.
For the planted models studied here, we use the following working hypothesis from the low-degree literature \citep{hopkins2018statistical,kunisky2019notes,wein2025computational}.

\begin{conjecture}[Working low-degree hypothesis]\label{conjecture: low-degree polynomial conjecture}
Consider the sequences of planted and null distributions $\{\bbP_{\nu_n}\}$ and $\{\bbQ_n\}$ defined above.
If there exist $\epsilon>0$ and degrees $D_n\geq(\log n)^{1+\epsilon}$ such that
\begin{align*}
\bbE_{\bfY\sim\bbQ_n}\big[L_{\nu_n}^{\leq D_n}(\bfY)^2\big]=O(1),
\end{align*}
then no polynomial-time test achieves strong detection.
\end{conjecture}

This hypothesis should not be interpreted as a valid general principle for arbitrary testing problems. 
Earlier formulations admit counterexamples \citep{holmgren2021counterexamples}, the quasi-polynomial-time version is false in general \citep{buhai2025quasi}, and \citet{mao2026polynomial} gives a counterexample to a permutation-invariant polynomial-time formulation even after independent noise is added.
Recent results obtain rigorous consequences from low-degree indistinguishability under additional structural or algorithmic assumptions \citep{hsieh2026rigorous,chen2026rigorous}.
 Accordingly, our low-degree lower bounds apply to polynomial statistics of degree at most $D$, in the sense of strong separation, while their interpretation as a lower bound against all polynomial-time tests is conditional on the working hypothesis above.
The code-based counterexamples of \citet{buhai2025quasi} and \citet{mao2026polynomial}, which target a specific formulation in \citet[Conjecture~2.2.4]{hopkins2018statistical}, do not directly fit the one-sided Gaussian model studied here, whose alternative mixes over latent supports.
Indeed, consider the natural attempted embedding in which the $1$-entries of a counterexample form the support $S$ and the observations are $Y_\ell=\mu\mathds{1}\{\ell\in S\}+Z_\ell$, where $Z_\ell\stackrel{\mathrm{iid}}{\sim}\mathrm N(0,1)$.
Exact low-order uniformity of these constructions gives $\Pr\{\ell\in S\}=1/2$ and pairwise independence, so that $\bbE[|S|]=n/2$ and $|S|/n\to1/2$ in probability.
The degree-one likelihood-ratio projection under this embedding satisfies $\bbE_{\bfY\sim\bbQ_n}[L_{\nu_n}^{\leq 1}(\bfY)^2]=1+n\mu^2/4$.
If $n\mu^2\to\infty$, the average statistic strongly detects, whereas if $n\mu^2=O(1)$, convexity of relative entropy gives $D_{\mathrm{KL}}(\bbP_{\nu_n}\,\|\,\bbQ_n)\leq\bbE_{S\sim\nu_n}[D_{\mathrm{KL}}(\bbP_{S,n}\,\|\,\bbQ_n)]=\mu^2\bbE[|S|]/2=n\mu^2/4=O(1)$, and hence strong detection is statistically impossible.
Thus, this natural embedding is either detectable by a degree-one statistic or statistically indistinguishable, and cannot reproduce the counterexamples' combination of low-degree indistinguishability and efficient detection.
The sign-based Gaussian extension of \citet{buhai2025quasi} uses codeword signs rather than one-sided support shifts, while their spectral construction uses a nonproduct null, so both fall outside the present model.

\begin{remark}[Conditional second moment method]\label{rem:conditional-second-moment-method}
Boundedness of the full and low-degree likelihood-ratio second moments is sufficient, but not necessary, to rule out strong detection and low-degree strong separation, respectively (\Cref{lemma: second moment method,lemma: low-degree second moment method}). 
Rare events under the alternative can make the second moment diverge even when strong detection remains impossible.
The conditional second moment method addresses this issue by restricting the alternative to a high-probability event that excludes such rare configurations before applying the second moment argument \citep{banks2017information,perry2018optimality,perry2020statistical,bandeira2022franz,dhawan2025detection}. In our analysis, we primarily use the conditional second moment method to transfer lower bounds from fixed-size models to their random-size counterparts through high-probability truncation; see, e.g., \Cref{corollary: random-size clique problem,theorem: multiple cliques problem,theorem: thick cycle problem} and their proofs.

\end{remark}

\subsection{Statistical and computational thresholds}\label{sec:gaussian-testing-setup}
To characterize the statistical and conjectural computational thresholds for the Gaussian combinatorial testing problem, we combine upper bounds on the testing risk with lower bounds from the full and low-degree likelihood-ratio second moments. Both second moments admit exact representations in terms of the overlap $|S\cap S'|$ between two independent draws $S,S'\sim\nu_n$. The following proposition records these identities; see \citet{arias2008searching,banks2017information} for the full second moment and \citet{hopkins2018statistical,kunisky2019notes} for its low-degree counterpart.

\begin{proposition}\label{proposition: gaussian overlapping pairs}
Consider the hypothesis testing problem defined in \Cref{equation: gaussian graph testing problem} with a prior distribution $\nu_n$ for $S$. Then,
\[
\bbE_{\bfY\sim \bbQ_n}\big[L_{\nu_n}(\bfY)^2\big]
= \bbE_{S,S'\sim\nu_n}\Big[\exp\big(\mu^2|S\cap S'|\big)\Big]
\]
and 
\[
\bbE_{\bfY\sim \bbQ_n}\big[L_{\nu_n}^{\leq D}(\bfY)^2\big] =  \bbE_{S,S'\sim\nu_n}\bigg[\sum_{d=0}^D \frac{\mu^{2d}|S\cap S'|^d}{d!}\bigg].
\]
\end{proposition}
This proposition reduces the derivation of statistical and computational lower bounds to the analysis of the overlap distribution induced by the alternative. In \Cref{sec:results}, we carry out this analysis for several combinatorial classes $\mathcal C$. 

For the upper-bound analysis, we consider three tests: the \emph{outlier test}, the \emph{average test}, and the \emph{scan test}, each obtained by thresholding the corresponding test statistic
\begin{equation}\label{equation: gaussian test statistics}
	\begin{aligned}
		T_{\mathrm{out}}(\bfy):=\max_{(i,j) \in [m]_2} y_{ij}, \qquad T_{\mathrm{avg}}(\bfy):=\sum_{(i,j) \in [m]_2} y_{ij}, \qquad \text{and} \qquad  T_{\mathrm{scan}}(\bfy):=\max_{s\in \mathcal C}\sum_{(i,j)\in s} y_{ij},
	\end{aligned}
\end{equation}
respectively. These tests, or variants of them, are widely used to establish upper bounds in signal detection problems; see, e.g., \citep{donoho2004higher,arias2008searching,arias2011detection,addario2010combinatorial,butucea2013detection,arias2014community,ma2015computational,sharpnack2016detecting,dadon2024detection,elimelech2025detecting}.
The outlier test uses the fact that the maximum of $n$ standard normal random variables is approximately $\sqrt{2\log n}$ and can therefore detect anomalous coordinates regardless of the size or shape of the latent structure $S$.
The average test uses the accumulated signal $\mu|S|$: the statistic $T_{\mathrm{avg}}(\bfY)$ follows $\mathrm N(0,n)$ under the null and $\mathrm N(\mu|S|,n)$ conditional on $S$ under the alternative.
Finally, the scan test exploits the combinatorial structure of $S$ by maximizing the total observed signal over all candidates in $\mathcal C$.
The outlier and average tests have computational complexity of order $m^2$ regardless of the alternative.
A direct implementation of the scan test instead requires enumerating $\mathcal C$ and is therefore superpolynomial for many examples considered in this paper.
 
The next result characterizes the statistical performance of these tests; see, e.g., \citet{addario2010combinatorial} for related bounds and \Cref{proof: gaussian test upper bound} for a proof.
For simplicity, the proposition assumes that $|S|=K$ deterministically; extensions to the random-size models considered here are justified in their respective proofs.
\begin{proposition}[Outlier, average, and scan tests]\label{proposition: gaussian test upper bound}
Consider the test statistics defined in \Cref{equation: gaussian test statistics} for the hypothesis testing problem defined in \Cref{equation: gaussian graph testing problem}, where every $s\in\mathcal C$ satisfies $|s|=K\geq1$, and let $\nu$ be a prior distribution on $\mathcal C$. Let $\delta \in (0,1)$ be a target risk level.
 Then the following holds:
\begin{itemize}
	\item Consider the threshold $t_{\mathrm{out}}:=\sqrt{2\log(2n/\delta)}$. The average risk of the outlier test defined by $\Delta_{\mathrm{out}}(\bfy):=\mathds1\{T_{\mathrm{out}}(\bfy) > t_{\mathrm{out}}\}$ satisfies $R_\nu(\Delta_{\mathrm{out}}) \leq \delta$ if
	\begin{align*}
	\mu \geq \sqrt{2\log\frac{2n}{\delta}}+\sqrt{2\log\frac{2}{\delta}}.
	\end{align*}
	\item Consider the threshold $t_{\mathrm{avg}}:=\mu K/2$. The average risk of the average test defined by $\Delta_{\mathrm{avg}}(\bfy):=\mathds1\{T_{\mathrm{avg}}(\bfy) > t_{\mathrm{avg}}\}$ satisfies $R_\nu(\Delta_{\mathrm{avg}}) \leq \delta$ if
	\begin{align*}
	\mu \geq \sqrt{\frac{8n}{K^2}\log\frac{2}{\delta}}.
	\end{align*}
	\item Consider the threshold $t_{\mathrm{scan}}:=\big(\mu K+\sqrt{2K\log|\mathcal C|}\big)/2$. The average risk of the scan test defined by $\Delta_{\mathrm{scan}}(\bfy):=\mathds1\{T_{\mathrm{scan}}(\bfy) > t_{\mathrm{scan}}\}$ satisfies $R_\nu(\Delta_{\mathrm{scan}}) \leq \delta$ if
	\begin{align*}
	\mu \geq \sqrt{\frac{2}{K}\log|\mathcal C|}+2\sqrt{\frac{2}{K}\log\frac{2}{\delta}}.
	\end{align*}
\end{itemize}
\end{proposition}

\begin{remark}[Degree test]\label{remark: gaussian degree test}
	Another commonly used computationally efficient procedure, also with complexity of order $m^2$, is the degree test \citep{arias2014community,hajek2018submatrix,elimelech2025detecting}, based on the statistic
    \[
    T_{\mathrm{deg}}(\bfy):=\max_{i\in[m]}T_i^{\mathrm{deg}}(\bfy),
    \]
    where $T_i^{\mathrm{deg}}(\bfy):=\sum_{\{e\in [m]_2\,:\,i\in e\}} y_e$ is the weighted degree of vertex $i\in[m]$.
	For $s\subseteq[m]_2$, let $\mathfrak{d}_s(i):=|\{e\in s\,:\,i\in e\}|$ denote the degree of vertex $i$ in the graph with vertex set $[m]$ and edge set $s$.
	Suppose that $S\sim \nu$ satisfies $\max_{i\in[m]}\mathfrak{d}_S(i)=\mathfrak{d}_{\max}\geq 1$ almost surely. The following bound summarizes the performance of the degree test; see \Cref{proof: gaussian test upper bound} for a proof.
	\begin{itemize}
	\item Let $\delta \in (0,1)$ be a target risk level. Consider the threshold $ t_{\mathrm{deg}}:=\mathfrak{d}_{\max}\mu/2.$ Then, the average risk of the degree test defined by $\Delta_{\mathrm{deg}}(\bfy):=\mathds1\{T_{\mathrm{deg}}(\bfy)>t_{\mathrm{deg}}\}$ satisfies $R_\nu(\Delta_{\mathrm{deg}})\leq\delta$ if 
	$$
	\mu\geq\sqrt{\frac{8(m-1)}{\mathfrak d_{\max}^2}\log\frac{2m}{\delta}}.
	$$
	\end{itemize} 
	The degree test is not used for the problems considered in \Cref{sec:results}, as it yields weaker upper bounds than the other tests in those settings. It can nevertheless provide a sharper upper bound for other structures, including the planted star studied in \Cref{appendix: k-size star}.

\end{remark}

To derive upper bounds for the statistical and computational thresholds, we use an appropriate combination of the tests in \Cref{proposition: gaussian test upper bound}, depending on the scaling of $\mu$.
Since the scan test can require superpolynomial time in $m$, the statistical and computational upper bounds may differ in regimes where scanning is more powerful than the available polynomial-time tests.
We now formalize the two thresholds. Because $\ll$ and $\gg$ denote polynomial separation, these definitions identify thresholds at polynomial order; unless stated otherwise, the main results therefore suppress polylogarithmic and other subpolynomial factors.
\begin{definition}[Statistical and computational thresholds]\label{definition: statistical and computational thresholds}
Consider a sequence of Gaussian graph testing problems defined in \Cref{equation: gaussian graph testing problem}, with prior $\nu=\nu_n$. A positive sequence $\mu_{\mathrm{stat}}^\ast$ is called a statistical threshold if the following holds for every signal strength sequence $\mu=\mu_n$:
\begin{itemize}
	\item if $\mu \gg \mu_{\mathrm{stat}}^\ast$, then there exists a sequence of tests $\Delta_n$ that achieves strong detection;
	\item if $\mu \ll \mu_{\mathrm{stat}}^\ast$, then no sequence of tests can achieve strong detection.
\end{itemize}
Similarly, a positive sequence $\mu_{\mathrm{comp}}^\ast$ is called a computational threshold if
\begin{itemize}
	\item if $\mu \gg \mu_{\mathrm{comp}}^\ast$, then there exists a sequence of polynomial-time tests $\Delta_n$ that achieves strong detection;
	\item if $\mu \ll \mu_{\mathrm{comp}}^\ast$, then, under \Cref{conjecture: low-degree polynomial conjecture}, no sequence of polynomial-time tests can achieve strong detection.
\end{itemize}
\end{definition}

Establishing either threshold requires matching upper and lower bounds. An upper bound $\mu^{\mathrm{UB}}$ follows by demonstrating a test that achieves strong detection whenever $\mu\gg\mu^{\mathrm{UB}}$. For a statistical lower bound $\mu^{\mathrm{LB}}$, it is sufficient to show that the full likelihood-ratio second moment remains bounded whenever $\mu\ll\mu^{\mathrm{LB}}$. For a computational lower bound, we instead bound the low-degree likelihood-ratio second moment and invoke \Cref{conjecture: low-degree polynomial conjecture} to interpret the resulting obstruction as evidence against all polynomial-time tests. If $\mu^{\mathrm{UB}}/\mu^{\mathrm{LB}}=n^{o(1)}$, either bound provides a valid representative of the polynomial-order threshold. We say that a statistical-computational gap occurs when $\mu_{\mathrm{stat}}^\ast\ll\mu_{\mathrm{comp}}^\ast$, with the computational lower bound understood in the context of \Cref{conjecture: low-degree polynomial conjecture}.

\section{Statistical-computational gaps in Gaussian graph testing}\label{sec:results}

We establish statistical-computational gaps, under \Cref{conjecture: low-degree polynomial conjecture}, for three classes of Gaussian combinatorial testing problems: planted cliques (\Cref{section: planted clique}), multiple planted cliques and hypercliques (\Cref{section: multiple planted cliques}), and thick cycles and paths (\Cref{section: thick cycle and path}).
Throughout, we consider symmetric alternatives induced by combinatorial classes of edge or hyperedge sets, with latent structures of either fixed or random size.
The appendix, beginning with \Cref{appendix: disjoint sets}, analyzes additional classes for which no gap occurs: $K$-disjoint sets, $K$-sets, perfect matchings, stars, spanning trees, $k$-stars, and paths and cycles of length $k$; \citet{addario2010combinatorial} derived statistical thresholds for the first five classes.
All results are summarized in \Cref{tab:summary-of-problems}.

\subsection{A canonical example: the planted clique}\label{section: planted clique}
We begin with the case in which the latent structure $S$ is a clique, a canonical benchmark for studying the limitations of polynomial-time algorithms in high-dimensional statistics and theoretical computer science.
Whereas the Erd\H{o}s--R\'enyi model has been studied extensively, its Gaussian counterpart is less well understood.
We therefore investigate the statistical-computational gap for planted clique detection under Gaussian observations.

Formally, we call a subset $s$ of the edges of $K_m$ a \emph{$k$-vertex clique} if $s$ is the edge set of a complete subgraph on $k$ vertices. We consider two scenarios: (i) the fixed-size case, in which $S$ is uniformly distributed over all $k$-vertex cliques in an $m$-vertex graph, and (ii) the random-size case, in which each of the $m$ vertices is included independently with probability $k/m$. As indicated in \Cref{sec:gaussian-testing}, the fixed-size model can be formulated using the auxiliary random vector $\dot \bfA = (\dot A_1,\ldots,\dot A_m)$, a uniformly random permutation of $[m]$. The random structure $S$ is then defined by
\begin{align}\label{equation: clique latent structure}
	S:= \{(i,j)\in[m]_2: \mathds{1}(\dot A_i\leq k)\mathds{1}(\dot A_j\leq k)=1\}
\end{align}
This construction induces the prior $\dot \nu$ on the set of all $k$-vertex cliques.
To model a random-size clique with expected vertex count $k$, we define $S$ analogously using the random vector $\tilde \bfA = (\tilde A_1, \ldots, \tilde A_m)$, whose coordinates are i.i.d.\ uniform random variables on $[0,m]$.

We start with the fixed-size case.
For this canonical benchmark, we state the statistical and low-degree lower bounds separately in \Cref{lemma: statistical lower bound for fixed-size clique problem,lemma: computational lower bound for fixed-size clique problem} and then combine them with the tests in \Cref{proposition: gaussian test upper bound} to obtain \Cref{theorem: fixed-size clique problem}.
For the subsequent models, we streamline the main text by stating the resulting threshold theorems directly and deferring the analogous intermediate bounds and proofs to the appendix.

Our first result gives a sufficient condition for the likelihood-ratio second moment to remain bounded, yielding a lower bound $\mu^{\textrm{LB}}_{\mathrm{stat}}$ for the statistical threshold $\mu_{\mathrm{stat}}^\ast$. The analysis follows the approach of \citet{addario2010combinatorial} but uses Finner's inequality \citep{finner1992generalization}; see \Cref{lemma: Finner's inequality} in \Cref{sec:app-technical-lemmas}. This inequality is well suited to our setting and extends naturally to the multiple-clique problem considered later.

\begin{lemma}[Statistical lower bound]\label{lemma: statistical lower bound for fixed-size clique problem}
    Fix any constant $\delta\in(0,0.1)$. The second moment of the likelihood ratio is bounded as long as $\mu$ satisfies: (i) if $k^3m^{-2} \leq m^{-\delta}$, then $\mu \leq \sqrt{\delta k^{-1}\log m}$; (ii) if $k^3m^{-2} =\Omega(1)$, then $\mu = O(mk^{-2})$.
\end{lemma}

\textit{Proof sketch.} The overlap $|S\cap S'|$ has the same distribution as $\binom H2$, where $H$ is a hypergeometric random variable with parameters $(m,k,k)$. By the convex ordering of the hypergeometric distribution by the binomial distribution \citep{hoeffding1963probability}, it suffices to work with $B\sim\mathrm{Bin}(k,k/m)$ instead of $H$. Thus, the lower-bound analysis reduces to showing that $\bbE\big[\exp(\mu^2\binom{B}{2})\big]$ is bounded. Writing $B \stackrel{d}{=} \sum_{i=1}^k B_i$, where the $B_i$ are i.i.d.\ Bernoulli random variables with success probability $k/m$, we have $\binom B2 \stackrel{d}{=} \sum_{i<j} B_i B_j$. Applying Finner's inequality (\Cref{lemma: Finner's inequality}) gives
\begin{align*}
	\ee_B\bigg[\exp\bigg(\mu^2\sum_{i<j} B_i B_j\bigg)\bigg] &\leq \prod_{i<j} \ee\big[\exp(\mu^2B_i B_j)^{k-1}\big]^{\frac{1}{k-1}}\\
	&= \ee\big[\exp\{\mu^2(k-1)B_1 B_2\}\big]^{k/2}\\
	&\leq \exp\bigg(\frac{k^3}{2m^2}\big(\exp(\mu^2k)-1\big)\bigg),
\end{align*}
where the last inequality follows from the moment generating function of a Bernoulli random variable $B_1B_2$ with success probability $(k/m)^2$ and the fact that $1+x\leq e^x$ for any $x\in\mathbb R$. Analyzing the asymptotic behavior of the final expression gives the conditions in \Cref{lemma: statistical lower bound for fixed-size clique problem}; the details appear in \Cref{proof: statistical lower bound for fixed-size clique problem}. \qed

The preceding lower bound suggests two regimes for the statistical threshold, separated at $k=m^{2/3}$. We next give matching upper bounds. Recall that $n=\binom{m}{2}$, $K=\binom{k}{2}$, and $|\dot {\mathcal C}|=\binom{m}{k}$ for the fixed-size planted clique problem. By \Cref{proposition: gaussian test upper bound}, the scan and average tests achieve strong detection when $\mu=\omega(\sqrt{k^{-1}\log m})$ and $\mu=\omega(mk^{-2})$, respectively. The scan test attains the matching upper bound in the first regime, but its direct implementation searches all $\binom{m}{k}$ candidates in $\dot{\mathcal C}$ and therefore requires superpolynomial time. The average test runs in polynomial time and attains the matching upper bound in the second regime.

The following lemma provides a lower bound $\mu^{\textrm{LB}}_{\mathrm{comp}}$ for the computational threshold $\mu_{\mathrm{comp}}^\ast$ by analyzing the second moment of the low-degree likelihood ratio. The argument uses moment bounds for sub-Poissonian random variables due to \citet{ahle2022sharp}.
\begin{lemma}[Computational lower bound]\label{lemma: computational lower bound for fixed-size clique problem}
    Fix any constant $\delta\in(0,0.1)$. For the sequence $D_m=\log^{1+\epsilon} m$ with $\delta>\epsilon>0$, the second moment of the low-degree likelihood ratio is bounded as long as $\mu$ satisfies: (i) if $k^2m^{-1} \leq m^{-\delta}$, then $\mu \leq C (\log m)^{(1-\epsilon)/2}$ for some sufficiently small constant $C>0$; (ii) if $k^2m^{-1} \geq \log^{1+\delta} m$, then $\mu = O\big(mk^{-2}\big)$.
\end{lemma}

\textit{Proof sketch.} The convex ordering of the hypergeometric distribution by the binomial distribution gives
$$
\ee_{S,S'\sim\nu}\big[|S\cap S'|^d\big] =\ee\bigg[\binom{H}{2}^d\bigg] \leq \ee\bigg[\bigg(\frac{H^2}{2}\bigg)^d\bigg] \leq \ee\bigg[\bigg(\frac{B^2}{2}\bigg)^d\bigg],
$$
where $B\sim\mathrm{Bin}(k,k/m)$. Thus, the moments of the binomial random variable $B$ up to order $2D$ must be controlled carefully to bound the second moment of the low-degree likelihood ratio.
Applying the moment bound of \citet{ahle2022sharp}, we obtain
\begin{align*}
	\bbE_{S,S'\sim\nu}\bigg[\sum_{d=0}^D \frac{\mu^{2d}|S\cap S'|^d}{d!}\bigg]
	\leq \bbE\bigg[\sum_{d=0}^D \frac{\mu^{2d}B^{2d}}{2^d d!}\bigg]
	\leq 1+ \sum_{d=1}^D \frac{2^d\mu^{2d}}{d!} \bigg(\frac{d}{\log(1+2dm/k^2)}\bigg)^{2d}.
\end{align*}
It remains to determine when this expression is bounded as a function of $\mu$, $k$, $m$, and $D$.
The details are given in \Cref{proof: computational lower bound for fixed-size clique problem}. \qed

Together with matching polynomial-time upper bounds, this result identifies the polynomial order of $\mu_{\mathrm{comp}}^\ast$ in both regimes of $k$.
Specifically, the outlier test achieves strong detection if $\mu=\omega(\sqrt{\log m})$, while the average test achieves strong detection if $\mu=\omega(mk^{-2})$.

\begin{figure}[t]
\centering
\includegraphics[width=0.36\textwidth]{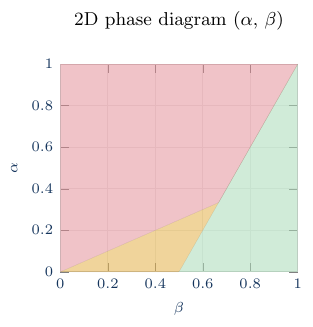}
\caption{Phase diagram for the planted clique problem under the parameterization $\mu=m^{-\alpha}$ and $k=m^{\beta}$.
Green, yellow, and pink denote the easy, hard, and impossible regions, respectively.}\label{fig:planted-clique-phase-diagram}
\end{figure}

With \Cref{lemma: statistical lower bound for fixed-size clique problem,lemma: computational lower bound for fixed-size clique problem} in place, we now state the statistical and computational thresholds for the fixed-size planted clique problem. Although the preceding lower and upper bounds retain subpolynomial factors, we focus here on the setting in which $k$ has positive polynomial order in $m$ and report the resulting threshold scales up to subpolynomial factors. This convention keeps the statement transparent while preserving the polynomial-order thresholds. The proof is given in \Cref{proof: fixed-size clique problem}, and \Cref{fig:planted-clique-phase-diagram} visualizes the result.
\begin{theorem}\label{theorem: fixed-size clique problem}
Consider the Gaussian graph testing problem with the latent structure $S$ defined in \Cref{equation: clique latent structure} and prior distribution $\dot \nu$. Suppose that $k$ is of positive polynomial order in $m$.
Then, the statistical and computational thresholds are 
\begin{align*}
\mu_{\mathrm{stat}}^\ast=\left\{
\begin{array}{@{}c@{\quad}l}
\sqrt{\frac{1}{k}}, & k \ll m^{2/3}, \\
\frac{m}{k^2}, & k \gg m^{2/3}
\end{array}
\right.
\qquad \text{and} \qquad 
\mu_{\mathrm{comp}}^\ast=\left\{
\begin{array}{@{}c@{\quad}l}
\vphantom{\sqrt{\frac{1}{k}}}1, & k \ll m^{1/2}, \\
\frac{m}{k^2}, & k \gg m^{1/2}.
\end{array}
\right.
\end{align*}
\end{theorem}

The statistical-computational phase diagram in \Cref{fig:planted-clique-phase-diagram} is not new, but is consistent with existing results in the literature on Gaussian submatrix detection.
In particular, \citet{ma2015computational} established the corresponding computational lower bounds through reductions from planted clique under the planted-clique conjecture (see, e.g., \citealt[Conjecture 4.13]{alon2007testing} and \citealt[Hypothesis $\mathsf B_{\mathsf {PC}}$]{berthet2013complexity}), and the subsequent work of \citet{brennan2019universality} extended these reduction-based bounds to submatrix detection with matching row and column supports.
The latter results imply the same computational lower bounds for the Gaussian planted clique model by restricting the observations to the strictly upper-triangular entries of the matrix.
In contrast to these reduction-based approaches, our analysis recovers equivalent computational lower bounds more directly within the low-degree polynomial framework by studying the second moment of the low-degree likelihood ratio.
This provides a complementary perspective on computational hardness and a unified approach that extends naturally to the more general combinatorial structures considered below.

We next consider the random-size case, which has the same statistical and computational thresholds as the fixed-size model.
The key fact is that the number of planted vertices concentrates around $k$ with high probability.
The detailed proof of the following corollary is given in \Cref{proof: random-size clique problem}.
\begin{corollary}\label{corollary: random-size clique problem}
	Consider the Gaussian graph testing problem with the latent structure $S$ defined in \Cref{equation: clique latent structure} and prior distribution $\tilde \nu$. Suppose that $k$ is of positive polynomial order in $m$. Then, the statistical and computational thresholds $\mu_{\mathrm{stat}}^\ast$ and $\mu_{\mathrm{comp}}^\ast$ are identical to those in \Cref{theorem: fixed-size clique problem}.
\end{corollary}
The fixed-size guarantees for the average and scan tests do not apply directly when the clique size is random, but the tests can be modified to detect a clique on $ck$ vertices for any fixed constant $c\in(0,1)$.
The planted clique has at least $ck$ vertices with high probability, and every larger clique contains a clique on $ck$ vertices.
The complementary event, on which the clique has fewer than $ck$ vertices, contributes only a negligible amount to the risk.
	
For the statistical lower bound, the ordinary second-moment analysis used in the fixed-size case yields a suboptimal result: on a rare event, the overlap is large enough to make $\ee[\exp(\mu^2 |S\cap S'|)]$ diverge even when strong detection is impossible. See \Cref{proof: random-size clique problem} for details. To bypass this limitation, we use the conditional second moment method (\Cref{rem:conditional-second-moment-method}) and condition on a high-probability event that excludes realizations with a disproportionate contribution to the second moment. We then show through stochastic dominance that the random-size clique detection problem is at least as hard as the fixed-size problem.
In particular, we truncate and renormalize the prior $\tilde \nu$ to the event that the clique size is at most $Ck$ for a fixed constant $C>1$. The resulting loss in risk is negligible because the induced alternative distribution changes by a vanishing amount in total variation distance. Under the modified prior, we show that $|S\cap S'|$ is stochastically dominated by $\binom H2$, where $H$ is a hypergeometric random variable with parameters $(m,Ck,Ck)$. This is the overlap distribution in the fixed-size model with clique size $Ck$. Intuitively, two random subsets of size at most $Ck$ tend to share fewer vertices than two subsets of fixed size $Ck$.

The same truncation argument yields the matching computational lower bound for the random-size model.
For this low-degree bound, however, truncation is not necessary: direct analysis of the ordinary low-degree likelihood-ratio second moment gives the same result.

\subsection{Multiple planted cliques and hypercliques}\label{section: multiple planted cliques}

\begin{figure}[t]
	\includegraphics[width=\textwidth]{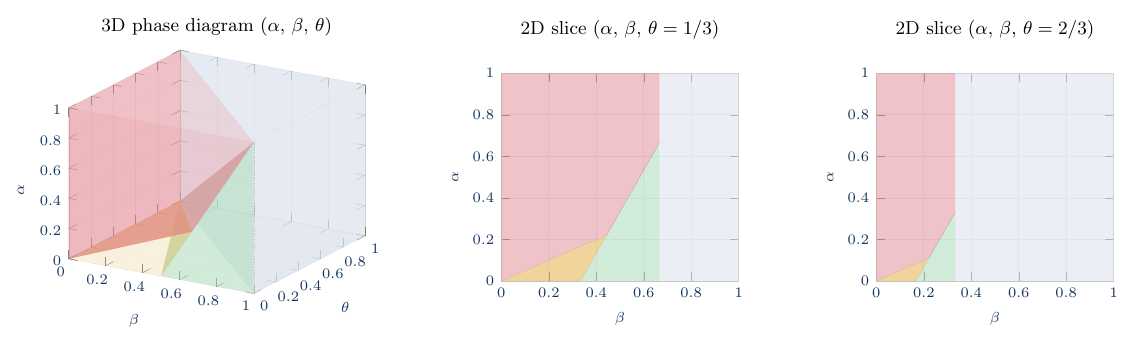}
	\caption{Phase diagram for the multiple planted cliques problem under the parameterization $\mu=m^{-\alpha}$, $k=m^{\beta}$, and $\gamma=m^{\theta}$.
The left panel shows the three-dimensional phase diagram, and the two right panels show two-dimensional slices at $\theta=1/3$ and $\theta=2/3$.
Green, yellow, and pink denote the easy, hard, and impossible regions, respectively; blue-gray denotes parameter values outside the admissible range.}\label{fig:multiple-cliques-phase-diagram}
\end{figure}

We treat the single-clique problem separately as a canonical benchmark that introduces the main proof ideas before extending them to several interacting planted components.
We now generalize the planted clique setting to alternatives consisting of multiple disjoint planted cliques.
Let $m,k,\gamma$
be positive integers satisfying $\gamma k\le m$, and $S$ be the union of $\gamma$ vertex-disjoint cliques, each with either a fixed vertex count $k$ or an expected vertex count $k$.
For the fixed-size model, let $\dot \bfA = (\dot A_1,\ldots, \dot A_m)$ be a uniformly random permutation of $[m]$ and define
\begin{align}\label{equation: multiple cliques latent structure}
	S:= \bigcup_{r=1}^{\gamma}\bigg\{(i,j)\in[m]_2:(r-1)k < \dot A_i, \dot A_j \le r k
\bigg\}.
\end{align}
A closely related model is the stochastic block model, whose statistical-computational gap has been studied extensively under the Erd\H{o}s--R\'enyi observation model; see, e.g., \citet{hopkins2017efficient,hopkins2018statistical,bandeira2021spectral,carpentier2025low,ding2025lowdegree,sohn2025sharp}. When $\gamma k=m$, the latent structure partitions all vertices into $\gamma$ communities, as in the stochastic block model. The following theorem analyzes the statistical and computational thresholds in both the fixed-size and random-size cases. In the random-size model, $S$ is defined as in \Cref{equation: multiple cliques latent structure}, with $\dot \bfA$ replaced by $\tilde \bfA = (\tilde A_1,\ldots, \tilde A_m)$, whose entries are i.i.d.\ uniform random variables on $[0,m]$.
The thresholds match those in \Cref{theorem: fixed-size clique problem} after replacing the ambient vertex count $m$ by $m/\gamma$, and the resulting phase diagram appears in \Cref{fig:multiple-cliques-phase-diagram}.
\begin{theorem}\label{theorem: multiple cliques problem}
	Consider the Gaussian graph testing problem with the latent structure $S$ defined in \Cref{equation: multiple cliques latent structure} and prior distribution $\dot \nu$ or $\tilde \nu$. Suppose that $k$ is of positive polynomial order in $m$. Then, in both cases, the statistical and computational thresholds are 
	\begin{align*}
		\mu_{\mathrm{stat}}^\ast=\left\{
		\begin{array}{@{}c@{\quad}l}
		\sqrt{\frac{1}{k}}, & \gamma^2 k^3 \ll m^2, \\
		\frac{m}{\gamma k^2}, & \gamma^2 k^3 \gg m^2
		\end{array}
		\right.
		\qquad \text{and} \qquad 
		\mu_{\mathrm{comp}}^\ast=\left\{
		\begin{array}{@{}c@{\quad}l}
		\vphantom{\sqrt{\frac{1}{k}}}1, & \gamma k^2 \ll m, \\
		\frac{m}{\gamma k^2}, & \gamma k^2 \gg m.
		\end{array}
		\right.
	\end{align*}
\end{theorem}

We next sketch the argument; the full proof is given in \Cref{proof: multiple cliques problem}.
For the statistical lower bound, a key step is to prove that the likelihood-ratio second moment in the multiple-clique setting is bounded by the $\gamma^2$ power of an appropriate binomial moment, using the conditional negative association of the clique overlaps. In the fixed-size case, let $S_r$ and $S'_r$ denote the $r$th planted cliques in $S$ and $S'$, respectively. We show that 
\begin{align*}
\bbE_{S,S'\sim\dot\nu}\Big[\exp\big(\mu^2|S\cap S'|\big)\Big]
&=\bbE_{S,S'\sim\dot\nu}\bigg[\exp\bigg(\mu^2\sum_{a,b=1}^\gamma |S_a\cap S'_b|\bigg)\bigg]\\
&\leq\bbE_B\bigg[\exp\bigg\{\mu^2\binom B2\bigg\}\bigg]^{\gamma^2}
\end{align*}
where $B\sim\mathrm{Binomial}(k,q_{\gamma,k})$ and $q_{\gamma,k}:=k/(m-(\gamma-1)k)$. 
We bound the last expectation using Finner's inequality, following the same argument as in the proof for the single-clique case (\Cref{lemma: statistical lower bound for fixed-size clique problem}).
For the random-size model, as in the proof of \Cref{corollary: random-size clique problem}, we apply the conditional second moment method after truncating the prior to the event that all $\gamma$ cliques have size at most $Ck$ for some constant $C>1$.

The computational lower bound requires a more delicate comparison.
One possible route is to represent the overlap as a U-statistic and control its moments using decoupling inequalities \citep{delapena1999decoupling,gine2000exponential}, as in the analysis of a different testing problem by \citet{mao2025information}.
We provide a proof based on this approach in \Cref{sec:alternative-proof}.
In our main analysis, however,
we use a direct graph-theoretic argument.
We first analyze the random-size model, where the overlap is the number of edges in a random $\gamma^2$-multiple-clique model, and bound its moments by those of an auxiliary single-clique model.
Specifically, we show that for every integer $d\geq 1$,
$$
\mathbb E_{S,S'\sim\tilde\nu}[|S\cap S'|^d]
\leq\mathbb E\bigg[\binom B2^d\bigg],
$$
where $B\sim\operatorname{Binomial}(
m,\gamma k^2/m^2)$.
This comparison lets us apply the tools developed for the single-clique problem, and a truncation argument then transfers the resulting bound to the fixed-size model.

\paragraph{Multiple planted hypercliques.}
We further generalize the multiple planted clique model to $t$-uniform planted hypercliques.
Formally, a $t$-uniform \emph{hyperclique} is the edge set of a complete $t$-uniform hypergraph. Each hyperedge is represented by an increasing $t$-tuple rather than a pair of vertices.
For $\bfI=(i_1,\ldots,i_t)\in[m]_t$, define the latent structure $S$ by
\begin{align*}
	S:= \bigcup_{r=1}^{\gamma}\bigg\{\bfI\in[m]_t: (r-1)k < \dot A_{i_1},\ldots, \dot A_{i_t} \le r k
\bigg\},
\end{align*}
where, as usual, $(\dot A_1,\ldots,\dot A_m)$ is a random permutation of $[m]$. The random-size case is defined analogously by replacing $\dot{\bfA}$ with $\tilde{\bfA}=(\tilde A_1,\ldots,\tilde A_m)$.

\begin{figure}[t]
\includegraphics[width=\textwidth]{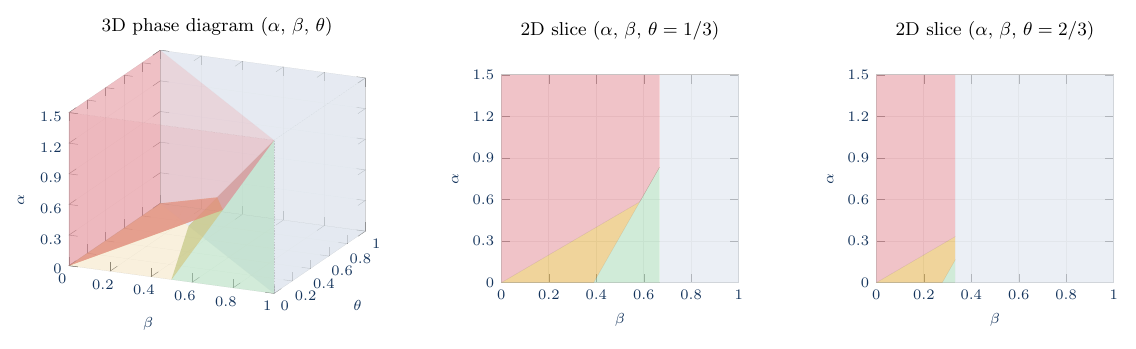}
\caption{Phase diagram for the multiple planted hypercliques problem with $t=3$ under the parameterization $\mu=m^{-\alpha}$, $k=m^{\beta}$, and $\gamma=m^{\theta}$.
The left panel shows the three-dimensional phase diagram, and the two right panels show two-dimensional slices at $\theta=1/3$ and $\theta=2/3$.
Green, yellow, and pink denote the easy, hard, and impossible regions, respectively; blue-gray denotes parameter values outside the admissible range.}\label{fig:multiple-hypercliques-phase-diagram}
\end{figure}

The proof strategy for multiple planted cliques extends to this hypergraph setting, yielding the following corollary.
The detailed proof is given in \Cref{proof: multiple planted hypercliques problem}.
For $t=3$, the resulting phase diagram is shown in \Cref{fig:multiple-hypercliques-phase-diagram}.
\begin{corollary}\label{corollary: multiple planted hypercliques problem}
	Consider the Gaussian hypergraph testing problem with the latent structure $S$ defined above and prior distribution $\dot \nu$ or $\tilde \nu$. Suppose that $t\geq 2$ is a fixed integer and $k$ is of positive polynomial order in $m$. Then, in both cases, the statistical and computational thresholds are 
	\begin{align*}
		\mu_{\mathrm{stat}}^\ast=\left\{
		\begin{array}{@{}c@{\quad}l}
		\sqrt{\frac{1}{k^{t-1}}}, & \gamma^2 k^{t+1} \ll m^t, \\
		\frac{m^{t/2}}{\gamma k^t}, & \gamma^2 k^{t+1} \gg m^t
		\end{array}
		\right.
		\qquad \text{and} \qquad 
		\mu_{\mathrm{comp}}^\ast=\left\{
		\begin{array}{@{}c@{\quad}l}
		\vphantom{\sqrt{\frac{1}{k^{t-1}}}}1, & \gamma^2 k^{2t} \ll m^t, \\
		\frac{m^{t/2}}{\gamma k^t}, & \gamma^2 k^{2t} \gg m^t.
		\end{array}
		\right.
	\end{align*}
\end{corollary}

\begin{remark}[Multiple $p$-partite hypercliques]\label{remark: multiple p partite hypercliques}
Another possible generalization is the detection of multiple $p$-partite hypercliques in $p$-partite hypergraphs, for any fixed integer $p\geq 2$. The case $p=2$ is equivalent to detecting multiple equal-sized submatrices with disjoint row supports and disjoint column supports, while $p\geq 3$ corresponds to multiple tensor detection. The submatrix detection problem and its statistical-computational gap have been studied extensively; see, e.g., \citet{butucea2013detection,ma2015computational,brennan2019universality,dadon2024detection}. It turns out that the statistical and computational thresholds for this problem coincide with those in \Cref{corollary: multiple planted hypercliques problem} with $t$ replaced by $p$, since the key arguments used to establish that result also apply to the $p$-partite setting; see \Cref{proof: multiple p-partite planted hypercliques problem} for details.
The closest work is by \citet{dadon2024detection}, who study multiple fixed-size Gaussian submatrix detection and derive statistical and low-degree computational strong detection thresholds.
Their statistical proof invokes negative association of the full block-overlap arrays, a property that need not hold; see \Cref{eq:overlap-positive-covariance}.
Their low-degree proof uses a bound on overlap moments whose justification in the cited preprint makes the same assumption.
Our proofs avoid this global assumption: the statistical lower bound uses conditional negative association followed by iterated conditioning, while the computational lower bound uses a separate moment comparison.
\end{remark}

\subsection{Thick cycles and paths}\label{section: thick cycle and path}
We next consider \emph{thick cycles} and \emph{thick paths}, which connect nearby vertices along an underlying cycle or path.
Specifically, two vertices are adjacent when their distance in the underlying ordering is at most $\Delta$.
For positive integers $\ell$ and $\Delta$ satisfying $2\Delta<\ell\le m$, define a length-$\ell$, $\Delta$-thick cycle by the edge set
\begin{align}\label{eq:think.cycle.definition}
	S:=\bigg\{(i,j)\in[m]_2: \dot A_i, \dot A_j\in(0,\ell], d^{(\ell)}_{cyc}(\dot A_{i},\dot A_{j})\leq \Delta\bigg\},
\end{align}
where $d^{(\ell)}_{cyc}(a,b):=\min\{|a-b|, \ell-|a-b|\}$ is the cyclic distance on a circle of circumference $\ell$, and $(\dot A_1,\ldots, \dot A_m)$ is a random permutation of $[m]$. We denote the distribution of $S$ by $\dot \nu$. A length-$\ell$, $\Delta$-thick path is defined analogously by replacing the cyclic distance with the linear distance $d^{(\ell)}_{lin}(a,b):=|a-b|$.
For random-size thick cycles and paths, we replace the coordinates of the random permutation in \Cref{eq:think.cycle.definition} by i.i.d.\ uniform random variables on $[0,m]$.
We write $\tilde \nu$ for the resulting distribution of $S$.

Related testing problems have been studied in \citet{carpentier2025low,mao2025information}, primarily under Erd\H{o}s--R\'enyi observations or with recovery rather than detection as the objective.
We begin by analyzing the statistical and computational thresholds for thick cycles.
\begin{theorem}\label{theorem: thick cycle problem}
	Consider the Gaussian graph testing problem with the latent structure $S$ defined above and prior distribution $\dot \nu$ or $\tilde \nu$. Suppose that $\Delta$ is of positive polynomial order in $m$. Then, the statistical and computational thresholds are 
	\begin{align*}
		\mu_{\mathrm{stat}}^\ast=\left\{
		\begin{array}{@{}c@{\quad}l}
		\sqrt{\frac{1}{\Delta}}, & \ell\sqrt{\Delta} \ll m, \\
		\frac{m}{\ell\Delta}, & \ell\sqrt{\Delta} \gg m
		\end{array}
		\right.
		\qquad \text{and} \qquad
		\mu_{\mathrm{comp}}^\ast=\left\{
		\begin{array}{@{}c@{\quad}l}
		\vphantom{\sqrt{\frac{1}{\Delta}}}1, & \ell\Delta \ll m, \\
		\frac{m}{\ell\Delta}, & \ell\Delta \gg m.
		\end{array}
		\right.
	\end{align*}
\end{theorem}

\begin{figure}[t]
	\includegraphics[width=\textwidth]{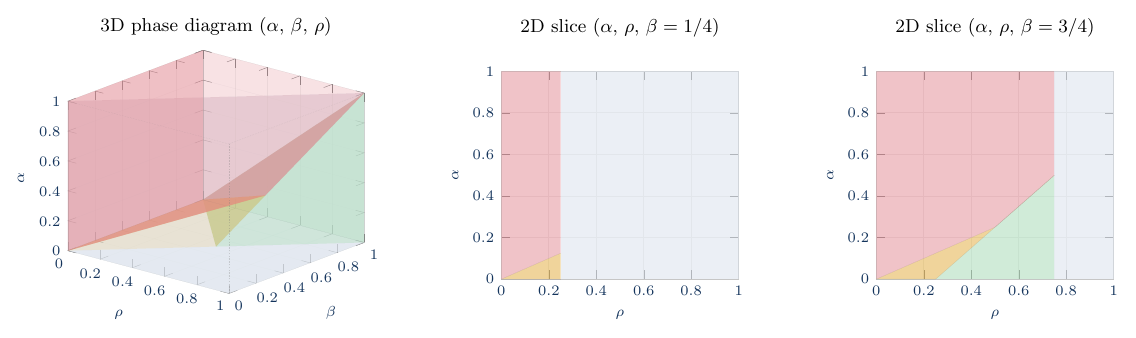}
	\caption{Phase diagram for the thick cycle problem under the parameterization $\mu=m^{-\alpha}$, $\ell=m^{\beta}$, and $\Delta=m^{\rho}$.
The left panel shows the three-dimensional phase diagram, and the two right panels show two-dimensional slices at $\beta=1/4$ and $\beta=3/4$.
Green, yellow, and pink denote the easy, hard, and impossible regions, respectively; blue-gray denotes parameter values outside the admissible range.}\label{fig:thick-cycle-phase-diagram}
\end{figure}

The phase diagram associated to the thresholds described in \Cref{theorem: thick cycle problem} is shown in  \Cref{fig:thick-cycle-phase-diagram}.

We briefly describe the proof strategy for the fixed-size model; the full proof is given in \Cref{proof: thick cycle problem}.
The random-size analysis then follows because the numbers of neighbors on either side of each vertex concentrate around $\Delta$ with high probability.
The upper bounds follow from the fact that a fixed-size thick cycle has $\ell\Delta$ edges and contains a clique on $\Delta+1$ consecutive vertices.
The average and outlier tests, together with a scan over all $(\Delta+1)$-vertex cliques, therefore yield the claimed upper bounds.

\begin{figure}[t]
    \centering
    \includegraphics[width=0.8\textwidth]{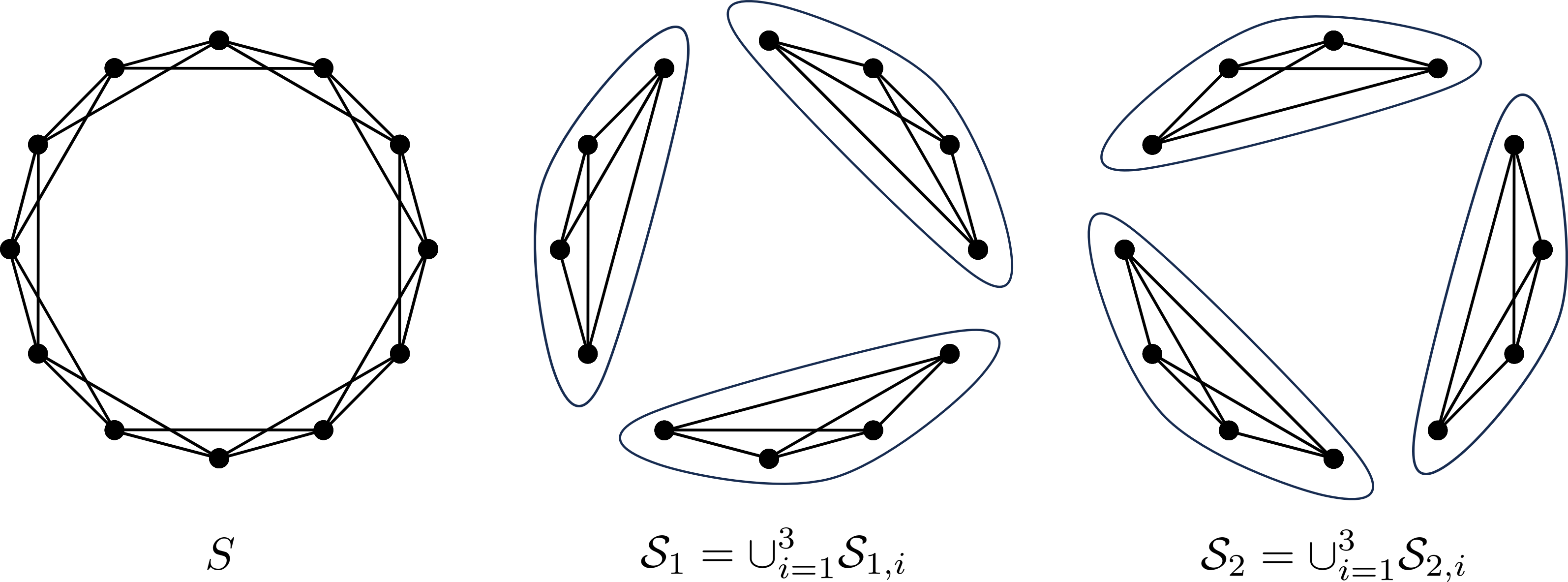}
    \caption{Illustration of the patchwork cover of a thick cycle. The left panel shows the thick cycle $S$ for $\ell=12$ and $\Delta=2$. The middle and right panels show the two patchwork components $\mathcal S_1$ and $\mathcal S_2$, respectively, whose union $\mathcal S=\mathcal S_1\cup\mathcal S_2$ covers $S$. Each patch $\mathcal S_{1,i}$ or $\mathcal S_{2,i}$ is a clique on $2\Delta$ vertices.}\label{fig:thick-cycle-patchwork}
\end{figure}

To establish the lower bounds, we control $|S\cap S'|$ using the overlap between patchwork covers of $S$ and $S'$.
For illustration, assume that $\ell$ is divisible by $2\Delta$. 
In the fixed-size model, set $q:=\ell/(2\Delta)$ and let $\mathcal S=\mathcal S_1\cup\mathcal S_2$ be a patchwork with two components
	\begin{align*}
	\mathcal S_1=\cup_{i=1}^q \mathcal S_{1,i} \quad\text{and}\quad \mathcal S_2=\cup_{i=1}^q \mathcal S_{2,i},
	\end{align*}
	each of which is a union of $\gamma=\ell/(2\Delta)$ disjoint cliques of size $k=2\Delta$; see \Cref{fig:thick-cycle-patchwork} for an illustration. Since $S$ is a subset of $\mathcal S$, the overlap can be upper bounded as
	\begin{equation}\label{equation: patchwork bound}
        \begin{aligned}
        |S\cap S'| &\leq |(\mathcal S_1\cup \mathcal S_2)\cap (\mathcal S'_1\cup \mathcal S'_2)|\\
        &\leq |\mathcal S_1\cap \mathcal S'_1| + |\mathcal S_1\cap \mathcal S'_2| + |\mathcal S_2\cap \mathcal S'_1| + |\mathcal S_2\cap \mathcal S'_2|,
        \end{aligned}
    \end{equation}
	where $\mathcal S'=\mathcal S'_1\cup\mathcal S'_2$ is the patchwork for $S'$ defined in the same way as $\mathcal S$. Each term on the right-hand side of \Cref{equation: patchwork bound} is the overlap between two independent copies of the multiple-clique model analyzed in \Cref{section: multiple planted cliques}.
It remains to combine the four contributions.
For the statistical lower bound, H\"older's inequality yields
	\begin{align*}
        \ee\bigg[\prod_{i,j=1}^2\exp\big(\mu^2 \big|\mathcal S_i\cap \mathcal S'_j\big|\big)\bigg]\leq \prod_{i,j=1}^2\Big(\ee\big[\exp\big(4\mu^2 \big|\mathcal S_i\cap \mathcal S'_j\big|\big)\big]\Big)^{1/4}=\ee\big[\exp\big(4\mu^2 \big|\mathcal S_1\cap \mathcal S'_1\big|\big)\big],
	\end{align*}
	where the last equality holds because $\mathcal S_i\cap \mathcal S'_j$, for $i,j\in\{1,2\}$, are identically distributed.
For the computational lower bound, the inequality $(\sum_{i=1}^4 x_i)^d \leq 4^{d-1} \sum_{i=1}^4 x_i^d$ for $x_i\geq 0$ gives
	\begin{align*}
		\ee\bigg[\sum^D_{d=1} \frac{\mu^{2d}}{d!}\bigg\{\sum_{i,j=1}^2 \big|\mathcal S_i\cap \mathcal S'_j\big|\bigg\}^d\bigg]\leq \ee\bigg[\sum^D_{d=1} \frac{4^{d-1}\mu^{2d}}{d!}\sum_{i,j=1}^2 \big|\mathcal S_i\cap \mathcal S'_j\big|^d\bigg]= \ee\bigg[\sum^D_{d=1} \frac{4^d\mu^{2d}}{d!} \big|\mathcal S_1\cap \mathcal S'_1\big|^d\bigg].
	\end{align*}
	Applying the bounds on overlap moments proved for the multiple-clique model (\Cref{theorem: multiple cliques problem}), with $\gamma=\ell/(2\Delta)$ and $k=2\Delta$,  yields the desired lower bounds.

We conclude by analyzing the statistical and computational thresholds for thick paths.
Write $\bfA=\dot\bfA$ in the fixed-size model and $\bfA=\tilde\bfA$ in the random-size model. A length-$\ell$, $\Delta$-thick path is defined by the edge set
\begin{align*}
	S:=\bigg\{(i,j)\in[m]_2: A_i,A_j\in(0,\ell], d^{(\ell)}_{lin}(A_{i},A_{j})\leq \Delta\bigg\},
\end{align*}
where $d^{(\ell)}_{lin}(a,b):=|a-b|$ is the linear distance on a line of length $\ell$. The upper- and lower-bound arguments for thick cycles extend to thick paths, with minor endpoint adjustments for the upper bounds, yielding the same statistical and computational thresholds.
The following corollary records the result; its proof is given in \Cref{proof: thick path problem}.
\begin{corollary}\label{corollary: thick path problem}
	Consider the Gaussian graph testing problem with the latent structure $S$ defined above and prior distribution $\dot \nu$ or $\tilde \nu$. Suppose that $\Delta$ is of positive polynomial order in $m$. Then, in both cases, the statistical and computational thresholds $\mu_{\mathrm{stat}}^\ast$ and $\mu_{\mathrm{comp}}^\ast$ are identical to those in \Cref{theorem: thick cycle problem}.
\end{corollary}

\section{Equivalence across natural exponential family models}\label{sec:nef-equivalence}
In this section, we extend the Gaussian framework for graph testing to natural exponential families (NEFs).
Under appropriate conditions, a canonical rescaling of the signal strength transfers statistical and computational thresholds from one NEF model to another.

We consider a NEF model with a nondegenerate reference measure $\rho$ and exponentially tilted distribution $P_\theta$ with natural parameter $\theta$.  The density of $P_\theta$ with respect to $\rho$ is given by
$$
dP_\theta(y) = \exp(\theta y - \psi(\theta))\rho(dy), \quad \theta \in \Theta
$$
Here, $\psi(\theta)=\log\int\exp(\theta y)\rho(dy)$ is the log-partition function, and $\Theta = \mathrm{int}\{ \theta \in \mathbb{R} \colon \psi(\theta) < \infty\}$ is the natural parameter space.
Representative NEF models and their natural parameterizations are listed in \Cref{tab:nef-chi-square}.
We consider an $n$-dimensional observation vector $\bfY=(Y_1,\ldots,Y_n)$ and a latent structure $S\sim\nu_n$, where $\nu_n$ is supported on a class $\mathcal C\subseteq 2^{[n]}$; the elements of $\mathcal C$ need not correspond to edge sets of subgraphs of $K_m$. Given a triplet $(\theta_0,\theta_1,\rho)$ with $\theta_0,\theta_1\in\Theta$, we consider the testing problem:
\begin{align}\label{equation:nef-hypothesis-testing}
H_0: Y_l \stackrel{iid}{\sim} P_{\theta_0}\quad \text{versus} \quad H_1: Y_l\,|\,S
 \stackrel{ind}{\sim}
\begin{cases}
P_{\theta_1}, & \text{if } l \in S, \\
P_{\theta_0}, & \text{if } l \notin S.
\end{cases}
\end{align}
The latent structure may have fixed or random size; fixed-size assumptions are stated explicitly where needed for the upper bounds.
The exponential tilts $\theta_0$ and $\theta_1$, corresponding to the null and alternative hypotheses, may depend on $n$.
In the Gaussian graph testing problem of \Cref{equation: gaussian graph testing problem}, the signal strength $\mu$ is the mean difference between $P_{\theta_0}=\mathrm N(0,1)$ and $P_{\theta_1}=\mathrm N(\mu,1)$.
For the NEF model in \Cref{equation:nef-hypothesis-testing}, the analogous quantity is the \emph{effective signal strength}
\begin{equation}\label{eq:eff.signal.strength}
\tau_n := \frac{\psi'(\theta_1)-\psi'(\theta_0)}{\sqrt{\psi''(\theta_0)}},
\end{equation}
which is the standardized mean difference between $P_{\theta_0}$ and $P_{\theta_1}$ because $\psi'(\theta)$ and $\psi''(\theta)$ are the mean and variance of $P_\theta$, respectively.
The following subsections develop lower bounds through the one-coordinate chi-square divergence and upper bounds in terms of the effective signal strength $\tau_n$.

\subsection{Lower bounds}\label{section: lower bound for NEF models}
Let $\chi^2(P_{\theta_1}\|P_{\theta_0})=\int(\frac{dP_{\theta_1}}{dP_{\theta_0}}-1)^2 dP_{\theta_0}$ denote the chi-square divergence between $P_{\theta_0}$ and $P_{\theta_1}$. In the following theorem, we show that the asymptotic behavior of the full and low-degree likelihood-ratio second moments for NEF models can be characterized by this one-coordinate divergence, without further dependence on the specific form of the distributions under the stated degree condition. This result allows us to transfer statistical and computational lower bounds for the Gaussian graph testing problem to other NEF models by matching the chi-square divergence between the two models. The proof is given in \Cref{proof: equivalence between NEF models}.
\begin{theorem}[$\chi^2$-characterization of lower bounds for NEF models]\label{theorem: equivalence between NEF models}
Consider the testing problem in \Cref{equation:nef-hypothesis-testing} with $\theta_0,\theta_1\in \Theta$ satisfying $\theta_0\neq\theta_1$ and $\psi''(\theta_0)>0$. Assume that $\psi(2\theta_1-\theta_0)<\infty$,
so that the chi-square divergence
$\chi_n^2:=\chi^2(P_{\theta_1}\|P_{\theta_0})$
is finite. Then, the second moment of the likelihood ratio $L_{\nu_n}(\bfY)$ is given by
\begin{align*}
	\bbE_{\bfY\sim \bbQ_n}\big[L_{\nu_n}(\bfY)^2\big]
	=
	\bbE_{S,S'\sim\nu_n}\Big[\big(1+\chi_n^2\big)^{|S\cap S'|}\Big].
\end{align*}
Moreover, for any sequence of degrees $D=D_n$ satisfying $D_n\log(\chi^2_n/\tau_{n}^2)=O(1)$, the second moment of the low-degree likelihood ratio $L_{\nu_n}^{\leq D}(\bfY)$ satisfies the asymptotic equivalence
\begin{align*}
	\bbE_{\bfY\sim \bbQ_n}\big[L_{\nu_n}^{\leq D}(\bfY)^2\big]
	\asymp
	\bbE_{S,S'\sim\nu_n}\Bigg[\sum^{D\wedge|S\cap S'|}_{d=0}\binom{|S\cap S'|}{d}(\chi_n^2)^d \Bigg].
\end{align*}
The implicit constants depend only on the upper bound for $D_n\log(\chi^2_n/\tau_{n}^2)$.
\end{theorem}

The theorem implies that, for two different NEF models defined by $(\theta^A_{0,n}, \theta^A_{1,n}, \rho^A)$ and $(\theta^B_{0,n}, \theta^B_{1,n}, \rho^B)$ with the same latent prior $\nu_n$ and matching chi-square divergences $\chi_n^2$, the second moments of their likelihood ratios coincide exactly.
Furthermore, if both models satisfy $D_n\log(\chi_n^2/\tau_n^2)=O(1)$, their low-degree likelihood-ratio second moments are equivalent up to constant factors.
Consequently, the resulting statistical and computational lower bounds obtained from second moment methods (\Cref{lemma: second moment method,lemma: low-degree second moment method}) for one model can be transferred to the other model by matching the chi-square divergences between the two models.
The degree condition $D_n\log(\chi^2_n/\tau_{n}^2)=O(1)$ can be parsed explicitly for specific NEF models.
For example, in a Gaussian location model with mean difference $\mu$ and known variance $\sigma^2$, we have $\chi^2_n=\exp(\mu^2/\sigma^2)-1$ and $\tau_n=\mu/\sigma$.
Since $\log((e^x-1)/x)\sim x/2$ for $x\to 0$, the degree condition is equivalent to $D_n\mu^2=O(1)$ for fixed $\sigma^2>0$.
In the Bernoulli model with success probabilities $p_0$ and $p_1$, we have $\chi_n^2=(p_1-p_0)^2/{p_0(1-p_0)}=\tau_n^2$, so the degree condition is always satisfied.
More generally, for any NEF model in the local alternative regime where $\theta_0$ is fixed and $\theta_1\to\theta_0$, a simple sufficient condition is $D_n(\theta_1-\theta_0)^2=O(1)$.
We formalize this observation in \Cref{proposition: local alternative regime}.
\begin{remark}[Compatibility with conditional second moment methods]\label{remark: compatibility with conditional second moment methods}
An astute reader may wonder whether the lower bounds
for the random-size problems in \Cref{sec:results}
can be transferred through \Cref{theorem: equivalence between NEF models} directly, since their proofs rely
not merely on the standard second moment method,
but on the conditional second moment method
(\Cref{rem:conditional-second-moment-method})
involving a truncation argument. Indeed, such conditional arguments, particularly those based on high-probability events involving the latent variables that define $S$, are compatible with \Cref{theorem: equivalence between NEF models},
since the $\chi^2$-characterization continues to hold
for the second moments under the truncated prior.
More specifically, given a prior $\nu$, let $\nu_\calE(\cdot):=\nu(\cdot\,|\,\calE)$ be a truncated prior defined by a high-probability event $\calE$ for $S$ under $\nu$.
Then, \Cref{theorem: equivalence between NEF models} ensures that the second moment of the likelihood ratio between the null distribution $\bbQ$ and the alternative distribution $\bbP_{\nu_\calE}:=\bbE_{S\sim\nu_\calE}[\bbP_S]$ is given by $\bbE_{S,S'\sim\nu_\calE}\big[\big(1+\chi_n^2\big)^{|S\cap S'|}\big]$, simply by replacing $\nu$ with $\nu_\calE$.
Since this expectation involves only the overlap structure under $\nu_\calE$ and the chi-square divergence $\chi_n^2$, the full likelihood-ratio second moments coincide across NEF models with the same conditional prior $\nu_\calE$ and matching chi-square divergence $\chi_n^2$.
The same argument yields equivalence up to constant factors for the low-degree second moments, and the degree condition $D_n\log(\chi^2_n/\tau_{n}^2)=O(1)$ is independent of the truncation event $\calE$.
Hence, under the degree condition, the lower bounds for the random-size problems in \Cref{sec:results} can also be transferred across NEF models by matching $\chi^2$.

 On the other hand, when the conditional second moment method involves conditioning on a high-probability event defined in terms of the observations $\bfY$ (see e.g., \citet{arias2014community,dhawan2025detection}) rather than the latent structure $S$, the conditioning event may depend on model-specific features of the observation distribution. Consequently, transferring the resulting bounds through \Cref{theorem: equivalence between NEF models} may not be applicable.
\end{remark}
One caveat is that the lower bounds obtained from \Cref{theorem: equivalence between NEF models} are not necessarily sharp across different NEF models unless matching upper bounds are available for each model; the full and low-degree (conditional) second moment methods provide only sufficient conditions for ruling out strong detection and low-degree strong separation, respectively.
Sharper lower bounds may therefore be obtainable through more refined arguments, such as the observation-dependent conditional second moment methods mentioned in \Cref{remark: compatibility with conditional second moment methods}.
A concrete example of this is the planted clique problem under the Erd\H{o}s--R\'enyi model, where the parameters $p_0$ and $p_1$ both go to zero as $n\to\infty$ at different rates.
In this problem, transferring the corresponding result for the Gaussian model via chi-square matching only yields loose lower bounds.
We will revisit this example in \Cref{section: equivalence on planted clique}.

To establish equivalence of the statistical and computational thresholds, we must also obtain matching upper bounds for the NEF models.
In the next subsection, we will show that the performance of the tests in \Cref{proposition: gaussian test upper bound} for general NEF models can be characterized in terms of the effective signal strength $\tau_n$ under appropriate regularity conditions.
Importantly, when $D_n\to\infty$, the degree condition $D_n\log(\chi_n^2/\tau_n^2)=O(1)$ implies $\chi_n^2/\tau_n^2\to1$.
Therefore, for two NEF models with matching chi-square divergences that both satisfy the degree condition, the magnitudes of their effective signal strengths are asymptotically equivalent, and the resulting upper bounds for the two models are of the same order.
Consequently, given $S\sim\nu$, if the lower bounds from second moment methods and upper bounds achieved by such tests for one model are sharp in the sense that they match up to a constant factor, then the transferred lower bounds and the corresponding upper bounds also match up to a constant factor under the degree condition and regularity assumptions on the upper bound.
We discuss when these assumptions hold in \Cref{sec: examples}.

\begin{table}[t]
\centering
\caption{Representative regular natural exponential families under convenient
reference measures $\rho$. The reference measure specifies the expression of $\theta$ in terms of the parameters of the family and is not unique.
For the Gaussian, Binomial, negative binomial, Gamma, and inverse Gaussian families, respectively, $\sigma^2$, $N$, $r$, $\alpha$, and $\lambda$ are held fixed. Here, $\beta$ denotes the scale parameter for the Gamma and exponential
families. For finite chi-square divergences, $p_1^2<p_0$ is required for the negative binomial family, $2\beta_0>\beta_1$ for the Gamma and exponential families, and $2\mu_0^2\geq\mu_1^2$ for the inverse Gaussian family.}
\label{tab:nef-chi-square}
\renewcommand{\arraystretch}{1.3}
\resizebox{\linewidth}{!}{
\begin{tabular}{cccccc}
\toprule
\textbf{Name} & \textbf{Reference measure} $\boldsymbol{\rho}$ & \textbf{Support} & \textbf{Exponential tilt} & $\boldsymbol{\Theta}$ & $\boldsymbol{\chi^2(P_{\theta_1}\|P_{\theta_0})}$ \\
\midrule

Gaussian & $\displaystyle \frac{1}{\sqrt{2\pi\sigma^2}}\exp\!\left(-\frac{x^2}{2\sigma^2}\right)dx$ & $x\in\mathbb R$ & $\displaystyle \theta=\frac{\mu}{\sigma^2}$ & $\mathbb R$ & $\displaystyle \exp\!\left(\frac{(\mu_1-\mu_0)^2}{\sigma^2}\right)-1$ \\[6pt]

Bernoulli & $\displaystyle \delta_0+\delta_1$ & $x\in\{0,1\}$ & $\displaystyle \theta=\log\frac{p}{1-p}$ & $\mathbb R$ & $\displaystyle \frac{(p_1-p_0)^2}{p_0(1-p_0)}$ \\[6pt]

Binomial & $\displaystyle \sum_{x=0}^{N}\binom{N}{x}\delta_x$ & $x\in\{0,\ldots,N\}$ & $\displaystyle \theta=\log\frac{p}{1-p}$ & $\mathbb R$ & $\displaystyle \left(1+\frac{(p_1-p_0)^2}{p_0(1-p_0)}\right)^N-1$ \\[6pt]

Poisson & $\displaystyle \sum_{x=0}^{\infty}\frac{1}{x!}\delta_x$ & $x\in\mathbb N_0$ & $\displaystyle \theta=\log\lambda$ & $\mathbb R$ & $\displaystyle \exp\!\left(\frac{(\lambda_1-\lambda_0)^2}{\lambda_0}\right)-1$ \\[6pt]

Negative binomial & $\displaystyle \sum_{x=0}^{\infty}\binom{x+r-1}{x}\delta_x$ & $x\in\mathbb N_0$ & $\displaystyle \theta=\log p$ & $(-\infty,0)$ & $\displaystyle \left[\frac{p_0(1-p_1)^2}{(1-p_0)(p_0-p_1^2)}\right]^r-1$ \\[6pt]

Gamma & $\displaystyle \frac{x^{\alpha-1}}{\Gamma(\alpha)}\,dx$ & $x>0$ & $\displaystyle \theta=-\frac{1}{\beta}$ & $(-\infty,0)$ & $\displaystyle \left[ \frac{\beta_0^2}{\beta_1(2\beta_0-\beta_1)}\right]^\alpha-1$ \\[6pt]

Exponential & $\displaystyle dx$ & $x>0$ & $\displaystyle \theta=-\frac{1}{\beta}$ & $(-\infty,0)$ & $\displaystyle \frac{\beta_0^2}{\beta_1(2\beta_0-\beta_1)}-1$ \\[6pt]

Inverse Gaussian & $\displaystyle \sqrt{\frac{\lambda}{2\pi x^3}}\exp\!\left(-\frac{\lambda}{2x}\right)dx$ & $x>0$ & $\displaystyle \theta=-\frac{\lambda}{2\mu^2}$ & $(-\infty,0)$ & $\displaystyle \exp\!\left\{\lambda\left(\frac{2}{\mu_1}-\frac{1}{\mu_0}-\sqrt{\frac{2}{\mu_1^2}-\frac{1}{\mu_0^2}}\right)\right\}-1$ \\

\bottomrule
\end{tabular}
}
\end{table}

\subsection{Upper bounds}
We now revisit the average and scan tests from \Cref{proposition: gaussian test upper bound} and analyze their performance under NEF models in terms of the effective signal strength.
The corresponding test statistics are
\begin{equation}\label{equation: test statistics}
	\begin{aligned}
		T_{\mathrm{avg}}(\bfy):=\sum^n_{i=1}y_i, \qquad \text{and} \qquad T_{\mathrm{scan}}(\bfy):=\max_{s\in \mathcal C}\sum_{l\in s} y_l,
	\end{aligned}
\end{equation}
respectively.
The risks of these tests under the Gaussian model are controlled using standard Gaussian tail bounds.
For a regular NEF, the centered moment generating function is finite in a neighborhood of the origin, so the observations are sub-exponential; more specifically, the following identity yields the local quadratic bound needed for our tests:
\begin{align*}
\log\big(\ee_{Y\sim P_\theta}\big[\exp\big\{t[Y-\psi'(\theta)]\big\}\big]\big)&=\psi(\theta+t)-\psi(\theta)-t\psi'(\theta)\\
&=t^2\int_0^1 (1-s)\psi''(\theta+st)ds,
\end{align*}
where the last equality follows from Taylor's expansion with integral remainder.
Therefore, if $\psi''$ is uniformly bounded near $\theta$, the centered log moment generating function has a quadratic upper bound on a local interval; equivalently, the random variable satisfies a local sub-Gaussian concentration bound while remaining sub-exponential in general.
This motivates a local regularity condition requiring the curvature of $\psi$ near $\theta_0$ and $\theta_1$ to be bounded by a constant multiple of $v_0:=\psi''(\theta_0)$.
\begin{definition}[Local regularity of the NEF models]\label{def: local regularity of the NEF models}
	Let $r_n>0$. We say that the NEF models $P_{\theta_0}$ and $P_{\theta_1}$ are locally regular at scale $r_n$ if $[\theta_j-r_n,\theta_j+r_n]\subseteq\Theta$ for each $j\in\{0,1\}$ and there exists a constant $C\geq1$ such that
		$$ \sup_{\substack{j\in\{0,1\}\\|u|\leq r_n}}\frac{\psi''(\theta_j+u)}{\psi''(\theta_0)}\leq C,$$
		for all sufficiently large $n$.
\end{definition}
Under local regularity, each coordinate $Y_l$, conditional on the latent structure, is locally sub-Gaussian with the common variance proxy $Cv_0$, whether its distribution is $P_{\theta_0}$ or $P_{\theta_1}$.
Consequently, the average and scan tests attain risk guarantees analogous to those obtained under the Gaussian model.
When the latent structure $S$ is a subgraph of the complete graph $K_m$ and each observation is indexed by an edge, we may also use the degree test introduced in \Cref{remark: gaussian degree test}.
The weighted degree of vertex $i$ and the maximum-degree statistic are, respectively,
$$
T_i^{\mathrm{deg}}(\bfy):=\sum_{\{e\in [m]_2:i\in e\}}y_e \qquad \text{and} \qquad T_{\mathrm{deg}}(\bfy):=\max_{i\in[m]}T_i^{\mathrm{deg}}(\bfy)
$$
Conditional on the latent structure, each weighted degree is a sum of independent observations, so the same concentration argument gives an analogous guarantee under NEF models.
The next proposition collects these bounds; its proof is given in \Cref{proof: upper bound for NEF models}.

\begin{proposition}[Upper bound for NEF models]\label{proposition: upper bound for NEF models}
	Consider the test statistics defined in \Cref{equation: test statistics}. Let $\delta=\delta_n\in(0,1)$ be a target risk, and let $d_\theta:=\psi'(\theta_1)-\psi'(\theta_0)$, $v_0:=\psi''(\theta_0)$, so that the effective signal strength is given by $\tau_{n}:=d_\theta/\sqrt{v_0}$. Suppose that every $s\in\mathcal C$ satisfies $|s|=K\geq1$, and let $\nu$ be a prior distribution on $\mathcal C$. Then,
	\begin{itemize}
		\item Assume that the NEF models $P_{\theta_0}$ and $P_{\theta_1}$ are locally regular at scale $r_{\textrm{avg}}:=\sqrt{2(nv_0)^{-1}\log(2/\delta)}$ with a constant $C_{\textrm{avg}}\geq 1$.
		For the threshold $t_{\mathrm{avg}}:=n\psi'(\theta_0)+Kd_\theta/2$, the average test $\Delta_{\mathrm{avg}}(\bfy):=\mathds1\{T_{\mathrm{avg}}(\bfy)>t_{\mathrm{avg}}\}$ satisfies $R_\nu(\Delta_{\mathrm{avg}})\leq\delta$ if
		$$
		\tau_{n}\geq\sqrt{\frac{8C_{\textrm{avg}}n}{K^2}\log\frac{2}{\delta}}.
		$$
		\item Assume that the NEF models $P_{\theta_0}$ and $P_{\theta_1}$ are locally regular at scale $r_{\textrm{scan}}:=\sqrt{2(Kv_0)^{-1}\log(2|\mathcal C|/\delta)}$ with a constant $C_{\textrm{scan}}\geq 1$.
		For the threshold $t_{\mathrm{scan}}:=K\psi'(\theta_0)+(Kd_\theta+\sqrt{2KC_{\textrm{scan}}v_0\log|\mathcal C|})/2$, the scan test $\Delta_{\mathrm{scan}}(\bfy):=\mathds1\{T_{\mathrm{scan}}(\bfy)>t_{\mathrm{scan}}\}$ satisfies $R_\nu(\Delta_{\mathrm{scan}})\leq\delta$ if
		$$
		\tau_{n}\geq\sqrt{\frac{2C_{\textrm{scan}}}{K}}\bigg(2\sqrt{\log\frac2\delta}+\sqrt{\log|\mathcal C|}\bigg).
		$$
	\end{itemize}
	Furthermore, suppose that the latent structure $S$ is a subgraph of the complete graph $K_m$, so that $n=\binom{m}{2}$ and each observation is indexed by an edge $e\in[m]_2$ connecting two vertices.
	Let $T_{\mathrm{deg}}(\bfy):=\max_{i\in[m]}T_i^{\mathrm{deg}}(\bfy)$, and suppose that the maximum number of signal edges incident to any vertex satisfies $\max_{i\in[m]}|\{e\in S:i\in e\}|=\mathfrak d_{\max}\geq1$ almost surely under $S\sim\nu$.
	Then,
	\begin{itemize}
		\item Assume that the NEF models $P_{\theta_0}$ and $P_{\theta_1}$ are locally regular with constant $C_{\mathrm{deg}}\geq1$ at scale
		\[
		r_{\mathrm{deg}}:=\sqrt{2\{(m-1)v_0\}^{-1}\log(2m/\delta)}.
		\]
		For the threshold $t_{\mathrm{deg}}:=(m-1)\psi'(\theta_0)+\mathfrak d_{\max}d_\theta/2$, the degree test $\Delta_{\mathrm{deg}}(\bfy):=\mathds1\{T_{\mathrm{deg}}(\bfy)>t_{\mathrm{deg}}\}$ satisfies $R_\nu(\Delta_{\mathrm{deg}})\leq\delta$ if
		\[
		\tau_n\geq\sqrt{\frac{8C_{\mathrm{deg}}(m-1)}{\mathfrak d_{\max}^2}\log\frac{2m}{\delta}}.
		\]
	\end{itemize}
\end{proposition}

Applying the proposition generally requires verifying local regularity for the particular NEF model.
For the Gaussian model with fixed variance, the condition holds globally with $C=1$ because $\psi''$ is constant.
For the Bernoulli model, if $p_0<p_1$ and $p_1/p_0=O(1)$, local regularity holds with some $C>1$ for any sequence $r_n=o(1)$.
If instead $p_1/p_0\to\infty$, local regularity may fail because $\psi''(\theta_1)/\psi''(\theta_0)=p_1(1-p_1)/\{p_0(1-p_0)\}$ may diverge to infinity.
Without local regularity, sums of observations are not guaranteed to satisfy the Gaussian-type tail bounds used above, so a different analysis may be needed.

	\begin{remark}[Generalized scan test]\label{remark: generalized scan test}
		The scan test in \Cref{proposition: upper bound for NEF models} can be extended to random-size alternatives and scans targeting substructures of the planted support. Specifically, suppose that there exists a class $\mathcal C_\star$ of subsets of $[n]$, each of cardinality $K_\star$, such that
		$$
		\Pr_{S\sim\nu_n}\left(\exists s_\star\in\mathcal C_\star:s_\star\subseteq S\right)\to1.
		$$
		Then, under the corresponding local regularity condition, the scan-test bound continues to hold with $(\mathcal C,K)$ replaced by $(\mathcal C_\star,K_\star)$, with an additional $o(1)$ term in the risk bound. Several scan tests considered in this paper are instances
		of this generalization. For example, in the random-size $k$-vertex clique problem, the planted support contains a clique on $ck$ vertices with high probability for any fixed $c\in(0,1)$, whereas in the fixed-size $\Delta$-thick cycle problem, the planted support always contains a clique on $\Delta+1$ vertices. In both cases, the scan test targets these cliques rather than the entire planted support; see \Cref{proof: random-size clique problem,proof: thick cycle problem}. Similarly, the average and degree tests can be extended to random-size alternatives using high-probability lower bounds on the support size and maximum degree, respectively.
	\end{remark}

We next identify a regime in which local regularity holds automatically and need not be checked separately for each model.

\subsection{Universality in the local alternative regime}\label{sec: examples}

We now identify a natural regime in which statistical and computational thresholds established by matching the preceding lower and upper bounds for one NEF model remain valid across other NEF models.
Let
$$
\theta_0=c \quad \text{and} \qquad \theta_1=c+h_n,
$$
where $c$ is a fixed point in the interior of the natural parameter space $\Theta$ and $\{h_n\}$ is a vanishing positive sequence.
We use $h_n>0$ to keep the one-sided tests and the effective signal strength $\tau_n$ oriented consistently; the case $h_n<0$ is analogous after reversing the tests and replacing $\tau_n$ by $|\tau_n|$.
We refer to this setting as the \textit{local alternative regime}.
Here, ``local'' means that the one-coordinate alternative $P_{\theta_1}$ approaches the null $P_{\theta_0}$; it does not prescribe a particular contiguity scale for the full structured testing problem.
For an $n$-fold product experiment in which every coordinate is shifted, differentiability in quadratic mean yields mutual contiguity of the sequences of product laws $P_{\theta_0}^{\otimes n}$ and $P_{\theta_1}^{\otimes n}$ when $h_n=O(n^{-1/2})$ \citep[Chapters~6--7]{vandervaart1998asymptotic},  whereas strong detection
requires $h_n=\omega(n^{-1/2})$.
In our setting, however, only the coordinates in the random structure $S$ are shifted, so the relevant detection scale depends on the size and geometry of $S$ and on its prior.
Nevertheless, convexity of the relative entropy, conditional independence given $S$, and $|S|\leq n$ give
\[
D_{\mathrm{KL}}(\bbP_{\nu_n}\,\|\,\bbQ_n)
\leq n D_{\mathrm{KL}}(P_{c+h_n}\,\|\,P_c)
=O(nh_n^2).
\]
The last equality follows from a second-order expansion at the fixed interior parameter $c$, and \Cref{equation: risk lower bound via BH} shows that any bounded subsequence of $nh_n^2$ prevents the risk from tending to zero.
Although $h_n=\omega(n^{-1/2})$ is still necessary for strong detection, this condition alone does not guarantee distinguishability in our structured testing problems. We therefore define the local alternative regime simply through the vanishing condition $h_n\to 0$ and characterize distinguishability scales through the statistical and computational thresholds derived below.

We define $\tau_{\mathrm{stat}}^\ast$ and $\tau_{\mathrm{comp}}^\ast$, analogously to $\mu_{\mathrm{stat}}^\ast$ and $\mu_{\mathrm{comp}}^\ast$ in \Cref{definition: statistical and computational thresholds}, as the statistical and computational thresholds in the local alternative regime.
Throughout this subsection, we restrict attention to polynomially vanishing effective signal strength sequences $\tau_n$, by which we mean sequences satisfying $\limsup_{n\to\infty}\log_n\tau_n<0$.

\begin{definition}[Statistical and computational thresholds in the local alternative regime]
\label{definition: statistical and computational thresholds in NEF models}

Consider a sequence of testing problems defined in \Cref{equation:nef-hypothesis-testing} in the local alternative regime, and let $\tau_n$ denote the corresponding effective signal strength. A positive sequence $\tau_{\mathrm{stat}}^\ast$ is called a statistical threshold in the local alternative regime if the following holds for every effective signal strength sequence $\tau_n$:
\begin{itemize}
\item if $\tau_n\gg\tau_{\mathrm{stat}}^\ast$, then there exists a sequence of tests $\Delta_n$ that achieves strong detection;
\item if $\tau_n\ll\tau_{\mathrm{stat}}^\ast$, then no sequence of tests can achieve strong detection.
\end{itemize}
Similarly, a positive sequence $\tau_{\mathrm{comp}}^\ast$ is called a computational threshold in the local alternative regime if
\begin{itemize}
\item if $\tau_n\gg\tau_{\mathrm{comp}}^\ast$, then there exists a sequence of polynomial-time tests $\Delta_n$ that achieves strong detection;
\item if $\tau_n\ll\tau_{\mathrm{comp}}^\ast$, then, under \Cref{conjecture: low-degree polynomial conjecture}, no sequence of polynomial-time tests can achieve strong detection.
\end{itemize}
If strong detection is impossible for every polynomially vanishing effective signal sequence in the local alternative regime using arbitrary tests (respectively, polynomial-time tests), we adopt the convention $\tau_{\mathrm{stat}}^\ast\asymp1$ (respectively, $\tau_{\mathrm{comp}}^\ast\asymp1$), for which the achievability direction is vacuously satisfied.
\end{definition}

We now state the main result of this subsection.
\begin{theorem}[Universality of tightness in the local alternative regime via chi-square matching]\label{theorem: universality of tightness in the local alternative regime via chi-square matching}
	For $M\in\{A,B\}$, fix a nondegenerate NEF $\{P_\theta^M:\theta\in\Theta_M\}$ with reference measure $\rho_M$ and a null parameter $c_M\in\mathrm{int}(\Theta_M)$.
	Consider the corresponding families of testing problems in \Cref{equation:nef-hypothesis-testing}, with
	$$
	\theta_{0}^M=c_M,\quad \theta_{1}^M=c_M+h_n^M,\quad h_n^M>0,\quad h_n^M\to 0,
	$$
	where $h_n^M$ ranges over all positive sequences tending to zero for which the effective signal strength $\tau_n^M$ vanishes polynomially. Also suppose that every $s \in \mathcal{C}$ satisfies $|s| = K$.
	Both families use the same sequence of latent priors $\nu_n$, fixed independently of the signal strength.
	Suppose that, for model $A$, the lower bounds obtained from the likelihood-ratio and degree-$D_n$ likelihood-ratio second moments match the upper bounds in \Cref{proposition: upper bound for NEF models}, achieved by the average and scan tests, together with the degree test when applicable, to establish statistical and computational thresholds $\tau^{A,\ast}_{\textrm{stat}}$ and $\tau^{A,\ast}_{\textrm{comp}}$ in the sense of \Cref{definition: statistical and computational thresholds in NEF models}.
	If the scan test provides the upper bound, also assume that $\log|\mathcal C|=o(K)$.
	Then $\tau^{A,\ast}_{\textrm{stat}}$ and $\tau^{A,\ast}_{\textrm{comp}}$ are also valid statistical and computational thresholds, respectively, for model $B$, in the sense of \Cref{definition: statistical and computational thresholds in NEF models}.
\end{theorem}

Although \Cref{theorem: universality of tightness in the local alternative regime via chi-square matching} is stated for fixed-size alternatives and relies on the upper bounds in \Cref{proposition: upper bound for NEF models}, the same transfer argument extends to random-size alternatives
through the generalized tests discussed in \Cref{remark: generalized scan test}. It also accommodates upper bounds obtained from substructure scans, provided that the corresponding upper and lower bounds match.

We prove \Cref{theorem: universality of tightness in the local alternative regime via chi-square matching} by simplifying the low-degree condition and verifying that local regularity holds automatically in the local alternative regime.

\paragraph{Lower bounds in the local alternative regime.}
The key to simplifying the lower-bound condition $D_n\log(\chi_n^2/\tau_n^2)=O(1)$ is the relation $\chi_n^2/\tau_n^2=1+O(h_n^2)$.
Indeed, the cumulant generating function $\log\ee_{Y\sim P_{\theta_0}}[\exp(tY)]=\psi(\theta_0+t)-\psi(\theta_0)$ gives
\begin{align*}
1+\chi^2(\theta_0,\theta_1)&=\ee_{Y\sim P_{\theta_0}}\big[\exp\big\{2(\theta_1-\theta_0)Y-2(\psi(\theta_1)-\psi(\theta_0))\big\}\big]\\
&=\exp\big\{\psi(\theta_0+2h_n)-2\psi(\theta_0+h_n)+\psi(\theta_0)\big\},
\end{align*}
while the effective signal strength is $\tau_n:=(\psi'(\theta_1)-\psi'(\theta_0))/\psi''(\theta_0)^{1/2}$.
Taylor expansions of $\psi$ and its derivatives therefore yield
$$
\frac{\chi_n^2}{\tau_n^2}=\frac{v_0h_n^2+c_3h_n^3+O(h_n^4)}{v_0h_n^2+c_3h_n^3+O(h_n^4)}=1+O(h_n^2),
$$
where $v_0:=\psi''(\theta_0)$ and $c_3:=\psi^{(3)}(\theta_0)$.
Combining this relation with $\log(1+x)=O(x)$ as $x\to0$ proves the following proposition; further details appear in \Cref{proof: local alternative regime}.

\begin{proposition}[Condition for the local alternative regime]\label{proposition: local alternative regime}
Let $\theta_0=c$ for some constant $c\in \Theta$, and let $\theta_1-\theta_0=:h_n\to 0$. Then, the condition $D_n\log(\chi^2_n/\tau_{n}^2)=O(1)$ is satisfied as long as $D_nh_n^2=O(1)$, or equivalently, $D_n\chi_n^2=O(1)$ or $D_n\tau_n^2=O(1)$.
\end{proposition}

For two NEF models $A$ and $B$ with matching chi-square divergences, it is therefore enough to verify $D_n\chi_n^2=O(1)$.
By \Cref{proposition: local alternative regime}, this is equivalent to verifying any one of $D_n\tau_{A,n}^2=O(1)$, $D_n\tau_{B,n}^2=O(1)$, $D_nh_{A,n}^2=O(1)$, or $D_nh_{B,n}^2=O(1)$.
Since we consider polylogarithmic degrees $D_n$, the condition $D_n\tau_n^2=O(1)$ is automatically satisfied under the assumption that the effective signal strength $\tau_n$ vanishes polynomially with $n$.

\paragraph{Upper bounds in the local alternative regime.} A second useful feature of the local alternative regime is that every NEF model automatically satisfies the local regularity condition in \Cref{def: local regularity of the NEF models}.
The following proposition formalizes this observation and allows us to apply the upper bounds in \Cref{proposition: upper bound for NEF models}.

\begin{proposition}[Local regularity in the local alternative regime]\label{proposition: upper bound for NEF models in local alternative regime}
	Suppose that $\theta_0\in\Theta$ is fixed and $\theta_1\to\theta_0$. Then, for any sequence $r_n>0$ satisfying $r_n\to 0$, 
	$$
	\limsup_{n\to\infty}\sup_{\substack{j\in\{0,1\}\\|u|\leq r_n}}\frac{\psi''(\theta_j+u)}{\psi''(\theta_0)}= 1.
	$$
	In particular, for any fixed $C>1$, the NEF models $P_{\theta_0}$ and $P_{\theta_1}$ are locally regular at scale $r_n$ with constant $C$ for sufficiently large $n$.
\end{proposition}
For fixed $\theta_0\in\Theta$, continuity of $\psi''$ and the limits $\theta_1\to\theta_0$ and $r_n\to0$ prove the proposition.
For the tests in \Cref{proposition: upper bound for NEF models}, local regularity holds whenever $\log(2/\delta)=o(n)$ for the average test, $\log(2|\mathcal C|/\delta)=o(K)$ for the scan test, and $\log(2m/\delta)=o(m)$ for the degree test.
Because the target risk $\delta$ may vanish arbitrarily slowly in the strong-detection problem, the average- and degree-test conditions hold automatically, while the scan-test condition reduces to $\log|\mathcal C|=o(K)$.
The fixed-size scanning classes used for the models in \Cref{sec:results}, including the substructure scans described in \Cref{remark: generalized scan test}, satisfy this condition when $k$ or $\Delta$ is polynomial in $n$.

Together, \Cref{proposition: local alternative regime,proposition: upper bound for NEF models in local alternative regime} show that, after matching chi-square divergences, both lower and upper bounds transfer across NEF models in the local alternative regime without model-specific analysis.
Thus, thresholds that are tight for one NEF model remain tight for any other such model.
The remaining technical details for \Cref{theorem: universality of tightness in the local alternative regime via chi-square matching} are deferred to \Cref{proof: universality of tightness in the local alternative regime via chi-square matching}.

We conclude by relating $\tau_n$ to the usual parameters of two representative NEF models.
For the Bernoulli model with success probabilities $p_0<p_1$, we have $\tau_n=(p_1-p_0)/\sqrt{p_0(1-p_0)}$.
For the Gamma model with fixed shape parameter $\alpha$ and scale parameters $\beta_0<\beta_1$, we have $\tau_n=\sqrt{\alpha}(\beta_1-\beta_0)/\beta_0$.
In the local alternative regime, $p_0$ and $\beta_0$ are fixed, so $\tau_n\asymp p_1-p_0$ for the Bernoulli model and $\tau_n\asymp\beta_1-\beta_0$ for the Gamma model.

\subsubsection{\texorpdfstring{Applications to Erd\H{o}s--R\'enyi models}{Applications to Erdos--Renyi models}}\label{section: applications to Erdos--Renyi models}
In this subsection, we apply the preceding equivalence results to Erd\H{o}s--R\'enyi observations, first for planted clique detection and then for arbitrary planted subgraphs.

\paragraph{Planted clique detection}\label{section: equivalence on planted clique}
We revisit the planted clique problem in \Cref{section: planted clique} to provide an explicit example of how our universality result operates in the local alternative regime, as well as of the possibilities and limitations of transferring lower and upper bounds beyond the local alternative regime. Specifically, we consider the Gaussian model in \Cref{sec:gaussian-testing} and Erd\H{o}s--R\'{e}nyi models and examine how the applicability of our equivalence results depends on the asymptotic scaling of the success probability parameters in the Erd\H{o}s--R\'{e}nyi model. Consider the latent structure $S$ defined in \Cref{equation: clique latent structure} and the prior distribution $\dot \nu$ for fixed-size case. Let model $A$ and $B$ be the standard Gaussian model and be the Erd\H{o}s--R\'enyi, respectively. More specifically, we consider the model A
\begin{align*}
	H_0^A\,:\,Y_l \stackrel{iid}{\sim} \mathrm{N}(0,1)\quad\mathrm{v.s.}\quad H_1^A\,:\,Y_l\,|\,S \stackrel{ind}{\sim} \begin{cases}\mathrm{N}(\mu,1), & \text{if } l \in S,\\ \mathrm{N}(0,1), & \text{if } l \notin S,\end{cases}
\end{align*}
with $\mu>0$, and consider the model B
\begin{align*}
	H_0^B\,:\,Y_l \stackrel{iid}{\sim} \mathrm{Ber}(p_0)\quad\mathrm{v.s.}\quad H_1^B\,:\,Y_l\,|\,S \stackrel{ind}{\sim} \begin{cases}\mathrm{Ber}(p_1), & \text{if } l \in S,\\ \mathrm{Ber}(p_0), & \text{if } l \notin S,\end{cases}
\end{align*}
where $1>p_1>p_0>0$. Note that the two models can be expressed with exponential tilts with $(\theta^A_0,\theta^A_1)=(0,\mu)$ and $(\theta^B_0,\theta^B_1)=(\log(p_0/(1-p_0)),\log(p_1/(1-p_1)))$. The chi-square divergences and effective signal strengths of the two models are given by
\begin{align*}
	\chi^2(P^A_{\theta_1}\|P^A_{\theta_0}) &= \exp(\mu^2)-1,\quad\chi^2(P^B_{\theta_1}\|P^B_{\theta_0}) = \frac{(p_1-p_0)^2}{p_0(1-p_0)}, \quad \tau^A_n = \mu, \quad \tau^B_n = \frac{p_1-p_0}{\sqrt{p_0(1-p_0)}}.
\end{align*}
This motivates the scaling of the signal strength between the two models, $\mu=\sqrt{\log(1+(p_1-p_0)^2/(p_0(1-p_0)))}$, which guarantees identical chi-square divergences, and thus the same $\bbE_{\bfY\sim \bbQ^M}\big[L_{\dot\nu}^M(\bfY)^2\big]$ and equivalent $\bbE_{\bfY\sim \bbQ^M}\big[L_{\dot\nu}^{M,\leq D}(\bfY)^2\big]$ for $M\in\{A,B\}$ under the assumption $D_n\chi^2_n=O(1)$. We study how the equivalence result can be applied to the planted clique problem in the following three scenarios, depending on the polynomial orders of $p_0$ and $p_1$ in the Erd\H{o}s--R\'enyi model. In particular, we analyze statistical and computational distinguishability in terms of the difference between the probability parameters $p_1-p_0$ and define the corresponding thresholds $(p_1-p_0)_{\mathrm{stat}}^\ast$ and $(p_1-p_0)_{\mathrm{comp}}^\ast$ analogously to \Cref{definition: statistical and computational thresholds}. Depending on the scenario, the model may fall within the local alternative regime or exhibit a polynomially vanishing effective signal strength. In such cases, a vacuous threshold may arise under the convention in \Cref{definition: statistical and computational thresholds in NEF models} or its extension. We explicitly indicate such cases when presenting the results.

\begin{itemize}
	\item \textbf{Scenario 1: $p_0=c$ and $p_1=c+m^{-\beta}$ for some fixed constants $c\in(0,1)$ and $\beta>0$.} Note that this case corresponds to the local alternative regime. Although the assumptions for the lower and upper bounds can be verified directly for the Gaussian and Bernoulli models, we examine them through the local alternative framework to illustrate the universality of the tightness of the lower and upper bounds. 
	For the lower bound, the condition $D_n\log(\chi^2_n/\tau_n^2)=O(1)$ boils down to $D_n\chi^2_n=O(1)$ in \Cref{proposition: local alternative regime}, which is satisfied for sequence $D_n=(\log m)^{1+\epsilon}$ for any $\epsilon>0$, since $\chi^2_n\asymp m^{-2\beta}$. For the upper bound, the local regularity condition in \Cref{proposition: upper bound for NEF models in local alternative regime} is automatically satisfied, and the condition $\log(|\mathcal C|)=o(K)$ is satisfied for the planted clique problem, since $|\mathcal C|=\binom{m}{k}\leq m^k$ and $\log m=o(k)$. Provided that $\chi^2/\tau_n^2\to 1$ from $D_n\log(\chi^2_n/\tau_n^2)=O(1)$,  $\tau^A_n/\tau^B_n\to 1$ is guaranteed, and we replace $\tau^A_n=\mu$ with $\tau^B_n=(p_1-p_0)/{\sqrt{p_0(1-p_0)}}$ in the statistical and computational thresholds for the Gaussian model displayed in \Cref{theorem: fixed-size clique problem} to obtain the corresponding thresholds for the Erd\H{o}s--R\'enyi model:
	\begin{align*}
		(p_1-p_0)_{\mathrm{stat}}^\ast \asymp\left\{
		\begin{array}{@{}c@{\quad}l}
		\sqrt{\frac{1}{k}}, & k \ll m^{2/3}, \\
		\frac{m}{k^2}, & k \gg m^{2/3},
		\end{array}
		\right.
		\qquad
		(p_1-p_0)_{\mathrm{comp}}^\ast\asymp\left\{
		\begin{array}{@{}c@{\quad}l}
		\vphantom{\sqrt{\frac{1}{k}}}1, & k\ll m^{1/2}\quad(\text{vacuous}),\\
		\frac{m}{k^2}, & k\gg m^{1/2},
		\end{array}
		\right.
	\end{align*}
	provided that $\tau^B_n\asymp p_1-p_0$.
	We note that the constant branch of the computational threshold for $k\ll m^{1/2}$ follows from our threshold convention in \Cref{definition: statistical and computational thresholds in NEF models}.

	\item \textbf{Scenario 2: $p_0=m^{-\alpha}$ and $p_1=c\cdot m^{-\alpha}$ for some fixed constants $c>1$ and $\alpha\in(0,2)$.} First, we note that the case $\alpha>2$ is not considered, because the expected number of edges in the planted clique under the alternative is $(n-K)p_0+Kp_1\leq cnp_0=o(1)$. Consequently, for every fixed structure $s$, Markov's inequality gives $\bbP_s^B\{\bfY=\mathbf{0}\}=1-o(1)$. Similarly, $\bbQ^B\{\bfY=\mathbf{0}\}=1-o(1)$, and averaging the fixed-structure bound over $S\sim\dot\nu$ gives $\bbP_{\dot\nu}^B\{\bfY=\mathbf{0}\}=1-o(1)$. It follows that the total variation distance between the null and alternative distributions converges to zero, so strong detection is impossible in this regime; hence, we restrict our attention to $\alpha\in(0,2)$.
	
	While this case does not fall into the local alternative regime, the model-specific analysis can be carried out to verify the assumptions for the lower and upper bounds. As pointed out earlier in \Cref{section: lower bound for NEF models}, the condition $D_n\log(\chi^2_n/\tau_n^2)=O(1)$ in \Cref{theorem: equivalence between NEF models} is automatically satisfied for the Bernoulli model, and also is satisfied for the Gaussian model with $\mu^2=O(D^{-1})$. 
	Note that the chi-square matching between the two models yields $\mu^2\asymp m^{-\alpha}$, and thus the condition $D_n\log(\chi^2_n/\tau_n^2)=O(1)$ is satisfied for sequence $D_n=(\log m)^{1+\epsilon}$ for any $\epsilon>0$. For the upper bound, the local regularity condition in \Cref{proposition: upper bound for NEF models} is required. Note that the curvature $\psi''(\theta)$ is the variance of the model, and thus the Gaussian model always satisfies the local regularity condition with constant $C=1$. For the Bernoulli model, we claim that for arbitrary sequence $r_n=O(1)$, the local regularity condition is satisfied with some constant $C>1$. To see this, from $p(\theta)=e^\theta/(1+e^\theta)$ and $\psi''(\theta)=p(\theta)(1-p(\theta))$, observe that
	$$
	\frac{\psi''(\theta_j+u)}{\psi''(\theta_j)}=\frac{e^u}{(1-p(\theta_j)+p(\theta_j)e^u)^2}\leq e^{|u|},
	$$
	for $j\in\{0,1\}$ and $|u|\leq r_n$. Consequently, for any $R>0$ such that $r_n\leq R$ for sufficiently large $n$, we have
\begin{align*}
	\sup_{\substack{j\in\{0,1\}\\|u|\leq r_n}}\frac{\psi''(\theta_j+u)}{\psi''(\theta_0)}= \sup_{\substack{j\in\{0,1\}\\|u|\leq r_n}}\bigg(\frac{\psi''(\theta_j+u)}{\psi''(\theta_j)}\frac{\psi''(\theta_j)}{\psi''(\theta_0)}\bigg)\leq \max\bigg\{1,\frac{p(\theta_1)(1-p(\theta_1))}{p(\theta_0)(1-p(\theta_0))}\bigg\}\sup_{\substack{j\in\{0,1\}\\|u|\leq r_n}}\frac{\psi''(\theta_j+u)}{\psi''(\theta_j)}.
\end{align*}
Therefore, if $p(\theta_1)/p(\theta_0)=c>1$ for some constant $c$, the local regularity condition holds with $C=ce^R$ for sufficiently large $n$.

    It remains to verify $r_{\textrm{avg}}=O(1)$ and $r_{\textrm{scan}}=O(1)$. Note that $n\asymp m^2$ and $v_0\asymp m^{-\alpha}$, and thus $r_{\textrm{avg}}=\sqrt{2(nv_0)^{-1}\log(2/\delta)}\asymp m^{-1+\alpha/2}\sqrt{\log(2/\delta)}$, and we can find $\delta\to0$ which satisfies $r_{\textrm{avg}}=O(1)$, provided that $\alpha<2$. On the other hand, we have $$
	r^2_{\textrm{scan}}=2(Kv_0)^{-1}\log(2|\mathcal C|/\delta)\asymp \frac{\log\binom{m}{k}+\log(1/\delta)}{k^2m^{-\alpha}}.
	$$
	Therefore, $r_{\textrm{scan}}=O(1)$ is satisfied if $\log\binom{m}{k}=O(k^2m^{-\alpha})$ and $\delta\to0$ sufficiently slowly.
	Substituting $\mu$ with $\tau_n^B\asymp\sqrt{p_1-p_0}$, the Gaussian statistical and computational detection bounds translate into the following thresholds for the Erd\H{o}s--R\'enyi model
	\begin{align*}
		(p_1-p_0)_{\mathrm{stat}}^\ast\asymp\left\{
		\begin{array}{@{}c@{\quad}l}
		\frac{1}{k}, & k \ll m^{2/3}, \\
		\frac{m^{2}}{k^4}, & k \gg m^{2/3},
		\end{array}
		\right.
		\qquad
		(p_1-p_0)_{\mathrm{comp}}^\ast\asymp\left\{
		\begin{array}{@{}c@{\quad}l}
		\vphantom{\frac{m^{2}}{k^4}}1, & k \ll m^{1/2}\quad(\text{vacuous}), \\
		\frac{m^{2}}{k^4}, & k \gg m^{1/2}.
		\end{array}
		\right.
	\end{align*}
	We note that the additional condition $\log\binom{m}{k}=O(k^2m^{-\alpha})$ required for local regularity is automatically satisfied whenever the scan test achieves strong detection at polynomial order. Thus, it imposes no further restriction on the scan-test upper bound. We also remark that, although this scenario does not fall within the local alternative regime, the effective signal strength $\tau_n^B$ always vanishes polynomially. When $k\ll m^{1/2}$, the transferred low-degree lower bound rules out computational strong detection for every such signal in this scenario. The constant branch of the computational threshold represents this impossibility under our convention.

	\item \textbf{Scenario 3: $p_0=m^{-\alpha}$ and $p_1=m^{-\beta}$ for some fixed  $\alpha>\beta>0$.} To illustrate how our equivalence result transfers and to point out its limitations, we focus on the regime $\alpha<1$, where both lower and upper bounds have been carefully studied in the existing literature. For instance, the tight statistical lower and upper bounds are established in \citet{arias2014community,verzelen2015community}, and the tight computational lower and upper bounds are studied in \citet{dhawan2025detection}. We rewrite their results in terms of $p_1-p_0$, with our notation:
	\begin{align*}
		(p_1-p_0)_{\mathrm{stat}}^\ast\asymp\left\{
		\begin{array}{@{}c@{\quad}l}
		\frac{1}{k}, & k \ll m^{1-\alpha/2}, \\
		\frac{m^{1-\alpha/2}}{k^2}, & k \gg m^{1-\alpha/2},
		\end{array}
		\right.
		\qquad
		(p_1-p_0)_{\mathrm{comp}}^\ast\asymp\left\{
		\begin{array}{@{}c@{\quad}l}
		\vphantom{\frac{m^{1-\alpha/2}}{k^2}}\frac{1}{k^\alpha}, & k \ll m^{1/2}, \\
		\frac{m^{1-\alpha/2}}{k^2}, & k \gg m^{1/2}.
		\end{array}
		\right.
	\end{align*}
	With these existing thresholds, we observe below that the lower- and upper-bound transfers provided by \Cref{theorem: equivalence between NEF models} and \Cref{proposition: upper bound for NEF models} are either loose or inapplicable in this scenario.

	We first restrict our attention to the case $\alpha<2\beta$. In this case, the chi-square matching between the two models requires that $\mu^2\asymp \log(1+m^{\alpha-2\beta})$, which is of the polynomially vanishing order $m^{\alpha-2\beta}$, and thus the degree condition $\mu^2=O(D^{-1})$ for the Gaussian model is satisfied for $D_n=(\log m)^{1+\epsilon}$. Consequently, both the lower bounds from the Gaussian model are transferrable to the Erd\H{o}s--R\'enyi model through \Cref{theorem: equivalence between NEF models}; the second moments of the likelihood ratio (LR) and the low-degree likelihood ratio (LDLR) are bounded if
	\begin{align*}
		\mathrm{LR}\,:\,(p_1-p_0)\ll\left\{
		\begin{array}{@{}c@{\quad}l}
		\frac{m^{-\alpha/2}}{\sqrt{k}}, & k \ll m^{2/3}, \\
		\frac{m^{1-\alpha/2}}{k^2}, & k \gg m^{2/3},
		\end{array}
		\right. \qquad \mathrm{LDLR}\,:\,(p_1-p_0)\ll\left\{
		\begin{array}{@{}c@{\quad}l}
		\frac1{m^{\alpha/2}}, & k \ll m^{1/2}, \\
		\frac{m^{1-\alpha/2}}{k^2}, & k \gg m^{1/2},
		\end{array}
		\right.
	\end{align*}
	respectively. Observe that the first branches of the transferred lower bounds are loose compared to the tight statistical and computational thresholds above. On the other hand, when $\alpha>2\beta$, the chi-square matching yields $\mu^2\asymp\log m$, causing the full likelihood ratio second moment of the matched Gaussian model to diverge. Consequently, although the second-moment equivalence remains valid, it provides no nontrivial statistical lower bound in this regime. In addition, the degree condition $\mu^2=O(D^{-1})$ is now violated for any $D_n=(\log m)^{1+\epsilon}$, and thus the low-degree equivalence cannot be applied, and transferring the computational lower bound from the Gaussian model through \Cref{theorem: equivalence between NEF models} is not applicable.

	This observation helps clarify the scope of applicability of \Cref{theorem: equivalence between NEF models}: its lower-bound transfer applies to bounds obtained through the second moment method or conditional second moment arguments based on the latent structure $S\sim\nu$, rather than to arbitrary lower bounds. Indeed, the tight statistical and computational lower bounds derived in the literature for the Bernoulli model rely on model-specific truncation arguments that exploit properties of the observed graph, which are not covered by \Cref{theorem: equivalence between NEF models}. For the statistical lower bound, \citet{arias2014community,verzelen2015community} use a truncated likelihood ratio based on a high-probability event that excludes atypical configurations in which some subset of the planted community contains unusually many edges. Similarly, for the computational lower bound, \citet{dhawan2025detection} obtains a sharp bound by excluding rare realizations of $\bfY$ in which the planted subgraph contains an unusually dense small subgraph.

	For the upper bounds, the local regularity condition fails for the Bernoulli model, regardless of the choice of $r_n=O(1)$, since $\psi''(\theta_1)/\psi''(\theta_0)=p_1(1-p_1)/p_0(1-p_0)\asymp m^{\alpha-\beta}\to \infty$.
	Consequently, the Gaussian upper bounds cannot be directly transferred to the Bernoulli model. Nevertheless, statistical upper bounds can be established directly by analyzing the average and scan tests using Bernoulli tail bounds rather than sub-Gaussian approximations. Such bounds have been established in \citet{arias2014community,verzelen2015community,dhawan2025detection}. In particular, the sufficient conditions for these tests to achieve strong detection are
	\begin{align*}
		(p_1-p_0)\gg\left\{
		\begin{array}{@{}c@{\quad}l}
		\frac{1}{k}, & k \ll m^{1-\alpha/2}, \\
		\frac{m^{1-\alpha/2}}{k^2}, & k \gg m^{1-\alpha/2}.
		\end{array}
		\right.
	\end{align*}
	where the first branch is achieved by the scan test and the second by the average test. These upper bounds match the statistical thresholds stated above. However, the average test alone does not attain the sharp computational threshold in the regime $k\ll m^{1/2}$. To close this gap, \citet{dhawan2025detection} develop a polynomial-time algorithm based on counting copies of a suitably chosen constant-size balanced subgraph. In this regime, their test achieves vanishing risk whenever $(p_1-p_0)\gg k^{-\alpha}$, thereby attaining the computational threshold.

\end{itemize}

\paragraph{Detecting arbitrary planted subgraphs.}

For Erd\H{o}s--R\'enyi models, the statistical and computational thresholds in the fixed-size case are studied in \citet{elimelech2025detecting} in terms of some graph-theoretic properties such as maximum degree, maximum subgraph density and the number of edges or vertices.
Notably, in the local alternative regime under the assumptions of their Theorem 16, the computational threshold at polynomial order is determined by the detection boundaries of the average and degree tests.
We remark that the computational lower bounds in their work are established using the low-degree second moment method\footnote{However, we note that their lower bounds (Theorem 16) assume $D=O(\log m)$ instead of $D=\log^{1+\epsilon} m$.}.
Moreover, their upper-bound analysis yields the same asymptotic polynomial-order signal requirements as our analysis.
By the transfer argument underlying \Cref{theorem: universality of tightness in the local alternative regime via chi-square matching}, with computational thresholds interpreted under a logarithmic-degree version of the low-degree conjecture, the corresponding lower and upper bounds can be transferred to arbitrary NEF models through chi-square matching. Consequently, the corresponding computational threshold in the local alternative regime is characterized as
\begin{align}\label{equation: generic computational threshold for NEF models}
\tau^\ast_{\textrm{comp}}
\asymp
1\wedge \frac{m}{K}\wedge \frac{\sqrt m}{\mathfrak{d}_{\textrm{max}}(S)}
\end{align}
up to subpolynomial factors, under a version of the low-degree conjecture that allows $D=O(\log m)$.

Now we shift our focus to the statistical threshold. For the upper bound, their scan statistic differs from our naive scan statistic in that it scans copies of a maximum-density subgraph $H^\star\subseteq S$, rather than copies of the entire planted graph $S$. Nevertheless, the performance of this modified scan statistic can be analyzed by essentially the same argument as in \Cref{proof: scan test for NEF models}. Let $\mathcal C_{H^\star}$ denote the collection of edge sets of all copies of $H^\star$ in the ambient complete graph, and let $K_\star=|H^\star|$. The proof of our NEF scan bound applies to this modified statistic with $(\mathcal C,K)$ replaced by $(\mathcal C_{H^\star},K_\star)$. The resulting sufficient condition for the modified scan test is then
$\tau_n\geq \sqrt{2C_{scan}K_\star^{-1}}\big(2\sqrt{\log(2/\delta)}+\sqrt{\log|\mathcal C_{H^\star}|}\big).$
This condition is satisfied if $\liminf_{m\to \infty}|V(H^\star)|^{-1}K_\star(\log m)^{-1}\tau_n^2>2$, provided that $|\mathcal C_{H^\star}|\le m^{|V(H^\star)|}$ and $\delta$ vanishes sufficiently slowly. This is of the same asymptotic order as the scan detection boundary $\liminf_{m\to\infty}|V(H^\star)|^{-1}K_\star(\log m)^{-1}D_{\mathrm{KL}}(\mathrm{Ber}(p_1)\,\|\,\mathrm{Ber}(p_0))>1$ in \citealp{elimelech2025detecting}, as can be verified from the fact that the Taylor expansion of the KL divergence in the local alternative regime yields $D_{\mathrm{KL}}(\mathrm{Ber}(p_1)\,\|\,\mathrm{Ber}(p_0))=\tau_n^2(1+o(1))/2.$ Moreover, since $\tau_n\to0$ in the local alternative regime, the above sufficient condition necessarily implies
$|V(H^\star)|K_\star^{-1}\log m\to0.$
Together with $|\mathcal C_{H^\star}|\leq m^{|V(H^\star)|}$, the local regularity condition $\log|\mathcal C_{H^\star}|=o(K_\star)$ required for our NEF scan bound is automatically satisfied.

Therefore, the upper bound of \citealp{elimelech2025detecting} which consists of the average, maximum-degree, and modified scan tests can be transferred to any other NEF model through chi-square matching, and this implies that strong detection is possible if
$$
\tau_n\gg\frac{m}{K}\wedge \frac{\sqrt m}{\mathfrak{d}_{\textrm{max}}(S)} \wedge \sqrt{\frac{|V(H^\star)|}{K_\star}}$$
in the local alternative regime.

On the other hand, the lower bound in their work is obtained by controlling the second moment of the likelihood ratio, and hence can again be transferred to another NEF model through chi-square matching. In our notation, their statistical lower bound implies that statistical detection is impossible when
$$
\tau_n\ll\begin{cases}
\displaystyle
\frac{m}{K}\wedge\frac{\sqrt m}{\mathfrak d_{\mathrm{max}}(S)}\sqrt{\frac{|V(H^\star)|}{K_\star}}\wedge\sqrt{\frac{|V(H^\star)|}{K_\star}},
&\displaystyle\lim_{m\to\infty}\frac{\log\!\left(K_\star/|V(H^\star)|\right)}{\log |V(S)|}<1,
\\
\displaystyle\frac{m}{K}\wedge\frac{\sqrt m}{\mathfrak d_{\mathrm{max}}(S)}\wedge\sqrt{\frac{|V(H^\star)|}{K_\star}},&\displaystyle\lim_{m\to\infty}
\frac{\log\!\left(K_\star/|V(H^\star)|\right)}
{\log |V(S)|}=1.
\end{cases}
$$
Since these lower bounds are obtained through the likelihood-ratio second moment, our equivalence result extends the same statistical impossibility regions to arbitrary NEF models in the local alternative regime through chi-square matching. 

While these results provide a useful structure-agnostic characterization of the statistical threshold, the corresponding lower and upper bounds do not match in all regimes, as noted by \citet{elimelech2025detecting}. Our multiple-clique problem provides a concrete example in which these general statistical lower and upper bounds do not match. Consider $\gamma$ disjoint $k$-cliques with $k=m^{3/5}$ and $\gamma=m^{1/20}$. Then
$$
K\asymp \gamma k^2=m^{5/4},\qquad\mathfrak d_{\mathrm{max}}(S)\asymp k=m^{3/5},
$$
and we may take $H^\star$ to be a single $k$-clique, for which $
|V(H^\star)|=k$ and $K_\star\asymp k^2.$
Moreover, since
$$
\frac{\log(K_\star/|V(H^\star)|)}{\log |V(S)|}\to \frac{12}{13}<1,
$$
the first case of the statistical lower bound applies. The general upper bound is then of order
$$
m^{-1/4}\wedge m^{-1/10}\wedge m^{-3/10}=m^{-3/10},
$$
whereas the general lower bound is only of order
$$
m^{-1/4}\wedge m^{-2/5}\wedge m^{-3/10}=m^{-2/5}.
$$
Thus, the general lower and upper bounds leave a polynomial gap in this example.
In contrast, our analysis of the lower bound captures the correct polynomial order and thereby closes this gap.
We also point out that, while the statistical and computational thresholds in our main text are stated only up to polynomial order for clarity of presentation, the corresponding appendices retain sharper signal scales including logarithmic factors.
These more refined rates do not follow directly from the general bounds of \citealp{elimelech2025detecting}, particularly in regimes where those bounds do not match.
Finally, we emphasize that the framework of \citealp{elimelech2025detecting} is formulated for ordinary graphs, where each edge is a pair of vertices, and does not directly cover hypergraph settings such as \Cref{corollary: multiple planted hypercliques problem}.
Accordingly, the statistical and computational thresholds for NEF models obtained above by transferring their results should likewise be understood as applying to the ordinary graph setting.

\section{Discussion}\label{sec:discussion}

Our results identify the overlap distribution of two independent latent structures as a common basis for analyzing statistical and computational thresholds in graph combinatorial testing.
In the Gaussian model, both the full and low-degree likelihood-ratio second moments depend only on this overlap, yielding a unified analysis of planted cliques, multiple cliques and hypercliques, and thick cycles and paths.
Comparing the resulting lower bounds with the scan, average, and coordinatewise maximum tests identifies regimes with statistical-computational gaps.
Spectral methods provide another important approach to planted-structure detection, and, under suitable conditions on the spectral gap and moments, matrix powers or trace moments connect such methods to low-degree polynomial statistics \citep{kunisky2019notes,wein2025computational}.
We do not analyze spectral tests separately because our benchmark polynomial-time tests already match the low-degree lower bounds at polynomial order for the Gaussian models considered.
The low-degree lower bounds rule out strong separation by polynomial statistics of polylogarithmic degree; their interpretation as lower bounds for all polynomial-time tests relies on the working low-degree hypothesis in \Cref{conjecture: low-degree polynomial conjecture}.

The natural exponential family results separate the role of the latent structure from that of the observation distribution.
Under the conditions of \Cref{theorem: equivalence between NEF models}, matching the one-coordinate chi-square divergence preserves the full likelihood-ratio second moments and makes the low-degree second moments equivalent up to constant factors.
The upper bounds depend on the observation model through an effective signal strength and a local regularity condition.
In the local alternative regime, chi-square matching also matches the effective signal strengths asymptotically, and local regularity holds automatically; under its stated assumptions, \Cref{theorem: universality of tightness in the local alternative regime via chi-square matching} therefore transfers tightness of these bounds across natural exponential families.
Outside this regime, chi-square matching alone need not yield sharp thresholds: the relevant signal scales may differ, local regularity may fail, and sharp lower bounds may require observation-dependent truncations.
The sparse Erd\H{o}s--R\'enyi planted-clique example in \Cref{section: equivalence on planted clique} illustrates these limitations.

Several questions remain beyond the polynomial-order thresholds established here, including sharp logarithmic factors, optimal constants, and behavior on critical phase boundaries.
Our sharp results focus on strong detection with observations that are independent conditional on the latent structure, known noise parameters, and symmetric priors; weak detection, recovery, adaptation, and models with conditionally dependent or heterogeneous observations may exhibit different behavior.
Specialized polynomial-time statistics may also improve on the benchmark tests outside the local alternative regime.

One direction for further work is to characterize broader classes of graphs and hypergraphs through structural properties of their overlap distributions.
Another direction is to study mean-matched nulls, whose coordinate means agree with those of the mixture alternative.
This removes the mean shift exploited by our average test and calls for a separate analysis of detection thresholds and distributional transfer.
A transfer principle for observation-dependent conditional second moments could help extend the universality theory beyond the local alternative regime.
Rigorous average-case reductions or other complexity-theoretic lower bounds would help determine when our low-degree predictions describe the limits of all polynomial-time algorithms.

\section*{Acknowledgments}

We thank Wasim Huleihel and Matey Neykov for helpful discussions. Ikjun Choi and Alessandro Rinaldo acknowledge partial support from the NIH Award RF1NS121913.

\bibliographystyle{plainnat}
\bibliography{ref}

\appendix
\clearpage

\begin{center}
\Large
{\bf Appendix}
\end{center}

The appendix contains supporting technical results, complete proofs, and additional examples of latent structures.
\Cref{sec:app-technical-lemmas} collects the technical lemmas used throughout the paper, and \Cref{sec:app-proofs} contains the proofs.
\Cref{sec:app-additional-examples} studies further examples that complement the problems considered in the main text.
A detailed outline is provided below.

\section*{Outline}

\addcontentsline{toc}{section}{Appendix} 
\startcontents
\printcontents{}{1}{\setcounter{tocdepth}{2}}

\section{Technical lemmas}\label{sec:app-technical-lemmas}

This section collects the moment and decoupling inequalities used in the proofs.

The following lemma provides a moment bound for sub-Poissonian random variables.
Its assumptions hold for Poisson and binomial random variables, as well as for random variables dominated by them in the convex order.
A recurring example in this paper is a hypergeometric random variable, which can be represented as a sum of negatively associated Bernoulli random variables.

\begin{lemma}[{\citealt[Corollary 1]{ahle2022sharp}}]\label{lemma: moment bound of sub-Poissonian distribution}
	Let $X$ be a nonnegative random variable with mean $\lambda>0$ whose moment generating function satisfies $\ee[\exp(tX)]\leq\exp\{\lambda(e^t-1)\}$ for every $t>0$.
	Then, for every $d>0$,
	$$
	\ee[(X/\lambda)^d]\leq \bigg(\frac{d/\lambda}{\log(1+d /\lambda)}\bigg)^d\leq \bigg(1+\frac{d}{2\lambda}\bigg)^d\leq\exp\bigg(\frac{d^2}{2\lambda}\bigg).
	$$
\end{lemma}

Finner's inequality generalizes H\"older's inequality by accounting for the coordinates on which each function depends in a product measure space.

\begin{lemma}[{\citealt[Theorem 2.1]{finner1992generalization}}]\label{lemma: Finner's inequality}
Let $I$ and $J$ be finite index sets, and let $\emptyset\neq T_j\subseteq I$ and $p_j\geq1$ for each $j\in J$.
For each $i\in I$, define $J_i:=\{j\in J:i\in T_j\}$, and assume that $\sum_{j\in J_i}p_j^{-1}=1$.
Let $(\Omega_i,\mathcal A_i,\mu_i)$, $i\in I$, be measure spaces, and for each $T\subseteq I$, write
\[
(\Omega_T,\mathcal A_T,\mu_T):=\bigg(\prod_{i\in T}\Omega_i,\bigotimes_{i\in T}\mathcal A_i,\bigotimes_{i\in T}\mu_i\bigg).
\]
Suppose that $f_j\in L_{p_j}(\Omega_{T_j},\mathcal A_{T_j},\mu_{T_j})$ for each $j\in J$, and view $f_j$ as a function on $\Omega_I$ depending only on the coordinates in $T_j$.
Then $\prod_{j\in J}|f_j|\in L_1(\Omega_I,\mathcal A_I,\mu_I)$ and
$$
\int \prod_{j\in J}|f_j|\,d\mu_I
\leq
\prod_{j\in J}
\left(
\int |f_j|^{p_j}\,d\mu_{T_j}
\right)^{1/p_j}.
$$

\end{lemma}

The next inequality bounds moments of a U-statistic by those of its decoupled counterpart; see \citet{delapena1999decoupling,mao2025information}.
\begin{lemma}\label{lemma: decoupling of U-statistics}
	Let $X_1,\ldots,X_m$ be i.i.d.\ random variables taking values in a measurable space $(\mathcal X,\mathcal A)$, and let $(X'_1,\ldots,X'_m)$ be an independent copy of this sequence.
	For each ordered pair $(i,j)$ with $i\neq j$, let $f_{ij}:\mathcal X^2\to\mathbb R$ be a measurable kernel satisfying $\ee[|f_{ij}(X_i,X_j)|]<\infty$.
	There exists an absolute constant $C>0$ such that, for every convex nondecreasing function $\Phi:[0,\infty)\to[0,\infty)$ with $\ee[\Phi(|f_{ij}(X_i,X_j)|)]<\infty$ for all $i\neq j$,
	\begin{align*}
		\ee\Phi\bigg(\bigg|\sum_{\substack{i,j\in[m]\\i\neq j}} f_{ij}(X_i,X_j)\bigg|\bigg) &\leq \ee\Phi\bigg(C\bigg|\sum_{\substack{i,j\in[m]\\i\neq j}} f_{ij}(X_i,X'_j)\bigg|\bigg).
	\end{align*}
The same bound applies to sums over $i<j$ by setting $f_{ij}=0$ for $i>j$.
\end{lemma}

The following version applies to U-statistics of any fixed order.
\begin{lemma}[{\citealt[Theorem 3.1.1]{delapena1999decoupling}}]\label{lemma: decoupling of t order U-statistics}
	Let $m\geq t\geq1$ be integers, and let $X_1,\ldots,X_m$ be independent random variables taking values in a measurable space $(\mathcal X,\mathcal A)$.
	For $j\in[t]$, let $(X_i^{(j)})_{i=1}^m$ be copies of this sequence that are mutually independent and independent of the original sequence.
	For each $\bfi=(i_1,\ldots,i_t)\in[m]^t$, let $H_{\bfi}:\mathcal X^t\to\mathbb R$ be a measurable kernel satisfying $\ee[|H_{\bfi}(X_{i_1},\ldots,X_{i_t})|]<\infty$, and assume that $H_{\bfi}\equiv0$ whenever $\bfi$ has repeated indices.
	Then, for every convex nondecreasing function $\Phi:[0,\infty)\to[0,\infty)$ with $\ee[\Phi(|H_{\bfi}(X_{i_1},\ldots,X_{i_t})|)]<\infty$ for all $\bfi\in[m]^t$,
	\begin{align*}
		\ee\bigg[\Phi\bigg(\bigg|\sum_{\bfi\in [m]^t} H_{\bfi}(X_{i_1},\dots,X_{i_t})\bigg|\bigg)\bigg] \leq \ee\bigg[\Phi\bigg(C_t\bigg|\sum_{\bfi\in [m]^t} H_{\bfi}(X^{(1)}_{i_1},\dots,X^{(t)}_{i_t})\bigg|\bigg)\bigg],
	\end{align*}
	where $C_t=2^t\prod_{r=2}^t(r^r-1)$ depends only on $t$.
\end{lemma}

Finally, we record a moment bound for decoupled U-statistics with nonnegative kernels.
\begin{lemma}[{\citealt[Proposition 2.1]{gine2000exponential}}]\label{lemma: moment bound of nonnegative U-statistics}
	Let $m,t\in\mathbb N$, $d>1$, and let $\{X_i^{(j)}:i\in[m],\ j\in[t]\}$ be independent random variables.
	For each $\bfi=(i_1,\ldots,i_t)\in[m]^t$, let $H_{\bfi}$ be a nonnegative measurable function satisfying $\ee[H_{\bfi}(\mathbf X_{\bfi})^d]<\infty$, where $\mathbf X_{\bfi}:=(X^{(1)}_{i_1},\ldots,X^{(t)}_{i_t})$.
	For $J\subseteq[t]$, write $\bfi_J:=(i_j)_{j\in J}$, let $J^c:=[t]\setminus J$, and define $\mathcal F_J:=\sigma(X_i^{(j)}:j\in J,\ i\in[m])$.
	Let $\ee_J$ denote expectation over all variables $X_i^{(j)}$ with $j\in J$, holding the remaining variables fixed.
	There exists a universal constant $C>0$ such that
	$$
	\ee\bigg[\bigg(\sum_{\bfi \in [m]^t} H_{\bfi}(\mathbf X_{\bfi})\bigg)^d\bigg] \leq C^{td} \sum_{J \subseteq [t]} \Bigg[d^{|J|d}\sum_{\bfi_J \in [m]^{|J|}}\ee_{J}\Bigg\{\Bigg(\sum_{\bfi_{J^c} \in [m]^{t-|J|}}\ee_{{J^c}}\big[H_{\bfi}(\mathbf X_{\bfi})\,\big|\,\mathcal F_J\big]\Bigg)^d\Bigg\}\Bigg],
	$$
	By convention, $\sum_{\bfi_{\emptyset}}a=a$ and $\ee_\emptyset[a]=a$.
\end{lemma}

\section{Proofs}\label{sec:app-proofs}

\subsection{\texorpdfstring{Proof of \Cref{lemma: second moment method}}{Proof of Lemma \getrefnumber{lemma: second moment method}}}\label{proof: second moment method}
The variational characterization of total variation and the Bretagnolle--Huber inequality \citep[Section 2.4.1]{tsybakov2009introduction} give, for every test $\Delta_n$,
\begin{equation}\label{equation: risk lower bound via BH}
\begin{aligned}
R_{\nu_n}(\Delta_n)
&\geq 1-\operatorname{TV}(\bbP_{\nu_n},\bbQ_n)\\
&\geq 1-\sqrt{1-\exp\left\{-D_{\mathrm{KL}}(\bbP_{\nu_n}\,\|\,\bbQ_n)\right\}}.
\end{aligned}
\end{equation}
Jensen's inequality under $\bbP_{\nu_n}$ yields
\begin{align*}
D_{\mathrm{KL}}(\bbP_{\nu_n}\,\|\,\bbQ_n)
&=\bbE_{\bfY\sim\bbP_{\nu_n}}\big[\log L_{\nu_n}(\bfY)\big]\\
&\leq \log\bbE_{\bfY\sim\bbP_{\nu_n}}\big[L_{\nu_n}(\bfY)\big]\\
&=\log\bbE_{\bfY\sim\bbQ_n}\big[L_{\nu_n}(\bfY)^2\big]
\leq\log C.
\end{align*}
Substitution into \Cref{equation: risk lower bound via BH} proves the claim. An alternative proof based on classical Le Cam contiguity theory \citep{lecam1960locally} is given, for example, by \citet[Lemma~2]{montanari2017limitation}.

\subsection{\texorpdfstring{Proof of \Cref{proposition: gaussian test upper bound}}{Proof of Proposition \getrefnumber{proposition: gaussian test upper bound}}}\label{proof: gaussian test upper bound}

We prove the three bounds in \Cref{proposition: gaussian test upper bound} and the degree-test bound in \Cref{remark: gaussian degree test}.

\paragraph{Outlier test.}
Recall that $t_{\mathrm{out}}:=\sqrt{2\log(2n/\delta)}$.
Under $\bbQ$, each coordinate is standard normal, so the union bound and the Gaussian tail bound give
\[
\bbQ\{\Delta_{\mathrm{out}}(\bfY)=1\}
\leq n\exp(-t_{\mathrm{out}}^2/2)
=\frac{\delta}{2}.
\]

For the type II error, fix any $s$ in the support of $\nu$, so that $|s|=K\geq1$.
Under $\bbP_s$, the variables $\{Y_\ell:\ell\in s\}$ are independent and distributed as $\mathrm N(\mu,1)$.
Writing $\Phi$ for the standard normal distribution function, we have
\begin{align*}
\bbP_s\{\Delta_{\mathrm{out}}(\bfY)=0\}
&=\bbP_s\bigg\{\max_{\ell\in[n]}Y_\ell\leq t_{\mathrm{out}}\bigg\}\\
&\leq\bbP_s\bigg\{\max_{\ell\in s}Y_\ell\leq t_{\mathrm{out}}\bigg\}\\
&=\prod_{\ell\in s}\bbP_s\{Y_\ell\leq t_{\mathrm{out}}\}\\
&=\Phi(t_{\mathrm{out}}-\mu)^{|s|}.
\end{align*}
For $\mu\geq t_{\mathrm{out}}$, the Gaussian tail bound therefore gives
\begin{align*}
\bbP_s\{\Delta_{\mathrm{out}}(\bfY)=0\}
&\leq\exp\bigg(-\frac{|s|(\mu-t_{\mathrm{out}})^2}{2}\bigg)\\*
&\leq\exp\bigg(-\frac{(\mu-t_{\mathrm{out}})^2}{2}\bigg).
\end{align*}
The last inequality uses $|s|\geq1$.
This bound is at most $\delta/2$ whenever $\mu\geq t_{\mathrm{out}}+\sqrt{2\log(2/\delta)}$.
Since it holds uniformly over $s$ in the support of $\nu$, averaging the type II error and adding the type I error yields $R_\nu(\Delta_{\mathrm{out}})\leq\delta$.

\paragraph{Average test.}
Under $\bbQ$, the average statistic $T_{\mathrm{avg}}(\bfY)$ has distribution $\mathrm N(0,n)$, while under $\bbP_s$ it has distribution $\mathrm N(\mu K,n)$ for every $s$ in the support of $\nu$.
At the threshold $t_{\mathrm{avg}}:=\mu K/2$, the Gaussian tail bound therefore gives
\[
R_\nu(\Delta_{\mathrm{avg}})\leq2\exp\bigg(-\frac{\mu^2K^2}{8n}\bigg).
\]
The right-hand side is at most $\delta$ whenever $\mu\geq\sqrt{(8n/K^2)\log(2/\delta)}$, as claimed.

\paragraph{Scan test.}
For $t_{\mathrm{scan}}:=(\mu K+\sqrt{2K\log|\mathcal C|})/2$, the signal condition in \Cref{proposition: gaussian test upper bound} implies
\begin{align*}
t_{\mathrm{scan}}
&\geq\sqrt{2K\log|\mathcal C|}+\sqrt{2K\log(2/\delta)}
\geq\sqrt{2K\log(2|\mathcal C|/\delta)},\\
\mu K-t_{\mathrm{scan}}
&\geq\sqrt{2K\log(2/\delta)}.
\end{align*}
Under $\bbQ$, the sum over any fixed candidate $s\in\mathcal C$ has distribution $\mathrm N(0,K)$.
The union bound therefore yields
\[
\bbQ\{\Delta_{\mathrm{scan}}(\bfY)=1\}
\leq|\mathcal C|\exp\bigg(-\frac{t_{\mathrm{scan}}^2}{2K}\bigg)
\leq\frac{\delta}{2}.
\]
For the type II error, fix any $s$ in the support of $\nu$.
Under $\bbP_s$, the sum $\sum_{\ell\in s}Y_\ell$ has distribution $\mathrm N(\mu K,K)$ and is bounded above by $T_{\mathrm{scan}}(\bfY)$.
Hence,
\begin{align*}
\bbP_s\{\Delta_{\mathrm{scan}}(\bfY)=0\}
&\leq\bbP_s\bigg\{\sum_{\ell\in s}Y_\ell\leq t_{\mathrm{scan}}\bigg\}\\
&\leq\exp\bigg(-\frac{(\mu K-t_{\mathrm{scan}})^2}{2K}\bigg)
\leq\frac{\delta}{2}.
\end{align*}
Averaging over $s$ and adding the type I error gives $R_\nu(\Delta_{\mathrm{scan}})\leq\delta$.

\paragraph{Degree test.}
We now prove the bound in \Cref{remark: gaussian degree test}, using the threshold $t_{\mathrm{deg}}:=\mu\mathfrak d_{\max}/2$.
Under $\bbQ$, each weighted degree satisfies
$$
T_i^{\mathrm{deg}}(\bfY)=
\sum_{\substack{e\in[m]_2\\ i\in e}}Y_e\sim \mathrm{N}(0,m-1),
$$
because there are $m-1$ edges incident to each vertex.
The union bound and the Gaussian tail bound give
\begin{align*}
\bbQ\big\{\Delta_{\mathrm{deg}}(\bfY)=1\big\}&=\bbQ\big\{T_{\mathrm{deg}}(\bfY)>t_{\mathrm{deg}}\big\}\\
&\leq\sum_{i=1}^m\bbQ\bigg\{T_i^{\mathrm{deg}}(\bfY)>\frac{\mathfrak d_{\max}\mu}{2}\bigg\}\\
&\leq m\exp\bigg(-\frac{\mathfrak d_{\max}^2\mu^2}{8(m-1)}\bigg).
\end{align*}

For the type II error, fix any $s$ in the support of $\nu$ and choose a vertex $i_s\in[m]$ with $\mathfrak d_s(i_s)=\mathfrak d_{\max}$.
Under $\bbP_s$, we have $T_{i_s}^{\mathrm{deg}}(\bfY)\sim\mathrm N(\mu\mathfrak d_{\max},m-1)$.
Since $T_{\mathrm{deg}}(\bfY)\geq T_{i_s}^{\mathrm{deg}}(\bfY)$, we have
\begin{align*}
\bbP_s\big\{\Delta_{\mathrm{deg}}(\bfY)=0\big\}&=\bbP_s\big\{T_{\mathrm{deg}}(\bfY)\leq t_{\mathrm{deg}}\big\}\\
&\leq\bbP_s\bigg\{T_{i_s}^{\mathrm{deg}}(\bfY)\leq\frac{\mathfrak d_{\max}\mu}{2}\bigg\}\\
&\leq\exp\bigg(-\frac{\mathfrak d_{\max}^2\mu^2}{8(m-1)}\bigg),
\end{align*}
where the last inequality follows from the Gaussian tail bound.
The bound holds uniformly over $s$ in the support of $\nu$, so averaging the type II error and adding the type I error gives
\begin{align*}
R_\nu(\Delta_{\mathrm{deg}})&\leq(m+1)\exp\bigg(-\frac{\mathfrak d_{\max}^2\mu^2}{8(m-1)}\bigg)\\
&\leq2m\exp\bigg(-\frac{\mathfrak d_{\max}^2\mu^2}{8(m-1)}\bigg).
\end{align*}
Consequently, $R_\nu(\Delta_{\mathrm{deg}})\leq\delta$ whenever
$$
\mu\geq\sqrt{\frac{8(m-1)}{\mathfrak d_{\max}^2}\log\frac{2m}{\delta}}.
$$

\subsection{\texorpdfstring{Proof of \Cref{lemma: statistical lower bound for fixed-size clique problem}}{Proof of Lemma \getrefnumber{lemma: statistical lower bound for fixed-size clique problem}}}\label{proof: statistical lower bound for fixed-size clique problem}
Let $\dot V:=\{i\in[m]:\dot A_i\leq k\}$ be the vertex set of $S$, and define $\dot V'$ using an independent copy $\dot\bfA'$ of $\dot\bfA$.
Then $H:=|\dot V\cap\dot V'|$ has the hypergeometric distribution with parameters $(m,k,k)$, and $|S\cap S'|=\binom H2$.
The function $x\mapsto\exp\{\mu^2x(x-1)/2\}$ is convex, so the comparison between sampling with and without replacement in \citet{hoeffding1963probability} gives
\begin{align*}
\ee_{S,S'\sim\dot\nu}\big[\exp\{\mu^2|S\cap S'|\}\big]
&=\ee_H\bigg[\exp\bigg\{\mu^2\binom H2\bigg\}\bigg]\\
&\leq\ee_B\bigg[\exp\bigg\{\mu^2\binom B2\bigg\}\bigg],
\end{align*}
where $B\sim\mathrm{Binomial}(k,k/m)$.
The case $k=1$ is trivial because the overlap is empty, so assume $k\geq2$.
Take $B_1,\ldots,B_k\stackrel{\mathrm{iid}}{\sim}\mathrm{Bernoulli}(k/m)$ and set $B:=\sum_{i=1}^k B_i$.
Since $\binom B2=\sum_{i<j}B_iB_j$, we have
\begin{equation}
    \label{eq:Bi}
\ee_B\bigg[\exp\bigg\{\mu^2\binom B2\bigg\}\bigg]
=\ee\bigg[\exp\bigg(\mu^2\sum_{i<j}B_iB_j\bigg)\bigg].
\end{equation}

We apply Finner's inequality (\Cref{lemma: Finner's inequality}) with $I=[k]$, $J=[k]_2$, and $T_{(i,j)}=\{i,j\}$.
Each coordinate $B_i$ appears in $k-1$ factors, so choosing $p_{ij}=k-1$ makes the sum of the reciprocal exponents equal to one at every coordinate.
Since $B_iB_j\sim\mathrm{Bernoulli}(k^2/m^2)$ for $i<j$, the inequality yields
\begin{equation}\label{eq:Finner inequality application}
\begin{aligned}
\ee\bigg[\exp\bigg(\mu^2\sum_{i<j}B_iB_j\bigg)\bigg]
&\leq\prod_{i<j}\ee\big[\exp\{\mu^2(k-1)B_iB_j\}\big]^{1/(k-1)}\\
&=\ee\big[\exp\{\mu^2(k-1)B_1B_2\}\big]^{k/2}\\
&=\bigg(1+\frac{k^2}{m^2}\big(\exp\{\mu^2(k-1)\}-1\big)\bigg)^{k/2}\\
&\leq\exp\bigg(\frac{k^3}{2m^2}\big(\exp(\mu^2k)-1\big)\bigg),
\end{aligned}
\end{equation}
where the last inequality uses $\log(1+x)\leq x$ for $x\geq0$ and $k-1\leq k$.
Thus, the likelihood-ratio second moment is bounded whenever $\frac{k^3}{m^2}(\exp(\mu^2k)-1)=O(1)$.
More explicitly, for any fixed constant $C>0$, the condition
\[
\mu\leq\sqrt{\frac1k\log\bigg(1+C\frac{m^2}{k^3}\bigg)}
\]
ensures that the second moment is at most $e^{C/2}$.
We verify this condition in the two regimes of \Cref{lemma: statistical lower bound for fixed-size clique problem}.

First, if $k^3m^{-2}\leq m^{-\delta}$ and $\mu\leq\sqrt{\delta k^{-1}\log m}$, then
\[
\frac{k^3}{2m^2}\big(\exp(\mu^2k)-1\big)
\leq\frac{m^{-\delta}}{2}(m^\delta-1)
\leq\frac12.
\]
Consequently, the second moment is at most $e^{1/2}$.

Second, suppose that $k^3m^{-2}=\Omega(1)$ and $\mu=O(mk^{-2})$.
Choose constants $c_0,C_0>0$ such that $k^3/m^2\geq c_0$ and $\mu\leq C_0m/k^2$ for all sufficiently large $m$, so that $\mu^2k\leq C_0^2/c_0$.
Using $e^x-1\leq xe^x$ for $x\geq0$, we obtain
\[
\frac{k^3}{2m^2}\big(\exp(\mu^2k)-1\big)
\leq\frac{\mu^2k^4}{2m^2}\exp(\mu^2k)
\leq\frac{C_0^2}{2}\exp(C_0^2/c_0),
\]
which is bounded independently of $m$ and proves the second assertion.

\subsection{\texorpdfstring{Proof of \Cref{lemma: computational lower bound for fixed-size clique problem}}{Proof of Lemma \getrefnumber{lemma: computational lower bound for fixed-size clique problem}}}\label{proof: computational lower bound for fixed-size clique problem}

Fix constants $\delta\in(0,0.1)$ and $\epsilon\in(0,\delta)$ as in the lemma, and set $D:=\lfloor\log^{1+\epsilon}m\rfloor$.
By \Cref{proposition: gaussian overlapping pairs}, the low-degree likelihood-ratio second moment is
\[
\bbE_{\mathbf Y\sim\bbQ}\big[L_{\dot\nu}^{\leq D}(\mathbf Y)^2\big]
=\ee_{S,S'\sim\dot\nu}\bigg[\sum_{d=0}^D\frac{\mu^{2d}}{d!}|S\cap S'|^d\bigg].
\]
As in \Cref{proof: statistical lower bound for fixed-size clique problem}, let $H:=|\dot V\cap\dot V'|$, so that $H$ has the hypergeometric distribution with parameters $(m,k,k)$ and $|S\cap S'|=\binom H2\leq H^2/2$.
For each integer $d\geq1$, the function $x\mapsto x^{2d}$ is convex on $[0,\infty)$, so the comparison with $B\sim\mathrm{Binomial}(k,k/m)$ gives $\ee_H[H^{2d}]\leq\ee_B[B^{2d}]$.
Combining these bounds with \Cref{lemma: moment bound of sub-Poissonian distribution}, due to \citet[Corollary 1]{ahle2022sharp}, yields
\begin{equation}\label{equation: computational lower bound moment bound}
\begin{aligned}
\bbE_{\mathbf Y\sim\bbQ}\big[L_{\dot\nu}^{\leq D}(\mathbf Y)^2\big]
&\leq1+\sum_{d=1}^D\frac{\mu^{2d}}{2^dd!}\bigg(\frac{2d}{\log(1+2dm/k^2)}\bigg)^{2d}\\
&=1+\sum_{d=1}^D\frac{2^d\mu^{2d}}{d!}\bigg(\frac{d}{\log(1+2dm/k^2)}\bigg)^{2d}.
\end{aligned}
\end{equation}
It remains to bound the series in the two regimes of the lemma.
\begin{itemize}
\item \textbf{Case 1.} $k^2/m\leq m^{-\delta}$.

Since $\log(1+2dm/k^2)\geq\log(m/k^2)$ for $d\geq1$, the series in \Cref{equation: computational lower bound moment bound} is at most $\sum_{d=1}^D T_d$, where
\[
T_d:=\frac{2^d\mu^{2d}}{d!}\bigg(\frac{d}{\log(m/k^2)}\bigg)^{2d}.
\]
For $1\leq d<D$, we have
\begin{align*}
\frac{T_{d+1}}{T_d}
&=\frac{2\mu^2(d+1)}{\log^2(m/k^2)}\bigg(1+\frac1d\bigg)^{2d}\\
&\leq\frac{2e^2\mu^2D}{\log^2(m/k^2)},
\end{align*}
where we used $(1+1/d)^d\leq e$.
Thus, if $\mu\leq\log(m/k^2)/(2e\sqrt D)$, then $T_{d+1}\leq T_d/2$ and
\[
\sum_{d=1}^D T_d\leq2T_1
=\frac{4\mu^2}{\log^2(m/k^2)}
\leq\frac{1}{e^2D}.
\]
The assumption $k^2/m\leq m^{-\delta}$ gives $\log(m/k^2)\geq\delta\log m$.
Since $D\leq\log^{1+\epsilon}m$, the required condition holds whenever $\mu\leq C(\log m)^{(1-\epsilon)/2}$ with $0<C\leq\delta/(2e)$.
This proves the first assertion.

\item \textbf{Case 2.} $k^2/m\geq\log^{1+\delta}m$.

Uniformly over $d\in[D]$, we have $0<dm/k^2\leq\log^{\epsilon-\delta}m\to0$, since $\epsilon<\delta$.
Using $\log(1+2x)\geq x$ for $0\leq x\leq1$, we obtain, for all sufficiently large $m$,
\[
\log\bigg(1+\frac{2dm}{k^2}\bigg)\geq\frac{dm}{k^2},\qquad d\in[D].
\]
Consequently, the series in \Cref{equation: computational lower bound moment bound} satisfies
\begin{align*}
\sum_{d=1}^D\frac{2^d\mu^{2d}}{d!}\bigg(\frac{d}{\log(1+2dm/k^2)}\bigg)^{2d}
&\leq\sum_{d=1}^D\frac{1}{d!}\bigg(\frac{2\mu^2k^4}{m^2}\bigg)^d\\
&\leq\exp\bigg(\frac{2\mu^2k^4}{m^2}\bigg)-1.
\end{align*}
This bound is $O(1)$ whenever $\mu=O(mk^{-2})$, proving the second assertion.
\end{itemize}
By \Cref{lemma: low-degree second moment method}, the bounded second moment rules out strong separation by polynomials of degree at most $D$.
Its interpretation as a lower bound against all polynomial-time tests is conditional on \Cref{conjecture: low-degree polynomial conjecture}.

\subsection{\texorpdfstring{Proof of \Cref{theorem: fixed-size clique problem}}{Proof of Theorem \getrefnumber{theorem: fixed-size clique problem}}}\label{proof: fixed-size clique problem}

For the fixed-size planted clique problem, $n=\binom m2$, $K=\binom k2$, and $|\dot{\mathcal C}|=\binom mk$.
We combine the upper bounds in \Cref{proposition: gaussian test upper bound} with \Cref{lemma: statistical lower bound for fixed-size clique problem,lemma: computational lower bound for fixed-size clique problem} to identify the thresholds at polynomial order.

\paragraph{Upper bounds.}
Applying \Cref{proposition: gaussian test upper bound} with prior $\dot\nu$ gives the following guarantees.
\begin{itemize}
\item Taking the target risk to be $\delta_m=1/\log m$, the outlier test satisfies $R_{\dot\nu}(\Delta_{\mathrm{out}})\to0$ whenever
\[
\mu\geq C\sqrt{\log m}
\]
for a sufficiently large constant $C>0$.

\item The average-test risk bound proved in \Cref{proof: gaussian test upper bound} gives
\[
R_{\dot\nu}(\Delta_{\mathrm{avg}})
\leq2\exp\bigg(-\frac{\mu^2K^2}{8n}\bigg).
\]
Since $\sqrt n/K\asymp m/k^2$, the risk tends to zero whenever $\mu=\omega(mk^{-2})$.

\item For any fixed $c\in(0,1)$, the scan test satisfies $R_{\dot\nu}(\Delta_{\mathrm{scan}})\to0$ whenever $k\leq(1-c)m$ and
\[
\mu\geq C\sqrt{\frac1k\log\frac{m}{k}}
\]
for a sufficiently large constant $C>0$ depending only on $c$.
This follows by taking $\delta_m=1/\log m$ and using $\log\binom mk\leq k\log(em/k)$, together with the assumption that $k$ grows at least polynomially in $m$.
\end{itemize}
The outlier and average tests run in $O(m^2)$ time.
In the regime $k\ll m^{2/3}$, the condition $k\leq(1-c)m$ holds eventually, and a direct scan over all $\binom mk$ candidates requires superpolynomial time because $k$ grows at least polynomially in $m$.

\paragraph{Statistical threshold.}
If $k\ll m^{2/3}$, then $k^3/m^2\leq m^{-\delta}$ eventually for some constant $\delta\in(0,0.1)$.
Thus, \Cref{lemma: statistical lower bound for fixed-size clique problem}(i) and \Cref{lemma: second moment method} rule out strong detection whenever $\mu\ll k^{-1/2}$, while the scan test achieves strong detection whenever $\mu\gg k^{-1/2}$.
If $k\gg m^{2/3}$, then \Cref{lemma: statistical lower bound for fixed-size clique problem}(ii) rules out strong detection whenever $\mu\ll mk^{-2}$, while the average test achieves strong detection whenever $\mu\gg mk^{-2}$.
These bounds give the two branches of $\mu_{\mathrm{stat}}^\ast$ in \Cref{theorem: fixed-size clique problem}.

\paragraph{Computational threshold.}
If $k\ll m^{1/2}$, then $k^2/m\leq m^{-\delta}$ eventually for some constant $\delta\in(0,0.1)$.
Choose any $\epsilon\in(0,\delta)$.
By \Cref{lemma: computational lower bound for fixed-size clique problem}(i), the low-degree likelihood-ratio second moment is bounded for degrees $D_m=\log^{1+\epsilon}m$ whenever $\mu\ll1$, while the outlier test achieves strong detection whenever $\mu\gg1$.
If $k\gg m^{1/2}$, then $k^2/m\geq\log^{1+\delta}m$ eventually for any fixed $\delta\in(0,0.1)$, so \Cref{lemma: computational lower bound for fixed-size clique problem}(ii) bounds the low-degree second moment whenever $\mu\ll mk^{-2}$.
The average test achieves strong detection whenever $\mu\gg mk^{-2}$.
Under \Cref{conjecture: low-degree polynomial conjecture}, these low-degree bounds give the corresponding lower bounds against all polynomial-time tests.
Since logarithmic factors are subpolynomial in $m$, \Cref{definition: statistical and computational thresholds} yields the stated two branches of $\mu_{\mathrm{comp}}^\ast$ and completes the proof.

\subsection{\texorpdfstring{Proof of \Cref{corollary: random-size clique problem}}{Proof of Corollary \getrefnumber{corollary: random-size clique problem}}}\label{proof: random-size clique problem}

In the random-size model, $\tilde{\mathcal C}$ is the collection of clique edge sets of arbitrary size in the complete graph with $m$ vertices.
For the latent vector $\tilde \bfA$, define the vertex set $ \tilde V:=\{l\in[m]\,:\,\mathds1(\tilde A_l\leq k)=1\}$ of the planted clique. We define $\tilde V'$ similarly for a random vector $\tilde \bfA'$. With a slight abuse of notation, we use $\tilde\nu$ for both the distribution of $\tilde V$ (and also of $\tilde V'$) and the induced distribution of $S$ on $\tilde{\mathcal C}$.

\paragraph{Conditioning and total variation.}
We first compare the total variation distance between the alternative distribution and the alternative induced by a conditional mixing measure.
Given a mixing measure $\tilde\nu$, consider an event $\mathcal E$ and define a new mixing measure $\tilde\nu_{\mathcal{E}}$ as the conditional distribution of $\tilde\nu$ given $\mathcal E$, i.e., for any $s\in \tilde{\mathcal C},$ $\tilde\nu_{\mathcal{E}}(s) = \tilde\nu(s \,|\, \mathcal E)$.
Here and below, conditioning on $\mathcal E$ refers to the joint construction of the latent structure and observations, and $\tilde\nu(\mathcal E)$ denotes the probability of the latent event.
In the density calculation below, we use the same symbols for probability laws and their Lebesgue densities.
The corresponding alternative distribution $\bbP_{\tilde\nu_{\mathcal{E}}}$ satisfies
$$
\bbP_{\tilde\nu_{\mathcal{E}}}(\cdot) = \ee_{S\sim \tilde\nu_{\mathcal{E}}}[\bbP_{S}(\cdot)] = \ee_{S\sim \tilde\nu}[\bbP_{S}(\cdot) \,|\, \mathcal E] : = \bbP_{\tilde\nu}(\cdot \,|\, \mathcal E).
$$
We compare $\bbP_{\tilde\nu}$ and $\bbP_{\tilde\nu_{\mathcal{E}}}$ in terms of total variation distance. Note that the distribution $\bbP_{\tilde\nu}$ can be decomposed as
\begin{align*}
	\bbP_{\tilde\nu}(\cdot) &= \bbP_{\tilde\nu}(\cdot \,|\, \mathcal E) \tilde\nu(\mathcal E) + \bbP_{\tilde\nu}(\cdot \,|\, \mathcal E^c) (1-\tilde\nu(\mathcal E))\\
&= \bbP_{\tilde\nu_{\mathcal{E}}}(\cdot) \tilde\nu(\mathcal E) + \bbP_{\tilde\nu}(\cdot \,|\, \mathcal E^c) \tilde\nu(\mathcal E^c).
\end{align*}
Therefore, we have
$$
\bbP_{\tilde\nu}(\cdot) - \bbP_{\tilde\nu_{\mathcal{E}}}(\cdot) = \tilde\nu(\mathcal E^c)\big(\bbP_{\tilde\nu}(\cdot \,|\, \mathcal E^c) - \bbP_{\tilde\nu_{\mathcal{E}}}(\cdot)\big),
$$
which implies that 
\begin{equation}\label{eq:tv.bound}
\begin{aligned}
	\mathrm{TV}(\bbP_{\tilde\nu},\bbP_{\tilde\nu_{\mathcal{E}}}) &= \frac{1}{2}\int |\bbP_{\tilde\nu}(\mathbf y) - \bbP_{\tilde\nu_{\mathcal{E}}}(\mathbf y)| d\mathbf y\\
	&= \tilde\nu(\mathcal E^c) \cdot \frac{1}{2}\int |\bbP_{\tilde\nu}(\mathbf y \,|\, \mathcal E^c) - \bbP_{\tilde\nu_{\mathcal{E}}}(\mathbf y)| d\mathbf y\\
	&\leq \tilde\nu(\mathcal E^c).
\end{aligned}
\end{equation}

\paragraph{Upper bounds.}
To apply the tests introduced in \Cref{proposition: gaussian test upper bound}, we need to address the randomness in the size of the clique. At a vanishing cost in total variation, we modify the mixing measure $\tilde\nu$ so that, under the alternative hypothesis, the number of clique vertices is at least $c\cdot k$ for some constant $c\in (0,1).$ To be precise, consider the event defined by $\underline {\mathcal E} = \underline {\mathcal E}(c) :=\{\tilde V\in 2^{[m]}: |\tilde V| \geq c\cdot k\},$ where $c\in (0,1)$ is a constant. Define a new mixing measure $\underline{\tilde\nu}$ as the conditional distribution of $\tilde\nu$ given $\underline {\mathcal E}$, i.e., for any $s\in \tilde{\mathcal C},$ $\underline{\tilde\nu}(s) = \tilde\nu(s \,|\, \underline {\mathcal E})$. For the mixing measures $\tilde\nu$ and $\underline{\tilde\nu}$, the corresponding risks $R_{\tilde\nu}(\Delta)$ and $R_{\underline{\tilde\nu}}(\Delta)$ of a test $\Delta$ satisfy
\begin{equation}\label{equation: risk lower bound via TV with conditioning for upper bound}
\begin{aligned}
	R_{\tilde\nu}(\Delta)&= \bbQ\{\Delta(\mathbf Y)=1\} + \bbP_{\tilde\nu}\{\Delta(\mathbf Y)=0\}\\
	&= \bbQ\{\Delta(\mathbf Y)=1\}+\bbP_{\underline{\tilde\nu}}\{\Delta(\mathbf Y)=0\}+\big(\bbP_{\tilde\nu}\{\Delta(\mathbf Y)=0\}-\bbP_{\underline{\tilde\nu}}\{\Delta(\mathbf Y)=0\} \big)\\
	&\leq R_{\underline{\tilde\nu}}(\Delta)+\mathrm{TV}(\bbP_{\tilde\nu},\bbP_{\underline{\tilde\nu}})\\
	&\leq R_{\underline{\tilde\nu}}(\Delta)+\tilde\nu(\underline {\mathcal E}^c),
\end{aligned}
\end{equation}
where the last step follows from \Cref{eq:tv.bound}.
Since $\mathds 1\{l\in \tilde V\}$ for $l\in[m]$ are i.i.d.\ $\mathrm{Bernoulli}(k/m)$ random variables and thus $|\tilde V|=\sum^m_{l=1}\mathds 1\{l\in \tilde V\}$ follows a binomial distribution with parameters $m$ and $k/m$, by the Chernoff bound, we have
$$
\tilde\nu(\underline {\mathcal E}^c) = \Pr\big\{|\tilde V| < c\cdot k\big\}\leq \exp\bigg(-\frac{(1-c)^2}{2} k\bigg),
$$
which converges to zero as long as $k\to \infty.$ Therefore, it suffices to show $R_{\underline{\tilde\nu}}(\Delta)\to 0$ to conclude that $R_{\tilde\nu}(\Delta)\to 0.$

When $ck$ and $Ck$ specify vertex counts below, we round them down and up, respectively; the Chernoff bounds remain valid under these roundings.
We apply the three tests in \Cref{proposition: gaussian test upper bound} with the modified mixing measure $\underline{\tilde\nu}$. Note that every $S$ in the support of $\underline{\tilde\nu}$ contains a clique on $c\cdot k$ vertices. Thus, the theoretical guarantees of the tests we analyzed in \Cref{proof: fixed-size clique problem} still hold with $k$ replaced by $c\cdot k$.
More precisely, the average test uses threshold $\mu\binom{\lfloor ck\rfloor}{2}/2$, and the scan searches all $\lfloor ck\rfloor$-vertex cliques.
The same error bounds hold for larger planted cliques because the average statistic has mean at least $\mu\binom{\lfloor ck\rfloor}{2}$ and the scan includes a fully signaled candidate.
Therefore, we have the following upper bounds on the critical signal strength $\mu^*$:
\begin{itemize}
	\item For some constant $C>0$, the outlier test achieves $R_{\underline{\tilde\nu}}(\Delta_{\textrm{out}})\to 0$ whenever
	$\mu \geq C \sqrt{\log m}.$
	\item The average test achieves $R_{\underline{\tilde\nu}}(\Delta_{\textrm{avg}})\to 0$ as long as $\mu =\omega(\frac{m}{k^2}).$
	\item For some constant $C>0$, the scan test achieves $R_{\underline{\tilde\nu}}(\Delta_{\textrm{scan}})\to 0$ if $k\leq(1-c_0)m$ for a fixed constant $c_0\in(0,1)$ and
	$$
	\mu \geq C \sqrt{\frac{1}{k}\log\frac{m}{k}}.
	$$
\end{itemize}
The restriction on $k$ for the scan bound holds in the regime $k\ll m^{2/3}$ where that bound is needed.

\paragraph{Statistical lower bound.}
We now derive a statistical lower bound on the critical signal strength $\mu^*$. We first consider the case $k\leq \kappa m$ for some constant $\kappa\in (0,1)$ and address the case $k>\kappa m$ separately at the end of the statistical lower-bound argument.

When $k\leq \kappa m$, we again condition the mixing measure, this time to obtain a lower bound on the risk.
Consider the event $\overline {\mathcal E} = \overline{\mathcal E}(C) :=\{\tilde V\in2^{[m]}: |\tilde V| \leq C\cdot k\},$ where $C\in(1,1/\kappa)$ is a constant. Define a new mixing measure $\overline{\tilde\nu}$ as the conditional distribution of $\tilde\nu$ given $\overline {\mathcal E}$, i.e., for any $s\in \tilde{\mathcal C},$ $\overline{\tilde\nu}(s) = \tilde\nu(s \,|\, \overline {\mathcal E})$. By Chernoff's bound, we have
$$
\tilde\nu(\overline {\mathcal E}^c) = \Pr\big\{|\tilde V| > C\cdot k\big\}\leq \exp\bigg(-\frac{(C-1)^2}{2C} k\bigg) \to 0
$$
 as $k\to \infty.$ 

 For any test $\Delta$, by \Cref{eq:tv.bound} and \Cref{equation: risk lower bound via BH}, together with the argument that follows, we have
\begin{equation}\label{equation: risk lower bound via TV with conditioning}
\begin{aligned}
	R_{\tilde\nu}(\Delta)&\geq 1-\mathrm{TV}(\bbQ,\bbP_{\tilde\nu})\\
	&\geq 1-\big(\mathrm{TV}(\bbQ,\bbP_{\overline{\tilde\nu}})+\mathrm{TV}(\bbP_{\tilde\nu},\bbP_{\overline{\tilde\nu}})\big)\\
	&\geq 1-\sqrt{1-\frac{1}{\bbE_{\mathbf Y\sim \bbQ_n}\Big[\Big(\frac{d\bbP_{\overline{\tilde\nu}}(\bfY)}{d\bbQ(\bfY)}\Big)^2\Big]}}-\tilde\nu(\overline{\mathcal E}^c).
\end{aligned}
\end{equation}
Since $\tilde\nu(\overline{\mathcal E}^c) \to  0$ and  $\bbE_{\mathbf Y\sim \bbQ_n}\Big[\Big(\frac{d\bbP_{\overline{\tilde\nu}}(\bfY)}{d\bbQ(\bfY)}\Big)^2\Big] > 1$, to keep $R_{\tilde\nu}(\Delta)$ bounded away from zero asymptotically, it suffices to establish that
\begin{align*}
\bbE_{\mathbf Y\sim \bbQ_n}\Big[\Big(\frac{d\bbP_{\overline{\tilde\nu}}(\bfY)}{d\bbQ(\bfY)}\Big)^2\Big] &= \ee_{S,S'\sim \overline{\tilde\nu}}[\exp\{\mu^2 |S\cap S'|\}]\\
&=\ee_{\tilde V,\tilde V'\sim \overline{\tilde\nu}}\bigg[\exp\bigg\{\mu^2\binom{|\tilde V\cap\tilde V'|}{2}\bigg\}\bigg]=O(1).
\end{align*}
We claim that, under $\overline{\tilde\nu}$, $|\tilde V\cap\tilde V'|$  is stochastically dominated by a hypergeometric distribution with parameters $(m,Ck,Ck).$ The following coupling proves the claim: by construction, $\tilde V$ and $\tilde V'$ are independent samples from $\overline{\tilde\nu}.$ Add $Ck-|\tilde V|$ (resp.\ $Ck-|\tilde V'|$) elements to $\tilde V$ (resp.\ $\tilde V'$), drawn uniformly at random from $[m]\setminus\tilde V$ (resp.\ $[m]\setminus\tilde V'$), thus producing a new set $\dot V$ (resp.\ $\dot V'$) of cardinality exactly $Ck$. By the permutation invariance of  $\overline{\tilde\nu}$, $\dot V$ and $\dot V'$ are independent and uniformly distributed over the set $\big\{a\in 2^{[m]}: \sum_{l=1}^m \mathds1\{l\in a\} = Ck\big\}$.
Under this coupling, we have $|\tilde V\cap\tilde V'| \leq |\dot V\cap\dot V'|$ almost surely, which implies that
\begin{equation}
\label{eq:stoc.dominance}
\Pr\{|\tilde V\cap\tilde V'| \geq t\} \leq \Pr\{|\dot V\cap\dot V'| \geq t\}, \quad \forall t \geq 0.
\end{equation}
	Therefore,  we have
	\begin{equation}\label{equation: truncated likelihood ratio second moment}
\begin{aligned}
	\ee_{S,S'\sim \overline{\tilde\nu}}[\exp\{\mu^2 |S\cap S'|\}] &= \ee_{\tilde V,\tilde V'}\bigg[\exp\bigg\{\mu^2\binom{|\tilde V\cap\tilde V'|}{2}\bigg\}\bigg]\\
&\stackrel{(a)}\leq \ee_{\dot V,\dot V'}\bigg[\exp\bigg\{\mu^2\binom{|\dot V\cap \dot V'|}{2}\bigg\}\bigg]\\
&\stackrel{(b)}\leq \ee_B\bigg[\exp\bigg\{\mu^2\binom B2\bigg\}\bigg]\\
&\stackrel{(c)}\leq \exp\bigg(\frac{(Ck)^3}{2m^2}\big(\exp(\mu^2Ck)-1\big)\bigg),
\end{aligned}
\end{equation}
where $ B\sim \mathrm{Binomial}(Ck,Ck/m)$.
Above, step (a) follows from the stochastic dominance relation in \Cref{eq:stoc.dominance} and step (b) from the fact that the hypergeometric distribution with parameters $(m,Ck,Ck)$ is dominated in the convex ordering by the binomial distribution with parameters $(Ck, Ck/m)$ \citep{hoeffding1963probability}, applied to the convex function $x\mapsto\exp\{\mu^2x(x-1)/2\}$.
Step (c) is \Cref{eq:Finner inequality application} with $k$ replaced by $Ck$.
We can therefore apply the same analysis in \Cref{proof: statistical lower bound for fixed-size clique problem}, replacing $k$ by $Ck,$ to obtain the same statistical lower bound.

When $k>\kappa m$, we can directly bound the second moment without truncation. Since $|S\cap S'|\leq \binom{m}{2}$ almost surely, we have
$$
\mathbb E_{\mathbf Y\sim\mathbb Q_n}[L_{\tilde\nu}(\mathbf Y)^2]
=\mathbb E_{S,S'\sim\tilde\nu}\big[e^{\mu^2|S\cap S'|}\big]
\leq \exp\bigg(\mu^2\binom{m}{2}\bigg).
$$
Since $k>\kappa m$, the right-hand side is bounded whenever $\mu=O(m/k^2)$, yielding the same statistical lower bound as in the case $k\leq \kappa m$.

\paragraph{Computational lower bound.}
To derive a computational lower bound, we present two approaches: (i) working directly with the original mixing measure (without truncation), and (ii) working with a truncated mixing measure, as in the statistical lower bound analysis.

\textbf{(i) Original mixing measure.} 
 The starting point is the following bound for any integer degree cutoff $D\geq0$:
\begin{equation}\label{equation: low degree likelihood ratio random size}
\begin{aligned}
\ee_{\bfY\sim \bbQ}[L_{\tilde\nu}^{\leq D}(\mathbf Y)^2]&=\ee_{S,S'\sim \tilde\nu}\bigg[\sum_{d=0}^D \frac{\mu^{2d}}{d!}|S\cap S'|^d\bigg]\\
&=\ee_{\tilde V,\tilde V'\sim \tilde\nu}\bigg[\sum_{d=0}^D \frac{\mu^{2d}}{d!}\binom{|\tilde V\cap\tilde V'|}{2}^d\bigg]\\
&\leq \sum_{d=0}^D \frac{\mu^{2d}}{2^d d!}\ee_{\tilde V,\tilde V'\sim \tilde\nu}\Big[|\tilde V\cap\tilde V'|^{2d}\Big]\\
&= 1+ \sum_{d=1}^D \frac{\mu^{2d}}{2^d d!}\ee_{\tilde V,\tilde V'\sim \tilde\nu}\Big[|\tilde V\cap\tilde V'|^{2d}\Big].
\end{aligned}
\end{equation}
Letting $D=D_m=\lfloor\log^{1+\epsilon}m\rfloor$ with $0<\epsilon<\delta<0.1$ as in \Cref{lemma: computational lower bound for fixed-size clique problem},
\Cref{conjecture: low-degree polynomial conjecture} implies that a sufficient condition for computational indistinguishability is the boundedness of the last term above.

Observe that $\mathds 1\{i\in\tilde V\},\mathds 1\{i\in\tilde V'\}$ are i.i.d.\ $\mathrm{Bernoulli}(k/m)$ random variables for $i\in[m]$, and thus the vertex overlap $|\tilde V\cap\tilde V'|=\sum^{m}_i\mathds 1\{i\in\tilde V\}\mathds 1\{i\in\tilde V'\}$ follows $\mathrm{Binomial}(m,(k/m)^2),$ whose mean is $k^2/m$. Applying the sub-Poissonian moment bound (see \Cref{lemma: moment bound of sub-Poissonian distribution}) gives
\begin{align*}
	1+\sum_{d=1}^D \frac{\mu^{2d}}{2^d d!}\ee_{\tilde V,\tilde V'\sim \tilde\nu}\Big[|\tilde V\cap\tilde V'|^{2d}\Big]&\leq 1 + \sum_{d=1}^D \frac{\mu^{2d}}{2^d d!} \bigg(\frac{2d}{\log(1+2dm/k^2)}\bigg)^{2d}\\*
	&=1+ \sum_{d=1}^D \frac{2^d\mu^{2d}}{d!} \bigg(\frac{d}{\log(1+2dm/k^2)}\bigg)^{2d}.
\end{align*}
Note that the last equation is identical to the bound for the fixed-size case (see \Cref{equation: computational lower bound moment bound}). Repeating that analysis proves the bound for the original mixing measure.

\textbf{(ii) Truncated mixing measure $\overline{\tilde \nu}$.} We now use the truncation approach from the statistical lower-bound analysis to obtain the same computational lower bound. In what follows, we assume that $k\leq \kappa m$ for some constant $\kappa\in(0,1)$. The case $k>\kappa m$ can be handled as in the statistical lower bound analysis, using the fact that $|S\cap S'|\leq\binom{m}{2}$ almost surely and that the second moment of the low-degree likelihood ratio is bounded above by that of the full likelihood ratio.

We first compare the risks:
\begin{align*}
|R_{\tilde\nu}(\Delta)-R_{\overline{\tilde\nu}}(\Delta)| &= \big|\bbP_{\tilde\nu}\{\Delta(\mathbf Y)=0\}-\bbP_{\overline{\tilde\nu}}\{\Delta(\mathbf Y)=0\}\big|\\
&\leq \mathrm{TV}(\bbP_{\tilde\nu},\bbP_{\overline{\tilde\nu}})\\
&\leq \tilde\nu(\overline {\mathcal E}^c).
\end{align*}
This bound can be established in the same way as the bound in \Cref{equation: risk lower bound via TV with conditioning for upper bound}. Using this bound, and as long as $\tilde\nu(\overline {\mathcal E}^c)\to 0,$ we immediately obtain that any polynomial-time test under the truncated mixing measure $\overline {\tilde\nu}$ achieves strong detection if and only if it also achieves strong detection under the original mixing measure $\tilde\nu$.
Thus, under the low-degree polynomial conjecture (\Cref{conjecture: low-degree polynomial conjecture}), it suffices to find a scaling of $\mu$ such that $\ee_{\bfY\sim\bbQ}[L_{\overline{\tilde\nu}}^{\leq D}(\mathbf Y)^2]$ remains bounded for $D=D_m=\lfloor\log^{1+\epsilon}m\rfloor$ with $\epsilon$ as above, as this will imply that $R_{\overline{\tilde\nu}}(\Delta)$ will not vanish for any polynomial-time test $\Delta$ and, by the equivalence above, the same will hold for $R_{\tilde\nu}(\Delta)$ as well.

Repeating the steps leading to \Cref{equation: low degree likelihood ratio random size} with the truncated mixing measure $\overline {\tilde\nu}$, we obtain
\begin{align*}
\ee_{\bfY\sim\bbQ}[L_{\overline{\tilde\nu}}^{\leq D}(\mathbf Y)^2] &= \ee_{S,S'\sim \overline{\tilde\nu}}\bigg[\sum_{d=0}^D \frac{\mu^{2d}}{d!}|S\cap S'|^d\bigg]\\
&\leq 1+ \sum_{d=1}^D \frac{\mu^{2d}}{2^d d!}\ee_{\tilde V,\tilde V'\sim \overline{\tilde\nu}}\Big[|\tilde V \cap \tilde V'|^{2d}\Big],
\end{align*}
where $\tilde V$ and $\tilde V'$ are i.i.d.\ samples from the distribution on $2^{[m]}$ corresponding to the conditional distribution of $m$ i.i.d.\ $\mathrm{Bernoulli}(k/m)$ given that their sum is no larger than $ C\cdot k$, for some arbitrary constant $C\in(1,1/\kappa)$. %
As in the coupling leading to \Cref{equation: truncated likelihood ratio second moment}, enlarge $\tilde V$ and $\tilde V'$ independently to supersets $\dot V$ and $\dot V'$ of size $\lceil Ck\rceil$ by adding vertices uniformly without replacement from their respective complements in $[m]$.
By permutation invariance, $\dot V$ and $\dot V'$ are independent and uniformly distributed over $\{v\subseteq[m]:|v|=\lceil Ck\rceil\}$.
The coupling gives $|\tilde V\cap\tilde V'|\leq|\dot V\cap\dot V'|$ almost surely, so the former overlap is stochastically dominated by the latter, which has the hypergeometric distribution with parameters $(m,\lceil Ck\rceil,\lceil Ck\rceil)$.
Now we repeat the same analysis in \Cref{proof: computational lower bound for fixed-size clique problem} with $k$ replaced by $\lceil Ck\rceil$.
This substitution changes only the constants, not the asymptotic rate, and gives the same computational lower bound as above, completing the proof.

\subsection{\texorpdfstring{Proof of \Cref{theorem: multiple cliques problem}}{Proof of Theorem \getrefnumber{theorem: multiple cliques problem}}}\label{proof: multiple cliques problem}
In this section, we study the upper bounds provided by three tests for the multiple planted clique problem and prove \Cref{theorem: multiple cliques problem} by establishing matching statistical and low-degree lower bounds.
The computational interpretation of the latter uses \Cref{conjecture: low-degree polynomial conjecture}.
A formal statement of the lower-bound results we prove is as follows:
\begin{itemize}
	\item (Statistical lower bound) Fix any constant $\delta\in(0,0.1)$.
For either prior distribution $\dot\nu$ or $\tilde\nu$, strong detection is impossible under either of the following conditions: (i) $\gamma^2k^3m^{-2}\leq m^{-\delta}$ and $\mu\leq C_\delta\sqrt{k^{-1}\log m}$ for a sufficiently small constant $C_\delta>0$; or (ii) $\gamma^2k^3m^{-2}=\Omega(1)$ and $\mu=O\big(\gamma^{-1}k^{-2}m\big)$.
	\item (Computational lower bound) Fix any constants $\delta\in(0,0.1)$ and $\epsilon\in(0,\delta)$, and let $D_m=\lfloor\log^{1+\epsilon}m\rfloor$.
For either prior distribution $\dot\nu$ or $\tilde\nu$, no polynomial of degree at most $D_m$ achieves strong separation under either of the following conditions: (i) $\gamma k^2m^{-1}\leq m^{-\delta}$ and $\mu\leq C(\log m)^{(1-\epsilon)/2}$ for a sufficiently small constant $C>0$; or (ii) $\gamma k^2m^{-1}\geq\log^{1+\delta}m$ and $\mu=O\big(\gamma^{-1}k^{-2}m\big)$.
\end{itemize}

We first establish the claimed lower bounds under the restriction $\gamma k \leq \kappa m$ for an arbitrary fixed constant $\kappa\in(0,1)$, treating the fixed-size case and the random-size case in \Cref{proof: multiple fixed-size cliques problem} and \Cref{proof: random-size multiple cliques}, respectively.
We then remove this restriction by the grouping arguments below.

Throughout the following analysis, we assume that $k=\omega(\log m)$, which is automatically satisfied in the regime of interest where $k$ is polynomial in $m$.
The case $\gamma=1$ follows from \Cref{theorem: fixed-size clique problem,corollary: random-size clique problem}, so we assume $\gamma\geq2$ below.
When $ck$ or $Ck$ denotes a fixed vertex count, we round it down; this does not affect the asymptotic rates since $k\to\infty$.

\subsubsection{Fixed-size case}\label{proof: multiple fixed-size cliques problem}
We first consider the case where $\dot{\mathcal C}$ consists of all subgraphs that are the union of $\gamma$ disjoint cliques of size $k$. In this case, we have $n=\binom m2$, $K=\gamma\cdot\binom k2$ and $|\dot{\mathcal C}|=\frac{1}{\gamma!}\binom m k \cdot \binom{m-k}k \cdots \binom{m-(\gamma-1)k}k.$

\paragraph{Upper bounds.}
We analyze the upper bounds on the critical signal strength $\mu^*$ obtained by the three tests in \Cref{proposition: gaussian test upper bound}:
\begin{itemize}
	\item For some constant $C>0$, the outlier test achieves $R_{\dot \nu}(\Delta_{\textrm{out}})\to 0$ whenever
	$\mu \geq C \sqrt{\log m}.$
	\item The average test achieves $R_{\dot \nu}(\Delta_{\textrm{avg}})\to 0$ when $\mu =\omega(\frac{m}{\gamma k^2}).$
	\item Note that $|\dot{\mathcal C}|\leq \binom mk^\gamma$, and thus $\log|\dot{\mathcal C}|\leq \gamma k\log\frac{em}{k}$.
Therefore, for some constant $C>0$, the scan test achieves $R_{\dot \nu}(\Delta_{\textrm{scan}})\to 0$ whenever
	$$
	\mu \geq C \sqrt{\frac{1}{k}\log\frac{m}{k}}.
	$$
\end{itemize}
These guarantees follow by taking $\delta_m=1/\log m$ in \Cref{proposition: gaussian test upper bound}, except that for the average test we use its risk bound directly, as in \Cref{proof: fixed-size clique problem}.
Since $\gamma\geq2$ and $\gamma k\leq m$, we have $k\leq m/2$, so $\log(em/k)\asymp\log(m/k)$.

\paragraph{Statistical lower bound.} Fix $\kappa\in(0,1)$, and assume that $\gamma k \leq \kappa m$. Denote the $\gamma$ disjoint cliques in $S$ as $S_1,...,S_\gamma,$ and let $\dot V_1,\dots,\dot V_\gamma$ be their corresponding vertex membership sets. Formally, for $a \in [\gamma]$, define $\dot V_a:=\{l\in[m]\,:\,(a-1)k<\dot A_l\leq ak\}$ and let $S_a=\{(i,j)\in[m]_2: i,j\in \dot V_a\}.$ Then, $S = \bigcup_{a \in [\gamma]} S_a$. Let $\dot\bfA'$ be an independent copy of $\dot\bfA$, and define $S'_1,...,S'_\gamma$,  $S'$ and $\dot V'_1,...,\dot V'_\gamma$ similarly. Then,
\[
|S\cap S'| = \sum_{a,b} |S_a\cap S'_b|=\sum_{a,b} \binom{|\dot V_a\cap \dot V'_b|}{2}.
\]
For brevity, write $f(x):=\exp\{\mu^2x(x-1)/2\}$.
This function is nondecreasing on the nonnegative integers and convex on $\mathbb R$.
Then
\begin{equation}\label{equation: conditional expectation argument for multiple fixed size cliques}
\begin{aligned}
	\ee[\exp\{\mu^2 |S\cap S'|\}] &= \ee\bigg[\exp\bigg\{\mu^2 \sum_{a=1}^\gamma \sum_{b=1}^\gamma \binom{|\dot V_a\cap \dot V'_b|}{2}\bigg\}\bigg] \\
	&= \ee\bigg[\prod_{a=1}^\gamma \prod_{b=1}^\gamma \exp\bigg\{\mu^2 \binom{|\dot V_a\cap \dot V'_b|}{2}\bigg\}\bigg] \\
    & = \ee\bigg[\prod_{a=1}^\gamma \prod_{b=1}^\gamma f(|\dot V_a\cap \dot V'_b|)\bigg]  \\
	&=\ee\bigg[\ee\bigg\{\prod_{a=1}^\gamma \prod_{b=1}^\gamma  f\big(|\dot V_a\cap \dot V'_b|\big)\,\bigg|\,\dot V_1,\dots,\dot V_\gamma,\dot V'_1,\dots,\dot V'_{\gamma-1}\bigg\}\bigg]\\
	&= \ee\bigg[\prod_{a=1}^\gamma \prod_{b=1}^{\gamma-1}f\big(|\dot V_a\cap \dot V'_b|\big) \ee\bigg\{\prod_{a=1}^\gamma f\big(|\dot V_a\cap \dot V'_\gamma|\big)\,\bigg|\,\dot V_1,\dots,\dot V_\gamma,\dot V'_1,\dots,\dot V'_{\gamma-1}\bigg\}\bigg].
\end{aligned}
\end{equation}
Conditionally on $\dot V_1,\dots,\dot V_\gamma,\dot V'_1,\dots,\dot V'_{\gamma-1},$ the random variables $\big\{|\dot V_a\cap \dot V'_\gamma|\big\}_{a=1}^\gamma$ are negatively associated. This is because $\dot V'_\gamma$ is a uniformly random subset of size $k$ from the set $\mathcal{V}$ of $m-(\gamma-1)k$ nodes not contained in any of the first $\gamma -1$ node selections $\dot V'_1, \ldots, \dot V'_{\gamma- 1}$; that is, $\mathcal{V} = [m] \setminus \cup_{a=1}^{\gamma-1} \dot V'_a$. In turn, the set  $\mathcal{V}$ contains up to $\gamma$ disjoint (possibly empty) sets of the form $\mathcal{V} \cap \dot V_a$, for $a =1, \ldots, \gamma$. As a result, the random integer-valued vector
\[
\bigg( \underbrace{|\dot V_1 \cap \dot V'_\gamma|, \; \dots, \;|\dot V_\gamma \cap \dot V'_\gamma|}_{\gamma \text{ times}}, \; \Big| \left(\mathcal{V}  \setminus \cup_{a=1}^\gamma \dot V_a\right) \cap \dot V'_\gamma \Big| \bigg)
\]
has a multivariate hypergeometric distribution with parameters $k$ and 
\[
\bigg( \underbrace{|\mathcal{V} \cap \dot V_1 |, \; \dots, \;| \mathcal{V} \cap \dot V_\gamma|}_{\gamma \text{ times}}, \; \Big|\mathcal{V} \setminus \cup_{a=1}^\gamma \dot V_a \Big| \bigg)
\]
describing the outcome of sampling without replacement $k$ out of $m - (\gamma -1 ) k$ nodes, split into $\gamma + 1$ groups  $\dot V_1 \cap \mathcal{V}, \; \dots, \; \dot V_\gamma \cap \mathcal{V}, \; \mathcal{V} \setminus \cup_{a=1}^\gamma \dot V_a$, with each group  $\dot V_a \cap \mathcal{V}$, for $a=1,\ldots,\gamma$, possibly empty. Thus, as coordinates of a multivariate hypergeometric distribution, $\big\{|\dot V_a\cap \dot V'_\gamma|\big\}_{a=1}^\gamma$ are conditionally negatively associated.

Therefore,
\begin{equation}\label{equation: negative association argument for multiple fixed size cliques}
\begin{aligned}
	&\ee\bigg[\prod_{a=1}^\gamma \prod_{b=1}^{\gamma-1}f\big(|\dot V_a\cap \dot V'_b|\big) \ee\bigg\{\prod_{a=1}^\gamma f\big(|\dot V_a\cap \dot V'_\gamma|\big)\,\bigg|\,\dot V_1,\dots,\dot V_\gamma,\dot V'_1,\dots,\dot V'_{\gamma-1}\bigg\}\bigg]\\
	\leq &\ee\bigg[\prod_{a=1}^\gamma \prod_{b=1}^{\gamma-1}f\big(|\dot V_a\cap \dot V'_b|\big) \prod_{a=1}^\gamma \ee\bigg\{f\big(|\dot V_a\cap \dot V'_\gamma|\big)\,\bigg|\,\dot V_1,\dots,\dot V_\gamma,\dot V'_1,\dots,\dot V'_{\gamma-1}\bigg\}\bigg].
\end{aligned}
\end{equation}
The same argument shows that, for each $b=1,\dots,\gamma-1$, the random variables $\big\{|\dot V_a\cap\dot V'_b|\big\}_{a=1}^\gamma$ are negatively associated conditionally on $\dot V_1,\dots,\dot V_\gamma,\dot V'_1,\dots,\dot V'_{b-1}$.

Next,  for any $a,b=1,\dots,\gamma,$ it holds that
\begin{align}\label{equation: bound for conditional expectation of f for multiple fixed size cliques}
	&\ee\big[f\big(|\dot V_a\cap \dot V'_b|\big)\,\big|\,\dot V_1,\dots,\dot V_\gamma,\dot V'_1,\dots,\dot V'_{b-1}\big] \leq \ee[f(B)]
\end{align}
where $B\sim \mathrm{Binomial}(k,p_{\gamma,k}),$ and $p_{\gamma,k}:={k}/\{{m-(\gamma-1)k}\}$.
We prove this inequality below.
Using it, the right-hand side of \Cref{equation: negative association argument for multiple fixed size cliques} is bounded above by
\begin{align*}
    \ee\bigg[\prod_{a=1}^\gamma \prod_{b=1}^{\gamma-1}f\big(|\dot V_a\cap \dot V'_b|\big) \prod_{a=1}^\gamma \ee[f(B)]\bigg]
	= \ee[f(B)]^\gamma \cdot \ee\bigg[\prod_{a=1}^\gamma \prod_{b=1}^{\gamma-1}f\big(|\dot V_a\cap \dot V'_b|\big)\bigg].
\end{align*}
Plugging this bound into \Cref{equation: conditional expectation argument for multiple fixed size cliques} we obtain
\begin{equation}\label{eq:intermediate}
\begin{aligned}
&\ee[\exp\{\mu^2 |S\cap S'|\}]\\  
\leq &\ee[f(B)]^\gamma \cdot \ee\bigg[\prod_{a=1}^\gamma \prod_{b=1}^{\gamma-1}f\big(|\dot V_a\cap \dot V'_b|\big)\bigg] \\
 = &\ee[f(B)]^\gamma\cdot\ee\bigg[\prod_{a=1}^\gamma \prod_{b=1}^{\gamma-2}f\big(|\dot V_a\cap \dot V'_b|\big) \ee\bigg\{\prod_{a=1}^\gamma f\big(|\dot V_a\cap \dot V'_{\gamma-1}|\big)\,\bigg|\,\dot V_1,\dots,\dot V_\gamma,\dot V'_1,\dots,\dot V'_{\gamma-2}\bigg\}\bigg].
\end{aligned}
\end{equation}
By repeating the same arguments that led to \Cref{equation: negative association argument for multiple fixed size cliques}, 
\[
\begin{aligned}
&\ee\bigg[\prod_{a=1}^\gamma \prod_{b=1}^{\gamma-2}f\big(|\dot V_a\cap \dot V'_b|\big) \ee\bigg\{\prod_{a=1}^\gamma f\big(|\dot V_a\cap \dot V'_{\gamma-1}|\big)\,\bigg|\,\dot V_1,\dots,\dot V_\gamma,\dot V'_1,\dots,\dot V'_{\gamma-2}\bigg\}\bigg]\\
&\leq\ee[f(B)]^\gamma \ee\bigg[\prod_{a=1}^\gamma \prod_{b=1}^{\gamma-2}f\big(|\dot V_a\cap \dot V'_b|\big)\bigg].
\end{aligned}
\]
Plugging this bound into the right-hand side of \Cref{eq:intermediate},
\[
\ee[\exp\{\mu^2 |S\cap S'|\}] \leq \ee[f(B)]^{2\gamma} \cdot \ee\bigg[\prod_{a=1}^\gamma \prod_{b=1}^{\gamma-2}f\big(|\dot V_a\cap \dot V'_b|\big)\bigg].
\]
Iterating over $b = \gamma -2,\ldots, 1$ finally gives
\begin{equation}\label{equation: iterating conditional expectation argument for multiple fixed size cliques}
\ee[\exp\{\mu^2 |S\cap S'|\}] \leq \ee[f(B)]^{\gamma^2} = \ee\bigg[\exp\bigg(\mu^2 \binom B2\bigg)\bigg]^{\gamma^2}.
\end{equation}

\noindent {\bf Proof of \Cref{equation: bound for conditional expectation of f for multiple fixed size cliques}.}  %
Throughout, we will condition on $\dot V_1,\dots,\dot V_\gamma,\dot V'_1,\dots,\dot V'_{b-1}$. Then, $\dot V'_b$ is a uniformly random subset of size $k$ drawn without replacement
from the $M_b:=m-(b-1)k$ nodes not contained in any of the previously selected sets $\dot V'_1,\ldots,\dot V'_{b-1}$. Thus
$$
|\dot V_a\cap \dot V'_b|\,\big|\,\dot V_a,\dot V'_1,\dots,\dot V'_{b-1} \sim \mathrm{Hypergeometric}\bigg(M_b,\bigg|\dot V_a\setminus\bigcup_{l=1}^{b-1} \dot V'_{l}\bigg|,k\bigg).
$$
Since $\big|\dot V_a\setminus\bigcup_{l=1}^{b-1} \dot V'_{l}\big|\leq k$ and the hypergeometric distribution is stochastically increasing in the number of successes, %
$$
|\dot V_a\cap \dot V'_b|\,\big|\,\dot V_a,\dot V'_1,\dots,\dot V'_{b-1} \preceq_{\mathrm{st}} H,
$$
where $ H\sim \mathrm{Hypergeometric}(M_b,k,k)$
almost surely. In turn, this implies that 
$$
\ee\big[f(|\dot V_a\cap \dot V'_b|)\,\big|\,\dot V_a,\dot V'_1,\dots,\dot V'_{b-1}\big] \leq \ee[f(H)],
$$
since $ f(x)=\exp\!\big(\mu^2\binom{x}{2}\big)$ is nondecreasing on the nonnegative integers.
Furthermore, as the hypergeometric distribution is dominated by the corresponding binomial distribution in the convex ordering \citep{hoeffding1963probability}, we have
\begin{align*}
\mathrm{Hypergeometric}(M_b,k,k) &\preceq_{\mathrm{cx}} \mathrm{Binomial}\bigg(k,\frac{k}{M_b}\bigg)\\*
&\preceq_{\mathrm{st}} \mathrm{Binomial}\big(k,p_{\gamma,k}\big),
\end{align*}
where $p_{\gamma,k}:=k/(m-(\gamma-1)k)$ and the last stochastic ordering follows from the fact that $M_b\ge m-(\gamma-1)k.$
Therefore, we have $\ee[f(H)] \leq \ee[f(B)]$, where $B\sim \mathrm{Binomial}(k,p_{\gamma,k}).$ 
Since $f$ is nondecreasing on the nonnegative integers and convex, we obtain
$$
\mathbb E\!\left[f(|\dot V_a\cap \dot V'_b|)\,\middle|\,\dot V_1,\dots,\dot V_\gamma,\dot V'_1,\dots,\dot V'_{b-1}\right]
\le \mathbb E[f(B)],
$$
which proves \Cref{equation: bound for conditional expectation of f for multiple fixed size cliques}.

Getting back to \Cref{equation: iterating conditional expectation argument for multiple fixed size cliques}, 
by Finner's inequality (see the arguments leading to \Cref{eq:Finner inequality application} for details), 
\begin{equation}\label{equation: Finner inequality application for multiple fixed size cliques}
\begin{aligned}
	\ee\bigg[\exp\bigg\{\mu^2 \binom B2\bigg\}\bigg]\leq \exp\bigg\{\frac{kp_{\gamma,k}^2}{2}\big(\exp(\mu^2k)-1\big)\bigg\},
\end{aligned}
\end{equation}
which implies
\begin{align}\label{equation: final statistical bound for multiple fixed size cliques}
\ee[\exp\{\mu^2 |S\cap S'|\}]\leq \ee\bigg[\exp\bigg\{\mu^2 \binom B2\bigg\}\bigg]^{\gamma^2} \leq \exp\bigg\{\frac{\gamma^2 k^3}{2(1-\kappa)^2m^2}\big(\exp(\mu^2k)-1\big)\bigg\},
\end{align}
provided that $m-(\gamma-1)k\geq (1-\kappa)m$.
Thus, $\ee[\exp\{\mu^2 |S\cap S'|\}] = O(1)$ if $\frac{\gamma^2 k^3}{2(1-\kappa)^2m^2}\big(\exp(\mu^2k)-1\big)$ is at most $C$ for a fixed constant $C>0$. Equivalently,
$$
\mu \leq \sqrt{\frac{1}{k}\log\bigg(1+C\frac{2(1-\kappa)^2m^2}{\gamma^2 k^3}\bigg)}.
$$
In the regime $\gamma^2k^3=\omega(m^2)$, we have $C\frac{(1-\kappa)^2m^2}{\gamma^2k^3}\leq\log(1+C\frac{2(1-\kappa)^2m^2}{\gamma^2k^3})$ for all sufficiently large $m$. Therefore, the condition in the above equation holds as long as $\mu \leq {\sqrt C}(1-\kappa)\frac{m}{\gamma k^2}.$
The same rate holds under the weaker assumption $\gamma^2k^3/m^2=\Omega(1)$ stated above.
Indeed, choose fixed constants $c_0,C_0>0$ such that $\gamma^2k^3/m^2\geq c_0$ and $\mu\leq C_0m/(\gamma k^2)$ eventually.
Then $\mu^2k\leq C_0^2/c_0$, and $e^x-1\leq xe^x$ for $x\geq0$ gives
\[
\frac{\gamma^2k^3}{2(1-\kappa)^2m^2}\big(\exp(\mu^2k)-1\big)
\leq\frac{C_0^2}{2(1-\kappa)^2}\exp(C_0^2/c_0)=O(1).
\]
On the other hand, in the regime $\gamma^2k^3=O(m^{2-\delta})$ for any $\delta>0$, since $\log(1+C\frac{2(1-\kappa)^2m^2}{\gamma^2k^3})\geq\delta\log m+O(1)$, we have $\log(1+C\frac{2(1-\kappa)^2m^2}{\gamma^2k^3})\geq \frac{\delta}{2}\log m$ for sufficiently large $m$. Therefore, it suffices to have $\mu \leq \sqrt{\frac{\delta}{2k} \log m}.$

	We now remove the restriction $\gamma k\leq \kappa m$ for a fixed $\kappa\in(0,1)$ and extend the analysis to $\gamma k \leq m$. The case $\gamma=1$ reduces to the single planted clique problem considered above, so suppose that $\gamma\geq2$, and let $q:=\lceil \frac{\gamma}{2} \rceil$. Recall that $\dot V_a:=\{l\in[m]\,:\,(a-1)k<\dot A_l\leq ak\}$ and let $S_a=\{(i,j)\in[m]_2: i,j\in \dot V_a\}.$ Define
\begin{equation}\label{eq: multiple cliques grouping}
\mathcal S_1:=\bigcup_{r=1}^{q} S_r,\qquad \mathcal S_2:=\bigcup_{r=\gamma-q+1}^{\gamma} S_r,
\end{equation}
and define $\mathcal S_1',\mathcal S_2'$ analogously from the independent copy $\dot\bfA'$. Then, observe that 
$$
S=\mathcal S_1\cup\mathcal S_2,\qquad S'=\mathcal S_1'\cup\mathcal S_2',
$$
which implies that
\begin{align}\label{eq: multiple cliques decomposition}
|S\cap S'|\leq|\mathcal S_1\cap\mathcal S_1'|+|\mathcal S_1\cap\mathcal S_2'|+|\mathcal S_2\cap\mathcal S_1'|+|\mathcal S_2\cap\mathcal S_2'|.
\end{align}
By the permutation invariance of $\dot \bfA$, each of $\mathcal S_1$ and $\mathcal S_2$ is marginally distributed as a fixed-size multiple clique model consisting of $q$ cliques of size exactly $k$, and similarly for $\mathcal S_1'$ and $\mathcal S_2'$. Moreover, since $\dot \bfA$ and $\dot\bfA'$ are independent, the cliques in $\mathcal S_1$ and $\mathcal S_2$ are independent of those in $\mathcal S_1'$ and $\mathcal S_2'$. Therefore, for $i,j\in\{1,2\}$, the random variables $|\mathcal S_i\cap\mathcal S_j'|$ have the same distribution as the overlap between two independent multiple clique models.

Based on this observation, we first analyze the resulting statistical lower bound.
Observe that
\begin{equation}\label{eq: statistical lower bound extension multiple cliques}
\begin{aligned}
\ee[\exp(\mu^2 |S\cap S'|)] &\stackrel{(a)}\leq \ee\bigg[\exp\bigg(\mu^2 \bigg\{\sum_{i,j=1}^2 \big|\mathcal S_i\cap \mathcal S'_j\big|\bigg\}\bigg)\bigg] \\
&=  \ee\bigg[\prod_{i,j=1}^2\exp\big(\mu^2 \big|\mathcal S_i\cap \mathcal S'_j\big|\big)\bigg],\\
&\stackrel{(b)}\leq \prod_{i,j=1}^2\Big(\ee\big[\exp\big(4\mu^2 \big|\mathcal S_i\cap \mathcal S'_j\big|\big)\big]\Big)^{1/4}\\
&\stackrel{(c)}=\ee_{\mathcal S_1,\mathcal S_1'\sim \dot \nu_{q,k}}\big[\exp\big(4\mu^2 \big|\mathcal S_1\cap \mathcal S'_1\big|\big)\big],
\end{aligned}
\end{equation}
where the inequality $(a)$ follows from \Cref{eq: multiple cliques decomposition}, the inequality $(b)$ follows from H\"older's inequality, and the equality $(c)$ follows from the fact that all the terms have the same distribution. We note that $\frac\gamma2\leq q=\lceil\frac\gamma2\rceil\leq \frac{2}{3}\gamma$, and thus $qk\leq \frac23\gamma k\leq \frac23m.$ Therefore, we take $\kappa=2/3$ and apply \Cref{equation: final statistical bound for multiple fixed size cliques} with $(\gamma,k,\mu)$ replaced by $(q,k,2\mu)$ to obtain
\begin{align*}
\ee_{\mathcal S_1,\mathcal S_1'\sim \dot \nu_{q,k}}\big[\exp\{4\mu^2 |\mathcal S_1\cap \mathcal S'_1|\}\big]&\leq \exp\bigg\{\frac{q^2 k^3}{2(1-\kappa)^2m^2}\big(\exp(4\mu^2k)-1\big)\bigg\}\\
&=\exp\bigg\{\frac{9q^2 k^3}{2m^2}\big(\exp(4\mu^2k)-1\big)\bigg\}\\
&\leq \exp\bigg\{\frac{2\gamma^2 k^3}{m^2}\big(\exp(4\mu^2k)-1\big)\bigg\},
\end{align*}
where the last inequality follows from $q\leq \frac{2}{3}\gamma$. The remaining technical steps follow the analysis subsequent to \Cref{equation: Finner inequality application for multiple fixed size cliques}, and this completes the proof.

\paragraph{Computational lower bound.} 
From \Cref{proposition: gaussian overlapping pairs}, the second moment of the low-degree likelihood ratio can be expressed as
\begin{align}\label{equation: low degree likelihood ratio multiple cliques}
\ee_{\bfY\sim\bbQ}[L_{\dot\nu}^{\leq D}(\mathbf Y)^2]=\ee_{S,S'\sim \dot\nu}\bigg[\sum^D_{d=0} \frac{\mu^{2d}}{d!}|S\cap S'|^d\bigg].
\end{align}
We reduce boundedness of this quantity to boundedness of its counterpart in the random-size model, with mixing measure $\tilde\nu$ instead of $\dot\nu$.

Specifically, recall that the latent structure $S$ for multiple planted cliques is defined as 
$$
S:=\bigcup^\gamma_{r=1}\bigg\{(i,j)\in[m]_2\,:\, (r-1)k<A_i,A_j\leq rk\bigg\},
$$
where $\bfA=(A_1,\dots,A_m)$ is either $\dot\bfA$ for the fixed-size case or $\tilde\bfA$ for the random-size case, and let us denote the corresponding distributions of $S$ given $\gamma$ and $k$ as $\dot\nu=\dot\nu_{\gamma,k}$ and $\tilde\nu=\tilde\nu_{\gamma,k}$, respectively. Fix an arbitrary constant $C\in(1,1/\kappa)$ so that $\gamma(Ck)\leq C\kappa m<m$, and consider $\tilde\nu_{\gamma,Ck}$, which corresponds to the random-size case with the expected clique size scaled by $C$. Define an event $\underline{\mathcal E}=\{|\tilde V_a|\geq k,~\forall a\in[\gamma]\},$ and observe that the union bound and Chernoff's inequality give
$$
\tilde\nu_{\gamma,Ck}(\underline{\mathcal E}) \geq 1-\gamma \exp\bigg(-\frac{(C-1)^2 k}{2C}\bigg)=1-o(1),
$$
provided that $k=\omega(\log m)$.
Let us define $\underline{\mathcal E}'=\{|\tilde V'_a|\geq k,~\forall a\in[\gamma]\}$ analogously for the independent copy.
Writing $\Pr$ for probability under the joint law of these two independent latent samples, for sufficiently large $m$, we have $\Pr(\underline{\mathcal E}\cap\underline{\mathcal E}')=\tilde\nu_{\gamma,Ck}(\underline{\mathcal E})^2 \geq \frac{1}{2}.$
Then, for such $m$ the following bound holds:
\begin{equation}\label{equation: reduction to random-size multiple cliques}
\begin{aligned}
\ee_{\bfY\sim\bbQ}[L_{\tilde \nu_{\gamma,Ck}}^{\leq D}(\mathbf Y)^2] 
=& \ee_{S,S'\sim \tilde \nu_{\gamma,Ck}}\bigg[\sum^D_{d=0} \frac{\mu^{2d}}{d!}|S\cap S'|^d\bigg]\\
=&\Pr(\underline{\mathcal E}\cap \underline{\mathcal E}') \ee_{S,S'\sim \tilde \nu_{\gamma,Ck}}\bigg[\sum^D_{d=0} \frac{\mu^{2d}}{d!}|S\cap S'|^d\bigg|\underline{\mathcal E}\cap \underline{\mathcal E}'\bigg] \\
&+ \Pr(\underline{\mathcal E}^c\cup {\underline{\mathcal E}'}^c) \ee_{S,S'\sim \tilde \nu_{\gamma,Ck}}\bigg[\sum^D_{d=0} \frac{\mu^{2d}}{d!}|S\cap S'|^d\bigg|\underline{\mathcal E}^c\cup {\underline{\mathcal E}'}^c\bigg]\\
\geq &\frac{1}{2} \ee_{S,S'\sim \tilde \nu_{\gamma,Ck}}\bigg[\sum^D_{d=0} \frac{\mu^{2d}}{d!}|S\cap S'|^d\bigg|\underline{\mathcal E}\cap \underline{\mathcal E}'\bigg]\\
\geq &\frac{1}{2} \ee_{S,S'\sim \dot\nu_{\gamma,k}}\bigg[\sum^D_{d=0} \frac{\mu^{2d}}{d!}|S\cap S'|^d\bigg]\\
=&\frac{1}{2} \ee_{\bfY\sim\bbQ}[L_{\dot\nu_{\gamma,k}}^{\leq D}(\mathbf Y)^2],
\end{aligned}
\end{equation}
where the last inequality follows from the fact that the distribution of $|S\cap S'|$ conditioned on $\underline{\mathcal E}\cap\underline{\mathcal E}'$ under $\tilde \nu_{\gamma,Ck}$ stochastically dominates the distribution of $|S\cap S'|$ under $\dot\nu_{\gamma,k}$, since on the event $\underline{\mathcal E}\cap\underline{\mathcal E}'$, each clique in the random-size model has size at least $k$, which is the size of the cliques in the fixed-size model.
To construct the coupling, select a uniform $k$-subset of each clique, independently in the two copies, conditional on these events.
By vertex-permutation invariance, the resulting disjoint vertex sets have the fixed-size distribution, and their edge sets are contained in those of the random-size model.
In summary, we obtain the following
$$
\ee_{\bfY\sim\bbQ}[L_{\dot\nu_{\gamma,k}}^{\leq D}(\mathbf Y)^2]\leq 2\ee_{\bfY\sim\bbQ}[L_{\tilde \nu_{\gamma,Ck}}^{\leq D}(\mathbf Y)^2],
$$
and thus if we can show that $\ee_{\bfY\sim\bbQ}[L_{\tilde \nu_{\gamma,Ck}}^{\leq D}(\mathbf Y)^2]$ is $O(1)$ for $D=D_m$ as specified above, then we can conclude that $\ee_{\bfY\sim\bbQ}[L_{\dot\nu_{\gamma,k}}^{\leq D}(\mathbf Y)^2]$ is also $O(1)$ for the same $D$. Therefore, the computational lower bound established for the random-size case in \Cref{proof: random-size multiple cliques} can be applied directly to the fixed-size case through the above comparison. Strictly speaking, the comparison requires a bound for the random-size model with mixing measure $\tilde\nu_{\gamma,Ck}$, whereas the analysis on the random-size case states the result for $\tilde\nu_{\gamma,k}$. This distinction is immaterial for the asymptotic rate: as can be seen from the resulting bound, replacing $k$ by $Ck$ changes the sufficient condition on $\mu$ only by a constant factor.

We now consider the full feasible region $\gamma k\leq m$, removing the restriction $\gamma k\leq\kappa m$.
Recall the definition of $\mathcal S_1$ and $\mathcal S_2$ from \Cref{eq: multiple cliques grouping} and the bound $
|S\cap S'|\leq|\mathcal S_1\cap\mathcal S_1'|+|\mathcal S_1\cap\mathcal S_2'|+|\mathcal S_2\cap\mathcal S_1'|+|\mathcal S_2\cap\mathcal S_2'|.$ Then, observe that
\begin{equation}\label{eq: multiple cliques grouping on computational lower bound}
\begin{aligned}
	\ee_{\bfY\sim\bbQ}[L_{\dot\nu}^{\leq D}(\bfY)^2] &= 1+\ee\bigg[\sum^D_{d=1} \frac{\mu^{2d}}{d!}|S\cap S'|^d\bigg]\\
	&\leq 1+\ee\bigg[\sum^D_{d=1} \frac{\mu^{2d}}{d!}\bigg\{\sum_{i,j=1}^2 \big|\mathcal S_i\cap \mathcal S'_j\big|\bigg\}^d\bigg]\\
	&\stackrel{(a)}\leq 1+\ee\bigg[\sum^D_{d=1} \frac{4^{d-1}\mu^{2d}}{d!}\sum_{i,j=1}^2 \big|\mathcal S_i\cap \mathcal S'_j\big|^d\bigg]\\
	&= 1+\ee_{\mathcal S_1,\mathcal S_1'\sim \dot \nu_{q,k}}\bigg[\sum^D_{d=1} \frac{(2\mu)^{2d}}{d!} \big|\mathcal S_1\cap \mathcal S'_1\big|^d\bigg]\\
	&\stackrel{(b)}\leq 1+ 2\ee_{\mathcal S_1,\mathcal S_1'\sim \tilde \nu_{q,Ck}}\bigg[\sum^D_{d=1} \frac{(2\mu)^{2d}}{d!} \big|\mathcal S_1\cap \mathcal S'_1\big|^d\bigg]
\end{aligned}
\end{equation}
where $(a)$ follows from the fact that $(\sum_{i=1}^4 x_i)^d \leq 4^{d-1} \sum_{i=1}^4 x_i^d$ for $x_i\geq 0$, and $(b)$ follows from the same comparison of nonnegative moments as in \Cref{equation: reduction to random-size multiple cliques} since $qk\leq\frac{2}{3}m$ and thus the argument holds with $\kappa=\frac{2}{3}$ and $C\in(1,1/\kappa)$. Observe that the expectation on the right-hand side has the same form as the nonconstant part of the low-degree likelihood-ratio second moment for the random-size multiple-clique model with $q$ cliques, clique-size parameter $Ck$, and signal strength $2\mu$. Note that $q(Ck)<m$, and $\frac{q(Ck)^2}{m}\asymp\frac{\gamma k^2}{m}$ since $\frac{C^2}{2}\frac{\gamma k^2}{m}\leq\frac{q(Ck)^2}{m}\leq\frac{2C^2}{3}\frac{\gamma k^2}{m}$. Therefore, the random-size computational lower-bound analysis in \Cref{proof: random-size multiple cliques} applies with $(\gamma,k,\mu)$ replaced by $(q,Ck,2\mu)$, and these substitutions affect the resulting conditions only by constant factors.
To absorb these fixed factors in the regime conditions, replace $\delta$ by any $\delta'\in(\epsilon,\delta)$ and decrease the signal constant if necessary.

\subsubsection{Random-size case}\label{proof: random-size multiple cliques}
We now consider the model with random clique sizes.
More specifically, each vertex $v\in[m]$ is independently assigned to at most one of the $\gamma$ cliques with equal probability $k/m$, so that the expected size of each clique is $k$. For $a\in[\gamma]$, define $\tilde V_a:=\{v\in[m]\,:\,(a-1)k<\tilde A_v\leq ak\}$ and let $S_a=\{(i,j)\in[m]_2: i,j\in \tilde V_a\}.$ Let $\tilde\bfA'$ be an independent copy of $\tilde\bfA$, and define $S'_1,...,S'_\gamma$ and $ \tilde V'_1,...,\tilde V'_\gamma$ analogously.
Let $S:=\bigcup_{a=1}^{\gamma}S_a,~S':=\bigcup_{b=1}^{\gamma}S'_b.$

\paragraph{Upper bounds.}
As in the upper bound analysis in the random-size clique problem (see \Cref{proof: random-size clique problem}), we consider a high-probability event $\underline{\mathcal E}$ on which the size of each clique is at least $ck$ for some constant $c\in(0,1)$ and control the risk of the tests conditioned on this event. Specifically, consider the event
\begin{align}\label{equation: lower truncation event for random-size multiple cliques}
\underline{\mathcal E}:=\big\{|\tilde V_a|\geq ck,~\forall a\in[\gamma]\big\},
\end{align}
and let $\underline{\tilde\nu}$ denote the corresponding conditional distribution, $\underline{\tilde\nu}(\cdot):=\tilde\nu(\cdot|\underline{\mathcal E}).$ Since each $|\tilde V_a|$ is marginally distributed as $\mathrm{Binomial}(m,k/m)$, the union bound and Chernoff's inequality give
\begin{align*}
	\tilde\nu(\underline{\mathcal E}^c)\leq \gamma\exp\bigg(-\frac{(1-c)^2}{2}k\bigg)=o(1),
\end{align*}
provided that $\log\gamma=o(k)$, which is satisfied in the polynomial regime considered here. As shown in \Cref{equation: risk lower bound via TV with conditioning for upper bound}, the risk difference between the original and conditional distributions vanishes:
\begin{equation*}
\begin{aligned}
	R_{\tilde\nu}(\Delta)\leq R_{\underline{\tilde\nu}}(\Delta)+\tilde\nu(\underline {\mathcal E}^c),
\end{aligned}
\end{equation*}
which implies that we have $R_{\tilde\nu}(\Delta)\to 0$ as long as $R_{\underline{\tilde\nu}}(\Delta)\to 0$. 

Observe that under the event $\underline{\mathcal E}$, each $S$ contains $\gamma$ disjoint cliques of size at least $ck$, and thus the tests for the fixed-size case can be applied by targeting the cliques of size $ck$.
In particular, the average test uses threshold $\mu\gamma\binom{\lfloor ck\rfloor}{2}/2$, and the scan searches over unions of $\gamma$ disjoint $\lfloor ck\rfloor$-vertex cliques.
Larger planted cliques only increase the mean of the average statistic, and the scan still includes a fully signaled candidate.
Replacing $k$ with $ck$ only affects the constant factors in the upper bound, and hence we obtain the same upper bound as in \Cref{proof: multiple fixed-size cliques problem}.

\paragraph{Statistical lower bound.}
We first assume that $\gamma k \leq \kappa m$ for some constant $\kappa\in(0,1)$, and we will remove this assumption later.
We follow the truncation argument detailed in \Cref{proof: random-size clique problem}, which we sketch below for the reader's convenience. Fix a constant $C\in(1,1/\kappa)$ and define the event
\begin{equation}\label{equation: truncation event for random-size multiple cliques}
\overline{\mathcal E}:=\big\{|\tilde V_a|\leq Ck,~\forall a\in[\gamma]\big\}.
\end{equation}
Let $\overline{\tilde\nu}$ denote the conditional distribution of $\tilde\nu$ given $\overline{\mathcal E}$. Since each $|\tilde V_a|$ is marginally distributed as $\mathrm{Binomial}(m,k/m)$, the union bound and Chernoff's inequality give
$$
\tilde\nu(\overline{\mathcal E}^c)\leq \gamma\exp\bigg(-\frac{(C-1)^2}{2C}k\bigg)=o(1),
$$
where the last equality holds when $\log\gamma=o(k)$, as in the polynomial regime considered here. Therefore, by the same argument as in \Cref{equation: risk lower bound via TV with conditioning}, for any test $\Delta$,
$$
R_{\tilde\nu}(\Delta)\geq 1-\sqrt{1-\frac{1}{\ee_{\bfY\sim\bbQ}[L_{\overline{\tilde\nu}}(\bfY)^2]}}-\tilde\nu(\overline{\mathcal E}^c),
$$
and hence it suffices to show that $\ee_{\bfY\sim\bbQ}[L_{\overline{\tilde\nu}}(\bfY)^2]=O(1)$.
To prove this, we compare the truncated random-size model with the fixed-size model. See \Cref{proof: random-size clique problem} for detailed arguments. Under $\overline{\tilde\nu}$, the sets $\tilde V_1,\dots,\tilde V_\gamma$ are disjoint and satisfy $|\tilde V_a|\leq Ck$ for all $a\in[\gamma]$. Since $C\gamma k\leq C\kappa m<m$, we can enlarge each $\tilde V_a$ to a set $\dot V_a\supseteq\tilde V_a$ of size exactly $Ck$ by adding vertices uniformly at random from the unassigned vertices, and construct $\dot V'_1,\dots,\dot V'_\gamma$ independently in the same way from $\tilde V'_1,\dots,\tilde V'_\gamma$. By permutation invariance of $\overline{\tilde\nu}$ and of the filling procedure, $(\dot V_1,\dots,\dot V_\gamma)$ and $(\dot V'_1,\dots,\dot V'_\gamma)$ are independent and each has the same distribution as the vertex sets in the fixed-size multiple-clique model with clique size $Ck$.

Let $\dot S$ and $\dot S'$ denote the corresponding unions of fixed-size cliques. Since $\tilde V_a\subseteq\dot V_a$ and $\tilde V'_a\subseteq\dot V'_a$ for every $a\in[\gamma]$, we have $S\subseteq\dot S$ and $S'\subseteq\dot S'$, and therefore $|S\cap S'|\leq|\dot S\cap\dot S'|$ almost surely. It follows that
\begin{align}
\ee_{\bfY\sim\bbQ}[L_{\overline{\tilde\nu}}(\bfY)^2]
&=\ee_{S,S'\sim\overline{\tilde\nu}}\big[\exp\{\mu^2|S\cap S'|\}\big]\leq \ee_{\dot S,\dot S'}\big[\exp\{\mu^2|\dot S\cap\dot S'|\}\big]\label{equation: comparison of random-size and fixed-size multiple cliques}\\
&\leq \exp\bigg\{\frac{C^3\gamma^2k^3}{2(1-C\kappa)^2m^2}\big(\exp(C\mu^2k)-1\big)\bigg\},
\end{align}
where the last inequality follows by repeating the fixed-size analysis above with $k$ replaced by $Ck$, and $\kappa$ replaced by $C\kappa$. Since $C\in(1,1/\kappa)$ is a fixed constant, this replacement only changes constant factors. Therefore, repeating the two cases in the fixed-size analysis yields the same statistical lower bound as in \Cref{theorem: multiple cliques problem}.

We next remove the restriction $\gamma k \leq \kappa m$ and consider the general case $\gamma k \leq m$.
We use the grouping argument from \Cref{proof: multiple fixed-size cliques problem}.
Let $q:=\lceil \frac{\gamma}{2} \rceil$ and define
$$
\mathcal S_1:=\bigcup_{r=1}^{q} S_r,\qquad \mathcal S_2:=\bigcup_{r=\gamma-q+1}^{\gamma} S_r,
$$
and define $\mathcal S_1',\mathcal S_2'$ analogously from the independent copy $\tilde\bfA'.$
We apply the filling procedure used in \Cref{equation: comparison of random-size and fixed-size multiple cliques}. Since $q\leq \frac{2\gamma}{3}$, it follows that
$qCk\leq\frac{2C}{3}\gamma k\leq\frac{2C}{3}m.$ Therefore, taking $\kappa=2/3$, we have $qk\leq \kappa m$, and for any $C\in(1,1/\kappa)$,
$qCk\leq C\kappa m<m.$
Thus, each of the $q$ cliques in $\mathcal S_1$ and $\mathcal S_1'$ can be enlarged to a clique of size $Ck$, exactly as in the argument leading to \Cref{equation: comparison of random-size and fixed-size multiple cliques}. Let us denote the corresponding unions of such cliques by $\dot{\mathcal S}_1$ and $\dot{\mathcal S}_1'$, respectively. Observe that $\mathcal S_1\subseteq\dot{\mathcal S}_1$ and $\mathcal S_1'\subseteq\dot{\mathcal S}_1'$. Moreover, by permutation invariance of the conditional distribution and of the filling procedure,  $\dot{\mathcal S}_1$ and $\dot{\mathcal S}_1'$ are independent and identically distributed with law $\dot\nu_{q,Ck}$. Following the analysis in \Cref{eq: statistical lower bound extension multiple cliques}, we obtain
\begin{align*}
\ee_{S,S'\sim\overline{\tilde\nu}_{\gamma,k}}\big[\exp\{\mu^2|S\cap S'|\}\big]&\leq\ee_{\dot{\mathcal S}_1,\dot{\mathcal S}_1'\sim \dot\nu_{q,Ck}}\big[\exp\big(4\mu^2 \big|\dot{\mathcal S}_1\cap \dot{\mathcal S}'_1\big|\big)\big]\\
&\leq \exp\bigg\{\frac{C^3q^2k^3}{2(1-2C/3)^2m^2}\big(\exp(4C\mu^2k)-1\big)\bigg\}\\
&\leq \exp\bigg\{\frac{2C^3\gamma^2k^3}{9(1-2C/3)^2m^2}\big(\exp(4C\mu^2k)-1\big)\bigg\},
\end{align*}
where the second inequality is obtained by applying the fixed-size multiple-clique bound in \Cref{equation: final statistical bound for multiple fixed size cliques} with $(\gamma,k,\mu)\mapsto(q,Ck,2\mu)$, and the last inequality follows from $q\leq\frac{2\gamma}{3}$, which gives $q^2/2\leq2\gamma^2/9$. Thus, the resulting bound on $\mu$ differs from the previous bound by a constant factor.

\paragraph{Computational lower bound.} 
Our objective is to bound
\begin{align*}
\bbE_{\mathbf Y\sim\bbQ}=\ee_{S,S'\sim \tilde\nu}\bigg[\sum^D_{d=0} \frac{\mu^{2d}}{d!}|S\cap S'|^d\bigg],
\end{align*}
and we focus on bounding $\ee_{S,S'\sim \tilde\nu}[|S\cap S'|^d]$ for each $d\in [D].$
Our approach is to bound these moments using an auxiliary clique model whose overlap moments are no smaller and are easier to analyze.
For each vertex $v\in[m]$, define
\begin{align}\label{eq:definition of W_v}
W_v=\begin{cases}
e_{(a,b)} & \text{if }v\in\widetilde V_a\cap\widetilde V'_b\text{ for some }a,b\in[\gamma],\\
0 & \text{otherwise},
\end{cases}
\end{align}
where $e_{(a,b)}\in\{0,1\}^{\gamma\times \gamma}$ is the basis matrix with a 1 in the $(a,b)$-th entry and 0 elsewhere.
Then $W_1,\ldots,W_m$ are independent, and
$$
W_v=\begin{cases}
e_{(a,b)} & \text{with probability }\rho,\quad a,b\in[\gamma],\\
0 & \text{with probability }1-\gamma^2\rho.
\end{cases}
$$
where
$$
\rho:=k^2/m^2.
$$
Observe that $|S\cap S'|$ can be expressed as
\begin{equation}\label{equation: overlap in terms of W_v multiple cliques}
\begin{aligned}
|S\cap S'|
&=\sum_{a,b=1}^{\gamma}|S_a\cap S'_b|\\
&=\sum_{a,b=1}^{\gamma}\sum_{(i,j)\in[m]_2}\mathds 1\big\{i,j\in\widetilde V_a\cap\widetilde V'_b\big\}\\
&=\sum_{(i,j)\in[m]_2}\sum_{a,b=1}^{\gamma}\mathds 1\big\{i,j\in\widetilde V_a\cap\widetilde V'_b\big\}\\
&=\sum_{(i,j)\in[m]_2}\sum_{a,b=1}^{\gamma}\mathds 1\{W_i=e_{(a,b)}\}\mathds 1\{W_j=e_{(a,b)}\}\\
&=\sum_{(i,j)\in[m]_2}\mathds 1\{W_i=W_j\neq 0\}.
\end{aligned}
\end{equation}
An intuitive interpretation of the above expression is that $|S\cap S'|$ can be seen as the number of edges of another random $\gamma^2$-multiple clique model; specifically, consider a model where each vertex $v\in[m]$ is assigned a label $e_{(a,b)}$ when $W_v=e_{(a,b)}$, with equal probability $\rho$ (or $0$ with probability $1-\gamma^2\rho$), and an edge is formed between two vertices having the same nonzero label. Note that the number of labels is $\gamma^2$.

For notational convenience, write the overlap as $Z:=|S\cap S'|$. Then, our approach is to bound $\ee[Z^d]$ by comparing it with the $d$-th moment of an auxiliary random variable $Z_{\ast}$, which is defined as follows. Let $B_1,\ldots,B_m\stackrel{\mathrm{i.i.d.}}\sim\mathrm{Bernoulli}(p),$ where $p$ is a parameter to be determined below, and let $B:=\sum_{i=1}^m B_i\sim\mathrm{Binomial}(m,p)$. Then, we set
\[
Z_{\ast}:=\sum_{(i,j)\in[m]_2}B_iB_j= \binom{B}{2}.
\]
Then, $B$ can be seen as the size of a random single clique in which each of the $m$ vertices belongs to the clique independently and with probability $p$, and $Z_{\ast}$ is the number of edges in the resulting random clique. We seek to find a suitable value for $p$ such that $\ee[Z^d]\leq \ee[Z_{\ast}^d]$ for all $d\in[D]$, while matching the first moments. We begin with the case of $d=1$ and compare the first moments of $Z$ and $Z_{\ast}$. Letting $I_{(i,j)}:=\mathds 1\{W_i=W_j\neq 0\}$, we have
\begin{align*}
\ee[Z]=\sum_{(i,j)\in[m]_2}\ee[I_{(i,j)}]=\sum_{(i,j)\in[m]_2}\Pr(W_i=W_j\neq 0)=\sum_{(i,j)\in[m]_2}\gamma^2\rho^2=\frac{m(m-1)}{2}\gamma^2\rho^2.
\end{align*}
On the other hand, denoting $J_{(i,j)}:=\mathds 1\{B_i=B_j=1\}$ for all $i \neq j$, the first moment of $Z_{\ast}$ is given by
\begin{align*}
\ee[Z_{\ast}]=\sum_{(i,j)\in[m]_2}\ee[J_{(i,j)}]=\sum_{(i,j)\in[m]_2}\Pr(B_i=B_j=1)=\sum_{(i,j)\in[m]_2}p^2=\frac{m(m-1)}{2}p^2.
\end{align*}
This motivates us to choose $p=\gamma\rho=\gamma k^2/m^2$ so that $\ee[Z]=\ee[Z_{\ast}]$.
This is a valid Bernoulli parameter because $\gamma k\leq m$ implies $p\leq k/m\leq1$.
Next, we will show that this choice of $p$ indeed guarantees that $\ee[Z^d]\leq \ee[Z_{\ast}^d]$ for all $d\geq 1$.

For any integer $d \geq 2$, write $(\bfi,\bfj):=\big((i_1,j_1),\ldots,(i_d,j_d)\big)\in([m]_2)^d,$ and observe that
\begin{align*}
\ee[Z^d]&=\sum_{(\bfi,\bfj)\in([m]_2)^d}\ee\bigg[\prod_{r=1}^d I_{(i_r,j_r)}\bigg]=\sum_{(\bfi,\bfj)\in([m]_2)^d}\Pr\bigg(\prod_{r=1}^d I_{(i_r,j_r)}=1\bigg).
\end{align*}
Similarly, the $d$-th moment of $Z_{\ast}$ can be expressed as
$$\ee[Z_{\ast}^d]=\sum_{(\bfi,\bfj)\in([m]_2)^d}\Pr\bigg(\prod_{r=1}^d J_{(i_r,j_r)}=1\bigg).$$
For each $(\bfi,\bfj)\in([m]_2)^d$, let $v=v(\bfi,\bfj)$ be the number of distinct vertices in $(\bfi,\bfj)$, and let $c=c(\bfi,\bfj)$ be the number of connected components in the graph formed by the edges $(i_1,j_1),\ldots,(i_d,j_d)$. Then, for $\prod_{r=1}^d I_{(i_r,j_r)}$ to be equal to 1, all vertices in each connected component must have the same nonzero label $e_{(a,b)}$ for some $a,b\in[\gamma]$. Since there are $\gamma^2$ possible labels, $\gamma^{2c}$ is the number of ways to assign a nonzero label to each of the $c$ connected components: each connected component requires its vertices to be associated with some pair $(a,b)$ of cliques, for a total of $(\gamma^2)^c$ possibilities. Then, enumerate the connected components in any arbitrary way and let $v_h$ be the number of vertices in the $h$th connected component, so $v = \sum_{h=1}^c v_h$. For any given assignment of pairs to connected components (as we saw, there are $\gamma^{2c}$ such assignments), the probability that all the vertices in the $h$th connected component are of the same given type, say $(a_h,b_h)$, is $\rho^{v_h}$, yielding $\rho^{\sum_h v_h} = \rho^v$ for the overall probability that each connected component has vertices of the same given types (by independence).
Therefore, we have
\begin{align}\label{equation: moment decomposition for random-size multiple cliques}
\Pr\bigg(\prod_{r=1}^d I_{(i_r,j_r)}=1\bigg)=\gamma^{2c}\rho^{v}.
\end{align}

On the other hand, for the auxiliary model, $\prod_{r=1}^d J_{(i_r,j_r)}$ is 1 if and only if $B_i=1$ for every vertex appearing in the selected edges, an event of probability
$$
\Pr\bigg(\prod_{r=1}^d J_{(i_r,j_r)}=1\bigg)=p^{v}=\gamma^v\rho^v.
$$
As every connected component contains at least two vertices, $v\geq 2c$. This implies that 
\[
\Pr\bigg(\prod_{r=1}^d I_{(i_r,j_r)}=1\bigg) =\gamma^{2c}\rho^{v}\leq \gamma^v\rho^v = \Pr\bigg(\prod_{r=1}^d J_{(i_r,j_r)}=1\bigg),
\]
which results in
$$
\ee[Z^d]=\sum_{(\bfi,\bfj)\in([m]_2)^d}\Pr\bigg(\prod_{r=1}^d I_{(i_r,j_r)}=1\bigg)\leq \sum_{(\bfi,\bfj)\in([m]_2)^d}\Pr\bigg(\prod_{r=1}^d J_{(i_r,j_r)}=1\bigg)=\ee[Z_{\ast}^d].
$$

In summary, we have shown that $\ee[|S\cap S'|^d]\leq \ee[Z_{\ast}^d]$ for all $d\geq 1$, which allows us to work with the sum $\sum_{d=0}^D\frac{\mu^{2d}}{d!}\ee[Z_{\ast}^d]$. Note that $Z_{\ast}=\binom{B}{2}\leq B^2/2$, and $B\sim\mathrm{Binomial}(m,p)$ with $p=\gamma k^2/m^2$. Following the same argument as in \Cref{proof: computational lower bound for fixed-size clique problem}, we apply the polynomial moment bound for binomial random variables \citep[Corollary 1]{ahle2022sharp} given in \Cref{lemma: moment bound of sub-Poissonian distribution} to obtain
\begin{align}\label{equation: computational lower bound moment bound for random-size multiple cliques}
\sum_{d=0}^D\frac{\mu^{2d}}{2^d d!}\ee\big[B^{2d}\big]\leq 1+\sum_{d=1}^D\frac{\mu^{2d}}{2^d d!}\bigg(\frac{2d}{\log(1+2dm/(\gamma k^2))}\bigg)^{2d}.
\end{align}
Note that the only difference from \Cref{equation: computational lower bound moment bound} is that the denominator inside the logarithm is now $\gamma k^2$ instead of $k^2$. Therefore, the rest of the proof can be completed by replacing $k^2$ with $\gamma k^2$ in the analysis following \Cref{equation: computational lower bound moment bound}. We remark that the argument applies throughout the feasible region where $\gamma k \leq m$.
Combining the upper and lower bounds proves the statistical thresholds in \Cref{theorem: multiple cliques problem} at polynomial order, as defined in \Cref{definition: statistical and computational thresholds}.
The computational thresholds follow in the same way under \Cref{conjecture: low-degree polynomial conjecture}.

\subsection{\texorpdfstring{Proof of \Cref{corollary: multiple planted hypercliques problem}}{Proof of Corollary \getrefnumber{corollary: multiple planted hypercliques problem}}}\label{proof: multiple planted hypercliques problem}

We extend the multiple-clique analysis to the hypergraph setting.
We first state the statistical and low-degree lower bounds to be established; the computational interpretation of the latter uses \Cref{conjecture: low-degree polynomial conjecture}:
\begin{itemize} 
	\item (Statistical lower bound) Fix any constant $\delta\in(0,0.1)$.
For either prior distribution $\dot\nu$ or $\tilde\nu$, strong detection is impossible under either of the following conditions: (i) $\gamma^2k^{t+1}m^{-t}\leq m^{-\delta}$ and $\mu\leq C_\delta\sqrt{k^{-(t-1)}\log m}$ for a sufficiently small constant $C_\delta>0$; or (ii) $\gamma^2k^{t+1}m^{-t}=\Omega(1)$ and $\mu=O\big(\gamma^{-1}k^{-t}m^{t/2}\big)$.
	\item (Computational lower bound) Fix any constants $\delta\in(0,0.1)$ and $\epsilon\in(0,\delta)$, and let $D_m=\lfloor\log^{1+\epsilon}m\rfloor$.
For either prior distribution $\dot\nu$ or $\tilde\nu$, no polynomial of degree at most $D_m$ achieves strong separation under either of the following conditions: (i) $\gamma^{2/t}k^2m^{-1}\leq m^{-\delta}$ and $\mu\leq C(\log m)^{\{1-(t-1)\epsilon\}/2}$ for a sufficiently small constant $C>0$; or (ii) $\gamma^{2/t}k^2m^{-1}\geq\log^{1+\delta}m$ and $\mu=O\big(\gamma^{-1}k^{-t}m^{t/2}\big)$.
\end{itemize}
Throughout, $t\geq2$ is fixed and $k$ has positive polynomial order in $m$, as in \Cref{corollary: multiple planted hypercliques problem}; in particular, $k=\omega(\log m)$ and $k\geq t$ eventually.
Constants may depend on $t$, and we take $D=D_m$ in the low-degree bounds.
Noninteger fixed vertex counts are rounded down, as in \Cref{proof: multiple cliques problem}.
In both the fixed-size and random-size cases, the proof consists of three pieces: the upper bound, the statistical lower bound, and the computational lower bound.
For ease of analysis, we assume that $\gamma k\leq \kappa m$ for some fixed constant $\kappa\in(0,1)$ if needed.
However, this can be relaxed to the general case $\gamma k\leq m$, following the same argument as in \Cref{proof: multiple cliques problem}.
\begin{remark}[Extension to $\gamma k\leq m$]
A key idea is that the multiple hyperclique $S$ with $\gamma k \leq m$ can be expressed as a union of two groups of multiple hypercliques, say $\mathcal S_1$ and $\mathcal S_2$, where each group consists of $q:=\lceil\gamma/2\rceil$ hypercliques. Then, the overlap $|S\cap S'|$ can be bounded by $\sum_{i,j=1}^2|\mathcal S_i\cap \mathcal S_j'|$. Since $q=\lceil\gamma/2\rceil\leq {2\gamma}/{3}$ for $\gamma \geq 2$ and thus $qk\leq {2m}/{3}$, we can interpret $|\mathcal S_i\cap \mathcal S_j'|$ as the overlap between two independent multiple hypercliques satisfying the restricted occupancy condition $qk\leq \kappa m$ with $\kappa=2/3$, and apply the lower-bound results established for this regime. See \Cref{proof: multiple cliques problem} for more details.
Specifically, H\"older's inequality controls the full second moment, while $(x_1+\cdots+x_4)^d\leq4^{d-1}\sum_{i=1}^4x_i^d$ controls the low-degree moments; both comparisons only rescale $\mu$ by a constant factor.
For $\gamma=1$ and $k>\kappa m$, use $|S\cap S'|\leq\binom mt$ directly: if $\mu=O(m^{t/2}/k^t)$, then $\mu^2\binom mt=O(1)$, which bounds both the full and low-degree second moments.
\end{remark}

\subsubsection{Fixed-size case}
We consider the case where $\dot{\mathcal C}$ consists of all subhypergraphs that are the union of $\gamma$ disjoint hypercliques of size $k$ and dimension $t$. In this case, we have $n=\binom mt$, $K=\gamma\binom kt$ and $|\dot{\mathcal C}|=\frac{1}{\gamma!}\binom m k \cdot \binom{m-k}k \cdots \binom{m-(\gamma-1)k}k.$ Define $\dot V_a:=\{v\in[m]\,:\,(a-1)k<\dot A_v\leq ak\}$ and let $S_a=\{\bfi=(i_1,\dots,i_t)\in[m]_t: i_1,\dots,i_t\in \dot V_a\}.$ Let $\dot\bfA'$ be an independent copy of $\dot\bfA$, and define $S'_1,...,S'_\gamma$ and $ \dot V'_1,...,\dot V'_\gamma$ analogously. Finally, set $S:=\bigcup_{a=1}^{\gamma}S_a$ and $S':=\bigcup_{b=1}^{\gamma}S'_b.$

\paragraph{Upper bounds.}
The upper bounds on the critical signal strength $\mu^*$ can be obtained by applying the three tests in \Cref{proposition: gaussian test upper bound}:
\begin{itemize}
	\item For some constant $C>0$, the outlier test achieves $R_{\dot \nu}(\Delta_{\textrm{out}})\to 0$ whenever
	$\mu \geq C \sqrt{\log m}.$
	\item The average test achieves $R_{\dot \nu}(\Delta_{\textrm{avg}})\to 0$ when $\mu =\omega(\frac{m^{t/2}}{\gamma k^t}).$
	\item Note that $|\dot{\mathcal C}|\leq \binom mk^\gamma$, and thus $\log|\dot{\mathcal C}|\leq \gamma k\log\frac{em}{k}$.
Therefore, for some constant $C>0$, the scan test achieves $R_{\dot \nu}(\Delta_{\textrm{scan}})\to 0$ provided that $k\leq(1-c_0)m$ for a fixed $c_0\in(0,1)$ and
	$$
	\mu \geq C \sqrt{\frac{1}{k^{t-1}}\log\frac{m}{k}}.
	$$
\end{itemize}
The restriction on $k$ holds in the regime $\gamma^2k^{t+1}\ll m^t$ where the scan bound is needed.
As in \Cref{proof: multiple fixed-size cliques problem}, take $\delta_m=1/\log m$ for the outlier and scan tests, and use the average-test risk bound directly.

\paragraph{Statistical lower bound.} 
Assume that $\gamma k\leq \kappa m$ for some fixed constant $\kappa\in(0,1)$.
We can write
\begin{equation}
\begin{aligned}
	\ee[\exp\{\mu^2 |S\cap S'|\}] &= \ee\bigg[\exp\bigg\{\mu^2 \sum_{a,b=1}^\gamma \binom{|\dot V_a\cap \dot V'_b|}{t}\bigg\}\bigg] \\
	&= \ee\bigg[\prod_{a,b=1}^\gamma f\big(|\dot V_a\cap \dot V'_b|\big)\bigg], 
\end{aligned}
\end{equation}
where we now set $f(x):=\exp\{\mu^2\binom xt\}$ for nonnegative integers $x$, and extend $f$ linearly between consecutive integers.
The sequence $\binom xt$ is nondecreasing and discretely convex, so exponentiation and linear interpolation give a nondecreasing convex function on $[0,\infty)$.
This extension leaves all expectations below unchanged, since their arguments are integer-valued.
The conditional negative-association and convex-order arguments from \Cref{equation: conditional expectation argument for multiple fixed size cliques} through \Cref{equation: iterating conditional expectation argument for multiple fixed size cliques} therefore apply and give
\begin{align}\label{equation: fixed-size multiple hypercliques moment upper bound}
\ee\bigg[\prod_{a,b=1}^\gamma f\big(|\dot V_a\cap \dot V'_b|\big)\bigg]\leq \ee[f(B)]^{\gamma^2},
\end{align}
where $B = \sum_{i=1}^k B_i$, with  $B_1,\ldots,B_k$ i.i.d.\ Bernoulli random variables with parameter $p_{\gamma,k}:=k/(m-(\gamma-1)k).$ Then $B \sim \mathrm{Binomial}(k,p_{\gamma,k})$ and 
$$
\ee[f(B)]=\ee\bigg[\exp\bigg\{\mu^2 \binom Bt\bigg\}\bigg]=\ee\bigg[\exp\bigg\{\mu^2\sum_{\bfi\in[k]_t}B_{\bfi}\bigg\}\bigg] = \ee \bigg[ \prod_{\bfi\in[k]_t} \exp ( \mu^2 B_{\bfi} ) \bigg]
$$
where, for any $\bfi = (i_1,\ldots,i_t) \in [k]_t,$ we set
$B_{\bfi}=B_{i_1}\cdots B_{i_t}$, and the second equality above follows from the easily verifiable identity
\begin{equation}
    \label{eq:Bt}
\binom Bt = \sum_{\bfi\in[k]_t}B_{\bfi}.
\end{equation}
Next we apply Finner's inequality as given in \Cref{lemma: Finner's inequality}.
In the notation of that lemma, take $I=[k]$, $J=[k]_t$, and $T_{\bfi}=\{i_1,\ldots,i_t\}$ for each $\bfi=(i_1,\ldots,i_t)\in[k]_t$.
Thus, there are $\binom kt$ factors, each depending on $t$ coordinates, and we index the corresponding exponents by $p_{\bfi}$.
In particular, since each $B_{i}$ appears in $\binom{k-1}{t-1}$ terms in the summation in \Cref{eq:Bt}, we choose $p_{\bfi}=\binom{k-1}{t-1}$ for all $\bfi\in [k]_t$.
The reciprocal exponents therefore sum to one at each coordinate, as required by Finner's inequality.
Then, Finner's inequality (\Cref{lemma: Finner's inequality}) yields
\begin{equation}\label{equation: fixed-size multiple hypercliques finner bound}
\begin{aligned}
\ee_B\bigg[\exp\bigg(\mu^2\sum_{\bfi\in[k]_t} B_{\bfi}\bigg)\bigg] &\leq \prod_{\bfi\in[k]_t} \ee_B\big[\exp(\mu^2B_{\bfi})^{p_{\bfi}}\big]^{\frac{1}{p_{\bfi}}}\\
&= \prod_{\bfi\in[k]_t} \ee_B\big[\exp(\mu^2B_{\bfi})^{\binom{k-1}{t-1}}\big]^{\frac{1}{\binom{k-1}{t-1}}}\\
&= \ee_B\bigg[\exp\bigg\{\mu^2\binom{k-1}{t-1}B_{\bfi}\bigg\}\bigg]^{k/t}\\
&\leq \ee_B\bigg[\exp\bigg\{\mu^2k^{t-1}B_{\bfi}\bigg\}\bigg]^{k/t}\\
&= \bigg\{1+p_{\gamma,k}^t\big(\exp(\mu^2k^{t-1})-1\big)\bigg\}^{k/t} \quad \big(\because B_{\bfi} \sim \mathrm{Bernoulli}(p_{\gamma,k}^t)\big)\\
&\stackrel{(a)}\leq \exp\bigg\{\frac{kp_{\gamma,k}^t}{t}\big(\exp(\mu^2k^{t-1})-1\big)\bigg\}\\
&\stackrel{(b)}\leq \exp\bigg\{\frac{k^{t+1}}{t(1-\kappa)^tm^t}\big(\exp(\mu^2k^{t-1})-1\big)\bigg\}
\end{aligned}
\end{equation}
where $(a)$ follows from the fact that $1+x \leq e^x$ for all $x\in \mathbb R,$ and $(b)$ follows from the fact that $p_{\gamma,k} = k/(m-(\gamma-1)k)$ and $m-(\gamma-1)k\geq(1-\kappa)m$.
 Therefore, 
\begin{align}\label{equation: fixed-size multiple hypercliques moment bound}
\ee_{S,S'\sim\dot\nu}[\exp\{\mu^2 |S\cap S'|\}] \leq \exp\bigg\{\frac{\gamma^2k^{t+1}}{t(1-\kappa)^tm^t}\big(\exp(\mu^2k^{t-1})-1\big)\bigg\},
\end{align}
A sufficient condition for $\bbE_{\mathbf Y\sim\bbQ}[L_{\dot\nu}(\mathbf Y)^2]\leq e^C$, where $C>0$ is fixed, is
$$
\mu\leq\sqrt{\frac{1}{k^{t-1}}\log\bigg(1+C\frac{t(1-\kappa)^tm^t}{\gamma^2 k^{t+1}}\bigg)}.
$$
On one hand, if $\gamma^2k^{t+1}=\omega(m^t)$, we have $C\frac{(1-\kappa)^tm^t}{\gamma^2k^{t+1}}\leq\log(1+C\frac{t(1-\kappa)^tm^t}{\gamma^2k^{t+1}})$ for all sufficiently large $m$. Therefore, the condition in the above equation holds when $\mu \leq {\sqrt C}(1-\kappa)^{t/2}\gamma^{-1} k^{-t} m^{t/2}.$
The same rate holds under the stated assumption $\gamma^2k^{t+1}/m^t=\Omega(1)$.
Indeed, if $\gamma^2k^{t+1}/m^t\geq c_1>0$ and $\mu\leq C_1m^{t/2}/(\gamma k^t)$ eventually, then $\mu^2k^{t-1}\leq C_1^2/c_1$.
Using $e^x-1\leq xe^x$ for $x\geq0$, we obtain
\[
\frac{\gamma^2k^{t+1}}{t(1-\kappa)^tm^t}\big(\exp(\mu^2k^{t-1})-1\big)
\leq\frac{C_1^2}{t(1-\kappa)^t}\exp(C_1^2/c_1)=O(1).
\]
On the other hand, in the regime $\gamma^2k^{t+1}=O(m^{t-\delta})$ for any $\delta>0$, since $\log(1+C\frac{t(1-\kappa)^tm^t}{\gamma^2k^{t+1}})\geq\delta\log m+O(1)$, we have $\log(1+C\frac{t(1-\kappa)^tm^t}{\gamma^2k^{t+1}})\geq \frac{\delta}{2}\log m$ for sufficiently large $m$. Therefore, it suffices to have $\mu \leq \sqrt{\frac{\delta}{2k^{t-1}} \log m}.$

\paragraph{Computational lower bound.}
The computational lower bound for the fixed-size case follows directly from the random-size result by the same stochastic dominance argument as in \Cref{proof: multiple cliques problem}.
Fix $c\in(0,1-\kappa)$ and set $k'=k/(1-c)$, so that $\gamma k'\leq\kappa m/(1-c)<m$.
Specifically, conditioning both independent copies of the random-size model on the high-probability event that every hyperclique has size at least $(1-c)k'$ yields
\begin{align*}
\ee_{\bfY\sim\bbQ}\big[L_{\dot\nu_{\gamma,(1-c)k'}}^{\leq D}(\bfY)^2\big]
\leq 2\ee_{\bfY\sim\bbQ}\big[L_{\tilde\nu_{\gamma,k'}}^{\leq D}(\bfY)^2\big].
\end{align*}
The joint conditioning event has probability at least $1-2\gamma\exp(-c^2k'/2)=1-o(1)$, and hence at least $1/2$ eventually.
Selecting a uniform $k$-subset of each hyperclique in each copy gives a coupling with the fixed-size model, with containment of the corresponding hyperedge sets.
This choice of $k'$ transfers the boundedness of the LDLR second moment from the random-size case, which we cover in \Cref{proof: multiple planted hypercliques problem random size}, to the fixed-size case, and thus gives the desired computational lower bound.
As in \Cref{proof: multiple fixed-size cliques problem}, fixed multiplicative factors in the regime conditions can be absorbed by replacing $\delta$ by any $\delta'\in(\epsilon,\delta)$ and decreasing the signal constant if necessary.

\subsubsection{Random-size case}\label{proof: multiple planted hypercliques problem random size}

We define $\tilde V_a:=\{v\in[m]\,:\,(a-1)k<\tilde A_v\leq ak\}$ and $S_a:=\{\bfi=(i_1,\dots,i_t)\in[m]_t: i_1,\dots,i_t\in \tilde V_a\}.$ Also, let $\tilde\bfA'$ be an independent copy of $\tilde\bfA$, and define $S'_1,...,S'_\gamma$ and $ \tilde V'_1,...,\tilde V'_\gamma$ analogously. Finally, we denote $S:=\bigcup_{a=1}^{\gamma}S_a$ and $S':=\bigcup_{b=1}^{\gamma}S'_b.$

\paragraph{Upper bounds.}
Note that the same conditioning approach as in \Cref{proof: random-size multiple cliques} can be directly applied.
In particular, we consider the event that each hyperclique has size at least $ck$ for some constant $c\in(0,1)$, and the risk difference between the original and conditional models is at most the probability of the complementary event.
By the union bound and Chernoff's inequality, this probability is at most $\gamma\exp\{-(1-c)^2k/2\}=o(1)$.
Then, we can derive the upper bounds by targeting the hypercliques of size $ck$, which differs from the fixed-size case only by constant factors.
The average test uses threshold $\mu\gamma\binom{\lfloor ck\rfloor}{t}/2$, and the scan searches over unions of $\gamma$ disjoint $\lfloor ck\rfloor$-vertex hypercliques.

\paragraph{Statistical lower bound.} 
The truncation argument in \Cref{proof: random-size multiple cliques} applies without modification to the present hyperclique setting. Indeed, the argument from \Cref{equation: truncation event for random-size multiple cliques} to \Cref{equation: comparison of random-size and fixed-size multiple cliques} only uses the sizes and inclusions of the underlying vertex sets, and does not depend on the hyperedge dimension $t$. In particular, conditioning on the event $\overline{\mathcal E}:=\{|\tilde V_a|\leq Ck,~\forall a\in[\gamma]\}$ for a fixed constant $C\in(1,1/\kappa)$, we may enlarge each $\tilde V_a$ to a set $\dot V_a\supseteq\tilde V_a$ of size exactly $Ck$. Note that this is indeed feasible, since $C\gamma k \leq C\kappa m <m$.
Here, $\overline{\tilde\nu}$ denotes the conditional prior, and $\tilde\nu(\overline{\mathcal E}^c)\leq\gamma\exp\{-(C-1)^2k/(2C)\}=o(1)$, so the truncation has vanishing cost in total variation.
 The resulting union of random-size hypercliques $S$ is then contained in the corresponding union of fixed-size hypercliques $\dot S$, and the same construction for an independent copy gives $|S\cap S'|\leq|\dot S\cap\dot S'|$ almost surely. Therefore,
\begin{align*}
\ee_{\bfY\sim\bbQ}[L_{\overline{\tilde\nu}}(\bfY)^2]&\leq \ee_{\dot S,\dot S'}\big[\exp{\mu^2|\dot S\cap\dot S'|}\big]\leq \exp\bigg\{\frac{C^{t+1}\gamma^2k^{t+1}}{t(1-C\kappa)^tm^t}\big(\exp(C^{t-1}\mu^2k^{t-1})-1\big)\bigg\},
\end{align*}
where the inequality follows from the fixed-size analysis (see \Cref{equation: fixed-size multiple hypercliques moment bound}) with $k$ replaced by $Ck$, and $\kappa$ replaced by $C\kappa$. Since $C\in(1,1/\kappa)$ is a fixed constant, this replacement only changes constant factors. Hence, repeating the same two cases as in the fixed-size analysis yields the statistical lower bound stated in \Cref{corollary: multiple planted hypercliques problem}.

\paragraph{Computational lower bound.} We take a similar approach as in \Cref{proof: random-size multiple cliques} to bound the LDLR second moment. In particular, we will show that the moments of the overlap $|S\cap S'|$ can be tightly bounded by the moments of an auxiliary random variable $Z_{\ast}$, which corresponds to the number of hyperedges in a random single hyperclique model.
Recall the definition of $W_v$ in \Cref{eq:definition of W_v}:
\begin{align}
W_v=\begin{cases}
e_{(a,b)} & \text{if }v\in\widetilde V_a\cap\widetilde V'_b\text{ for some }a,b\in[\gamma],\\
0 & \text{otherwise},
\end{cases}
\end{align}
where $\Pr(W_v=e_{(a,b)})=\rho=k^2/m^2$ for all $a,b\in[\gamma]$ and $\Pr(W_v=0)=1-\gamma^2\rho$. Denote $\bfi:=(i_1,\ldots,i_t)\in[m]_t$. Then, analogous to \Cref{equation: overlap in terms of W_v multiple cliques}, we can express the overlap $|S\cap S'|$ as
\begin{equation}
\begin{aligned}
|S\cap S'|
&=\sum_{a,b=1}^{\gamma}|S_a\cap S'_b|\\
&=\sum_{a,b=1}^{\gamma}\sum_{\bfi\in[m]_t}\mathds 1\big\{i_1,\ldots,i_t\in\widetilde V_a\cap\widetilde V'_b\big\}\\
&=\sum_{\bfi\in[m]_t}\sum_{a,b=1}^{\gamma}\prod_{l=1}^t\mathds 1\{W_{i_l}=e_{(a,b)}\}\\
&=\sum_{\bfi\in[m]_t}\mathds 1\{W_{i_1}=\cdots=W_{i_t}\neq 0\}\\
&=:Z,
\end{aligned}
\end{equation}
where $Z$ can be interpreted as the number of hyperedges of another random $\gamma^2$-multiple hyperclique model.
Let $B_1,\ldots,B_m\stackrel{\textrm{i.i.d.}}\sim\mathrm{Bernoulli}(p)$, where $p$ is the probability parameter we will specify later. Define $B:=\sum_{i=1}^m B_i\sim\mathrm{Binomial}(m,p)$, and let $Z_{\ast}:=\sum_{\bfi\in[m]_t}\prod_{l=1}^t B_{i_l}=\binom{B}{t}$. Let us denote $I_{\bfi}:=\mathds 1\{W_{i_1}=\cdots=W_{i_t}\neq 0\}$ and $J_{\bfi}:=\prod_{l=1}^t B_{i_l}$. Then, we have
\begin{align*}
\ee[Z]=\sum_{\bfi\in[m]_t}\ee[I_{\bfi}]=\sum_{\bfi\in[m]_t}\Pr(W_{i_1}=\cdots=W_{i_t}\neq 0)=\sum_{\bfi\in[m]_t}\gamma^2\rho^t=\binom mt \gamma^2\rho^t,
\end{align*}
and
\begin{align*}
\ee[Z_{\ast}]=\sum_{\bfi\in[m]_t}\ee[J_{\bfi}]=\sum_{\bfi\in[m]_t}\Pr(B_{i_1}=\cdots=B_{i_t}=1)=\sum_{\bfi\in[m]_t}p^t=\binom mt p^t.
\end{align*}
Choose $p=\gamma^{2/t}\rho$ so that $\ee[Z]=\ee[Z_{\ast}]$.
This is a valid Bernoulli parameter since $t\geq3$ and $\gamma k\leq m$ imply $p\leq\gamma\rho\leq k/m\leq1$.
Next, we claim that this choice of $p$ satisfies
$$
\ee[Z^d]=\sum_{\bfI\in([m]_t)^d}\Pr\bigg(\prod_{r=1}^d I_{\bfi_r}=1\bigg)\leq \sum_{\bfI\in([m]_t)^d}\Pr\bigg(\prod_{r=1}^d J_{\bfi_r}=1\bigg)=\ee[Z_{\ast}^d]
$$
for all $d\geq 1$, where $\bfI=(\bfi_1,\ldots,\bfi_d)$. We define $v=v(\bfI)$ as the number of distinct vertices in $\bfI$, and $c=c(\bfI)$ as the number of connected components in the hypergraph formed by the hyperedges $\bfi_1,\ldots,\bfi_d$. More specifically, we call two hyperedges $\bfi_r$ and $\bfi_s$ adjacent if they share at least one vertex, and we define a connected component as a maximal collection of hyperedges linked by chains of adjacent hyperedges. Then, for $\prod_{r=1}^d I_{\bfi_r}$ to be equal to 1, all vertices in each connected component must have the same nonzero label $e_{(a,b)}$ for some $a,b\in[\gamma]$. Since there are $\gamma^2$ possible labels, we have
$$
\Pr\bigg(\prod_{r=1}^d I_{\bfi_r}=1\bigg)=\gamma^{2c}\rho^{v}.
$$
See \Cref{equation: moment decomposition for random-size multiple cliques} and the arguments following it for a detailed explanation of the above equation.
On the other hand, for the auxiliary model, $\prod_{r=1}^d J_{\bfi_r}$ is 1 if and only if $B_i=1$ for every vertex in the selected hyperedges, an event of probability
$$
\Pr\bigg(\prod_{r=1}^d J_{\bfi_r}=1\bigg)=p^{v}=\gamma^{2v/t}\rho^v.
$$
Observe that every connected component contains at least one $t$-uniform hyperedge and therefore at least $t$ vertices. This implies that $v\geq tc$, and thus
$$
\Pr\bigg(\prod_{r=1}^d I_{\bfi_r}=1\bigg) =\gamma^{2c}\rho^{v}\leq \gamma^{2v/t}\rho^v = \Pr\bigg(\prod_{r=1}^d J_{\bfi_r}=1\bigg),
$$
which proves
$$
\ee[Z^d]=\sum_{\bfI\in([m]_t)^d}\Pr\bigg(\prod_{r=1}^d I_{\bfi_r}=1\bigg)\leq \sum_{\bfI\in([m]_t)^d}\Pr\bigg(\prod_{r=1}^d J_{\bfi_r}=1\bigg)=\ee[Z_{\ast}^d].
$$

Consequently, we have $\ee[|S\cap S'|^d]\leq \ee[Z_{\ast}^d]$ for all $d\geq 1$, and thus we can bound the LDLR second moment by
\begin{align*}
\ee_{\bfY\sim\bbQ}[L_{\tilde\nu}^{\leq D}(\bfY)^2]&=\ee_{S,S'\sim \tilde\nu}\bigg[\sum^D_{d=0} \frac{\mu^{2d}}{d!}|S\cap S'|^d\bigg]\leq 1+\sum_{d=1}^D\frac{\mu^{2d}}{d!}\ee[Z_{\ast}^d].
\end{align*}
Recall that $Z_{\ast}=\binom{B}{t}\leq B^t/t!$, and $B\sim\mathrm{Binomial}(m,p)$ with $p=\gamma^{2/t} k^2/m^2$. Since $t$ is a fixed constant, the problem reduces to bounding the $td$th moments of a binomial random variable, which can be done by applying the polynomial moment bound for binomial random variables \citep[Corollary 1]{ahle2022sharp} given in \Cref{lemma: moment bound of sub-Poissonian distribution}. The resulting bound is
\begin{equation}\label{equation: random size hyperclique final computational bound}
\begin{aligned}
1+\sum_{d=1}^D\frac{\mu^{2d}}{d!}\ee[Z_{\ast}^d] &\leq 1+\sum_{d=1}^D\frac{\mu^{2d}}{(t!)^dd!}\ee\big[B^{td}\big]\\
&\leq 1+\sum_{d=1}^D\frac{\mu^{2d}}{(t!)^dd!}\bigg(\frac{td}{\log(1+tdm/(\gamma^{2/t} k^2))}\bigg)^{td}.
\end{aligned}
\end{equation}
Let $\eta:=\gamma^{2/t} k^2$, and split the analysis into two cases:
\begin{itemize}
	\item \textbf{Case 1:} $\eta=O(m^{1-\delta})$ for some constant $\delta>0$.
	
	Note that 
	\begin{align*}
	\sum_{d=1}^D\frac{\mu^{2d}}{(t!)^dd!}\bigg(\frac{td}{\log(1+tdm/\eta)}\bigg)^{td}&\leq \sum_{d=1}^D\frac{t^{td}\mu^{2d}}{(t!)^dd!}\bigg(\frac{d}{\log(m/\eta)}\bigg)^{td}\\
	&=:\sum_{d=1}^D T_d,
	\end{align*}
	and this is bounded as long as (i) $T_1=O(1)$ and (ii) $T_{d+1}/T_d\leq1/2$ for all $d\in[D-1]$. The first condition is equivalent to $\mu=O(\log^{t/2}(m/\eta)).$ For the second condition, observe that the ratio between consecutive terms is
	\begin{align*}
	\frac{T_{d+1}}{T_d}&=\frac{t^{t(d+1)}\mu^{2(d+1)}}{(t!)^{d+1}(d+1)!}\bigg(\frac{d+1}{\log(m/\eta)}\bigg)^{t(d+1)}\cdot\frac{(t!)^dd!}{t^{td}\mu^{2d}}\bigg(\frac{\log(m/\eta)}{d}\bigg)^{td}\\
	&=\frac{t^t\mu^2}{(t!) (d+1)}\bigg(\frac{d+1}{d}\bigg)^{td}\bigg(\frac{d+1}{\log(m/\eta)}\bigg)^t\\
	&\leq \frac{t^t}{t!}\mu^2e^t\frac{D^{t-1}}{\log^t(m/\eta)},
	\end{align*}
	where the last inequality follows from the fact that $(1+1/d)^d\leq e$ for all $d\geq 1$.
This ratio is at most $1/2$ whenever
	$$
	\mu\leq \sqrt{\frac{t!}{2t^te^t}}\frac{\log^{t/2}(m/\eta)}{D^{(t-1)/2}}.
	$$
	Observe that the second condition is more stringent than the first one for any $D=\Omega(1)$. Also, since
	$\eta=O(m^{1-\delta})$ implies $\log(m/\eta)=\Theta(\log m)$, the stated choice $D=D_m$ makes the preceding condition hold whenever $\mu\leq C(\log m)^{\{1-(t-1)\epsilon\}/2}$ for a sufficiently small constant $C>0$ depending on $\delta$ and $t$.
Indeed, $\log(m/\eta)\geq(\delta/2)\log m$ eventually and $D_m\leq\log^{1+\epsilon}m$.
In particular, if $\mu\leq \log^{1/2-\xi}m$ for some constant $\xi>0$, any $\epsilon\in(0,\min\{\delta,2\xi/(t-1)\})$ makes the second condition hold eventually.

	\item \textbf{Case 2:} $\eta=\Omega(m\log^{1+\delta}m)$ for some constant $\delta>0$.
	Recall that $D=D_m=\lfloor\log^{1+\epsilon}m\rfloor$ with $0<\epsilon<\delta$.
Then, there exists $M>0$ such that for all $m\geq M$ and $d\in[D],$ we have $tdm/\eta <1$. Since $\log(1 + tx) \geq x$, for all $0 < x < 1$ and $t\geq2$, we have that, for all such $m \geq M$,
	$$
	\log\bigg(1+t\frac{dm}{\eta}\bigg) \geq \frac{dm}{\eta}.$$ 
    Using this bound, 
\begin{align*}
	\sum_{d=1}^D\frac{\mu^{2d}}{(t!)^dd!}\bigg(\frac{td}{\log(1+tdm/\eta)}\bigg)^{td}&\leq \sum_{d=1}^D \frac{t^{td}\mu^{2d}}{(t!)^dd!} \bigg(\frac{\eta}{m}\bigg)^{td}\\
	&=\sum_{d=1}^D \frac{1}{d!}\bigg(\frac{t^t\mu^{2}\eta^t}{t!m^t}\bigg)^d\\
	&\leq \exp\bigg(\frac{t^t\mu^{2}\eta^t}{t!m^t}\bigg).
\end{align*}
In turn, a sufficient condition for this series to remain bounded is
$$\mu = O\bigg(\frac{m^{t/2}}{\eta^{t/2}}\bigg)=O\bigg(\frac{m^{t/2}}{\gamma k^t}\bigg).$$
\end{itemize}
The upper and lower bounds yield the thresholds in \Cref{corollary: multiple planted hypercliques problem} at polynomial order, in the sense of \Cref{definition: statistical and computational thresholds}.
The computational conclusion is conditional on \Cref{conjecture: low-degree polynomial conjecture}.

\subsubsection{Proof of multiple $p$-partite planted hypercliques problem}\label{proof: multiple p-partite planted hypercliques problem}

In this section, we formally define the multiple $p$-partite planted hypercliques problem (\Cref{remark: multiple p partite hypercliques}) and verify that the statistical and computational thresholds are the same as those of the multiple planted hypercliques problem.

We consider a $p$-partite hypergraph with $p$ disjoint vertex classes, each containing $m$ vertices, so that the set of possible hyperedges is indexed by $[m]^p$. Let $\dot\bfA^{(1)},\ldots,\dot\bfA^{(p)}$ be independent random vectors, where each $\dot\bfA^{(j)}=(\dot A_1^{(j)},\ldots,\dot A_m^{(j)})$ is a uniformly random permutation of $[m]$ in the fixed-size case. 
For each $a\in[\gamma]$ and $j\in[p]$, define the vertex set
$$
\dot V_a^{(j)}:=\big\{v\in[m]:(a-1)k<\dot A_v^{(j)}\leq ak\big\}.
$$
The edge set of the $a$th planted $p$-partite hyperclique is defined by
$$
S_a:=\dot V_a^{(1)}\times\cdots\times \dot V_a^{(p)}=\big\{\bfi=(i_1,\ldots,i_p)\in[m]^p:i_j\in \dot V_a^{(j)}\text{ for all }j\in[p]\big\}.
$$
The latent structure is then given by
$$
S:=\bigcup_{a=1}^{\gamma}S_a.
$$
For the random-size case, we replace $\dot\bfA^{(1)},\ldots,\dot\bfA^{(p)}$ with $\tilde\bfA^{(1)},\ldots,\tilde\bfA^{(p)}$, where each $\tilde\bfA^{(j)}=(\tilde A_1^{(j)},\ldots,\tilde A_m^{(j)})$ is a vector of i.i.d.\ uniform random variables on the interval $[0,m]$.
We denote the induced distributions of $S$ in the fixed-size and random-size cases by $\dot\nu$ and $\tilde\nu$, respectively. In the fixed-size case, each $\dot V_a^{(j)}$ has cardinality $k$, and hence $|S|=\gamma k^p$.

Since the proof closely follows that of the original hyperclique problem (\Cref{proof: multiple planted hypercliques problem}), rather than repeating the entire argument from scratch, we focus on verifying that the key arguments remain valid under the present $p$-partite extension.
In the following, we will focus on the statistical lower bound for the fixed-size case and the computational lower bound for the random-size case, as the other cases can be handled analogously using truncation and stochastic dominance arguments similar to those in \Cref{proof: multiple planted hypercliques problem}. 

\paragraph{Statistical lower bound for the fixed-size case.}

Let $S'$ be an independent copy of $S$, and define $\dot V_a'^{(j)}$ and $S_a'$ analogously. The overlap between $S$ and $S'$ is given by
$$
|S\cap S'|=\sum_{a,b=1}^{\gamma}\prod_{j=1}^p|\dot V_a^{(j)}\cap \dot V_b'^{(j)}|.
$$
Our starting point is to confirm the following inequality, which is the counterpart of \Cref{equation: iterating conditional expectation argument for multiple fixed size cliques} and \Cref{equation: fixed-size multiple hypercliques moment upper bound}:
\begin{align}
\ee[\exp\{\mu^2 |S\cap S'|\}]\leq \ee\big[F(B_1,\ldots,B_p)\big]^{\gamma^2},
\end{align}
where $F(x_1,\ldots,x_p):=\exp\big(\mu^2\prod_{j=1}^p x_j\big)$ and $B_1,\ldots,B_p$ are independent binomial random variables with parameters $k$ and $q_{\gamma,k}:=k/(m-(\gamma-1)k)$. Throughout the proof, we assume that $\gamma k \leq \kappa m$ for some constant $\kappa\in(0,1)$, as the extension to the case $\gamma k = \Theta(m)$ follows from the covering argument in \Cref{proof: multiple planted hypercliques problem}.

To prove this inequality, we follow the argument in \Cref{equation: conditional expectation argument for multiple fixed size cliques} and the subsequent analysis.
Write $X_{ab}^{(j)}:=|\dot V_a^{(j)}\cap \dot V_b'^{(j)}|$ for $a,b\in[\gamma]$ and $j\in[p]$ for notational convenience.

The full array of these overlap counts need not be negatively associated, as the following calculation shows.
Fix $j\in[p]$ and suppose that $\gamma\geq2$, so that $2k\leq m$.
Conditional on the unprimed vertex sets, any two vertices $u\in\dot V_1^{(j)}$ and $v\in\dot V_2^{(j)}$ are distinct, and uniform sampling of the primed disjoint blocks gives
\[
\Pr\big\{u\in\dot V_1'^{(j)},\ v\in\dot V_2'^{(j)}\mid\dot\bfA^{(j)}\big\}
=\frac{k}{m}\frac{k}{m-1}.
\]
Summing over these vertices gives conditional moments that do not depend on the unprimed sets, so unconditionally
\[
\bbE[X_{11}^{(j)}]=\bbE[X_{22}^{(j)}]=\frac{k^2}{m},\qquad
\bbE[X_{11}^{(j)}X_{22}^{(j)}]=\frac{k^4}{m(m-1)}.
\]
Consequently,
\begin{equation}\label{eq:overlap-positive-covariance}
\operatorname{Cov}\big(X_{11}^{(j)},X_{22}^{(j)}\big)
=\frac{k^4}{m(m-1)}-\frac{k^4}{m^2}
=\frac{k^4}{m^2(m-1)}>0.
\end{equation}
This positive covariance rules out negative association of the full array, but does not affect the conditional negative association used below.

For each $b\in[\gamma]$, let $\mathcal F_{b-1}$ denote the sigma-field generated by $\{\dot V_a^{(j)}: a\in[\gamma], j\in[p]\}$ and $\{\dot V_l'^{(j)}: l\in[b-1], j\in[p]\}$. Conditionally on $\mathcal F_{b-1}$, each $V_b'^{(j)}$ is a uniformly random subset of size $k$ drawn without replacement from the remaining $M_b:=m-(b-1)k$ vertices. Since $\dot V_1^{(j)},\ldots,\dot V_\gamma^{(j)}$ are disjoint, the random variables $\{X_{ab}^{(j)}\}_{a=1}^\gamma$ are coordinate marginals of a multivariate hypergeometric distribution and hence are negatively associated. Moreover, the pairs of permutations used in different partitions are independent, and $\mathcal F_{b-1}$ records the revealed vertex sets separately within each partition. Consequently, conditional on $\mathcal F_{b-1}$, the $p$ random vectors $\bfX^{(j)}:=(X_{1b}^{(j)},\ldots,X_{\gamma b}^{(j)})$ for $j\in[p]$ are independent, while the coordinates within each vector are negatively associated. Based on these properties, for each $b\in[\gamma]$, we obtain
\begin{align*}
&\ee\bigg[\prod_{a=1}^{\gamma}F(X_{ab}^{(1)},\ldots,X_{ab}^{(p)})\,\bigg|\,\mathcal F_{b-1}\bigg]\\
\stackrel{(a)}=&\ee\bigg[\ee\bigg\{\prod_{a=1}^{\gamma}F(X_{ab}^{(1)},\ldots,X_{ab}^{(p)})\,\bigg|\,\mathcal F_{b-1},\bfX^{(1)},\ldots,\bfX^{(p-1)}\bigg\}\,\bigg|\,\mathcal F_{b-1}\bigg]\\
\stackrel{(b)}\leq&\ee\bigg[\prod_{a=1}^{\gamma}\ee\big[F(X_{ab}^{(1)},\ldots,X_{ab}^{(p)})\,\big|\,\mathcal F_{b-1},\bfX^{(1)},\ldots,\bfX^{(p-1)}\big]\,\bigg|\,\mathcal F_{b-1}\bigg]\\
\stackrel{(c)}=&\ee\bigg[\ee\bigg\{\prod_{a=1}^{\gamma}\ee[F(X_{ab}^{(1)},\ldots,X_{ab}^{(p)})\mid\mathcal F_{b-1},\bfX^{(1)},\ldots,\bfX^{(p-1)}]\,\bigg|\,\mathcal F_{b-1},\bfX^{(1)},\ldots,\bfX^{(p-2)}\bigg\}\,\bigg|\,\mathcal F_{b-1}\bigg]\\
\stackrel{(d)}\leq&\ee\bigg[\prod_{a=1}^{\gamma}\ee\big[F(X_{ab}^{(1)},\ldots,X_{ab}^{(p)})\,\big|\,\mathcal F_{b-1},\bfX^{(1)},\ldots,\bfX^{(p-2)}\big]\,\bigg|\,\mathcal F_{b-1}\bigg]\\
\leq&\cdots\leq\prod_{a=1}^{\gamma}\ee\big[F(X_{ab}^{(1)},\ldots,X_{ab}^{(p)})\,\big|\,\mathcal F_{b-1}\big],
\end{align*}
where each step is justified as follows. $(a)$ and $(c)$ are direct applications of the law of iterated expectations. For $(b)$, note that the function $F$ given $\mathcal F_{b-1}$ and $\bfX^{(1)},\ldots,\bfX^{(p-1)}$ is $F(X_{ab}^{(1)},\ldots,X_{ab}^{(p-1)},x)=\exp(\mu^2x\prod_{j=1}^{p-1}X_{ab}^{(j)})$, which is a nonnegative nondecreasing function of $x$. Then, due to the conditional independence of $\bfX^{(p)}$ from $\bfX^{(1)},\ldots,\bfX^{(p-1)}$ given $\mathcal F_{b-1}$, additional conditioning on $\bfX^{(1)},\ldots,\bfX^{(p-1)}$ does not change the law of $\bfX^{(p)}$, and hence the negative association of $\{X_{ab}^{(p)}\}_{a=1}^{\gamma}$ gives the corresponding inequality. 
For $(d)$, the argument is analogous to that for $(b)$, except that we need to verify that taking the conditional expectation of $F$ over its $p$th coordinate preserves the nonnegative and nondecreasing properties in its $(p-1)$th coordinate. 
Indeed, note that
\begin{align*}
&\ee\big[F(X_{ab}^{(1)},\ldots,X_{ab}^{(p)})\mid\mathcal F_{b-1},\bfX^{(1)},\ldots,\bfX^{(p-1)}\big]\\
&=\sum_{t=0}^{k}\exp\bigg(X_{ab}^{(p-1)}\mu^2\prod^{p-2}_{j=1}X_{ab}^{(j)}\cdot t\bigg)\Pr\big(X_{ab}^{(p)}=t\mid\mathcal F_{b-1},\bfX^{(1)},\ldots,\bfX^{(p-1)}\big)\\
&=\sum_{t=0}^{k}\exp\bigg(X_{ab}^{(p-1)}\mu^2\prod^{p-2}_{j=1}X_{ab}^{(j)}\cdot t\bigg)\Pr\big(X_{ab}^{(p)}=t\mid\mathcal F_{b-1}\big),
\end{align*}
where the last equality follows from the conditional independence of $\bfX^{(p)}$ from $\bfX^{(1)},\ldots,\bfX^{(p-1)}$ given $\mathcal F_{b-1}$. Observe that $\mu^2\prod^{p-2}_{j=1}X_{ab}^{(j)}$ is a fixed nonnegative constant when $\mathcal F_{b-1},\bfX^{(1)},\ldots,\bfX^{(p-2)}$ are given. Therefore, the last expression is a nonnegative linear combination of nondecreasing functions of $X_{ab}^{(p-1)}$, and hence is itself a nonnegative nondecreasing function of $X_{ab}^{(p-1)}$. This verifies the required property for $(d)$.

Next, for any $a,b\in[\gamma]$ and $j\in[p]$, from the proof of \Cref{equation: bound for conditional expectation of f for multiple fixed size cliques} we obtain:
for every nondecreasing convex function $\varphi$,
$$
\ee\big[\varphi(X_{ab}^{(j)})\,\big|\,\mathcal F_{b-1}\big]\leq\ee[\varphi(B_j)]
$$
almost surely. We take $B_1,\ldots,B_p$ to be mutually independent and independent of all the vertex selections. For fixed nonnegative values of the other coordinates, $F$ is an exponential function with a nonnegative coefficient in the remaining coordinate, and hence is nondecreasing and convex in that coordinate. By conditional independence across partitions, the preceding comparison can therefore be applied successively to each coordinate, yielding
$$
\begin{aligned}
\ee&\big[F(X_{ab}^{(1)},\ldots,X_{ab}^{(p)})\,\big|\,\mathcal F_{b-1}\big]\\
\leq&\ee\big[F(B_1,X_{ab}^{(2)},\ldots,X_{ab}^{(p)})\,\big|\,\mathcal F_{b-1}\big]\\
\leq&\cdots\leq\ee\big[F(B_1,\ldots,B_p)\,\big|\,\mathcal F_{b-1}\big]
=\ee[F(B_1,\ldots,B_p)].
\end{aligned}
$$

Combining this bound with the negative association argument, we conclude that, for every $b\in[\gamma]$,
$$
\ee\bigg[\prod_{a=1}^{\gamma}F(X_{ab}^{(1)},\ldots,X_{ab}^{(p)})\,\bigg|\,\mathcal F_{b-1}\bigg]\leq\prod_{a=1}^{\gamma}\ee\big[F(X_{ab}^{(1)},\ldots,X_{ab}^{(p)})\,\big|\,\mathcal F_{b-1}\big]\leq\ee[F(B_1,\ldots,B_p)]^\gamma.
$$
Therefore, it follows that
$$
\begin{aligned}
&\ee\bigg[\prod_{a=1}^{\gamma}\prod_{l=1}^{b}F(X_{al}^{(1)},\ldots,X_{al}^{(p)})\bigg]\\
=&\ee\bigg[\prod_{a=1}^{\gamma}\prod_{l=1}^{b-1}F(X_{al}^{(1)},\ldots,X_{al}^{(p)})\ee\bigg[\prod_{a=1}^{\gamma}F(X_{ab}^{(1)},\ldots,X_{ab}^{(p)})\,\bigg|\,\mathcal F_{b-1}\bigg]\bigg]\\
\leq&\ee[F(B_1,\ldots,B_p)]^\gamma\ee\bigg[\prod_{a=1}^{\gamma}\prod_{l=1}^{b-1}F(X_{al}^{(1)},\ldots,X_{al}^{(p)})\bigg].
\end{aligned}
$$
Iterating this inequality over $b=\gamma,\ldots,1$, with the empty product interpreted as $1$, gives
$$
\ee[\exp\{\mu^2|S\cap S'|\}]\leq\ee[F(B_1,\ldots,B_p)]^{\gamma^2},
$$
which proves the desired inequality.

Now, it remains to check whether applying Finner's inequality to the right-hand side of the above inequality yields a bound equivalent to that in \Cref{equation: fixed-size multiple hypercliques finner bound}.
Write $B_j=\sum_{i=1}^k Z_i^{(j)}$ for $j\in[p]$, where $\{Z_i^{(j)}:i\in[k],j\in[p]\}$ are i.i.d.\ $\mathrm{Bernoulli}(q_{\gamma,k})$ random variables. For each $\bfi=(i_1,\ldots,i_p)\in[k]^p$, define $Z_{\bfi}:=\prod_{j=1}^p Z_{i_j}^{(j)}$. Then,
$$
\ee[F(B_1,\ldots,B_p)]=\ee\bigg[\exp\bigg\{\mu^2\prod_{j=1}^p B_j\bigg\}\bigg]=\ee\bigg[\exp\bigg\{\mu^2\sum_{\bfi\in[k]^p}Z_{\bfi}\bigg\}\bigg].
$$
Since each $Z_i^{(j)}$ appears in exactly $k^{p-1}$ terms in the above summation, we choose the exponents $p_{\bfi}=k^{p-1}$ for all $\bfi\in[k]^p$. Applying Finner's inequality (\Cref{lemma: Finner's inequality}) gives
\begin{align*}
\ee[F(B_1,\ldots,B_p)]
&\leq\prod_{\bfi\in[k]^p}\ee\big[\exp\{\mu^2 k^{p-1}Z_{\bfi}\}\big]^\frac{1}{k^{p-1}}\\
&=\ee\big[\exp\{\mu^2 k^{p-1}Z_{\bfi}\}\big]^{k}\\
&=\big\{1+q_{\gamma,k}^p\big(\exp(\mu^2k^{p-1})-1\big)\big\}^{k}\\
&\leq\exp\big\{kq_{\gamma,k}^p\big(\exp(\mu^2k^{p-1})-1\big)\big\},
\end{align*}
where the second equality follows from $Z_{\bfi}\sim\mathrm{Bernoulli}(q_{\gamma,k}^p)$ and the last inequality follows from $1+x\leq e^x$.
Combining this bound with the preceding inequality and recalling that $q_{\gamma,k}=k/(m-(\gamma-1)k)$, we obtain
\begin{align}\label{equation: fixed-size multiple p-partite hypercliques moment bound}
\ee[\exp\{\mu^2|S\cap S'|\}]
&\leq\ee[F(B_1,\ldots,B_p)]^{\gamma^2}\nonumber\\
&\leq\exp\bigg\{\frac{\gamma^2k^{p+1}}{(1-\kappa)^pm^p}\big(\exp(\mu^2k^{p-1})-1\big)\bigg\},
\end{align}
provided that $\gamma k\leq\kappa m$ for some fixed $\kappa\in(0,1)$. Comparing this bound with \Cref{equation: fixed-size multiple hypercliques finner bound}, the resulting statistical lower bound for the multiple $p$-partite planted hypercliques problem equivalent to that of the multiple planted hypercliques problem can be obtained by following the remaining steps in \Cref{proof: multiple planted hypercliques problem}.

\paragraph{Computational lower bound for the random-size case.}

We take an approach similar to that in \Cref{proof: multiple planted hypercliques problem random size}. In particular, we interpret the overlap $|S\cap S'|$ as the number of hyperedges in a random $\gamma^2$-multiple $p$-partite hyperclique model and compare its moments with those of an auxiliary single $p$-partite hyperclique. Throughout this argument, $p\geq2$ is fixed and $\gamma k\leq m$.

Let $S'$ be an independent copy of $S$, and write $\tilde V_a^{(j)}$ and $\tilde V_b'^{(j)}$ for the corresponding vertex membership sets. For each $j\in[p]$ and $v\in[m]$, define
$$
W_v^{(j)}:=\begin{cases}
e_{(a,b)}, & \text{if }v\in\tilde V_a^{(j)}\cap\tilde V_b'^{(j)}\text{ for some }a,b\in[\gamma],\\
0, & \text{otherwise}.
\end{cases}
$$
Since the vertex assignments are independent across both vertices and partitions in the random-size model, the random variables $\{W_v^{(j)}:v\in[m],j\in[p]\}$ are independent, with
$$
\Pr(W_v^{(j)}=e_{(a,b)})=\rho:=\frac{k^2}{m^2},\qquad \Pr(W_v^{(j)}=0)=1-\gamma^2\rho.
$$
For a hyperedge $\bfi=(i_1,\ldots,i_p)\in[m]^p$, let $I_{\bfi}:=\mathds1\{W_{i_1}^{(1)}=\cdots=W_{i_p}^{(p)}\neq0\}$. Then,
\begin{align*}
|S\cap S'|
&=\sum_{a,b=1}^{\gamma}\sum_{\bfi\in[m]^p}\prod_{j=1}^p\mathds1\{i_j\in\tilde V_a^{(j)}\cap\tilde V_b'^{(j)}\}\\
&=\sum_{\bfi\in[m]^p}\mathds1\{W_{i_1}^{(1)}=\cdots=W_{i_p}^{(p)}\neq0\}\\
&=\sum_{\bfi\in[m]^p}I_{\bfi}\\
&=:Z.
\end{align*}
Thus, $Z$ can be seen as the number of hyperedges of another random $\gamma^2$-multiple $p$-partite hyperclique model.
For the auxiliary model, let $\{U_v^{(j)}:v\in[m],j\in[p]\}$ be i.i.d.\ $\mathrm{Bernoulli}(q_\ast)$ random variables, independent of all the vertex selections, where $q_\ast:=\gamma^{2/p}\rho$. Define
$$
B_j:=\sum_{v=1}^m U_v^{(j)},\qquad J_{\bfi}:=\prod_{j=1}^p U_{i_j}^{(j)},\qquad Z_\ast:=\sum_{\bfi\in[m]^p}J_{\bfi}=\prod_{j=1}^p B_j.
$$
Then $B_1,\ldots,B_p$ are independent $\mathrm{Binomial}(m,q_\ast)$ random variables, and $Z_\ast$ is the number of hyperedges in the resulting random single $p$-partite hyperclique. 

Next, we show that $\ee[Z^d]\leq\ee[Z_\ast^d]$ for every integer $d\geq1$. For an ordered tuple $\bfI=(\bfi_1,\ldots,\bfi_d)\in([m]^p)^d$, let $v=v(\bfI)$ be the number of distinct vertices appearing in the hyperedges $\bfi_1,\ldots,\bfi_d$ and let $c=c(\bfI)$ be the number of connected components in the hypergraph formed by them.  For $\prod_{r=1}^d I_{\bfi_r}$ to equal $1$, all vertices in each connected component must share the same nonzero label. There are $\gamma^2$ possible labels for each component, and different components may use the same label. By independence of the vertex labels, it follows that
$$
\Pr\bigg(\prod_{r=1}^d I_{\bfi_r}=1\bigg)=\gamma^{2c}\rho^v.
$$
On the other hand, $\prod_{r=1}^d J_{\bfi_r}=1$ precisely when all the $v$ distinct vertices are selected in the auxiliary model. Hence,
$$
\Pr\bigg(\prod_{r=1}^d J_{\bfi_r}=1\bigg)=q_\ast^v=\gamma^{2v/p}\rho^v.
$$
Every connected component contains at least one $p$-partite hyperedge and therefore at least $p$ vertices. Thus, $v\geq pc$, and consequently
$$
\Pr\bigg(\prod_{r=1}^d I_{\bfi_r}=1\bigg)=\gamma^{2c}\rho^v\leq\gamma^{2v/p}\rho^v=\Pr\bigg(\prod_{r=1}^d J_{\bfi_r}=1\bigg).
$$
Summing this comparison over all $\bfI\in([m]^p)^d$, we obtain
\begin{align}\label{equation: random-size multiple p-partite hypercliques moment comparison}
\ee[Z^d]
&=\sum_{\bfI\in([m]^p)^d}\Pr\bigg(\prod_{r=1}^d I_{\bfi_r}=1\bigg)\nonumber\\
&\leq\sum_{\bfI\in([m]^p)^d}\Pr\bigg(\prod_{r=1}^d J_{\bfi_r}=1\bigg)
=\ee[Z_\ast^d].
\end{align}

We now apply this moment comparison to the low-degree likelihood ratio. By \Cref{proposition: gaussian overlapping pairs},
\begin{align*}
\ee_{\bfY\sim\bbQ}\big[L_{\tilde\nu}^{\leq D}(\bfY)^2\big]
&=1+\sum_{d=1}^D\frac{\mu^{2d}}{d!}\ee[Z^d]\\
&\leq1+\sum_{d=1}^D\frac{\mu^{2d}}{d!}\ee[Z_\ast^d]
=1+\sum_{d=1}^D\frac{\mu^{2d}}{d!}\prod_{j=1}^p\ee[(B_j)^d],
\end{align*}
where the last equality follows from the independence of $B_1,\ldots,B_p$. Applying \Cref{lemma: moment bound of sub-Poissonian distribution} gives $\ee[(B_j)^d]\leq(d/{\log(1+dm/\gamma^{2/p}k^2)})^d,$ and consequently,
\begin{align*}
\ee_{\bfY\sim\bbQ}\big[L_{\tilde\nu}^{\leq D}(\bfY)^2\big]\leq1+\sum_{d=1}^D\frac{\mu^{2d}}{d!}\bigg(\frac{d}{\log(1+dm/\gamma^{2/p}k^2)}\bigg)^{pd}.
\end{align*}
The remaining technical steps to derive the computational lower bound from this inequality are identical to those in \Cref{equation: random size hyperclique final computational bound} and subsequent analysis, and hence we omit the details.

\subsection{\texorpdfstring{Proof of \Cref{theorem: thick cycle problem}}{Proof of Theorem \getrefnumber{theorem: thick cycle problem}}}\label{proof: thick cycle problem}
We first state the statistical and low-degree lower bounds to be established; the computational interpretation of the latter uses \Cref{conjecture: low-degree polynomial conjecture}.
\begin{itemize} 
	\item (Statistical lower bound) Fix any constant $\delta\in(0,0.1)$.
For either prior distribution $\dot\nu$ or $\tilde\nu$, strong detection is impossible under either of the following conditions: (i) $\ell^2\Delta m^{-2}\leq m^{-\delta}$ and $\mu\leq C_\delta\sqrt{\Delta^{-1}\log m}$ for a sufficiently small constant $C_\delta>0$; or (ii) $\ell^2\Delta m^{-2}=\Omega(1)$ and $\mu=O\big(\ell^{-1}\Delta^{-1}m\big)$.
	\item (Computational lower bound) Fix any constants $\delta\in(0,0.1)$ and $\epsilon\in(0,\delta)$, and let $D_m=\lfloor\log^{1+\epsilon}m\rfloor$.
For either prior distribution $\dot\nu$ or $\tilde\nu$, no polynomial of degree at most $D_m$ achieves strong separation under either of the following conditions: (i) $\ell\Delta m^{-1}\leq m^{-\delta}$ and $\mu\leq C(\log m)^{(1-\epsilon)/2}$ for a sufficiently small constant $C>0$; or (ii) $\ell\Delta m^{-1}\geq\log^{1+\delta}m$ and $\mu=O\big(\ell^{-1}\Delta^{-1}m\big)$.
\end{itemize}
We assume that $\Delta=\omega(\log m)$ throughout the proof, which is satisfied in the polynomial regime considered here.
We retain the standing assumption $2\Delta<\ell\leq m$ and take $D=D_m$ in the low-degree bounds.

Recall that our convention $[m]_2={(i,j):1\leq i<j\leq m}$ represents each undirected edge by listing its endpoints in increasing order to avoid double-counting. In the proofs below, with a slight abuse of notation, let us write $(i_a,i_b)$ to denote the undirected edge joining vertices $i_a$ and $i_b$ regardless of their ordering, instead of the formal representation $(\min\{i_a,i_b\},\max\{i_a,i_b\})$.

\subsubsection{Fixed-size case}

We consider the case where $\dot{\mathcal C}$ consists of all fixed-size thick cycles of length $\ell$ and thickness $\Delta$. In this case, we have $n=\binom m2$, $K=\ell\Delta$.

\paragraph{Upper bounds.}
The upper bounds on the critical signal strength $\mu^*$ follow from the tests listed in \Cref{proposition: gaussian test upper bound}:
\begin{itemize}
	\item For some constant $C>0$, the outlier test achieves $R_{\dot \nu}(\Delta_{\textrm{out}})\to 0$ whenever
	$\mu \geq C \sqrt{\log m}.$
	\item The average test achieves $R_{\dot \nu}(\Delta_{\textrm{avg}})\to 0$ when $\mu =\omega(\frac{m}{\ell\Delta}).$
	\item For the scan test, instead of detecting the entire thick cycle, we target a clique of size $\Delta+1$, which is a subgraph of the thick cycle. We apply the scan test with parameters $|\dot{\mathcal C}_{\text{clique}}|:=\binom m{\Delta+1}$ and $K_{\text{clique}}:=\binom{\Delta+1}2$. Then, for some constant $C>0$, the scan test achieves $R_{\dot \nu}(\Delta_{\textrm{scan}})\to 0$ when
	$$
	\mu \geq C \sqrt{\frac{1}{\Delta}\log\frac{m}{\Delta}}.
	$$
\end{itemize}
For the outlier and clique-scan guarantees, take $\delta_m=1/\log m$ in \Cref{proposition: gaussian test upper bound}; for the average test, use the risk bound from \Cref{proof: gaussian test upper bound} directly.

\paragraph{Argument for lower bounds.}
For both lower bounds, we cover the thick cycle by a carefully chosen collection of overlapping cliques.
Controlling their overlap allows us to apply the multiple-clique analysis with parameters $(\gamma,k)=(\ell/(2\Delta),2\Delta)$.

In detail, for the explicit construction below, assume that $q:=\ell/(2\Delta)$ is an integer. We will extend the analysis to the case where $\ell$ is not divisible by $2\Delta$ at the end of this subsection. Given $\dot A_1,...,\dot A_m$, sort the vertices in $[m]$ as $i_1,\ldots,i_m$ according to the values of $\dot A_1,...,\dot A_m$ in ascending order. Recalling the definition of the cycle distance $d_{\mathrm{cyc}}^{(\ell)}(a,b):=\min\{|a-b|,\ell-|a-b|\}$, the thick cycle $S$ is composed of edges whose endpoints have cyclic distance at most $\Delta$ along the base cycle $(i_1,i_2),\ldots,(i_{\ell-1},i_\ell),(i_\ell,i_1)$ with vertices $i_1,...,i_\ell$.
 Now consider two subsets of the vertices, both indexed by $r\in[q]$:
$$
V_{1,r}=\{i_{2\Delta(r-1)+1},\dots, i_{2\Delta r}\} \text{ for all } r\in [q],
$$
and 
$$
V_{2,r}=\{i_{2\Delta(r-1)+\Delta+1},\dots, i_{2\Delta r+\Delta}\} \text{ for all } r\in [q],
$$
where $i_{\ell+s}:=i_s$ for $s\in[\ell]$; that is, the indices are interpreted modulo $\ell$. 
To each vertex set $V_{1,r}$ and $V_{2,r}$ associate a corresponding clique $S_{1,r}=\{(i_a,i_b)\in[m]_2: i_a,i_b\in V_{1,r}\}$ and $S_{2,r}=\{(i_a,i_b)\in[m]_2: i_a,i_b\in V_{2,r}\}$, respectively, for all $r \in [q]$. We then form a collection of cliques $\mathcal S$, which we will refer to as the patchwork of $S$, as
\[
	\mathcal S:=\mathcal S_1\cup \mathcal S_2,\]
    where
    \[
    \mathcal S_1=\cup_{r=1}^q S_{1,r} \quad\text{and}\quad \mathcal S_2=\cup_{r=1}^q S_{2,r}.
\]
By construction, $\mathcal S_1$ and $\mathcal S_2$ each consist of $q$ disjoint cliques of size $2\Delta$ and the patchwork $\mathcal S$ covers the thick cycle $S$. Indeed, consider an edge $(i_a,i_b)\in S$.
Orient its shorter cyclic arc clockwise, exchanging the endpoint names if necessary and interpreting their indices cyclically.
If $i_a$ and $i_b$ are in the same vertex set, say $V_{1,r}$ for some $r\in[q]$, then $(i_a,i_b)$ is covered by the corresponding clique. Otherwise, $i_a$ and $i_b$ straddle the boundary of two consecutive vertex sets, say $V_{1,r}$ and $V_{1,r+1}$ for some $r\in[q]$ (where $V_{1,q+1}:=V_{1,1}$). Since $d_{\mathrm{cyc}}^{(\ell)}(a,b)\leq\Delta$, the two vertices are at most $\Delta$ positions apart along the base cycle. Therefore, this forces $i_a$ to be in the last $\Delta$ vertices of $V_{1,r}$ and $i_b$ to be in the first $\Delta$ vertices of $V_{1,r+1}$. In this case, both vertices $i_a$ and $i_b$ are contained in the vertex set $V_{2,r}$, and thus $(i_a,i_b)$ is covered by the corresponding clique $S_{2,r}$.

Now, let $S$ and $S'$ be independent copies, each uniformly distributed over the set of length-$\ell$, $\Delta$-thick cycles as defined in \Cref{eq:think.cycle.definition}, with corresponding patchworks $\mathcal S = \mathcal{S}_1 \cup \mathcal{S}_2$ and $\mathcal S' = \mathcal{S'}_1 \cup \mathcal{S'}_2$, respectively. Then, the overlap size $|S\cap S'|$ can be bounded as
\begin{equation}\label{equation: patchwork bound2}
\begin{aligned}
|S\cap S'| &\leq |(\mathcal S_1\cup \mathcal S_2)\cap (\mathcal S'_1\cup \mathcal S'_2)|\\
&\leq |\mathcal S_1\cap \mathcal S'_1| + |\mathcal S_1\cap \mathcal S'_2| + |\mathcal S_2\cap \mathcal S'_1| + |\mathcal S_2\cap \mathcal S'_2|.
\end{aligned}
\end{equation}
We remark that $\mathcal S_1$ and $\mathcal S_2$ are not independent, and neither are $\mathcal S'_1$ and $\mathcal S'_2$. Instead, each of them is marginally distributed as a uniform distribution over the set of all unions of $q$ disjoint cliques of size $2\Delta$, and $\mathcal S_1$ and $\mathcal S_2$ are independent of $\mathcal S'_1$ and $\mathcal S'_2$. Based on this observation, we analyze the statistical and computational thresholds of the thick cycle problem by bounding $|\mathcal S_i\cap \mathcal S'_j|$ for $i,j\in\{1,2\}$ using the results from the multiple planted cliques problem.

\paragraph{Statistical lower bound.} 
We have that
\begin{equation}\label{equation: statistical lower bound for thick cycle after Holder inequality}
\begin{aligned}
\ee[\exp(\mu^2 |S\cap S'|)] &\stackrel{(a)}\leq \ee\bigg[\exp\bigg(\mu^2 \bigg\{\sum_{i,j=1}^2 \big|\mathcal S_i\cap \mathcal S'_j\big|\bigg\}\bigg)\bigg] \\
&=  \ee\bigg[\prod_{i,j=1}^2\exp\big(\mu^2 \big|\mathcal S_i\cap \mathcal S'_j\big|\big)\bigg],\\
&\stackrel{(b)}\leq \prod_{i,j=1}^2\Big(\ee\big[\exp\big(4\mu^2 \big|\mathcal S_i\cap \mathcal S'_j\big|\big)\big]\Big)^{1/4}\\
&\stackrel{(c)}=\ee\big[\exp\big(4\mu^2 \big|\mathcal S_1\cap \mathcal S'_1\big|\big)\big],
\end{aligned}
\end{equation}
where (a) follows from \Cref{equation: patchwork bound2}, (b) from H\"older's inequality, and (c) from the fact that the overlaps $|\mathcal S_i\cap\mathcal S'_j|$ are identically distributed for all $i,j\in \{1,2\}.$
We leverage our analysis of the multiple cliques problem. For ease of analysis, let us assume that $\ell\leq \kappa m$ for some constant $\kappa\in(0,1)$. We will relax this restriction later.
 Then, we refer to \Cref{equation: iterating conditional expectation argument for multiple fixed size cliques} and \Cref{equation: Finner inequality application for multiple fixed size cliques} and obtain, with $B\sim\textrm{Binomial}(2\Delta,p_{q,2\Delta})$ where $p_{q,2\Delta}:=2\Delta/(m-(q-1)2\Delta)$,
\begin{equation}\label{equation: statistical lower bound moment bound for thick cycle}
\begin{aligned}
\ee[\exp(\mu^2 |S\cap S'|)] &\leq \ee\bigg[\exp\bigg\{4\mu^2\binom{B}{2}\bigg\}\bigg]^{q^2} \\
&\leq \exp\bigg\{\frac{q^2(2\Delta)^3}{2(1-\kappa)^2m^2}\big(\exp((2\mu)^2(2\Delta))-1\big)\bigg\}\\
&\leq \exp\bigg\{\frac{\ell^2\Delta}{(1-\kappa)^2m^2}\big(\exp(8\mu^2\Delta)-1\big)\bigg\},
\end{aligned}
\end{equation}
provided that $m-(q-1)2\Delta\geq (1-\kappa)m$.
Therefore, $\ee[\exp(\mu^2 |S\cap S'|)]\leq e^C$ for some $C>0$, as long as
$$
\mu \leq \sqrt{\frac{1}{8\Delta}\log\bigg(1+C\frac{(1-\kappa)^2m^2}{\ell^2\Delta}\bigg)}.
$$

In the regime $\ell^2\Delta=\omega(m^2)$, we have $C\frac{(1-\kappa)^2m^2}{2\ell^2\Delta}\leq\log(1+C\frac{(1-\kappa)^2m^2}{\ell^2\Delta})$ for all sufficiently large $m$. Therefore, the condition in the above equation holds as long as $\mu \leq \frac{\sqrt C}{4}\frac{(1-\kappa)m}{\ell\Delta}.$
The same rate holds under the stated assumption $\ell^2\Delta/m^2=\Omega(1)$.
Indeed, if $\ell^2\Delta/m^2\geq c_0>0$ and $\mu\leq C_0m/(\ell\Delta)$ eventually, then $8\mu^2\Delta\leq8C_0^2/c_0$.
Using $e^x-1\leq xe^x$ for $x\geq0$ gives
\[
\frac{\ell^2\Delta}{(1-\kappa)^2m^2}\big(\exp(8\mu^2\Delta)-1\big)
\leq\frac{8C_0^2}{(1-\kappa)^2}\exp(8C_0^2/c_0)=O(1).
\]
On the other hand, in the regime $\ell^2\Delta=O(m^{2-\delta})$ for any $\delta>0$, since $\log(1+C\frac{(1-\kappa)^2m^2}{\ell^2\Delta})\geq\delta\log m+O(1)$, we have $\log(1+C\frac{(1-\kappa)^2m^2}{\ell^2\Delta})\geq \frac{\delta}{2}\log m$ for sufficiently large $m$. Therefore, it suffices to have $\mu \leq \sqrt{\frac{\delta}{16\Delta} \log m}.$

For the general case $\ell\leq m$, we apply the grouping argument on $\mathcal S_1$ and $\mathcal S_1'$, as in \Cref{proof: multiple fixed-size cliques problem}. When $q=1$, $\mathcal S_1$ consists of a single clique of size $2\Delta$, so the statistical and computational moment bounds follow directly from the corresponding single-clique analyses. Hence, we assume $q\geq2$ in the grouping argument below.
In particular, note that $\mathcal S_1$ consists of $q$ cliques of size $2\Delta$, i.e., $\mathcal S_1 = \bigcup_{r=1}^q S_{1,r}$. Let $h=\lceil q/2 \rceil$. Define
\begin{align}\label{eq: multiple cliques grouping on thick cycle}
\mathcal T_{1}:=\bigcup_{r=1}^h S_{1,r},\qquad \mathcal T_{2} := \bigcup_{r=q-h+1}^q S_{1,r},
\end{align}
and define $\mathcal T_1'$ and $\mathcal T_2'$ similarly. This construction is analogous to that in \Cref{eq: multiple cliques grouping}.
Since $h\leq2q/3$, each grouped model satisfies $h(2\Delta)\leq2m/3$.
Therefore, following the subsequent analysis in \Cref{eq: multiple cliques decomposition} and \Cref{eq: statistical lower bound extension multiple cliques}, we obtain
\begin{equation}\label{eq: final statistical lower bound for fixed size thick cycle}
\begin{aligned}
	\ee_{\mathcal S_1,\mathcal S_1'\sim \dot \nu_{q,2\Delta}}\big[\exp\big(4\mu^2 \big|\mathcal S_1\cap \mathcal S'_1\big|\big)\big]&\leq \ee_{\mathcal T_1,\mathcal T_1'\sim \dot \nu_{h,2\Delta}}\big[\exp\big(16\mu^2 \big|\mathcal T_1\cap \mathcal T'_1\big|\big)\big]\\
	&\leq \exp\bigg\{\frac{2q^2 (2\Delta)^3}{m^2}\big(\exp(16\mu^2(2\Delta))-1\big)\bigg\}\\
	&= \exp\bigg\{\frac{4\ell^2\Delta}{m^2}\big(\exp(32\mu^2\Delta)-1\big)\bigg\}.
\end{aligned}
\end{equation}
The resulting bounds on $\mu$ are equivalent to those obtained in the previous analysis up to constant factors, and establish the statistical lower bound when $2\Delta$ divides $\ell$.

\paragraph{Computational lower bound.} Observe that 
\begin{equation}\label{equation: computational lower bound moment bound for thick cycle}
\begin{aligned}
	\ee_{\bfY\sim\bbQ}[L_{\dot\nu}^{\leq D}(\bfY)^2] &= 1+\ee\bigg[\sum^D_{d=1} \frac{\mu^{2d}}{d!}|S\cap S'|^d\bigg]\\
	&\stackrel{(a)}\leq 1+\ee\bigg[\sum^D_{d=1} \frac{\mu^{2d}}{d!}\bigg\{\sum_{i,j=1}^2 \big|\mathcal S_i\cap \mathcal S'_j\big|\bigg\}^d\bigg]\\
	&\stackrel{(b)}\leq 1+\ee\bigg[\sum^D_{d=1} \frac{4^{d-1}\mu^{2d}}{d!}\sum_{i,j=1}^2 \big|\mathcal S_i\cap \mathcal S'_j\big|^d\bigg]\\
	&\stackrel{(c)}= 1+\ee\bigg[\sum^D_{d=1} \frac{4^d\mu^{2d}}{d!} \big|\mathcal S_1\cap \mathcal S'_1\big|^d\bigg],
\end{aligned}
\end{equation}
where (a) follows from \Cref{equation: patchwork bound2}, (b) follows from the fact that $(\sum_{i=1}^4 x_i)^d \leq 4^{d-1} \sum_{i=1}^4 x_i^d$ for $x_i\geq 0$, and (c) follows from the fact that the overlaps $|\mathcal S_i\cap\mathcal S'_j|$ are identically distributed for all $i,j\in \{1,2\}.$ Suppose that $\ell\leq \kappa m$ for some constant $\kappa\in(0,1)$, and notice that the last term can be viewed as the second moment of the low-degree likelihood ratio of the fixed-size multiple cliques problem, which is \Cref{equation: low degree likelihood ratio multiple cliques}, with $2\mu$ instead of $\mu$, and with parameters $k=2\Delta$ and $\gamma=q=\ell/(2\Delta)$.
As shown in \Cref{equation: reduction to random-size multiple cliques}, the fixed-size bound transfers from the random-size model after inflating the clique-size parameter by a fixed factor.
Thus, up to a fixed multiplicative constant, the analysis reduces to a quantity of the form on the right-hand side of \Cref{equation: computational lower bound moment bound for random-size multiple cliques}, with $k=C(2\Delta)$, $\gamma=q=\ell/(2\Delta)$, and signal strength $2\mu$, where $C\in(1,1/\kappa)$ is fixed.
These substitutions change only constant factors in the sufficient signal conditions, and the analysis in \Cref{proof: random-size multiple cliques} applies.
We omit the details.

Similar to the statistical lower bound, for the general case $\ell\leq m$, we apply the grouping argument in \Cref{eq: multiple cliques grouping on thick cycle}. Following \Cref{eq: multiple cliques grouping on computational lower bound} and the subsequent analysis therein, we obtain a bound of the same order as in the preceding analysis.

\paragraph{Handling the case where $\ell$ is not divisible by $2\Delta$.}
Now, we extend the above analysis to the case where $\ell$ is not divisible by $2\Delta$. 
Let $t:=\ell-2\Delta q$ be the remainder, where $q:=\lfloor \ell/(2\Delta)\rfloor$ is the quotient. If we construct $\mathcal S_1$ and $\mathcal S_2$ as before, their union may fail to cover some edges of $S$ incident to the last $t$ vertices. Our main idea is to augment the patchwork by introducing two additional coverings, $\mathcal S_3$ and $\mathcal S_4$, constructed analogously to $\mathcal S_1$ and $\mathcal S_2$ but rotated to cover the remaining edges, so that
$$
S\subseteq\mathcal S:=\mathcal S_1\cup\mathcal S_2\cup\mathcal S_3\cup\mathcal S_4.
$$

Specifically, for each $r\in[q]$, recall the definition of $V_{1,r}$ and $V_{2,r}$ and additionally define $V_{3,r}$ and $V_{4,r}$ as
\begin{align*}
V_{1,r}&=\{i_{2\Delta(r-1)+1},\dots, i_{2\Delta r}\},\\
V_{2,r}&=\{i_{2\Delta(r-1)+\Delta+1},\dots, i_{2\Delta r+\Delta}\},\\
V_{3,r}&=\{i_{2\Delta q+2\Delta(r-1)+1},\dots,i_{2\Delta q+2\Delta r}\},\\
V_{4,r}&=\{i_{2\Delta q+2\Delta(r-1)+\Delta+1},\dots,i_{2\Delta q+2\Delta r+\Delta}\},
\end{align*}
where the indices are interpreted modulo $\ell$, as before.
For $a\in[4]$ and $r\in[q]$, define the corresponding clique
$$
S_{a,r}:=\{(u,v)\in[m]_2:u,v\in V_{a,r}\},
$$
and let $\mathcal S_a:=\bigcup_{r=1}^q S_{a,r}$ for $a\in[4]$.
By construction, each $\mathcal S_a$, $a\in[4]$, consists of $q$ disjoint cliques of size $2\Delta$. Also, we claim that $\mathcal S$ covers the thick cycle $S$. To see this, consider an arbitrary edge $(i_a,i_b)\in S$. Without loss of generality, assume that $i_a$ is the starting vertex, in the sense that the clockwise cyclic distance from $a$ to $b$ is at most $\Delta$. If $a$ lies among the first $2\Delta q$ vertices, following the same argument as before, $(i_a,i_b)$ is covered by $\mathcal S_1\cup\mathcal S_2$, as in the previous analysis. Otherwise, if $a$ lies between $2\Delta q+1$ and $\ell$, then $i_a$ should be in either $V_{3,1}$ or $V_{4,1}$, and again following the same argument as before, $(i_a,i_b)$ is covered by $\mathcal S_3\cup\mathcal S_4$. Thus, the resulting patchwork $\mathcal S=\bigcup_{a=1}^4\mathcal S_a$ covers the thick cycle $S$.

Now, let $S'$ be an independent copy of $S$, with corresponding patchwork $\mathcal S'=\bigcup_{a=1}^4\mathcal S'_a$. Then, we have
$$
|S\cap S'|\leq\sum_{a,b=1}^4|\mathcal S_a\cap\mathcal S'_b|.
$$
Observe that each $\mathcal S_a$ is marginally distributed as a uniform union of $q$ disjoint cliques of size $2\Delta$, and $\mathcal S_a$ and $\mathcal S'_b$ are independent for all $a,b\in[4]$. Therefore, the overlaps $|\mathcal S_a\cap\mathcal S'_b|$ are identically distributed for all $a,b\in[4]$.

For the statistical lower bound, following the same argument as before, we obtain
\begin{align*}
\ee\bigg[\exp\bigg(\mu^2\sum_{a,b=1}^4|\mathcal S_a\cap\mathcal S'_b|\bigg)\bigg]
\stackrel{(a)}\leq\prod_{a,b=1}^4\bigg(\ee\big[\exp\big(16\mu^2|\mathcal S_a\cap\mathcal S'_b|\big)\big]\bigg)^{1/16}
\stackrel{(b)}=\ee\big[\exp\big(16\mu^2|\mathcal S_1\cap\mathcal S'_1|\big)\big],
\end{align*}
where (a) follows from H\"older's inequality and (b) from the fact that the overlaps $|\mathcal S_a\cap\mathcal S'_b|$ are identically distributed for all pairs $(a,b)\in[4]^2$.
Similarly, for the computational lower bound, we have
\begin{align*}
\ee\bigg[\sum_{d=1}^D\frac{\mu^{2d}}{d!}
\bigg\{\sum_{a,b=1}^4|\mathcal S_a\cap\mathcal S'_b|\bigg\}^d\bigg]
\stackrel{(a)}\leq \ee\bigg[\sum_{d=1}^D\frac{16^{d-1}\mu^{2d}}{d!}
\sum_{a,b=1}^4|\mathcal S_a\cap\mathcal S'_b|^d\bigg]
\stackrel{(b)}=\ee\bigg[\sum_{d=1}^D\frac{16^d\mu^{2d}}{d!}
|\mathcal S_1\cap\mathcal S'_1|^d\bigg],
\end{align*}
where $(a)$ follows by applying the inequality
$(\sum_{i=1}^nx_i)^d\leq n^{d-1}\sum_{i=1}^nx_i^d$
for $x_i\geq0, n\in\mathbb N$ and $d\geq 1$, and $(b)$ follows from the identical distributions of the overlaps.

The preceding statistical and computational lower-bound analyses following \Cref{equation: statistical lower bound for thick cycle after Holder inequality} and \Cref{equation: computational lower bound moment bound for thick cycle} can then be applied to the new bounds, with only a change in the constant factors. This completes the extension to general $\ell$, which may not be divisible by $2\Delta$.

\subsubsection{Random-size case}

\paragraph{Upper bounds.} We introduce an event on which the number of active vertices in the random-size thick cycle is at least $c_1\ell$, and each active vertex has at least $c_2\Delta$ active neighbors in each direction, where $c_1,c_2\in(0,1)$ are fixed constants. Formally, we let $V:=\{i\in[m]:\tilde A_i\in(0,\ell]\}$ be the set of active vertices; then the number of active vertices satisfies $|V|\sim\mathrm{Binomial}(m,\ell/m)$. For each $i\in V$, we define
\begin{align*}
N^{+}_i&:=\sum_{j\neq i} \mathds 1\{\tilde A_{j}\in(0,\ell],~d_{\mathrm{cyc}}^{(\ell),+}(\tilde A_i,\tilde A_j)\leq \Delta\},\\
N^{-}_i&:=\sum_{j\neq i} \mathds 1\{\tilde A_{j}\in(0,\ell],~d_{\mathrm{cyc}}^{(\ell),-}(\tilde A_i,\tilde A_j)\leq \Delta\},
\end{align*}
where $d_{\mathrm{cyc}}^{(\ell),+}(a,b):=b-a~(\mathrm{mod}~\ell)$ and $d_{\mathrm{cyc}}^{(\ell),-}(a,b):=a-b~(\mathrm{mod}~\ell)$ are the clockwise and counterclockwise distances on a circle of circumference $\ell$ between $a$ and $b$, respectively.
We define an event
$$
\underline{\mathcal{E}}:=\Big\{|V|\geq c_1\ell,~\min_{i\in V} N^{+}_i\geq c_2\Delta,~\min_{i\in V} N^{-}_i\geq c_2\Delta\Big\},
$$
where $1>c_1,c_2>0$ are constants. Then, observe that the Chernoff bound guarantees the existence of a constant $c_3>0$ such that
$$
\Pr(|V|\geq c_1\ell)\geq 1-\exp(-c_3\ell),
$$
Conditional on $i\in V$ and its latent location, each of $N_i^+$ and $N_i^-$ has distribution $\mathrm{Binomial}(m-1,\Delta/m)$.
Hence, for some fixed $c_4>0$ and all sufficiently large $m$,
$$
\Pr(N^{+}_i\geq c_2\Delta\mid i\in V)\geq 1-\exp(-c_4\Delta),\quad \Pr(N^{-}_i\geq c_2\Delta\mid i\in V)\geq 1-\exp(-c_4\Delta).
$$
We apply the union bound to obtain
\begin{align*}
\Pr(\exists~i\in V: N^{+}_i< c_2\Delta) &\leq \sum^m_{i=1}\Pr(i\in V,~N^{+}_i< c_2\Delta)\\
&\leq m\frac{\ell}{m}\exp(-c_4\Delta)=\ell\exp(-c_4\Delta),
\end{align*}
and similarly, $\Pr(\exists~i\in V: N^{-}_i< c_2\Delta)\leq \ell\exp(-c_4\Delta)$. Therefore, we apply the union bound again to obtain
$$
\Pr(\underline{\mathcal{E}}^c)\leq \exp(-c_3\ell)+2\ell\exp(-c_4\Delta)=o(1),
$$
where the last equality follows from $\Delta=\omega(\log m)$.
Now, for an event $\underline{\mathcal{E}}$ defined in terms of the
auxiliary latent vector $\tilde{\mathbf A}$, 
define the induced conditional distribution of $S$ by
$$
\underline{\tilde\nu}(B):=\frac{\Pr\big(S(\tilde{\mathbf A})\in B,\underline{\mathcal{E}}\big)}{\Pr(\underline{\mathcal{E}})},
\qquad B\subseteq\mathcal C.
$$
The total variation distance between $\tilde\nu$ and $\underline{\tilde\nu}:=\tilde\nu(\cdot\mid\underline{\mathcal E})$ vanishes, as does the risk difference $|R_{\tilde\nu}(\Delta)-R_{\underline{\tilde\nu}}(\Delta)|$ for any test $\Delta$.
To apply the tests, set $h:=\lfloor c_2\Delta\rfloor$.
On $\underline{\mathcal E}$, each active vertex has at least $h$ neighbors in each of two disjoint oriented arcs, so $|S|\geq|V|h\geq c_1\ell h$.
Also, any active vertex and the next $h$ active vertices clockwise form a clique: they all lie in a clockwise arc of length at most $\Delta$.
Thus, the average test with threshold $\mu c_1\ell h/2$, the scan over $(h+1)$-vertex cliques, and the outlier test give the same upper bounds as in the fixed-size case, up to constant factors.

\paragraph{Truncation argument for lower bounds.}
For the lower bounds, we reduce the analysis of the random-size thick cycle case to the fixed-size thick cycle case via a truncation argument similar to that in \Cref{proof: random-size clique problem}. Specifically, let $V$ be the set of active vertices in the random-size thick cycle, i.e., $V:=\{i\in[m]:\tilde A_i\in(0,\ell]\}$. Then, observe that $|V|\sim\mathrm{Binomial}(m,\ell/m)$, and thus $\ee[|V|]=\ell$. Then, for any constant $C_1>1$, the Chernoff bound guarantees the existence of a constant $c_1>0$ such that
\begin{equation}\label{eq: truncation Chernoff first}
\Pr\big(|V|\geq C_1\ell\big)\leq \exp(-c_1\ell)=o(1),
\end{equation}
where $c_1$ only depends on $C_1$. Note that if $\ell$ is large, then $C_1\ell$ may exceed $m$. Therefore, for later use, let $r_{\mathrm{max}}:=\min\{\lfloor C_1\ell\rfloor,m\}$.

For each $i\in V$, let $N_i$ be the number of other active vertices within cyclic distance $\Delta$ of $i$, that is,
$$
N_i:=\sum_{j\neq i} \mathds 1\{\tilde A_{j}\in(0,\ell],~d_{\mathrm{cyc}}^{(\ell)}(\tilde A_i,\tilde A_j)\leq \Delta\},\quad i\in V,
$$
which is the degree of vertex $i$ in the thick cycle $S$. Since the total length of the two oriented arcs within distance $\Delta$ from $i$ is $2\Delta$, we have $N_i\mid\{i\in V\}\sim\mathrm{Binomial}(m-1,2\Delta/m)$, and thus $\ee[N_i\mid i\in V]=(m-1)2\Delta/m$. By the Chernoff bound, there exists a constant $c_2>0$ such that
$$
\Pr\big(N_i\geq 2C_2\Delta\,|\,i\in V\big)\leq \exp\big(-c_2\Delta\big),
$$
where $C_2>1$ is a constant that can be chosen arbitrarily. Taking a union bound over the $m$ possible vertex labels gives
\begin{equation}\label{eq: truncation Chernoff second}
\Pr\Big(\max_{i\in V} N_i\geq 2C_2\Delta\Big)\leq m\exp(-c_2\Delta)=o(1),
\end{equation}
where the last equality follows from $\Delta=\omega(\log m)$.
Define the event
$$
\mathcal{E}
:=\Big\{1\leq|V|\leq r_{\mathrm{max}},~\max_{i\in V} N_i\leq 2C_2\Delta\Big\}.
$$
The additional nonemptiness condition costs at most $\Pr(|V|=0)=(1-\ell/m)^m\leq e^{-\ell}=o(1)$.
Then, $\Pr(\mathcal{E})=1-o(1)$ by \Cref{eq: truncation Chernoff first} and \Cref{eq: truncation Chernoff second}. Define the conditional distribution of $S$ by
$$
\bar{\tilde \nu}(B):=\frac{\Pr\big(S(\tilde{\mathbf A})\in B,\mathcal{E}\big)}{\Pr(\mathcal{E})},
\qquad B\subseteq\mathcal C.
$$
Following the same truncation argument as in \Cref{proof: random-size clique problem}, the lower-bound analysis under the truncated mixing measure $\bar{\tilde \nu}$ is sufficient to establish the same lower bounds for the original mixing measure $\tilde \nu$.

We further refine the event $\mathcal{E}$ by conditioning on the number of active vertices $|V|=r$ for $r\in[r_{\mathrm{max}}]$.
All sums and suprema over these indices below are understood to include only conditional laws with positive mixing weight.
For $r\in[r_{\mathrm{max}}]$, consider the event $\mathcal{E}_r:=\{|V|=r,~\max_{i\in V} N_i\leq 2C_2\Delta\}$ and let $\bar{\tilde \nu}_r(\cdot):=\bar{\tilde \nu}(\cdot\,|\,\mathcal{E}_r)$ be the corresponding conditional mixing measure given $\mathcal{E}_r$. Also, let $S_r\sim\bar{\tilde \nu}_r$.
Given a conditional latent realization with edge set $s_r$, denote the active vertices as $v_1,...,v_r\in[m]$, arranged according to their latent locations on a circle of circumference $\ell$; that is, $0<\tilde A_{v_1}<\cdots<\tilde A_{v_r}\leq \ell$.
Write $b_\ast:=\lfloor2C_2\Delta\rfloor$ for the integer covering bandwidth.
As a covering of $s_r$, we consider the cyclic band graph
$$
\dot s_r:=\big\{(v_a,v_b)\in[m]_2\,:\,d_{\mathrm{cyc}}^{(r)}(a,b)\leq 2C_2\Delta\big\},
$$
where $d_{\mathrm{cyc}}^{(r)}(a,b):=\min\{|a-b|,~r-|a-b|\}$ is the cyclic distance on a circle of circumference $r$. We claim that $s_r\subseteq \dot s_r$. Indeed, for any edge $(v_a,v_b)\in s_r$, we have $d_{\mathrm{cyc}}^{(\ell)}(\tilde A_{v_a},\tilde A_{v_b})\leq \Delta$ by definition. Choose one of the two oriented arcs connecting $\tilde A_{v_a}$ and $\tilde A_{v_b}$ whose length is at most $\Delta$. Without loss of generality, assume that this arc from $\tilde A_{v_a}$ to $\tilde A_{v_b}$ is in the clockwise direction.
If $R$ is the number of cyclic rank steps along this oriented arc, then $d_{\mathrm{cyc}}^{(r)}(a,b)\leq R$; the opposite arc may have fewer rank steps.
Then, every active vertex along this arc after $v_a$ including $v_b$ must be within distance $\Delta$ from $v_a$, and thus $R\leq N_{v_a}$. Under the event $\mathcal{E}_r$, we have $N_{v_a}\leq 2C_2\Delta$, and thus $d_{\mathrm{cyc}}^{(r)}(a,b)\leq 2C_2\Delta$.

If $r>2b_\ast+1$, this cover is a fixed-size thick cycle with length $r$ and thickness $b_\ast$; if $r\leq2b_\ast+1$, it is the complete graph on the $r$ active vertices.
Applying the construction to the conditional latent variables defines the random cover $\dot S_r$.
Then, we have $S_r\subseteq \dot S_r$ almost surely.
It remains to verify that this cover has the uniform fixed-size law (or the uniform $r$-vertex clique law in the latter case).
To see this, it suffices to observe that conditionally on the event $\mathcal{E}_r$, the ordered active vertices $v_1,...,v_r$ are uniformly distributed over all possible ordered $r$-tuples of distinct vertices in $[m]$. Consider any two ordered $r$-tuples of distinct vertices $(i_1,...,i_r)$ and $(i'_1,...,i'_r)$ in $[m]$. Then, there exists a permutation $\pi$ on $[m]$ such that $\pi(i_j)=i'_j$ for all $j\in[r]$. Since the latent variables $\tilde A_1,...,\tilde A_m$ are exchangeable and the event $\mathcal {E}_r$ is invariant under vertex labeling, relabeling by $\pi$ gives
$$
\Pr((V_1,...,V_r)=(i_1,...,i_r))=\Pr((V_1,...,V_r)=(i'_1,...,i'_r)),
$$
where $V_1,...,V_r$ are the random indices defined by the ordered active vertices under the conditional latent law corresponding to $S_r\sim \bar{\tilde \nu}_r$.
Lastly, observe that $\dot S_r$ is a deterministic function of $(V_1,...,V_r)$, so it has the claimed uniform law.

\paragraph{Statistical lower bound.}
Let $S_r,S'_r\sim \bar{\tilde \nu}_r$ be independent, and let $\dot S_r,\dot S'_r$ be their corresponding covering fixed-size thick cycles. Then, $|S_r\cap S'_r|\leq |\dot S_r\cap \dot S'_r|$ almost surely, which implies 
\begin{align*}
	\ee[\exp(\mu^2 |S_r\cap S'_r|)] &\leq \ee[\exp(\mu^2 |\dot S_r\cap \dot S'_r|)].
\end{align*}
We next aggregate the above bound over $r\in[r_{\mathrm{max}}]$. Observe that the truncated mixing measure $\bar{\tilde \nu}$ can be expressed as a mixture of the conditional mixing measures $\bar{\tilde \nu}_r$, i.e., $\bar{\tilde \nu}(\cdot)=\sum_{r=1}^{r_{\mathrm{max}}}\Pr(|V|=r\,|\,\mathcal{E})\bar{\tilde \nu}_r(\cdot)$. Let $S,S'\sim \bar{\tilde \nu}$ be independent, and let $S_r,S'_r\sim \bar{\tilde \nu}_r$ be independent for each $r\in[r_{\mathrm{max}}]$. Then,
\begin{align}
	\ee_{\bf Y\sim \bbQ}\big[L_{\bar{\tilde \nu}}(\mathbf Y)^2\big] &= \ee_{\bfY\sim \bbQ}\bigg[\bigg(\sum_{r=1}^{r_{\mathrm{max}}}\Pr(|V|=r\,|\,\mathcal{E})L_{\bar{\tilde \nu}_r}(\mathbf Y)\bigg)^2\bigg]\notag\\
	&=\sum_{r=1}^{r_{\mathrm{max}}}\sum_{r'=1}^{r_{\mathrm{max}}}\Pr(|V|=r\,|\,\mathcal{E})\Pr(|V|=r'\,|\,\mathcal{E})\ee_{\bfY\sim \bbQ}\big[L_{\bar{\tilde \nu}_r}(\mathbf Y)L_{\bar{\tilde \nu}_{r'}}(\mathbf Y)\big]\notag\\
	&\stackrel{(\dagger)}\leq \sum_{r=1}^{r_{\mathrm{max}}}\sum_{r'=1}^{r_{\mathrm{max}}}\Pr(|V|=r\,|\,\mathcal{E})\Pr(|V|=r'\,|\,\mathcal{E})\sqrt{\ee_{\bfY\sim \bbQ}\big[L_{\bar{\tilde \nu}_r}(\mathbf Y)^2\big]\ee_{\bfY\sim \bbQ}\big[L_{\bar{\tilde \nu}_{r'}}(\mathbf Y)^2\big]}\notag\\
	&\leq \sum_{r=1}^{r_{\mathrm{max}}}\sum_{r'=1}^{r_{\mathrm{max}}}\Pr(|V|=r\,|\,\mathcal{E})\Pr(|V|=r'\,|\,\mathcal{E}) \sup_{r\in[r_{\mathrm{max}}]}\ee_{\bfY\sim \bbQ}\big[L_{\bar{\tilde \nu}_r}(\mathbf Y)^2\big]\notag\\
	&= \sup_{r\in[r_{\mathrm{max}}]}\ee_{\bfY\sim \bbQ}\big[L_{\bar{\tilde \nu}_r}(\mathbf Y)^2\big]\notag\\
	&= \sup_{r\in[r_{\mathrm{max}}]}\ee[\exp(\mu^2 |S_r\cap S'_r|)]\notag\\
	&\leq \sup_{r\in[r_{\mathrm{max}}]}\ee[\exp(\mu^2 |\dot S_r\cap \dot S'_r|)].\label{equation: statistical lower bound decomposition for random-size thick cycle}
\end{align}
where $(\dagger)$ follows from the Cauchy--Schwarz inequality. Therefore, $\ee_{\bf Y\sim \bbQ}\big[L_{\bar{\tilde \nu}}(\mathbf Y)^2\big]$ is bounded as long as $\ee[\exp(\mu^2 |\dot S_r\cap \dot S'_r|)]$ is uniformly bounded in $r\in[r_{\mathrm{max}}]$.
For $r>2b_\ast+1$, apply the four-component extension of \Cref{eq: final statistical lower bound for fixed size thick cycle} with length $r$ and thickness $b_\ast\leq2C_2\Delta$.
The H\"older factor increases from $4$ to $16$, so the exponential constant $32$ in that bound becomes $128$.
Here $q=\lfloor r/(2b_\ast)\rfloor$ satisfies $2b_\ast q\leq r$; the single-clique case $q=1$ is bounded by the same expression.
For $2\leq r\leq2b_\ast+1$, the single-clique bound in \Cref{eq:Finner inequality application} instead gives an upper bound $\exp\{r^3(\exp(\mu^2r)-1)/(2m^2)\}$, which is no larger than the bound below since $r\leq4C_2\Delta+1\leq5C_2\Delta$.
For $r=1$, the overlap is zero.
Thus, in each case we obtain the uniform estimate
\begin{align*}
\ee[\exp(\mu^2 |\dot S_r\cap \dot S'_r|)]&\leq \exp\bigg(\frac{8C_2r^2\Delta}{m^2}\big(\exp(256\mu^2C_2\Delta)-1\big)\bigg).
\end{align*}
Since the displayed upper bound is increasing in $r$, it is bounded above by its value at $r=r_{\mathrm{max}}$.
Therefore, we have
\begin{align*}
\sup_{r\in[r_{\mathrm{max}}]}\ee[\exp(\mu^2 |\dot S_r\cap \dot S'_r|)]&\leq \exp\bigg(\frac{8C_1^2C_2\ell^2\Delta}{m^2}\big(\exp(256C_2\mu^2\Delta)-1\big)\bigg),
\end{align*}
since $r_{\mathrm{max}}\leq C_1\ell$.
The analysis of this bound only differs from the analysis we conducted in the fixed-size case by constant factors. Thus, without any divisibility assumption, the statistical lower bound for the random-size thick cycle problem is identical to the fixed-size one.

\paragraph{Computational lower bound.}
Note that a decomposition analogous to that in \Cref{equation: statistical lower bound decomposition for random-size thick cycle} can be applied to the low-degree likelihood ratio, which yields
\begin{align*}
\ee_{\bfY\sim \bbQ}\big[L_{\bar{\tilde \nu}}^{\leq D}(\mathbf Y)^2\big] &\leq \sup_{r\in[r_{\mathrm{max}}]}\ee_{\bfY\sim \bbQ}\big[L_{\bar{\tilde \nu}_r}^{\leq D}(\mathbf Y)^2\big]\\
	&= 1+\sup_{r\in[r_{\mathrm{max}}]}\ee\bigg[\sum^D_{d=1} \frac{\mu^{2d}}{d!}|S_r\cap S'_r|^d\bigg]\\
	&\leq 1+\sup_{r\in[r_{\mathrm{max}}]}\ee\bigg[\sum^D_{d=1} \frac{\mu^{2d}}{d!}|\dot S_r\cap \dot S'_r|^d\bigg].
\end{align*}
We apply the four-component extension of \Cref{equation: computational lower bound moment bound for thick cycle} from the fixed-size analysis, allowing every conditional length $r$, to obtain
\begin{align*}
\sup_{r\in[r_{\mathrm{max}}]}\ee\bigg[\sum^D_{d=1} \frac{\mu^{2d}}{d!}|\dot S_r\cap \dot S'_r|^d\bigg] &\leq\sup_{r\in[r_{\mathrm{max}}]}\ee\bigg[\sum^D_{d=1} \frac{16^d\mu^{2d}}{d!}|\mathcal S_{1,r}\cap \mathcal S'_{1,r}|^d\bigg],
\end{align*}
where $\mathcal S_{1,r}$ and $\mathcal S'_{1,r}$ are the first patchwork components of $\dot S_r$ and $\dot S'_r$, respectively.
For $r\leq2b_\ast+1$, take these components to be the covering cliques themselves, so the displayed inequality still holds.
For the four-component construction with $r>2b_\ast+1$, each component is a union of $q=\lfloor r/(2b_\ast)\rfloor$ disjoint cliques of size $k=2b_\ast$, and hence has the fixed-size multiple-clique law with parameters $(\gamma,k)=(q,2b_\ast)$.
For fixed $k$, coupling these models using the first $q$ blocks of a common random permutation shows that the overlap moments are nondecreasing in the integer $q$.
Also, $qk^2\leq2b_\ast r\leq4C_1C_2\ell\Delta$, so the multiple-clique moment bound is controlled uniformly by the same binomial series as in the fixed-size analysis, with its overlap parameter changed only by a constant factor.
For the covering cliques when $r\leq2b_\ast+1$, use \Cref{equation: computational lower bound moment bound}; their corresponding parameter satisfies $r^2=O(\Delta^2)=O(\ell\Delta)$, giving the same uniform control.
The four-component construction therefore gives the same computational lower-bound rates for all conditional lengths $r$, without any divisibility assumption.

\subsection{\texorpdfstring{Proof of \Cref{corollary: thick path problem}}{Proof of Corollary \getrefnumber{corollary: thick path problem}}}\label{proof: thick path problem}
In general, the formal upper and lower bounds we aim to establish are the same as those for the thick cycle problem; see \Cref{proof: thick cycle problem}. In this section, we show how the thick path problem reduces to the thick cycle problem.
As in \Cref{proof: thick cycle problem}, assume $2\Delta<\ell\leq m$ and $\Delta=\omega(\log m)$ throughout.

\paragraph{Upper bounds.}
We note that the analysis of the upper bounds for the thick path problem is nearly identical to that of the thick cycle problem in \Cref{proof: thick cycle problem}, with minor modifications due to the absence of wrap-around edges. We address such modifications in the following.

For the fixed-size case, the outlier and clique-scan tests apply without change.
For the average test, the number of edges is $\ell\Delta-\frac{\Delta(\Delta+1)}{2}$ instead of $\ell\Delta$.
Under $2\Delta<\ell$, this count lies between $\ell\Delta/2$ and $\ell\Delta$, so the same upper-bound rate follows.

For the random-size case, we modify the truncation event to handle the vertices close to the endpoints of the interval $(0,\ell]$. Specifically, we modify the definition of $N_i^{+}$ and $N_i^{-}$ for each $i\in V$ as follows:
\begin{align*}
	N_i^{+}&:=\sum_{j\neq i} \mathds 1\{\tilde A_{j}\in(0,\ell],~0\leq d_{\mathrm{lin}}^{(\ell),+}(\tilde A_i,\tilde A_j)\leq \Delta\},\\
	N_i^{-}&:=\sum_{j\neq i} \mathds 1\{\tilde A_{j}\in(0,\ell],~0\leq d_{\mathrm{lin}}^{(\ell),-}(\tilde A_i,\tilde A_j)\leq \Delta\},
\end{align*}
where $d_{\mathrm{lin}}^{(\ell),+}(a,b):=b-a$ and $d_{\mathrm{lin}}^{(\ell),-}(a,b):=a-b$ are the signed displacements on the interval $(0,\ell]$ between $a$ and $b$ in the right and left directions, respectively. Also, to handle the vertices close to the boundary, define
$$
V^+:=\{i\in V\,:\,\tilde A_i\leq \ell-\Delta\},\quad V^-:=\{i\in V\,:\,\tilde A_i\geq \Delta\}.
$$
Consider the truncation event
$$
\underline{\mathcal{E}}:=\Big\{|V|\geq c_1\ell,~\min_{i\in V^+} N^{+}_i\geq c_2\Delta,~\min_{i\in V^-} N^{-}_i\geq c_2\Delta\Big\},
$$
where $0<c_2<c_1<1$ are fixed constants, and let $\underline{\tilde \nu}:=\tilde \nu(\cdot\,|\,\underline{\mathcal{E}})$.
We interpret a minimum over an empty set as $+\infty$.
We claim that under this event, $S$ contains a thick path with parameters $\lfloor c_1\ell\rfloor$ and $\lfloor c_2\Delta\rfloor$.
To see this, let $0<\tilde A_{v_1}<\cdots<\tilde A_{v_{|V|}}\leq \ell$ be the ordered latent locations of the active vertices for $S\sim\underline{\tilde \nu}$.
Write $h:=\lfloor c_2\Delta\rfloor$, and $r:=\lfloor c_1\ell\rfloor$.
Since $c_2<c_1$ and $2\Delta<\ell$, we have $r\geq2h$ and $h\to\infty$.
Then, observe that
$\dot S:=\{(v_a,v_b):1\leq a<b\leq r,~b-a\leq h\}$ 
is a thick path with parameters $r$ and $h$, and it suffices to show that $\dot S\subseteq S$ almost surely.
For $a<b$ satisfying $b-a\leq h$:
\begin{itemize}
	\item \textbf{Case 1:} $\tilde A_{v_a}\leq \ell-\Delta$. In this case, the event $\underline{\mathcal E}$ guarantees that there are at least $h$ other active vertices in $(\tilde A_{v_a},\tilde A_{v_a}+\Delta]$, and thus $\tilde A_{v_b}\leq \tilde A_{v_a}+\Delta$, provided that $b-a\leq h$. Therefore, $(v_a,v_b)\in S$.
	\item \textbf{Case 2:} $\tilde A_{v_b}\geq\Delta$. Similar to Case 1, there are at least $h$ other active vertices in $[\tilde A_{v_b}-\Delta,\tilde A_{v_b})$, and thus $\tilde A_{v_a}\geq \tilde A_{v_b}-\Delta$, provided that $b-a\leq h$. Therefore, $(v_a,v_b)\in S$.
	\item \textbf{Case 3:} $\tilde A_{v_a}>\ell-\Delta$ and $\tilde A_{v_b}<\Delta$. Note that this case can only occur if $\Delta>\ell/2$. In this case, $\tilde A_{v_b}-\tilde A_{v_a}<\Delta$, and thus $(v_a,v_b)\in S$.
Thus, this case does not arise under the standing assumption $2\Delta<\ell$.
\end{itemize}
Therefore, we have $\dot S\subseteq S$ almost surely, and the upper bound analysis for detecting $\dot S$ is identical to that of the fixed-size thick path problem up to constant factors.
Indeed, $\dot S$ has $rh-h(h+1)/2\geq rh/2$ edges and contains an $(h+1)$-vertex clique, with $r=\Theta(\ell)$ and $h=\Theta(\Delta)$.
Under the original latent law, conditional on $i\in V^+$ and its latent location, $N_i^+$ has distribution $\mathrm{Binomial}(m-1,\Delta/m)$; the analogous conditional law holds for $N_i^-$ when $i\in V^-$.
Chernoff bounds and a union bound therefore give, for some fixed $c_3,c_4>0$,
\[
\Pr(\underline{\mathcal E}^c)\leq\exp(-c_3\ell)+2m\exp(-c_4\Delta)=o(1).
\]
The comparison of risks using total variation in \Cref{proof: random-size clique problem} then transfers these upper bounds from the truncated prior to $\tilde\nu$.

\paragraph{Lower bounds.}
The lower bound follows from a coupling with the thick cycle problem. For either choice of the auxiliary vector $\dot{\mathbf A}$ or $\tilde{\mathbf A}$, a thick path of length $\ell$
$$
S_{\mathrm{path}}:=\big\{(i,j)\in[m]_2:A_i,A_j\in(0,\ell],~d_{\mathrm{lin}}^{(\ell)}(A_i,A_j)\leq \Delta\big\}
$$
is a subgraph of the corresponding thick cycle of length $\ell$
$$
S_{\mathrm{cycle}}:=\big\{(i,j)\in[m]_2:A_i,A_j\in(0,\ell],~d_{\mathrm{cyc}}^{(\ell)}(A_i,A_j)\leq \Delta\big\}.
$$
Indeed, recall that $d_{\mathrm{lin}}^{(\ell)}(a,b):=|a-b|$ and $d_{\mathrm{cyc}}^{(\ell)}(a,b):=\min\{|a-b|,~\ell-|a-b|\}$.
Therefore, for any $i,j\in[m]$, if $d_{\mathrm{lin}}^{(\ell)}(A_i,A_j)\leq \Delta$, then $d_{\mathrm{cyc}}^{(\ell)}(A_i,A_j)\leq \Delta$, which in turn guarantees that $S_{\mathrm{path}}\subseteq S_{\mathrm{cycle}}$.
Using an independent copy of the same latent vector to define the primed paths and cycles, we obtain
$$
|S_{\mathrm{path}}\cap S'_{\mathrm{path}}|\leq |S_{\mathrm{cycle}}\cap S'_{\mathrm{cycle}}|
$$
almost surely.
The full and low-degree likelihood-ratio second moments are expectations of increasing functions of this overlap, so the moment bounds in \Cref{proof: thick cycle problem} transfer to thick paths.
In the random-size case, the same coupling remains valid when both latent copies are conditioned on the high-probability event used for the thick-cycle lower bounds.
Thus, the statistical and computational lower bounds established for the thick cycle problem in \Cref{theorem: thick cycle problem} imply the same bounds for the thick path problem.

\subsection{\texorpdfstring{Proof of \Cref{proposition: upper bound for NEF models}}{Proof of Proposition \getrefnumber{proposition: upper bound for NEF models}}}\label{proof: upper bound for NEF models}
Throughout this section, we use the notation
\begin{align*}
	m_{0}:=\psi'(\theta_{0}),\quad m_{1}:=\psi'(\theta_{1}),\quad v_0:=\psi''(\theta_0)\mathrel{>}0,\\
	d_{\theta}:=m_1-m_0>0,\quad \tau_n:=d_\theta/\sqrt{v_0}.
\end{align*}
Let $\delta=\delta_n\in(0,1)$, and retain the assumption $|S|=K\geq1$ almost surely under $S\sim\nu$.
Recall the local regularity condition in \Cref{def: local regularity of the NEF models}: we say that the NEF models $P_{\theta_0}$ and $P_{\theta_1}$ are locally regular at scale $r_n$ if $ [\theta_{j}-r_n,\theta_{j}+r_n]\subseteq \Theta,~ j\in\{0,1\},$ and there exists a constant $C\geq 1$ such that
$$ \sup_{\substack{j\in\{0,1\}\\|u|\leq r_n}}\frac{\psi''(\theta_j+u)}{\psi''(\theta_0)}\leq C,$$
for sufficiently large $n$.
Note that $C$ is no less than 1, since the supremum includes $j=0$ and $u=0$.
In the following lemma, we show that the centered sum of independent samples from the NEF models satisfies a local sub-Gaussian tail bound under the local regularity condition.
\begin{lemma}[Local concentration for NEF models]\label{lemma: local sub-Gaussian concentration of the NEF models}
	Suppose that the NEF models $P_{\theta_0}$ and $P_{\theta_1}$ are locally regular at scale $r_n$ with constant $C\geq1$.
For $N=N_n\geq 1$, let $Y_1,\ldots,Y_N$ be independent samples such that $Y_i\sim P_{\theta^{(i)}}$ for $\theta^{(i)}\in\{\theta_0,\theta_1\}$.
Then, for all sufficiently large $n$ and every $0\leq x\leq CNv_0 r_n,$
	the centered sum satisfies the following tail bounds:
	\begin{align*}
		\Pr\bigg\{\sum^N_{i=1}\big(Y_i-\psi'(\theta^{(i)})\big)\geq x\bigg\}&\leq \exp\bigg(-\frac{x^2}{2CNv_0}\bigg),\\
		\Pr\bigg\{\sum^N_{i=1}\big(Y_i-\psi'(\theta^{(i)})\big)\leq -x\bigg\}&\leq \exp\bigg(-\frac{x^2}{2CNv_0}\bigg).
	\end{align*}
	
\end{lemma}
\proof{
	Note that the NEF models allow the following equality for the cumulant generating function:
$$
\log\big(\ee_{Y\sim P_\theta}\big[\exp\big\{t[Y-\psi'(\theta)]\big\}\big]\big)=\psi(\theta+t)-\psi(\theta)-t\psi'(\theta),\quad \theta,\theta+t\in\Theta,
$$
which can also be rewritten as the remainder term in the Taylor expansion:
$$
\psi(\theta+t)-\psi(\theta)-t\psi'(\theta)=t^2\int_0^1 (1-s)\psi''(\theta+st)ds.
$$
Under the local regularity condition in \Cref{def: local regularity of the NEF models}, we have an upper bound on the curvature of the log-partition function in the neighborhood of $\theta_0$ and $\theta_1$:
$$
\sup_{\substack{\theta\in\{\theta_0,\theta_1\}\\|u|\leq r_n}} \psi''(\theta+u)\leq C v_0.
$$
For $\theta\in\{\theta_0,\theta_1\}$, the centered variable $Y-\psi'(\theta)$ is sub-exponential with parameters $(Cv_0,1/r_n)$ in the following sense: for $|t|\leq r_n$, its log moment generating function satisfies
\begin{align*}
	\log\big(\ee_{Y\sim P_\theta}\big[\exp\big\{t[Y-\psi'(\theta)]\big\}\big]\big)&\mathrel{=} t^2\int_0^1 (1-s)\psi''(\theta+st)ds\\
	&\leq \frac{Cv_0 t^2}{2}.
\end{align*} 
By the independence of $Y_1,\ldots,Y_N$, the sum $\sum_{i=1}^N (Y_i-\psi'(\theta^{(i)}))$ is sub-exponential with parameters $(CNv_0,1/r_n)$; its log moment generating function is bounded by $CNv_0t^2/2$ for $|t|\leq r_n$.
For $0\leq x\leq CNv_0r_n$, choose $t=x/(CNv_0)\leq r_n$ in the Chernoff bound to obtain
\[
\Pr\bigg\{\sum_{i=1}^N\big(Y_i-\psi'(\theta^{(i)})\big)\geq x\bigg\}
\leq\exp\bigg(-tx+\frac{CNv_0t^2}{2}\bigg)
=\exp\bigg(-\frac{x^2}{2CNv_0}\bigg).
\]
Applying the same argument to the negative of the centered sum gives the lower-tail bound.

\qed
}

Thus, the proof procedure for the average, scan, and degree tests consists of the following two steps: (i) we specify the scale $r_n$ for the local regularity of the NEF models, and (ii) we check that each target deviation lies in the range $x\leq CNv_0r_n$ of \Cref{lemma: local sub-Gaussian concentration of the NEF models} and is no larger than the corresponding gap between the mean and the test threshold.
The gap itself may exceed the local concentration range; in that case, we use event inclusion to apply the lemma at the smaller target deviation.

\subsubsection{Average test}
For the average test, define
$$
r_{\textrm{avg},n}:=\sqrt{\frac{2}{nv_0}\log\frac{2}{\delta}},
$$
and assume that the NEF models are locally regular at scale $r_{\textrm{avg},n}$ with constant $C_{\textrm{avg}}$. Then, the average test $\Delta_{\mathrm{avg}}$ achieves $R_\nu(\Delta_{\mathrm{avg}})\leq \delta$ if
$$
\tau_n\geq \sqrt{\frac{8C_{\textrm{avg}}n}{K^2}\log\frac{2}{\delta}}.
$$

To prove the claim, we control the Type~I and Type~II errors of the average test separately by $\delta/2$.
Recall that the average test statistic and the threshold for the average test are defined as follows:
$$
T_{\mathrm{avg}}(\bfY):=\sum_{i=1}^n Y_i, \quad t_{\mathrm{avg}}:=nm_{0}+\frac{Kd_\theta}{2}.
$$

\paragraph{Type I error.}
Under the null hypothesis $H_0$, we have $Y_i\sim P_{\theta_0}$ independently for $i=1,\dots,n$, which implies that the test statistic $\sum_{i=1}^n Y_i$ has mean $nm_0$. Therefore,
$$
\bbQ\big\{T_{\mathrm{avg}}(\bfY)>t_{\mathrm{avg}}\big\}= \Pr_{Y_i\sim P_{\theta_0}}\bigg\{\sum_{i=1}^n (Y_i-m_0)> \frac{Kd_\theta}{2}\bigg\}.
$$ 
If the local regularity condition holds for the sequence $r_{\mathrm{avg},n}$ with a constant $C_{\mathrm{avg}}\geq 1$, then for any $0\leq x\leq C_{\mathrm{avg}}nv_0 r_{\mathrm{avg},n}$, the local sub-Gaussian concentration inequality in \Cref{lemma: local sub-Gaussian concentration of the NEF models} guarantees
$$
\Pr_{Y_i\sim P_{\theta_0}}\bigg\{\sum_{i=1}^n (Y_i-m_0)\geq x\bigg\}\leq \exp\bigg(-\frac{x^2}{2C_{\mathrm{avg}}nv_0}\bigg).
$$
To make the right-hand side less than or equal to $\delta/2$, define a target deviation $x_\delta:=\sqrt{2C_{\mathrm{avg}}nv_0\log(2/\delta)}$, so that $\exp(-x_\delta^2/(2C_{\mathrm{avg}}nv_0))=\delta/2$.
 Then, if 
$$
\frac{Kd_\theta}{2}\geq x_\delta,
$$
we have
$$
\bigg\{\sum_{i=1}^n (Y_i-m_0)> \frac{Kd_\theta}{2}\bigg\}\subseteq \bigg\{\sum_{i=1}^n (Y_i-m_0)> x_\delta\bigg\},
$$
which implies that the Type I error is bounded by $\delta/2$. The condition $\frac{Kd_\theta}{2}\geq x_\delta$ is equivalent to
$$
\tau_n\geq \sqrt{\frac{8C_{\textrm{avg}}n}{K^2}\log\frac{2}{\delta}}.
$$
Finally, it remains to check that the target deviation $x_\delta$ satisfies the condition $x_\delta\leq C_{\mathrm{avg}}nv_0 r_{\mathrm{avg},n}$ in \Cref{lemma: local sub-Gaussian concentration of the NEF models}, which holds as long as 
$$
r_{\mathrm{avg},n}\geq \sqrt{\frac{2}{C_{\mathrm{avg}}nv_0}\log\frac{2}{\delta}}.
$$
Since $C_{\mathrm{avg}}\geq 1$, by the definition of $r_{\mathrm{avg},n}=\sqrt{2/(nv_0)\log(2/\delta)}$, the condition is satisfied.

\paragraph{Type II error.}
Fix any $s$ in the support of $\nu$, and define $\theta_s^{(i)}:=\theta_1$ for $i\in s$ and $\theta_s^{(i)}:=\theta_0$ for $i\notin s$. Under $\bbP_s$, we have $\bbE_{\bfY\sim\bbP_s}[\sum_{i=1}^n Y_i]=\sum_{i=1}^n \psi'(\theta_s^{(i)})=(n-K)m_0+Km_1=nm_0+Kd_\theta$. Therefore, the Type~II error for this fixed $s$ satisfies
$$
\bbP_s\big\{T_{\mathrm{avg}}(\bfY)\leq t_{\mathrm{avg}}\big\}= \Pr_{Y_i\sim P_{\theta_s^{(i)}}}\bigg\{\sum_{i=1}^n (Y_i-\psi'(\theta_s^{(i)}))\leq -\frac{Kd_\theta}{2}\bigg\}.
$$
We repeat the same argument as in the Type I error case, and the only difference is that we apply the local sub-Gaussian concentration inequality in \Cref{lemma: local sub-Gaussian concentration of the NEF models} to the lower tail.
Consequently, the Type II error is also bounded by $\delta/2$ under the same condition $$
\tau_n\geq \sqrt{\frac{8C_{\textrm{avg}}n}{K^2}\log\frac{2}{\delta}}.
$$
Because the bound is uniform over all $s$ in the support of $\nu$,
$$
\bbE_{S\sim\nu}\big[\bbP_S\{\Delta_{\mathrm{avg}}(\bfY)=0\}\big]\leq\frac{\delta}{2}.
$$

Combining the Type I and Type II error bounds, we conclude that the average test achieves $R_\nu(\Delta_{\mathrm{avg}})\leq \delta$ under the stated condition on $\tau_n$.

\subsubsection{Scan test}\label{proof: scan test for NEF models}
For the scan test, define
$$
r_{\textrm{scan},n}:=\sqrt{\frac{2}{Kv_0}\log\frac{2|\mathcal C|}{\delta}},
$$
and assume that the NEF models are locally regular at scale $r_{\textrm{scan},n}$ with constant $C_{\textrm{scan}}$. Consider the scan statistic and the threshold defined as
$$
T_{\mathrm{scan}}(\mathbf{Y}) = \max_{s \in \mathcal{C}} \sum_{l \in s} Y_l,\quad t_{\mathrm{scan}} = Km_0 + \frac{Kd_\theta + \sqrt{2KC_{\textrm{scan}}v_0\log|\mathcal C|}}{2}.
$$
Then, the scan test $\Delta_{\mathrm{scan}}:=\mathds 1\{T_{\mathrm{scan}}(\mathbf{Y})>t_{\mathrm{scan}}\}$ achieves $R_\nu(\Delta_{\mathrm{scan}})\leq \delta$ if
$$
\tau_n\geq \sqrt{\frac{2C_{\textrm{scan}}}{K}}\bigg(2\sqrt{\log\frac{2}{\delta}}+\sqrt{\log|\mathcal C|}\bigg).
$$

We control the Type~I and Type~II errors separately as follows.

\paragraph{Type I error.} Note that for any fixed $s\in\mathcal C$, the sum $\sum_{l \in s} Y_l$ has mean $Km_0$ under the null hypothesis $H_0$. Therefore, by the union bound, we obtain
\begin{align*}
	\bbQ\big\{T_{\mathrm{scan}}(\mathbf{Y})>t_{\mathrm{scan}}\big\}&\leq \sum_{s\in\mathcal C} \bbQ\bigg\{\sum_{l \in s} Y_l > t_{\mathrm{scan}}\bigg\}\\
	&=\sum_{s\in\mathcal C} \Pr_{Y_l\sim P_{\theta_0}}\bigg\{\sum_{l \in s} (Y_l-m_0)> \frac{Kd_\theta + \sqrt{2KC_{\textrm{scan}}v_0\log|\mathcal C|}}{2}\bigg\}\\
	&=|\mathcal C|\cdot \Pr_{Y_l\sim P_{\theta_0}}\bigg\{\sum_{l \in s} (Y_l-m_0)> \frac{Kd_\theta + \sqrt{2KC_{\textrm{scan}}v_0\log|\mathcal C|}}{2}\bigg\}.
\end{align*} 
In the last line, $s$ is any fixed member of $\mathcal C$; all candidate sums have the same distribution under $\bbQ$ because each candidate has size $K$.
No independence across candidates is needed for the union bound.
If the local regularity condition holds for the sequence $r_{\mathrm{scan},n}$ with a constant $C_{\mathrm{scan}}\geq 1$, then for any $0\leq x\leq C_{\mathrm{scan}}Kv_0 r_{\mathrm{scan},n}$, the local sub-Gaussian concentration inequality in \Cref{lemma: local sub-Gaussian concentration of the NEF models} guarantees
$$
\Pr_{Y_l\sim P_{\theta_0}}\bigg\{\sum_{l \in s} (Y_l-m_0)> x\bigg\}\leq \exp\bigg(-\frac{x^2}{2C_{\mathrm{scan}}Kv_0}\bigg).
$$
We aim to make the right-hand side less than or equal to $\delta/(2|\mathcal C|)$, and thus define a target deviation $x_{\delta,\mathrm{I}}:=\sqrt{2C_{\mathrm{scan}}Kv_0\log(2|\mathcal C|/\delta)}$.
Then, it remains to verify two conditions.

First, the target deviation $x_{\delta,\mathrm{I}}$ should satisfy the condition $x_{\delta,\mathrm{I}}\leq C_{\mathrm{scan}}Kv_0 r_{\mathrm{scan},n},$ equivalently, we require
$$
r_{\mathrm{scan},n}\geq \sqrt{\frac{2}{C_{\mathrm{scan}}Kv_0}\log\frac{2|\mathcal C|}{\delta}}.
$$
By the definition of $r_{\mathrm{scan},n}=\sqrt{2/(Kv_0)\log(2|\mathcal C|/\delta)}$, the condition is satisfied since $C_{\mathrm{scan}}\geq 1$.

Second, we require 
\begin{align}\label{eq:scan-type-I-condition}
\frac{Kd_\theta + \sqrt{2KC_{\textrm{scan}}v_0\log|\mathcal C|}}{2}\geq x_{\delta,\mathrm{I}},
\end{align}
so that we can bound the Type I error as
\begin{align*}
	\bbQ\big\{T_{\mathrm{scan}}(\mathbf{Y})>t_{\mathrm{scan}}\big\}&\leq |\mathcal C|\cdot \Pr_{Y_l\sim P_{\theta_0}}\bigg\{\sum_{l \in s} (Y_l-m_0)> \frac{Kd_\theta + \sqrt{2KC_{\textrm{scan}}v_0\log|\mathcal C|}}{2}\bigg\}\\
	&\leq |\mathcal C|\cdot \Pr_{Y_l\sim P_{\theta_0}}\bigg\{\sum_{l \in s} (Y_l-m_0)> x_{\delta,\mathrm{I}}\bigg\}\\
	&\leq |\mathcal C|\cdot \frac{\delta}{2|\mathcal C|}=\frac{\delta}{2}.
\end{align*}
Substituting the definition of $x_{\delta,\mathrm{I}}$ into \Cref{eq:scan-type-I-condition} gives the equivalent condition
$$
Kd_\theta\geq \sqrt{2KC_{\textrm{scan}}v_0}\bigg(2\sqrt{\log\frac{2|\mathcal C|}{\delta}}-\sqrt{\log|\mathcal C|}\bigg),
$$
and since $d_\theta/\sqrt{v_0}=\tau_n$, this holds as long as
\begin{align}\label{eq:scan-type-I-condition-mu}
\tau_n\geq \sqrt{\frac{2C_{\textrm{scan}}}{K}}\bigg(2\sqrt{\log\frac{2|\mathcal C|}{\delta}}-\sqrt{\log|\mathcal C|}\bigg).
\end{align}

\paragraph{Type II error.} Fix any $s^\star$ in the support of $\nu$. Under $\bbP_{s^\star}$, the mean of $\sum_{l \in s^\star}Y_l$ is $Km_1=Km_0+Kd_\theta$. Therefore, the Type~II error for this fixed $s^\star$ satisfies
\begin{align*}
	\bbP_{s^\star}\big\{T_{\mathrm{scan}}(\mathbf{Y})\leq t_{\mathrm{scan}}\big\}&\leq \bbP_{s^\star}\bigg\{\sum_{l \in s^\star} Y_l \leq t_{\mathrm{scan}}\bigg\}\\
	&= \Pr_{Y_l\sim P_{\theta_1}}\bigg\{\sum_{l \in s^\star} (Y_l-m_1)\leq -\frac{Kd_\theta - \sqrt{2KC_{\textrm{scan}}v_0\log|\mathcal C|}}{2}\bigg\}.
\end{align*}
Again, if the local regularity condition holds for the sequence $r_{\mathrm{scan},n}$ with a constant $C_{\mathrm{scan}}\geq 1$, then for any $0\leq x\leq C_{\mathrm{scan}}Kv_0 r_{\mathrm{scan},n}$, the local sub-Gaussian concentration inequality in \Cref{lemma: local sub-Gaussian concentration of the NEF models} guarantees
$$
\Pr_{Y_l\sim P_{\theta_1}}\bigg\{\sum_{l \in s^\star} (Y_l-m_1)\leq -x\bigg\}\leq \exp\bigg(-\frac{x^2}{2C_{\mathrm{scan}}Kv_0}\bigg).
$$
To control the Type II error by $\delta/2$, the local sub-Gaussian concentration inequality suggests the target deviation $x_{\delta,\mathrm{II}}:=\sqrt{2C_{\mathrm{scan}}Kv_0\log(2/\delta)}$, so that $\exp(-x_{\delta,\mathrm{II}}^2/(2C_{\mathrm{scan}}Kv_0))=\delta/2$, and it remains to check the following two conditions.

First, we verify that the target deviation $x_{\delta,\mathrm{II}}$ lies in the local concentration regime, i.e., $x_{\delta,\mathrm{II}}\leq C_{\mathrm{scan}}Kv_0 r_{\mathrm{scan},n}$, which is equivalent to
$$
r_{\mathrm{scan},n}\geq \sqrt{\frac{2}{C_{\mathrm{scan}}Kv_0}\log\frac{2}{\delta}}.
$$
By the definition of $r_{\mathrm{scan},n}=\sqrt{2/(Kv_0)\log(2|\mathcal C|/\delta)}$, the condition is satisfied since $C_{\mathrm{scan}}\geq 1$ and $|\mathcal C|\ge 1$, hence the local sub-Gaussian concentration lemma applies.
Also, it remains to show that
$$
\frac{Kd_\theta - \sqrt{2KC_{\textrm{scan}}v_0\log|\mathcal C|}}{2}\geq x_{\delta,\mathrm{II}},
$$
so that the Type II error is bounded by $\delta/2$. 
Note that it suffices to have
$$
Kd_\theta\geq \sqrt{2KC_{\textrm{scan}}v_0}\bigg(2\sqrt{\log\frac{2}{\delta}}+\sqrt{\log|\mathcal C|}\bigg),
$$
which is equivalent to
\begin{align}\label{eq:scan-type-II-condition-mu}
\tau_n\geq \sqrt{\frac{2C_{\textrm{scan}}}{K}}\bigg(2\sqrt{\log\frac{2}{\delta}}+\sqrt{\log|\mathcal C|}\bigg).
\end{align}
Therefore, under this condition, the Type~II error is at most $\delta/2$ uniformly over all $s^\star$ in the support of $\nu$, and hence
$$
\bbE_{S\sim\nu}\big[\bbP_S\{\Delta_{\mathrm{scan}}(\bfY)=0\}\big]\leq\frac{\delta}{2}.
$$

Finally, we combine the Type I and Type II error controls; observe that the signal strength condition for the Type II error control in \Cref{eq:scan-type-II-condition-mu} is stronger than the condition for the Type I error control in \Cref{eq:scan-type-I-condition-mu}, since
\begin{align*}
	2\sqrt{\log\frac{2}{\delta}}+\sqrt{\log|\mathcal C|}&=2\bigg(\sqrt{\log\frac{2}{\delta}}+\sqrt{\log|\mathcal C|}\bigg)-\sqrt{\log|\mathcal C|}\\
	&\geq 2\sqrt{\log|\mathcal C|+\log\frac{2}{\delta}}-\sqrt{\log|\mathcal C|},
\end{align*}
where the last inequality follows from the fact that $\sqrt{a}+\sqrt{b}\geq \sqrt{a+b}$ for all $a,b\geq 0$. Therefore, if the condition in \Cref{eq:scan-type-II-condition-mu} is satisfied, both Type I and Type II errors are controlled by $\delta/2$, and we conclude that the scan test achieves $R_\nu(\Delta_{\mathrm{scan}})\leq \delta$ under the stated condition on $\tau_n$ and the local regularity condition.

\subsubsection{Degree test}\label{proof: degree test for NEF models} 
For the degree test, recall the scale and threshold from \Cref{proposition: upper bound for NEF models}:
\[
\begin{aligned}
r_{\mathrm{deg}}&:=\sqrt{\frac{2}{(m-1)v_0}\log\frac{2m}{\delta}},\\
t_{\mathrm{deg}}&:=(m-1)m_0+\frac{\mathfrak d_{\max}d_\theta}{2}.
\end{aligned}
\]
We assume the signal-strength condition in that proposition and retain its assumption that every $s$ in the support of $\nu$ has maximum signal degree $\mathfrak d_{\max}\geq1$.
Recall that, under the local regularity condition at scale $r_{\mathrm{deg}}$ with constant $C_{\mathrm{deg}}$, \Cref{lemma: local sub-Gaussian concentration of the NEF models} guarantees that for any sum of $m-1$ independent observations from $P_{\theta_0}$ and $P_{\theta_1}$ and every $0\leq x\leq C_{\mathrm{deg}}(m-1)v_0r_{\mathrm{deg}}$, $$ \Pr\bigg\{\sum_{j=1}^{m-1}\big(Y_j-\bbE[Y_j]\big)\geq x\bigg\}\leq\exp\bigg(-\frac{x^2}{2C_{\mathrm{deg}}(m-1)v_0}\bigg), $$ and the analogous bound holds for the lower tail. 
\paragraph{Type I error.} Under the null hypothesis, for every vertex $i\in[m]$, the weighted degree $T_i^{\mathrm{deg}}(\bfY)$ is a sum of $m-1$ independent observations from $P_{\theta_0}$ and thus has mean $(m-1)m_0$. Hence, by the union bound, 
\begin{align*} 
	\bbQ\big\{T_{\mathrm{deg}}(\bfY)>t_{\mathrm{deg}}\big\} &\leq\sum_{i=1}^m\bbQ\bigg\{T_i^{\mathrm{deg}}(\bfY)-(m-1)m_0>\frac{\mathfrak d_{\max} d_\theta}{2}\bigg\}.
\end{align*} 
Define the target deviation 
$$
x_{\delta,\mathrm{I}}:=\sqrt{2C_{\mathrm{deg}}(m-1)v_0\log\frac{2m}{\delta}}. 
$$ 
By the definition of $r_{\mathrm{deg}}$ and the fact that $C_{\mathrm{deg}}\geq1$, we have $x_{\delta,\mathrm{I}}\leq C_{\mathrm{deg}}(m-1)v_0r_{\mathrm{deg}},$ so the local concentration bound (\Cref{lemma: local sub-Gaussian concentration of the NEF models}) applies. Moreover, the assumed lower bound on $\tau_n$ is equivalent to 
$$
\frac{\mathfrak d_{\max} d_\theta}{2}\geq x_{\delta,\mathrm{I}}.
$$ 
Therefore, 
\begin{align*} 
	\sum_{i=1}^m\bbQ\bigg\{T_i^{\mathrm{deg}}(\bfY)-(m-1)m_0>\frac{\mathfrak d_{\max} d_\theta}{2}\bigg\} &\leq\sum_{i=1}^m\bbQ\big\{T_i^{\mathrm{deg}}(\bfY)-(m-1)m_0>x_{\delta,\mathrm{I}}\big\}\\
	&\leq m\exp\bigg(-\frac{x_{\delta,\mathrm{I}}^2}{2C_{\mathrm{deg}}(m-1)v_0}\bigg) =\frac{\delta}{2}. 
\end{align*} 
\paragraph{Type II error.} Fix any $s$ in the support of $\nu$. By assumption, there exists a vertex $i_s\in[m]$ such that $\mathfrak d_s(i_s)=\mathfrak d_{\max}$. Under $\bbP_s$, the weighted degree of this vertex has mean
$$ 
\bbE_{\bfY\sim\bbP_s}\big[T_{i_s}^{\mathrm{deg}}(\bfY)\big]=(m-1)m_0+\mathfrak d_{\max}d_\theta.
$$ 
Since $T_{\mathrm{deg}}(\bfY)\geq T_{i_s}^{\mathrm{deg}}(\bfY)$, the Type~II error for this fixed $s$ satisfies
\begin{align*} 
	\bbP_s\big\{T_{\mathrm{deg}}(\bfY)\leq t_{\mathrm{deg}}\big\} &\leq\bbP_s\big\{T_{i_s}^{\mathrm{deg}}(\bfY)\leq t_{\mathrm{deg}}\big\}\\
	&\leq\bbP_s\bigg\{T_{i_s}^{\mathrm{deg}}(\bfY)-\bbE_{\bfY\sim\bbP_s}[T_{i_s}^{\mathrm{deg}}(\bfY)]\leq-\frac{\mathfrak d_{\max} d_\theta}{2}\bigg\}.
\end{align*} 
Define 
$$ 
x_{\delta,\mathrm{II}}:=\sqrt{2C_{\mathrm{deg}}(m-1)v_0\log\frac{2}{\delta}}. 
$$ 
Since $\log(2/\delta)\leq\log(2m/\delta)$, we have 
$$ 
x_{\delta,\mathrm{II}}\leq x_{\delta,\mathrm{I}}\leq C_{\mathrm{deg}}(m-1)v_0r_{\mathrm{deg}},
$$ 
and the assumed lower bound on $\tau_n$ also implies $\mathfrak d_{\max} d_\theta/2\geq x_{\delta,\mathrm{II}}$. Therefore, by the lower-tail bound in \Cref{lemma: local sub-Gaussian concentration of the NEF models},
\begin{align*} 
	\bbP_s\bigg\{T_{i_s}^{\mathrm{deg}}(\bfY)-\bbE_{\bfY\sim\bbP_s}[T_{i_s}^{\mathrm{deg}}(\bfY)]\leq-\frac{\mathfrak d_{\max} d_\theta}{2}\bigg\} &\leq\bbP_s\bigg\{T_{i_s}^{\mathrm{deg}}(\bfY)-\bbE_{\bfY\sim\bbP_s}[T_{i_s}^{\mathrm{deg}}(\bfY)]\leq-x_{\delta,\mathrm{II}}\bigg\}\\
	&\leq\exp\bigg(-\frac{x_{\delta,\mathrm{II}}^2}{2C_{\mathrm{deg}}(m-1)v_0}\bigg) =\frac{\delta}{2}. 
\end{align*} 
The bound holds uniformly over every $s$ in the support of $\nu$, and hence
$$ 
\bbE_{S\sim\nu}\big[\bbP_S\{\Delta_{\mathrm{deg}}(\bfY)=0\}\big]\leq\frac{\delta}{2}.
$$ 
Combining the Type~I and Type~II error bounds yields $R_\nu(\Delta_{\mathrm{deg}})\leq\delta$.

\subsection{\texorpdfstring{Proof of \Cref{theorem: equivalence between NEF models}}{Proof of Theorem \getrefnumber{theorem: equivalence between NEF models}}}\label{proof: equivalence between NEF models}
\paragraph{Second moment of the likelihood ratio.}
We first derive the exact second-moment identity for the likelihood ratio.
Fix one of the two models and suppress the superscript $M\in \{A,B\}$ for notational convenience.
Write $\chi^2=\chi^2(\theta_0,\theta_1)=\chi_n^2$ and $D=D_n$ for brevity.
Write
\[
L_{\theta_0,\theta_1}(y)=\frac{dP_{\theta_1}(y)}{dP_{\theta_0}(y)}=\exp\big\{(\theta_1-\theta_0)y-(\psi(\theta_1)-\psi(\theta_0))\big\},
\]
for the one-coordinate likelihood ratio.
The assumption $\psi(2\theta_1-\theta_0)<\infty$ ensures that this ratio belongs to $L^2(P_{\theta_0})$.
Then, for any fixed $s\in\mathcal C_n$, the likelihood ratio of the $n$-dimensional observation $\mathbf Y=(Y_1,\dots,Y_n)$ is \(L_{s,n}(\mathbf y)=\prod_{l\in s}L_{\theta_0,\theta_1}(y_l)\).
Since the alternative is a mixture over $S\sim\nu_n$, the overall likelihood ratio is
\[
L_{\nu_n}(\mathbf y)=\ee_{S\sim\nu_n}[L_{S,n}(\mathbf y)].
\]
Let $S$ and $S'$ be independent copies with law $\nu_n$.
The second moment of the likelihood ratio under the null is therefore:
\begin{align*}
	\ee_{\mathbf Y\sim \bbQ_n}[L_{\nu_n}(\mathbf Y)^2]&= \ee_{S,S'\sim\nu_n}\Big[\ee_{\mathbf Y\sim \bbQ_n}\big[L_{S,n}(\mathbf Y)L_{S',n}(\mathbf Y)\big]\Big]\\
	&=\ee_{S,S'\sim\nu_n}\Bigg[\ee_{\mathbf Y\sim \bbQ_n}\bigg[\prod_{l\in S\cap S'}L_{\theta_0,\theta_1}(Y_l)^2\prod_{i\in S\triangle S'}L_{\theta_0,\theta_1}(Y_i)\bigg]\Bigg]\\
	&=\ee_{S,S'\sim\nu_n}\Bigg[\prod_{l\in S\cap S'}\ee_{Y_l\sim P_{\theta_0}}[L_{\theta_0,\theta_1}(Y_l)^2]\prod_{l\in S\triangle S'}\ee_{Y_l\sim P_{\theta_0}}[L_{\theta_0,\theta_1}(Y_l)]\Bigg]\\
	&=\ee_{S,S'\sim\nu_n}\Big[(1+\chi^2)^{|S\cap S'|}\Big],
\end{align*}
as $\ee_{Y\sim P_{\theta_0}}[L_{\theta_0,\theta_1}(Y)]=1$, and $\chi^2=\ee_{Y\sim P_{\theta_0}}[L_{\theta_0,\theta_1}(Y)^2]-1$.
Therefore, for two models with the same prior $\nu_n$, if $\chi^2_A=\chi^2_B$, then
\[
\ee_{\mathbf Y\sim \bbQ^A_n}[L_{\nu_n}^A(\mathbf Y)^2] = \ee_{\mathbf Y\sim \bbQ^B_n}[L_{\nu_n}^B(\mathbf Y)^2].
\]

\paragraph{Second moment of the low-degree likelihood ratio.}
We now study the low-degree likelihood ratio.
It suffices to consider integers $D\geq1$, since $D=0$ gives $L_{\nu_n}^{\leq0}\equiv1$.
Again we suppress the superscript $M\in \{A,B\}$ for notational convenience.
Let $\{\phi_d\}_{d\ge 0}$ be orthonormal polynomials in $L^2(P_{\theta_0})$, with $\phi_0\equiv 1$ and $\phi_1(y)=(y-\psi'(\theta_0))/\sqrt{\psi''(\theta_0)}$.
Note that if $P_{\theta_0}$ has finite support, then the sequence $\{\phi_d\}_{d\ge 0}$ terminates. %
For each available degree, define
\[
r_d(\theta_0,\theta_1):=\ee_{Y\sim P_{\theta_0}}[L_{\theta_0,\theta_1}(Y)\phi_d(Y)],\quad d\geq 1.
\]
In the finite-support case, extend $r_d$ by zero beyond the last available degree.
For an admissible multi-index $\boldsymbol\alpha=(\alpha_1,\dots,\alpha_n)\in\mathbb Z^n_{\ge 0}$, write
\[
\Phi_{\boldsymbol\alpha}(\mathbf y)=\prod_{l=1}^n\phi_{\alpha_l}(y_l),\quad |\boldsymbol\alpha|=\sum_{l=1}^n\alpha_l,
\]
and $\mathrm{supp}(\boldsymbol\alpha)=\{l:\alpha_l>0\}$.
In the finite-support case, retain only indices for which every $\phi_{\alpha_l}$ exists, and set coefficients at all other indices to zero in the sums below.
Then, the set of tensor products $\{\Phi_{\boldsymbol\alpha}: |\boldsymbol\alpha|\le D\}$ forms an orthonormal basis for the space of polynomials of degree at most $D$ in $L^2(\bbQ_n)$.
Let
\begin{align*}
R_{\boldsymbol\alpha}&:=\ee_{\mathbf Y\sim \bbQ_n}[L_{\nu_n}(\mathbf Y)\Phi_{\boldsymbol\alpha}(\mathbf Y)]\\
&=\ee_{S\sim\nu_n}\Big[\ee_{\mathbf Y\sim \bbQ_n}\big[L_{S,n}(\mathbf Y)\Phi_{\boldsymbol\alpha}(\mathbf Y)\big]\Big]\\
&=\ee_{S\sim\nu_n}\Bigg[\prod_{l=1}^n\ee_{Y\sim P_{\theta_0}}\big[L_{\theta_0,\theta_1}(Y)^{\mathds 1\{{l\in S}\}}\phi_{\alpha_l}(Y)\big]\Bigg].
\end{align*}
We now focus on the inner expectation.
For each $l \in [n],$ there are only three possible cases: (i) if $\alpha_l=0$, then $\ee_{Y\sim P_{\theta_0}}[L_{\theta_0,\theta_1}(Y)^{\mathds 1\{{l\in S}\}}\phi_{\alpha_l}(Y)]=1$; (ii) if $\alpha_l\geq1$ and $l\in S$, then the factor equals $\ee_{Y\sim P_{\theta_0}}[L_{\theta_0,\theta_1}(Y)\phi_{\alpha_l}(Y)]=r_{\alpha_l}(\theta_0,\theta_1)$; (iii) if $\alpha_l\geq1$ and $l\notin S$, then $\ee_{Y\sim P_{\theta_0}}[\phi_{\alpha_l}(Y)]=0$.
Therefore, the entire product is zero unless $\mathrm{supp}(\boldsymbol\alpha)\subseteq S$.
Consequently,
\begin{align*}
R_{\boldsymbol\alpha}&=\ee_{S\sim\nu_n}\Bigg[\mathds 1\{\mathrm{supp}(\boldsymbol\alpha)\subseteq S\}\prod_{l\in \mathrm{supp}(\boldsymbol\alpha)}r_{\alpha_l}(\theta_0,\theta_1)\Bigg]\\
&=\Pr_{S\sim\nu_n}\big[\mathrm{supp}(\boldsymbol\alpha)\subseteq S\big]\prod_{l\in \mathrm{supp}(\boldsymbol\alpha)}r_{\alpha_l}(\theta_0,\theta_1).
\end{align*}
By the same logic,
\begin{align*}
R_{\boldsymbol\alpha}^2=\Pr_{S,S'\sim\nu_n}\big[\mathrm{supp}(\boldsymbol\alpha)\subseteq S\cap S'\big]\prod_{l\in \mathrm{supp}(\boldsymbol\alpha)}r_{\alpha_l}(\theta_0,\theta_1)^2.
\end{align*}
Now, recall that the low-degree likelihood ratio is 
\[
L_{\nu_n}^{\leq D}(\mathbf Y)=\sum_{|\boldsymbol\alpha|\le D}R_{\boldsymbol\alpha}\Phi_{\boldsymbol\alpha}(\mathbf Y).
\]
Therefore, Parseval's identity gives
\begin{align*}
\ee_{\mathbf Y\sim \bbQ_n}[L_{\nu_n}^{\leq D}(\mathbf Y)^2]&=\sum_{|\boldsymbol\alpha|\le D}R_{\boldsymbol\alpha}^2\\
&=\sum_{|\boldsymbol\alpha|\le D}\Pr_{S,S'\sim\nu_n}\big[\mathrm{supp}(\boldsymbol\alpha)\subseteq S\cap S'\big]\prod_{l\in \mathrm{supp}(\boldsymbol\alpha)}r_{\alpha_l}(\theta_0,\theta_1)^2\\
&=\ee_{S,S'\sim\nu_n}\Bigg[\sum_{\substack{\boldsymbol\alpha:\mathrm{supp}(\boldsymbol\alpha)\subseteq S\cap S'\\|\boldsymbol\alpha|\le D}}\prod_{l\in \mathrm{supp}(\boldsymbol\alpha)}r_{\alpha_l}(\theta_0,\theta_1)^2\Bigg].
\end{align*}
We group this sum according to the support size $d=|\mathrm{supp}(\boldsymbol\alpha)|$, which ranges from $0$ to $D$, and obtain
\begin{align}
\ee_{\mathbf Y\sim \bbQ_n}[L_{\nu_n}^{\leq D}(\mathbf Y)^2]&=\ee_{S,S'\sim\nu_n}\Bigg[\sum_{d=0}^{D}\sum_{\substack{\boldsymbol\alpha:\mathrm{supp}(\boldsymbol\alpha)\subseteq S\cap S'\\|\mathrm{supp}(\boldsymbol\alpha)|=d\\|\boldsymbol\alpha|\leq D}}\prod_{l\in \mathrm{supp}(\boldsymbol\alpha)}r_{\alpha_l}(\theta_0,\theta_1)^2\Bigg] \nonumber\\
&=\ee_{S,S'\sim\nu_n}\Bigg[\sum_{d=0}^{D}\sum_{\substack{T\subseteq S\cap S'\\|T|=d}}\sum_{\substack{\boldsymbol\alpha:\mathrm{supp}(\boldsymbol\alpha)=T\\|\boldsymbol\alpha|\leq D}}\prod_{l\in T}r_{\alpha_l}(\theta_0,\theta_1)^2\Bigg]  \nonumber\\
&=\ee_{S,S'\sim\nu_n}\Bigg[\sum_{d=0}^{D\wedge |S\cap S'|}\binom{|S\cap S'|}{d}\sum_{\substack{\alpha_1,\dots,\alpha_d\geq 1\\\alpha_1+\cdots+\alpha_d\leq D}}~\prod_{i=1}^d r_{\alpha_i}(\theta_0,\theta_1)^2\Bigg],\label{eq:try}
\end{align}
The $d=0$ term corresponds to the empty support; its empty product and inner sum are both $1$.

Now, we bound this expression from above and below.
For $1\leq d\leq D$, the case $\alpha_1=\cdots=\alpha_d=1$ is always included in the sum, and thus we have
\begin{equation}
    \label{eq:r}
\sum_{\substack{\alpha_1,\alpha_2,\dots,\alpha_d\geq 1\\\alpha_1+\cdots+\alpha_d\leq D}}~\prod_{i=1}^d r_{\alpha_i}(\theta_0,\theta_1)^2\geq r_1(\theta_0,\theta_1)^{2d},
\end{equation}
which gives the lower bound
\[
\ee_{\mathbf Y\sim \bbQ_n}[L_{\nu_n}^{\leq D}(\mathbf Y)^2]\ge \ee_{S,S'\sim\nu_n}\Bigg[\sum_{d=0}^{D\wedge |S\cap S'|}\binom{|S\cap S'|}{d}r_1(\theta_0,\theta_1)^{2d}\Bigg].
\]
On the other hand, for each $d \geq 1$, 
\begin{align*}
\sum_{\substack{\alpha_1,\dots,\alpha_d\geq 1\\\alpha_1+\cdots+\alpha_d\leq D}}~\prod_{i=1}^d r_{\alpha_i}(\theta_0,\theta_1)^2&\leq \sum_{\alpha_1,\dots,\alpha_d\geq 1}~\prod_{i=1}^d r_{\alpha_i}(\theta_0,\theta_1)^2\\
&=\bigg(\sum_{\alpha\geq 1} r_{\alpha}(\theta_0,\theta_1)^2\bigg)^d\\
&\mathrel{\leq} \|L_{\theta_0,\theta_1}-1\|_{L^2(P_{\theta_0})}^{2d}\\
& =\chi^2(\theta_0,\theta_1)^d.
\end{align*}
The norm bound is Bessel's inequality applied to $L_{\theta_0,\theta_1}-1$ and the orthonormal polynomials of positive degree.
The above inequality, when combined with \Cref{eq:r}, yields
\begin{equation}\label{eq:r2}
    \chi^2(\theta_0,\theta_1)^d \geq r_1(\theta_0,\theta_1)^{2d},
\end{equation}
and, when combined with \Cref{eq:try},
\[
\ee_{\mathbf Y\sim \bbQ_n}[L_{\nu_n}^{\leq D}(\mathbf Y)^2]\le \ee_{S,S'\sim\nu_n}\Bigg[\sum_{d=0}^{D\wedge |S\cap S'|}\binom{|S\cap S'|}{d}\chi^2(\theta_0,\theta_1)^d\Bigg].
\]
Thus, we have shown that
\begin{equation}\label{eq:low-degree sandwiching}
	\begin{aligned}
\ee_{S,S'\sim\nu_n}\Bigg[\sum_{d=0}^{D\wedge |S\cap S'|}\binom{|S\cap S'|}{d}r_1(\theta_0,\theta_1)^{2d}\Bigg]&\leq \ee_{\mathbf Y\sim \bbQ_n}[L_{\nu_n}^{\leq D}(\mathbf Y)^2]\\
&\leq \ee_{S,S'\sim\nu_n}\Bigg[\sum_{d=0}^{D\wedge |S\cap S'|}\binom{|S\cap S'|}{d}\chi^2(\theta_0,\theta_1)^d\Bigg].
	\end{aligned}
\end{equation}
 
Now we compare $r_1(\theta_0,\theta_1)^2$ and $\chi^2(\theta_0,\theta_1)$. Note that $r_1(\theta_0,\theta_1)^2=\tau^2_{n}$, since
\begin{align*}
r_1(\theta_0,\theta_1)&=\ee_{Y\sim P_{\theta_0}}[L_{\theta_0,\theta_1}(Y)\phi_1(Y)]\\
&=\ee_{Y\sim P_{\theta_1}}[\phi_1(Y)]\\
&=\frac{\ee_{Y\sim P_{\theta_1}}[Y]-\ee_{Y\sim P_{\theta_0}}[Y]}{\sqrt{\vv_{Y\sim P_{\theta_0}}[Y]}}\\
&=\frac{\psi'(\theta_1)-\psi'(\theta_0)}{\psi''(\theta_0)^{1/2}}.
\end{align*}
Since the NEF is nondegenerate and $\theta_0\neq\theta_1$, we have $\tau_n\neq0$.
By \Cref{eq:r2},
\[
\frac{\chi^2(\theta_0,\theta_1)}{\tau_n^2}\geq1,
\qquad \log\frac{\chi^2(\theta_0,\theta_1)}{\tau_n^2}\geq0.
\]
Under the assumption $D\log(\chi^2(\theta_0,\theta_1)/\tau^2_{n})=O(1)$, for $0\leq d\leq D$ we have
\begin{align*}
r_1(\theta_0,\theta_1)^{2d}&=\exp\{-d\log(\chi^2(\theta_0,\theta_1)/\tau^2_{n})\}\chi^{2}(\theta_0,\theta_1)^{d}\\
&\geq \exp\{-D\log(\chi^2(\theta_0,\theta_1)/\tau^2_{n})\}\chi^{2}(\theta_0,\theta_1)^{d}\\
&\geq e^{-C}\chi^2(\theta_0,\theta_1)^d
\end{align*}
for some fixed $C>0$, uniformly over $0\leq d\leq D$ and all sufficiently large $n$.
Therefore, we obtain
\begin{align*}
	\ee_{S,S'\sim\nu_n}\Bigg[\sum_{d=0}^{D\wedge |S\cap S'|}\binom{|S\cap S'|}{d}r_1(\theta_0,\theta_1)^{2d}\Bigg]\geq e^{-C}\ee_{S,S'\sim\nu_n}\Bigg[\sum_{d=0}^{D\wedge |S\cap S'|}\binom{|S\cap S'|}{d}\chi^2(\theta_0,\theta_1)^d\Bigg],
\end{align*}
which implies that the two bounds in \Cref{eq:low-degree sandwiching} are equivalent up to a constant factor.
The comparison constants depend only on an eventual upper bound for $D\log(\chi^2/\tau_n^2)$.
This completes the proof.

\subsection{\texorpdfstring{Proof of \Cref{proposition: local alternative regime}}{Proof of Proposition \getrefnumber{proposition: local alternative regime}}}\label{proof: local alternative regime}
In this proof, we parse the condition $D\log(\chi^2(\theta_0,\theta_1)/\tau^2_{n})=O(1)$ and discuss its implications for the fixed-null setting.
For brevity, write $D=D_n$, $h=h_n$, and $\chi^2=\chi^2(\theta_0,\theta_1)=\chi_n^2$.
We first fix $\theta_0\in\Theta$ and let $\theta_1=\theta_0+h$, where $h\to0$ with $h\neq0$.
Define $v_0:=\psi^{(2)}(\theta_0)\mathrel{>}0$ and $c_3:=\psi^{(3)}(\theta_0)$.
Because $\theta_0$ is an interior point, the segments $\{\theta_0+uh:u\in[0,2]\}$ lie in a fixed compact subset of $\Theta$ for all sufficiently small $|h|$.
In particular, $\theta_0+2h\in\Theta$, so the chi-square divergence below is finite.
Then, we have
\begin{align*}
1+\chi^2&=\ee_{Y\sim P_{\theta_0}}[L_{\theta_0,\theta_1}(Y)^2]\\
&=\ee_{Y\sim P_{\theta_0}}[\exp\{2(\theta_1-\theta_0) Y-2(\psi(\theta_1)-\psi(\theta_0))\}]\\
&=\exp\{-2(\psi(\theta_1)-\psi(\theta_0))\}\ee_{Y\sim P_{\theta_0}}[\exp\{2(\theta_1-\theta_0) Y\}]\\
&\stackrel{(\dagger)}=\exp\{-2(\psi(\theta_1)-\psi(\theta_0))\}\exp\{\psi(2\theta_1-\theta_0)-\psi(\theta_0)\}\\
&=\exp\{\psi(2\theta_1-\theta_0)-2\psi(\theta_1)+\psi(\theta_0)\}\\
&=\exp\{\psi(\theta_0+2h)-2\psi(\theta_0+h)+\psi(\theta_0)\},
\end{align*}
where $(\dagger)$ follows from $\psi(\theta_0+t)-\psi(\theta_0)=\log \ee_{Y\sim P_{\theta_0}}[\exp(t Y)]$ at $t=2h$.
By Taylor expansion with Lagrange remainder, there exist $u_1, u_2 \in [0,1]$ such that
\[
\begin{aligned}
\psi(\theta_0+2h) &= \psi(\theta_0) + 2h \psi^{(1)}(\theta_0) + 4h^2 \frac{\psi^{(2)}(\theta_0) }{2} + 8h^3 \frac{\psi^{(3)}(\theta_0) }{6}\\
&\quad+16h^4 \frac{\psi^{(4)}(\theta_0 + 2u_1h)}{24}
\end{aligned}
\]
and 
\[
2 \psi(\theta_0+h) = 2 \psi(\theta_0) + 2h \psi^{(1)}(\theta_0) + 2 h^2 \frac{\psi^{(2)}(\theta_0) }{2} + 2 h^3  \frac{\psi^{(3)}(\theta_0) }{6} +  2 h^4 \frac{\psi^{(4)}(\theta_0 + u_2 h) }{24} .
\]
Using the above identities, we obtain 
\[
\psi(\theta_0+2h)-2\psi(\theta_0+h)+\psi(\theta_0) = v_0 h^2+c_3h^3+O(h^4),
\]
where the remainder is controlled by $\sup_{u\in[0,2]}|\psi^{(4)}(\theta_0+uh)|$ and is uniform for sufficiently small $|h|$.
Thus, as $h \to 0$, 
\[
\chi^2 = \exp\{v_0 h^2+c_3h^3+O(h^4)  \} - 1,
\]
and therefore, using the fact that $e^x-1=x+O(x^2)$ for $x\to0$, 
\begin{align}\label{eq:chi2.h2}
\chi^2=v_0 h^2+{c_3}h^3+O(h^4) = \Theta(h^2),
\end{align}
for $h \to 0$.
As for $\tau_n=\frac{\psi'(\theta_1)-\psi'(\theta_0)}{\psi''(\theta_0)^{1/2}}$, the Taylor expansion of $\psi^{(1)}(\theta_1)$ around $\theta_0$ yields, for some $u_3\in[0,1]$,
\begin{align}
\tau_n^2 & = \frac{ (v_0 h + c_3  h^2/2  +  \psi^{(4)}(\theta_0 + u_3 h) h^3 /6 )^2} {v_0} \nonumber \\
& = \frac{ v^2_0 h^2 +  v_0 c_3 h^3  +  O(h^4)} {v_0} \nonumber \\
& = v_0 h^2+c_3 h^3+O(h^4) = \Theta(h^2). \label{eq:mu2.h2}
\end{align}
For $|h|\leq1$, the final remainder $O(h^4)$ can be bounded in absolute value by $h^4$ times
\[
c_4+\frac{\max\{c_3^2,c_4^2,|c_3|c_4\}}{v_0},
\]
with $c_4:=\sup_{u\in[0,2]}|\psi^{(4)}(\theta_0+uh)|$, which is uniformly bounded for small $|h|$.
Therefore, as $h \to 0$, 
$$
\frac{\chi^2}{\tau_n^2}=\frac{v_0 h^2+c_3 h^3+O(h^4)}{v_0 h^2+c_3 h^3+O(h^4)}=\frac{1+\frac{c_3}{v_0}h+O(h^2)}{1+\frac{c_3}{v_0}h+O(h^2)}=1+O(h^2),
$$
which implies that, as $h \to 0$,
\[
\log\frac{\chi^2}{\tau_n^2} = \log(1+O(h^2)) = O(h^2).
\]
The last expression, combined with \Cref{eq:chi2.h2} and \Cref{eq:mu2.h2}, gives
\[
\log\frac{\chi^2}{\tau_n^2} = O(\chi^2) = O(\tau_n^2),
\]
as $h \to 0$.
Consequently, $D_nh_n^2=O(1)$ implies
\[
D_n\log\frac{\chi_n^2}{\tau_n^2}=O(D_nh_n^2)=O(1).
\]
Moreover, \Cref{eq:chi2.h2,eq:mu2.h2} show that $h_n^2$, $\chi_n^2$, and $\tau_n^2$ are comparable up to fixed positive constants, so the three sufficient degree conditions in \Cref{proposition: local alternative regime} are equivalent.

The fixed-null assumption can be relaxed when the Taylor expansions above hold uniformly.
In particular, let $\theta_{0,n}$ vary with $n$, set $\theta_{1,n}=\theta_{0,n}+h_n$ with $h_n\to0$, and suppose that the segments $\{\theta_{0,n}+u h_n:u\in[0,2]\}$ are contained in a fixed compact subset of $\operatorname{int}(\Theta)$.
Then $\psi''(\theta_{0,n})$ is bounded away from zero and the derivatives of $\psi$ appearing in the expansions are uniformly bounded, so the preceding calculation gives
\[
\frac{\chi_n^2}{\tau_n^2}=1+O(h_n^2)
\]
with an implicit constant independent of $n$.
Consequently, the conclusion of \Cref{proposition: local alternative regime} continues to hold for such varying null parameters.
Without this condition on the null parameters, uniform derivative bounds or a positive lower bound on $\psi''(\theta_{0,n})$ may fail, so the extension requires model-specific verification, as in the sparse Bernoulli regimes discussed in \Cref{section: equivalence on planted clique}.

\subsection{\texorpdfstring{Proof of \Cref{theorem: universality of tightness in the local alternative regime via chi-square matching}}{Proof of Theorem \getrefnumber{theorem: universality of tightness in the local alternative regime via chi-square matching}}}\label{proof: universality of tightness in the local alternative regime via chi-square matching}
Fix an arbitrary positive sequence $h_n^B\to0$ for which the effective signal strength $\tau_n^B$ vanishes polynomially.
We first construct a chi-square-matched sequence $h_n^A\to0$.
For $M\in\{A,B\}$ and sufficiently small $h\geq0$, define
\[
g_M(h):=\chi^2\big(P^M_{c_M+h}\,\|\,P^M_{c_M}\big)
=\exp\big\{\psi^M(c_M+2h)-2\psi^M(c_M+h)+\psi^M(c_M)\big\}-1.
\]
This function is continuous with $g_M(0)=0$, and for sufficiently small $h>0$,
\[
\frac{d}{dh}\log\big(1+g_M(h)\big)
=2\big\{(\psi^M)'(c_M+2h)-(\psi^M)'(c_M+h)\big\}>0,
\]
where the strict inequality follows from $(\psi^M)''>0$.
Also, the expansion in \Cref{eq:chi2.h2} gives $g_M(h)=(\psi^M)''(c_M)h^2(1+O(h))$.
Hence, for all sufficiently large $n$, continuity and strict monotonicity give a unique sufficiently small $h_n^A>0$ such that $g_A(h_n^A)=g_B(h_n^B)$, and $h_n^A\to0$.

For $M\in\{A,B\}$, write
$$
\chi_{M,n}^2:=\chi^2\big(P^M_{\theta_1^M}\,\|\,P^M_{\theta_0^M}\big).
$$
By the local alternative expansion in \Cref{proposition: local alternative regime},
$$
\frac{\chi_{M,n}^2}{(\tau_n^M)^2}=1+O\big((h_n^M)^2\big)=1+o(1).
$$
The two models have matching chi-square divergences, $\chi_{A,n}^2=\chi_{B,n}^2=\chi_n^2$, and both effective signal strengths are positive because $h_n^M>0$ and $(\psi^M)''>0$.
Therefore,
$$
\frac{\tau_n^A}{\tau_n^B}\to 1.
$$
In particular, the matched effective signal strength $\tau_n^A$ also vanishes polynomially, so the assumptions for model $A$ apply to this sequence.
Hence, multiplication by a $1+o(1)$ factor preserves the relations of polynomial separation $\ll$ and $\gg$.
In particular, for any positive sequence $a_n$,
$$
\tau_n^A\ll a_n \Longleftrightarrow \tau_n^B\ll a_n,\qquad \tau_n^A\gg a_n \Longleftrightarrow \tau_n^B\gg a_n.
$$

We first consider the statistical lower bound.
By \Cref{theorem: equivalence between NEF models}, the second moment of the likelihood ratio under model $M$ is
$$
\bbE_{\bfY\sim\bbQ_n^M}\big[L_{\nu_n}^M(\bfY)^2\big]
=
\bbE_{S,S'\sim\nu_n}\Big[\big(1+\chi_n^2\big)^{|S\cap S'|}\Big].
$$
Thus, under chi-square matching, the likelihood-ratio second moments of models $A$ and $B$ are identical.
Therefore, if $\tau_n^B\ll\tau_{\mathrm{stat}}^{A,\ast}$, then $\tau_n^A\ll\tau_{\mathrm{stat}}^{A,\ast}$, and boundedness of the second moment for model $A$ implies the same boundedness for model $B$.
Hence, by \Cref{lemma: second moment method}, the statistical lower bound at the scale $\tau_{\mathrm{stat}}^{A,\ast}$ transfers to model $B$.

Next, consider the computational lower bound.
Suppose that $\tau_n^B\ll\tau_{\mathrm{comp}}^{A,\ast}$.
Then $\tau_n^A\ll\tau_{\mathrm{comp}}^{A,\ast}$ as well.
By the standing assumption that the effective signal strengths vanish polynomially and that $D_n$ is polylogarithmic,
$$
D_n(\tau_n^M)^2=o(1),\qquad M\in\{A,B\}.
$$
It follows from \Cref{proposition: local alternative regime} that
$$
D_n\log\bigg(\frac{\chi_n^2}{(\tau_n^M)^2}\bigg)=O(1),\qquad M\in\{A,B\}.
$$
Therefore, the low-degree part of \Cref{theorem: equivalence between NEF models} applies to both models and gives
$$
\bbE_{\bfY\sim\bbQ_n^M}\big[L_{\nu_n}^{M,\leq D_n}(\bfY)^2\big]
\asymp
\bbE_{S,S'\sim\nu_n}\Bigg[\sum_{d=0}^{D_n\wedge|S\cap S'|}\binom{|S\cap S'|}{d}(\chi_n^2)^d\Bigg],\qquad M\in\{A,B\}.
$$
The expression on the right-hand side is common to the two models.
Thus, boundedness of the degree-$D_n$ likelihood-ratio second moment for model $A$ implies boundedness for model $B$.
By \Cref{lemma: low-degree second moment method}, this rules out strong separation by degree-$D_n$ polynomials under model $B$.
Under \Cref{conjecture: low-degree polynomial conjecture}, with the same degree requirement as for model $A$, the computational lower bound at the scale $\tau_{\mathrm{comp}}^{A,\ast}$ therefore transfers to model $B$.

It remains to establish the corresponding upper bounds.
Because each null parameter is a fixed interior point, the neighborhoods $[\theta_j^M-r_n,\theta_j^M+r_n]$, $j\in\{0,1\}$, are contained in $\Theta_M$ for all sufficiently large $n$ whenever $r_n\to0$.
Since both models are in the local alternative regime, \Cref{proposition: upper bound for NEF models in local alternative regime} gives, for every sequence $r_n\to0$,
$$
\sup_{\substack{j\in\{0,1\}\\|u|\leq r_n}}\frac{(\psi^M)''(\theta_j^M+u)}{(\psi^M)''(\theta_0^M)}=1+o(1),\qquad M\in\{A,B\}.
$$
Hence, the local regularity conditions required for the average, degree, and scan tests are satisfied with a constant arbitrarily close to one, provided that the corresponding radius converges to zero.
For the average test and, when applicable, the degree test, this holds by choosing the target risk $\delta_n\to0$ sufficiently slowly, so that $\log(2/\delta_n)=o(n)$ and $\log(2m/\delta_n)=o(m)$, respectively.
For the scan test, since $(\psi^M)''(\theta_0^M)$ is a fixed positive constant,
$$
r_{\mathrm{scan}}^2\asymp\frac{\log|\mathcal C|+\log(1/\delta_n)}{K},
$$
where the radius is evaluated separately for each model.
Whenever the scan-test sufficient condition gives strong detection in model $A$, the assumption $\tau_n^A\to0$ forces $K\to\infty$.
Thus, the assumption $\log|\mathcal C|=o(K)$ allows us to choose $\delta_n\to0$ sufficiently slowly so that $r_{\mathrm{scan}}\to0$ for both models.

By \Cref{proposition: upper bound for NEF models} and \Cref{remark: gaussian degree test}, the sufficient conditions for these tests depend on the NEF model only through its effective signal strength, up to the local regularity constant.
Therefore, the corresponding average, degree, or scan test that yields the matching upper bound for model $A$ yields an upper bound of the same order for model $B$.
Because the latent-structure class and the test statistic are unchanged, this transfer also preserves the computational complexity of the test.
Consequently, strong detection in model $B$ is achievable whenever $\tau_n^B\gg\tau_{\mathrm{stat}}^{A,\ast}$, and polynomial-time strong detection is achievable whenever $\tau_n^B\gg\tau_{\mathrm{comp}}^{A,\ast}$.

Since the starting admissible signal sequence in model $B$ was arbitrary, combining the transferred lower and upper bounds shows that $\tau_{\mathrm{stat}}^{A,\ast}$ and $\tau_{\mathrm{comp}}^{A,\ast}$ are also valid statistical and computational thresholds, respectively, for model $B$, in the sense of \Cref{definition: statistical and computational thresholds in NEF models}.
This proves the claim.

\subsection{\texorpdfstring{Alternative proofs of computational lower bounds for multiple cliques and hypercliques}{Alternative proofs of computational lower bounds for multiple cliques and hypercliques}}
\label{sec:alternative-proof}
In this subsection, we provide alternative proofs of the computational lower bounds for the random-size multiple clique problem (\Cref{theorem: multiple cliques problem}) and the random-size multiple planted hypercliques problem (\Cref{corollary: multiple planted hypercliques problem}).
We follow the approach of \citet{mao2025information}, based on $U$-statistics and their decoupled counterparts \citep{delapena1999decoupling,gine2000exponential}.
We note that the resulting bounds match those obtained from the original proofs in \Cref{proof: random-size multiple cliques,proof: multiple planted hypercliques problem random size}, up to logarithmic factors.

\subsubsection{Multiple random-size cliques}

Recall the definition in \Cref{eq:definition of W_v}:
\begin{align}
W_v=\begin{cases}
e_{(a,b)} & \text{if }v\in\widetilde V_a\cap\widetilde V'_b\text{ for some }a,b\in[\gamma],\\
0 & \text{otherwise},
\end{cases}
\end{align}
and note that $\Pr(W_v=e_{(a,b)})=\rho=k^2/m^2$ for all $a,b\in[\gamma]$.
The labels $W_1,\ldots,W_m$ are i.i.d., with $\Pr(W_v=0)=1-\gamma^2\rho$; the feasibility condition $\gamma k\leq m$ ensures that $\gamma^2\rho\leq1$.
Define the kernel function $h:(\{0,1\}^{\gamma \times \gamma})^2\to\mathbb R$ by
$
h(w_1,w_2):=\sum_{a,b=1}^{\gamma}\mathds 1\{w_1=e_{(a,b)}\}\mathds 1\{w_2=e_{(a,b)}\}.
$
From \Cref{equation: overlap in terms of W_v multiple cliques}, the overlap size $|S\cap S'|$ can be expressed as
\begin{equation}\label{eq:U-statistic representation of overlap size}
\begin{aligned}
|S\cap S'|=\sum_{(i,j)\in[m]_2}h(W_i,W_j).
\end{aligned}
\end{equation}
Hence, this can be seen as a $U$-statistic of order 2 with kernel $h$.
We bound the overlap size $|S\cap S'|$ using its decoupled version.
Consider an independent copy $W'_1,\ldots,W'_m$ of $W_1,\ldots,W_m$, and  the kernel 
$$
H_{i,j}(w_1,w_2):=\mathds1\{i\neq j\}h(w_1,w_2), \quad (i,j)\in[m]\times[m].
$$
Note that $H_{i,j}$ is zero whenever $i=j$, and is equal to $h(w_1,w_2)$ otherwise.
Therefore, since $h$ is symmetric, we obtain
\begin{align*}
|S\cap S'| &=\sum_{(i,j)\in [m]_2} h(W_{i},W_{j})\\&=\frac{1}{2}\sum_{(i,j) \in [m]^2} H_{i,j}(W_{i},W_{j})\\
&=:U(W_1,...,W_m).
\end{align*}
Define the decoupled version of $U(W_1,...,W_m)$ as
\begin{align*}
U'(W_1,...,W_m;W'_1,...,W'_m) =\frac{1}{2}\sum_{(i,j) \in [m]\times[m]} H_{i,j}(W_{i},W'_{j}).
\end{align*}
Then, by \Cref{lemma: decoupling of U-statistics}, 
there exists a constant $c>0$ such that
$$
\ee\bigg[\bigg(\sum_{(i,j) \in [m]^2} H_{i,j}(W_{i},W_{j})\bigg)^d\bigg]\leq c^d\ee\bigg[\bigg(\sum_{(i,j) \in [m]^2} H_{i,j}(W_{i},W'_{j})\bigg)^d\bigg],
$$
for every integer $d\geq1$, which implies that
\begin{align}\label{eq:decoupling of U-statistics for multiple cliques}
	\ee[|S\cap S'|^d] \leq {c'}^d \ee\bigg[\bigg(\sum_{(i,j) \in [m]^2} H_{i,j}(W_{i},W'_{j})\bigg)^d\bigg],
\end{align}
for $c' = c/2$.

For notational convenience, define
$$
\lambda:=\gamma^2\bigg(\frac{k^2}{m}\bigg)^2.
$$
Then, when $d=1$, 
\begin{equation}\label{eq:first moment bound for multiple cliques}
\begin{aligned}
\ee\bigg[\sum_{(i,j)\in[m]^2}H_{i,j}(W_i,W'_j)\bigg]&=\sum_{(i,j)\in[m]^2}\ee\left[H_{i,j}(W_i,W'_j)\right]\\
&=\sum_{i\neq j}\ee\left[H_{i,j}(W_i,W'_j)\right]\\
&=\sum_{i\neq j}\Pr\big(H_{i,j}(W_i,W'_j)=1\big)\\
&=\sum_{i\neq j}\sum_{a,b=1}^{\gamma}\Pr\big(W_i=e_{(a,b)}\big)\Pr\big(W'_j=e_{(a,b)}\big)\\
&=\sum_{i\neq j}\sum_{a,b=1}^{\gamma}\rho^2\\
&\leq m^2\gamma^2\rho^2\\
&=\lambda.
\end{aligned}
\end{equation}
For $d\geq 2$, \Cref{lemma: moment bound of nonnegative U-statistics} yields 
\begin{align*}
\ee\bigg[\bigg(\sum_{(i,j)\in[m]^2}H_{i,j}(W_i,W'_j)\bigg)^d\bigg]&\leq C_0^{2d}\sum_{J\subseteq[2]}d^{|J|d}\sum_{\bfi_J}\ee_J\bigg[\bigg(\sum_{\bfi_{J^c}}\ee_{J^c}[H_{i_1,i_2}(W_{i_1},W'_{i_2})\mid\mathcal F_J]\bigg)^d\bigg]\\
&=:C_0^{2d}\sum_{J\subseteq[2]}Q_J,
\end{align*}
where $C_0>0$ is a universal constant, $\bfi_J:=(i_l:l\in J)$, and $\bfi_{J^c}:=(i_l:l\in J^c)$.
Here, $\mathcal F_J$ and $\ee_J$ are defined as in \Cref{lemma: moment bound of nonnegative U-statistics}, with $(W_i)_{i=1}^m$ and $(W'_i)_{i=1}^m$ as the two independent coordinate blocks.
We bound the four terms $Q_J$ for $J=\emptyset,\{1\},\{2\},\{1,2\}$ separately.
When $J=\emptyset$,
\begin{equation}\label{eq:Q0 for multiple cliques}
\begin{aligned}
Q_\emptyset&=\bigg(\sum_{i=j}\ee\left[H_{i,j}(W_i,W'_j)\right]+\sum_{i\neq j}\ee\left[H_{i,j}(W_i,W'_j)\right]\bigg)^d\\
&=\bigg(\sum_{i\neq j}\ee\left[H_{i,j}(W_i,W'_j)\right]\bigg)^d\\
&\stackrel{(a)}=\bigg(\sum_{i\neq j}\Pr\left(H_{i,j}(W_i,W'_j)=1\right)\bigg)^d\\
&=\bigg(\sum_{i\neq j}\sum_{a,b=1}^{\gamma}\Pr\left(W_i=e_{(a,b)}\right)\Pr\left(W'_j=e_{(a,b)}\right)\bigg)^d\\
&=\bigg(\sum_{i\neq j}\sum_{a,b=1}^{\gamma}\rho^2\bigg)^d\\
&\leq (m^2\gamma^2\rho^2)^d\\
&=\lambda^d,
\end{aligned}
\end{equation}
where $(a)$ follows from the fact that $H_{i,j}$ takes values in $\{0,1\}$.

Next, consider $J=\{1\}$.
For ease of notation, for $i\in[m]$, let $Q_{\{1\}}(i):=\sum_{j\in[m]}\ee\left[H_{i,j}(W_i,W'_j)\mid W_i\right]$.
Then,
\begin{align*}
Q_{\{1\}}&=d^{d}\sum_{i\in [m]}\ee\bigg[\bigg(\sum_{j\in[m]}\ee\left[H_{i,j}(W_i,W'_j)\mid W_i\right]\bigg)^d\bigg]\\
&=d^{d}\sum_{i\in [m]}\ee\left[Q_{\{1\}}(i)^d\right].
\end{align*}
If $W_i=0$, then $Q_{\{1\}}(i)=0$. On the other hand, if $W_i=e_{(a,b)}$ for some $a,b\in[\gamma]$, then  $\Pr\left(W'_j=e_{(a,b)}\right)=\rho$ for every $j\neq i$.
Therefore
$$
Q_{\{1\}}(i)\leq m\rho\,\mathds{1}\{W_i\neq 0\}.
$$
As $\Pr(W_i\neq 0)=\gamma^2\rho$, we obtain
\begin{align*}
Q_{\{1\}}&=d^{d}\sum_{i\in [m]}\ee\left[Q_{\{1\}}(i)^d\right]\\
&\leq d^{d}\sum_{i=1}^m\ee\left[m^d\rho^d\,\mathds{1}\{W_i\neq 0\}\right]\\
&=d^{d}\sum_{i=1}^m(m^d\rho^d)\Pr(W_i\neq 0)\\
&=d^{d}m^{d+1}\gamma^2\rho^{d+1}\\
&=d^{d}\lambda \bigg(\frac{k^2}{m}\bigg)^{d-1}.
\end{align*}
By symmetry, the same bound holds for $J=\{2\}$. Therefore, the total contribution from $J=\{1\}$ and $J=\{2\}$ is at most
\begin{align}\label{eq:Q1Q2 for multiple cliques}
Q_{\{1\}}+Q_{\{2\}}\leq 2d^d\lambda \left( \frac{k^2}{m}\right)^{d-1}.
\end{align}

Finally, consider $J=\{1,2\}$. Since $H_{i,j}$ takes values in $\{0,1\}$, we have
$$
H_{i,j}(W_i,W'_j)^d=H_{i,j}(W_i,W'_j).
$$
Thus,
\begin{equation}\label{eq:Q12 for multiple cliques}
\begin{aligned}
Q_{\{1,2\}}&=d^{2d}\sum_{(i,j)\in[m]^2}\ee\left[H_{i,j}(W_i,W'_j)^d\right]\\
&=d^{2d}\sum_{(i,j)\in[m]^2}\ee\left[H_{i,j}(W_i,W'_j)\right]\\
&=d^{2d}\sum_{(i,j)\in[m]^2}\Pr\left(H_{i,j}(W_i,W'_j)=1\right)\\
&=d^{2d}\sum_{i\neq j}\Pr\left(H_{i,j}(W_i,W'_j)=1\right)\\
&=d^{2d}\sum_{i\neq j}\sum_{a,b=1}^{\gamma}\Pr\left(W_i=e_{(a,b)}\right)\Pr\left(W'_j=e_{(a,b)}\right)\\
&\leq d^{2d}m^2\gamma^2\rho^2\\
&=d^{2d}\lambda.
\end{aligned}
\end{equation}

Combining the bounds \Cref{eq:Q0 for multiple cliques,eq:Q1Q2 for multiple cliques,eq:Q12 for multiple cliques}, we obtain, for any integer $d \geq 2$,
$$
\ee\bigg[\bigg(\sum_{(i,j)\in[m]^2}H_{i,j}(W_i,W'_j)\bigg)^d\bigg]\leq C_0^{2d}\left[\lambda^d+2d^d\lambda \bigg(\frac{k^2}{m}\bigg)^{d-1}+d^{2d}\lambda\right],
$$
for some universal constant $C_0>0$. %

Combining the above bound with \Cref{eq:first moment bound for multiple cliques}, which deals with the case of $d=1$, and with \Cref{eq:decoupling of U-statistics for multiple cliques}, and taking a sufficiently large universal constant $C\geq1$ to absorb the constants from decoupling and the moment bounds, we arrive at the bound
\begin{align*}
\bbE_{\bfY\sim\bbQ}\left[L_{\tilde\nu}^{\leq D}(\bfY)^2\right] & \leq 1+\sum_{d=1}^D\frac{\mu^{2d}}{d!}C^d\left[\lambda^d+2d^d\lambda\left(\frac{k^2}{m}\right)^{d-1}+d^{2d}\lambda\right].
\end{align*}
We now handle the terms on the right-hand side of the last expression. 
For the first term, we have
$$
\sum_{d=1}^D\frac{\mu^{2d}}{d!}C^d\lambda^d\leq \exp(C\mu^2\lambda)-1,
$$
which is $O(1)$ whenever $\mu=O\big(m/(\gamma k^2)\big)$.
Next, we bound the quantity
$$
\sum_{d=1}^D\frac{\mu^{2d}}{d!}C^d\left[2d^d\lambda\bigg(\frac{k^2}{m}\bigg)^{d-1}+d^{2d}\lambda\right].
$$
We distinguish two cases.

\textbf{Case 1: $\lambda\leq 1$.}
Since $\lambda=\gamma^2(k^2/m)^2$ and $\gamma\geq 1$, we have $k^2/m\leq 1$. Therefore, for every $d\geq 1$,
$$
2d^d\lambda \bigg(\frac{k^2}{m}\bigg)^{d-1}+d^{2d}\lambda\leq 3\lambda d^{2d}.
$$
Using $d!\geq(d/e)^d$, we obtain
$$
\sum_{d=1}^D\frac{\mu^{2d}}{d!}C^d\left[2d^d\lambda\bigg(\frac{k^2}{m}\bigg)^{d-1}+d^{2d}\lambda\right]\leq 3\lambda\sum_{d=1}^D\left(eC\mu^2d\right)^d\leq 3\lambda\sum_{d=1}^D\left(eC\mu^2D\right)^d,
$$
which is of order $O(1)$ whenever $\mu^2D\leq c_1$ for a sufficiently small constant $c_1>0$.

\textbf{Case 2: $\lambda>1$.}
We first note that, for every $d\geq 1$,
\begin{equation}\label{eq:bounded.by.lambdad}
\lambda \bigg(\frac{k^2}{m}\bigg)^{d-1}\leq\lambda^d.
\end{equation}
When  $k^2/m\leq 1$, this bound is straightforward. On the other hand, if $k^2/m>1$, then
$$
\lambda \bigg(\frac{k^2}{m}\bigg)^{d-1}=\gamma^2\bigg(\frac{k^2}{m}\bigg)^{d+1}\leq\gamma^{2d}\bigg(\frac{k^2}{m}\bigg)^{2d}=\lambda^d,
$$
as claimed, where we used the facts that $d+1\leq 2d$ and $\gamma\geq 1$. 

Using \Cref{eq:bounded.by.lambdad} and since $\lambda>1$, 
$$
2d^d\lambda \bigg(\frac{k^2}{m}\bigg)^{d-1}+d^{2d}\lambda\leq 3d^{2d}\lambda^d.
$$
It follows that
\begin{align*}
\sum_{d=1}^D\frac{\mu^{2d}}{d!}C^d\left[2d^d\lambda\bigg(\frac{k^2}{m}\bigg)^{d-1}+d^{2d}\lambda\right]&\leq 3\sum_{d=1}^D\frac{(C\mu^2\lambda)^dd^{2d}}{d!}\\
&\leq 3\sum_{d=1}^D\left(eC\mu^2\lambda d\right)^d\\
&\leq 3\sum_{d=1}^D\left(eC\mu^2\lambda D\right)^d,
\end{align*}
where in the second inequality we use again the bound $d!\geq(d/e)^d$.
The last quantity is $O(1)$ whenever $\mu^2\lambda D\leq c_2$ for a sufficiently small constant $c_2>0$.

A sufficient condition for both cases is
$$
\mu\leq \frac{c}{\sqrt D}\left(\frac{m}{\gamma k^2}\wedge 1\right),
$$
where $c=\min\{\sqrt{c_1},\sqrt{c_2}\}$.
For $D\geq1$, this condition also gives $\mu^2\lambda=O(1)$, so the degree-$D$ likelihood-ratio second moment is $O(1)$.
By \Cref{lemma: low-degree second moment method}, no degree-$D$ polynomial strongly separates the two hypotheses in this regime.

\subsubsection{Multiple random-size hypercliques}
Throughout this part, $t\geq3$ is fixed, and constants may depend on $t$ but not on $m$, $\gamma$, $k$, $d$, or $D$.
We recall \Cref{eq:definition of W_v} and define the kernel function $h:(\{0,1\}^{\gamma \times \gamma})^t\to\mathbb R$ by
$$
h(w_1,...,w_t):=\sum_{a,b=1}^{\gamma}\prod_{l=1}^t\mathds 1\{w_l=e_{(a,b)}\}.
$$
Following the analysis in \Cref{eq:U-statistic representation of overlap size}, the overlap size $|S\cap S'|$ can be expressed as
\begin{equation}
	|S\cap S'|=\sum_{\bfi\in[m]_t}h(W_{i_1},...,W_{i_t}),
\end{equation}
which can be seen as a $U$-statistic of order $t$ with kernel $h$.
For the decoupled version, let us define the kernel
$$
H_{\bfi}(w_1,...,w_t):=\mathds1\{i_1,...,i_t\text{ are pairwise distinct}\}h(w_1,...,w_t),
$$
for each $\bfi\in[m]^t$.
Note that $H_{\bfi}$ is zero whenever any two indices coincide, and is equal to $h(w_1,...,w_t)$ otherwise.
Therefore, since $h$ is symmetric, we obtain
\begin{align*}
|S\cap S'| &=\sum_{\bfi\in [m]_t} h(W_{i_1},...,W_{i_t})\\
&=\frac{1}{t!}\sum_{\bfi \in [m]^t} H_{\bfi}(W_{i_1},...,W_{i_t})\\
&=:U(W_1,...,W_m).
\end{align*}
Let the decoupled version of $U$ be defined as
\begin{align*}
U'(W^{(1)}_1,\ldots,W^{(1)}_m;\ldots;W^{(t)}_1,\ldots,W^{(t)}_m) =\frac{1}{t!}\sum_{\bfi \in [m]^t} H_{\bfi}(W^{(1)}_{i_1},...,W^{(t)}_{i_t}),
\end{align*}
where the $t$ sequences $(W_i^{(j)})_{i=1}^m$, $j\in[t]$, are mutually independent copies of $(W_i)_{i=1}^m$ and are independent of the original sequence.
Then, by \Cref{lemma: decoupling of t order U-statistics}, there exists a constant $c>0$ depending only on $t$ such that
$$
\ee\bigg[\bigg(\sum_{\bfi \in [m]^t} H_{\bfi}(W_{i_1},...,W_{i_t})\bigg)^d\bigg]\leq c^d\ee\bigg[\bigg(\sum_{\bfi \in [m]^t} H_{\bfi}(W^{(1)}_{i_1},...,W^{(t)}_{i_t})\bigg)^d\bigg],
$$
for every integer $d\geq1$, which implies that
\begin{align}\label{eq:decoupling of U-statistics for multiple hypercliques}
	\ee[|S\cap S'|^d] \leq {c'}^d \ee\bigg[\bigg(\sum_{\bfi \in [m]^t} H_{\bfi}(\bfW_\bfi)\bigg)^d\bigg],
\end{align}
for some constant $c'>0$ depending only on $t$ and $\bfW_\bfi:=(W^{(1)}_{i_1},...,W^{(t)}_{i_t})$.
For $d=1$, the decoupled sum has expectation at most $m^t\gamma^2\rho^t$.
For integers $d\geq2$, we apply \Cref{lemma: moment bound of nonnegative U-statistics} to bound the $d$-th moment:
\begin{align*}
\ee\bigg[\bigg(\sum_{\bfi \in [m]^t} H_{\bfi}(\mathbf W_{\bfi})\bigg)^d\bigg] &\leq C^{td} \sum_{J \subseteq [t]} \Bigg[d^{|J|d}\sum_{\bfi_J}\ee_{J}\Bigg\{\Bigg(\sum_{\bfi_{J^c} \in [m]^{t-|J|}}\ee_{{J^c}}\big[H_{\bfi}(\mathbf W_{\bfi})\,\big|\,\mathcal F_J\big]\Bigg)^d\Bigg\}\Bigg]\\
	&= C^{td} \sum_{r=0}^t \sum_{\substack{J\subseteq[t]\\|J|=r}}\Bigg[d^{rd}\sum_{\bfi_J\in[m]^r}\ee_{J}\Bigg\{\Bigg(\sum_{\bfi_{J^c} \in [m]^{t-r}}\ee_{{J^c}}\big[H_{\bfi}(\mathbf W_{\bfi})\,\big|\,\mathcal F_J\big]\Bigg)^d\Bigg\}\Bigg],
\end{align*}
where $\mathcal F_J=\sigma \big(W^{(j)}_i:j\in J,~i\in[m]\big)$, $J^c$ is the complement of $J$ in $[t]$, $\sum_{\bfi_{\emptyset}}a=a$ and $\mathbb E_\emptyset[a]=a$ by convention.
Let
\begin{align*}
Q_J(\bfi_J) &:= \sum_{\bfi_{J^c} \in [m]^{t-r}} \ee_{{J^c}}\big[H_{\bfi}(\mathbf W_{\bfi})\,\big|\,\mathcal F_J\big]\\
&= \sum_{\bfi_{J^c} \in [m]^{t-r}} \Pr_{J^c}\big(H_{\bfi}(\mathbf W_{\bfi})=1\,\big|\,\mathcal F_J\big),
\end{align*}
where $r=|J|.$
For $J\neq\emptyset$, let $\mathcal E_{\bfi_J}$ be the event that the $r$ variables $(W^{(j)}_{i_j})_{j\in J}$ have the same nonzero label, i.e.,
$$
\mathcal E_{\bfi_J}=\big\{\exists (a,b)\in [\gamma]^2 \text{ such that } \forall j \in J,~ W^{(j)}_{i_{j}}=e_{(a,b)}\big\}.
$$ 
Then, under the event $\mathcal E_{\bfi_J}$, the only way to have $H_{\bfi}(\mathbf W_{\bfi})=1$ is to have $W^{(j)}_{i_{j}}=e_{(a,b)}$ for all $j\in J^c$ as well, where $(a,b)$ is the common label determined by the variables indexed by $J$.
By independence, these remaining labels match with probability $\rho^{t-r}$, where $\rho:=k^2/m^2$; repeated vertex indices can only make $H_{\bfi}$ vanish.
On the other hand, if $\mathcal E_{\bfi_J}$ does not hold, we have $H_{\bfi}(\mathbf W_{\bfi})=0$ almost surely.
Therefore, we have
\begin{align*}
Q_J(\bfi_J) = \sum_{\bfi_{J^c} \in [m]^{t-r}}\Pr_{J^c}\big(H_{\bfi}(\mathbf W_{\bfi})=1\,\big|\,\mathcal F_J\big) \leq m^{t-r}\rho^{t-r}\mathds 1\{\mathcal E_{\bfi_J}\}.
\end{align*}
Also note that $\Pr(\mathcal E_{\bfi_J}) =\gamma^2\rho^r$ for $r\geq1$.
Therefore, 
\begin{align*}
&\sum_{\substack{J\subseteq[t]\\|J|=r}}\sum_{\bfi_J\in[m]^r}\ee_{J}\Bigg\{\Bigg(\sum_{\bfi_{J^c} \in [m]^{t-r}}\ee_{{J^c}}\big[H_{\bfi}(\mathbf W_{\bfi})\,\big|\,\mathcal F_J\big]\Bigg)^d\Bigg\}\\
&\qquad= \sum_{\substack{J\subseteq[t]\\|J|=r}}\sum_{\bfi_J\in[m]^r}\ee_{J}\big[Q_J(\bfi_J)^d\big]\\
&\qquad\leq \sum_{\substack{J\subseteq[t]\\|J|=r}}\sum_{\bfi_J\in[m]^r} m^{(t-r)d}\rho^{(t-r)d} \Pr(\mathcal E_{\bfi_J})\\
&\qquad\leq \binom tr m^r(m\rho)^{(t-r)d}\gamma^2\rho^r\\
&\qquad=\binom tr\gamma^2 (m\rho)^{(t-r)d+r},
\end{align*}
for $r\geq 1$. For $J=\emptyset,$ we have
\begin{align*}
Q_\emptyset^d &= \bigg(\sum_{\bfi \in [m]^t} \Pr\big(H_{\bfi}(\mathbf W_{\bfi})=1\big)\bigg)^d\leq (m^t\cdot\gamma^2\rho^t)^d.
\end{align*}
Summarizing the above analysis, we obtain
\begin{align*}
\ee\bigg[\bigg(\sum_{\bfi \in [m]^t} H_{\bfi}(\mathbf W_{\bfi})\bigg)^d\bigg] &\leq C^{td} \sum_{r=0}^t \sum_{\substack{J\subseteq[t]\\|J|=r}}\Bigg[d^{rd}\sum_{\bfi_J\in[m]^r}\ee_{J}\Bigg\{\Bigg(\sum_{\bfi_{J^c} \in [m]^{t-r}}\ee_{{J^c}}\big[H_{\bfi}(\mathbf W_{\bfi})\,\big|\,\mathcal F_J\big]\Bigg)^d\Bigg\}\Bigg]\\
&\leq C^{td} \Bigg[(m^t\cdot\gamma^2\rho^t)^d + \sum_{r=1}^t \binom tr d^{rd}\gamma^2 (m\rho)^{(t-r)d+r}\Bigg]\\
&\leq C^{td} \Bigg[\gamma^{2d}\bigg(\frac{k^2}{m}\bigg)^{td} + \gamma^2\sum_{r=1}^t \binom trd^{rd} \bigg(\frac{k^2}{m}\bigg)^{(t-r)d+r}\Bigg]\\
&=C^{td} \Bigg[\lambda^d + \lambda\sum_{r=1}^t \binom trd^{rd} \bigg(\frac{k^2}{m}\bigg)^{(t-r)(d-1)}\Bigg]\\
&\leq C_1^{td} \Bigg[\lambda^d + \lambda\sum_{r=1}^t d^{rd} \bigg(\frac{k^2}{m}\bigg)^{(t-r)(d-1)}\Bigg],
\end{align*}
where $\lambda=\gamma^2(k^2/m)^t$ and $C_1=2C$, and the last inequality follows from the fact that $\binom tr \leq 2^t$ for all $r\in [t]$.
The same bound holds for $d=1$ after enlarging $C_1$ if necessary, by the first-moment bound above.
Therefore, combining this bound with \Cref{eq:decoupling of U-statistics for multiple hypercliques} and further enlarging $C_1$ to absorb the decoupling constant, we have
\begin{align*}
	\bbE_{\bfY\sim\bbQ} &\leq 1+\sum^D_{d=1} \frac{\mu^{2d}}{d!} C_1^{td} \Bigg[\lambda^d + \lambda\sum_{r=1}^t d^{rd} \bigg(\frac{k^2}{m}\bigg)^{(t-r)(d-1)}\Bigg]\\
	&\leq 1+\sum^D_{d=1} \frac{\mu^{2d}}{d!} C_1^{td} \lambda^d+\sum^D_{d=1} \frac{\mu^{2d}}{d!} C_1^{td} \lambda\sum_{r=1}^t d^{rd} \bigg(\frac{k^2}{m}\bigg)^{(t-r)(d-1)},
\end{align*}
and it is sufficient to show
\begin{align*}
(\mathrm{I}):&~\sum^D_{d=1} \frac{\mu^{2d}}{d!} C_1^{td} \lambda^d=O(1),\quad (\mathrm{II}):~\sum^D_{d=1} \frac{\mu^{2d}}{d!} C_1^{td} \lambda\sum_{r=1}^t d^{rd} \bigg(\frac{k^2}{m}\bigg)^{(t-r)(d-1)}=O(1).
\end{align*}
For (I), note that
$$
\sum^D_{d=1} \frac{\mu^{2d}}{d!} C_1^{td} \lambda^d \leq \exp\big(\mu^2 C_1^t \lambda\big)-1,
$$
and thus it suffices to show $\mu^2 C_1^t \lambda=O(1),$ which is equivalent to $\mu = O\big(\frac{m^{t/2}}{\gamma k^t}\big).$ 
For (II), we split into two cases.
\begin{itemize}
	\item \textbf{Case 1:} $\lambda \leq 1.$ In this case, observe that 
	\begin{align*}
	\frac{k^2}{m} \leq \lambda^{1/t} \leq 1,
	\end{align*}
	since $\gamma \geq 1.$ Therefore, for every $r\in [t]$ and $d\in \mathbb N,$ we have
	\begin{align*}
	\lambda\sum^t_{r=1} d^{rd} \bigg(\frac{k^2}{m}\bigg)^{(t-r)(d-1)} \leq \lambda\sum^t_{r=1}d^{td}\leq td^{td},
	\end{align*}
	as $\lambda \leq 1$ and $d^{rd} \leq d^{td}$. Therefore, note that (II) can be upper bounded by
	\begin{align*}
	\sum^D_{d=1} \frac{\mu^{2d}}{d!} C_1^{td} \lambda\sum_{r=1}^t d^{rd} \bigg(\frac{k^2}{m}\bigg)^{(t-r)(d-1)} &\leq \sum^D_{d=1} \frac{\mu^{2d}}{d!} C_1^{td} td^{td} \\
	&\leq t\sum^D_{d=1} {\mu^{2d}} e^dC_1^{td} d^{(t-1)d} \quad \big(\,\because d! \geq (d/e)^d\,\big)\\
	&\leq t\sum^D_{d=1} {\mu^{2d}} e^dC_1^{td} D^{(t-1)d}\\
	&= t\sum^D_{d=1} \big(e\mu^2 C_1^t D^{t-1}\big)^d,
	\end{align*}
	which is $O(1)$ when $\mu^2 D^{t-1}\leq c_1,$ for a sufficiently small constant $1>c_1>0.$ 
	\item \textbf{Case 2:} $\lambda > 1.$ If $\frac{k^2}{m} \leq 1,$ we follow the same argument as in Case 1 to obtain
    $$
    \lambda d^{rd} \bigg(\frac{k^2}{m}\bigg)^{(t-r)(d-1)} \leq \lambda d^{rd} \leq \lambda d^{td} \leq \lambda^d d^{td},
    $$
    as $\lambda > 1.$
    On the other hand, if $\frac{k^2}{m} > 1,$ then observe that for every $r\in [t]$ and $d\in \mathbb N,$
    $$
    (t-r)d+r \leq td,
    $$
    and thus
    \begin{align*}
        \lambda d^{rd}\bigg(\frac{k^2}{m}\bigg)^{(t-r)(d-1)}=d^{rd}\gamma^2\bigg(\frac{k^2}{m}\bigg)^{(t-r)d+r} \leq d^{td}\gamma^2\bigg(\frac{k^2}{m}\bigg)^{td}.
    \end{align*}
    Since
    $$
    \gamma^2\bigg(\frac{k^2}{m}\bigg)^{td} \leq \gamma^{2d}\bigg(\frac{k^2}{m}\bigg)^{td} =\lambda^d
    $$ we have
    \begin{align*}
        \lambda d^{rd} \bigg(\frac{k^2}{m}\bigg)^{(t-r)(d-1)} \leq d^{td}\lambda^d.
    \end{align*}
	Therefore, in either case, we have
	\begin{align*}
	\lambda\sum^t_{r=1} d^{rd} \bigg(\frac{k^2}{m}\bigg)^{(t-r)(d-1)} \leq t \lambda^d d^{td},
	\end{align*}
	and thus (II) can be upper bounded by
	\begin{align*}
	\sum^D_{d=1} \frac{\mu^{2d}}{d!} C_1^{td} \lambda\sum_{r=1}^t d^{rd} \bigg(\frac{k^2}{m}\bigg)^{(t-r)(d-1)} &\leq \sum^D_{d=1} \frac{\mu^{2d}}{d!} C_1^{td} t \lambda^d d^{td} \\
	&\leq \sum^D_{d=1} \mu^{2d}e^d C_1^{td} t \lambda^d d^{(t-1)d} \\
	&\mathrel{\leq} t\sum^D_{d=1} \big(e\mu^2 C_1^t \lambda D^{t-1}\big)^d,
	\end{align*}
	which is $O(1)$ when $\mu^2\lambda D^{t-1} \leq c_2,$ for a sufficiently small constant $c_2>0.$ 
\end{itemize}
Combining the above two cases, we conclude that a sufficient condition is
$$
\mu \leq \frac{c}{D^{(t-1)/2}}\bigg(\frac{m^{t/2}}{\gamma k^{t}}\wedge 1\bigg),
$$
for $c=\min\{\sqrt{c_1},\sqrt{c_2}\}$.
For $D\geq1$, this condition also bounds term (I), so the degree-$D$ likelihood-ratio second moment is $O(1)$ and \Cref{lemma: low-degree second moment method} rules out strong separation by degree-$D$ polynomials.
For fixed $t$ and polylogarithmic $D$, the factors $D^{-1/2}$ and $D^{-(t-1)/2}$ in the two bounds are logarithmic losses.
Thus, when $D$ satisfies the degree requirement in \Cref{conjecture: low-degree polynomial conjecture}, these estimates recover the computational lower bounds in \Cref{theorem: multiple cliques problem,corollary: multiple planted hypercliques problem}, up to logarithmic factors, under that conjecture.

\section{Additional examples of latent structures}\label{sec:app-additional-examples}
\subsection{Examples from \texorpdfstring{\citet{addario2010combinatorial}}{Addario-Berry et al. (2010)}}
In this subsection, we collect the results from other combinatorial testing examples discussed in \citet{addario2010combinatorial} for completeness.
First, note that the second moment of the low-degree likelihood ratio is always bounded above by that of the full likelihood ratio, i.e.,
\begin{align}\label{equation: low-degree LR vs full LR}
\ee_{\bfY\sim\bbQ}[L_\nu^{\leq D}(\bfY)^2]
=\ee_{S,S'\sim \nu}\Bigg[\sum_{d=0}^D \frac{\mu^{2d}}{d!}|S\cap S'|^d\Bigg]\leq \ee_{S,S'\sim \nu}\big[\exp(\mu^2 |S\cap S'|)\big]
= \ee_{\bfY\sim\bbQ}[L_\nu(\bfY)^2].
\end{align}
This implies that a statistical lower bound obtained by controlling the second moment of the full likelihood ratio also yields a computational lower bound through the low-degree framework.
Here, the lower bound against polynomial-time tests is unconditional, since statistical impossibility already rules out every test.
This observation is particularly useful for problems in which there is no statistical-computational gap.
Indeed, for the problems considered in this subsection, the statistical lower bounds established by \citet{addario2010combinatorial} are matched by upper bounds achieved by polynomial-time tests, so no separate analysis of computational lower bounds is necessary.

Accordingly, our goal in this subsection is not to establish new lower-bound results for these problems.
Rather, we aim to complement the analysis of \citet{addario2010combinatorial}, which relies on structural properties of the underlying set classes, such as symmetry and negative association, by providing alternative proofs based on more explicit characterizations of the overlap distribution $|S\cap S'|$ or problem-specific combinatorial arguments.

Throughout this subsection, we assume that the prior $\nu$ is uniform over $\mathcal C$ and that $|s|=K$ for all $s\in\mathcal C$.
Also, we denote $N=|\mathcal C|$.
In each overlap calculation, $S$ and $S'$ are independent draws from $\nu$, and $C>1$ is a fixed second-moment bound independent of the ambient dimension.
For each problem, we aim to derive a sufficient condition on $\mu$ under which $\ee_{\bfY\sim\bbQ}[L_\nu(\bfY)^2]$ remains bounded by a constant, thereby establishing the impossibility of strong detection via \Cref{lemma: second moment method}.
We also identify the test that achieves a matching upper bound.
	\subsubsection{Disjoint sets}\label{appendix: disjoint sets}
	Consider the case where $\mathcal C$ is the collection of $N$ disjoint sets of size $K$.
	In this case, $|S\cap S'|$ is a random variable that takes the value $K$ with probability $1/N$ and $0$ otherwise.
	Therefore, we have
	$$
	\ee_{\bfY\sim\bbQ}[L_\nu(\bfY)^2] = \ee_{S,S'\sim \nu}[\exp(\mu^2 |S\cap S'|)] = \frac{1}{N}\exp(\mu^2 K) + \frac{N-1}{N},$$
	and thus, if
	$$\mu \leq \sqrt{\frac{1}{K}\log\big(1+N(C-1)\big)},$$
	then $\ee_{\bfY\sim\bbQ}[L_\nu(\bfY)^2]$ is bounded by $C$.
	When $N\to\infty$, the matching upper bound (up to a constant factor) can be obtained by the scan test in \Cref{proposition: gaussian test upper bound}, which runs in $O(NK)$ time, hence $O(n)$ time, because the candidate sets are disjoint.
	If $N$ remains bounded, the corresponding condition for strong detection is $\mu\sqrt K\to\infty$.

	\subsubsection{\texorpdfstring{$K$-subsets}{K-subsets}}\label{appendix: K-sets}
	We consider the case where $\mathcal C$ consists of all $K$-subsets of $[n]$.
	In this case, the overlap size $|S\cap S'|$ between two randomly chosen $K$-sets $S$ and $S'$ from $\mathcal C$ follows a hypergeometric distribution with parameters $(n,K,K)$.
	By the convex-order comparison with a binomial random variable \citep[Theorem~4]{hoeffding1963probability}, we have
	\begin{align*}
		\ee_{\bfY\sim\bbQ}[L_\nu(\bfY)^2] &= \ee_{S,S'\sim \nu}[\exp(\mu^2 |S\cap S'|)] \\
		&\leq \ee[\exp(\mu^2 B)]\\
		&= \bigg(1+\frac{K}{n}\big(e^{\mu^2}-1\big)\bigg)^K\\
		&\leq \exp\bigg(\frac{K^2}{n}\big(e^{\mu^2}-1\big)\bigg),
	\end{align*}
	where $B\sim \mathrm{Binomial}(K,K/n).$ Therefore, if 
	$$
	\mu \leq \sqrt{\log\bigg(1+\frac{n}{K^2}\log C\bigg)},
	$$
	then $\ee_{\bfY\sim\bbQ}[L_\nu(\bfY)^2]$ is bounded by $C$.
	We highlight two regimes.
	$(i)$ If $K^2=O(n^{1-\delta})$ for some fixed $\delta\in(0,0.1)$, then $\log(1+\frac{n}{K^2}\log C)\geq \delta\log n +O(1)$, so $\log(1+\frac{n}{K^2}\log C)\geq \frac\delta 2\log n$ for sufficiently large $n$.
	Therefore, it suffices to have $\mu\leq \sqrt{\frac\delta2\log n}$.
	A matching upper bound up to a constant factor is achieved by the outlier test in \Cref{proposition: gaussian test upper bound}, which runs in $O(n)$ time.
	$(ii)$ If $K^2=\omega(n)$, then $\frac{n}{2K^2}\log C\leq \log(1+\frac{n}{K^2}\log C)$ for sufficiently large $n$.
	Thus, the sufficient condition is $\mu\leq \sqrt{\frac{1}{2}\frac{n}{K^2}\log C}$.
	A matching signal scale is achieved by the average test in \Cref{proposition: gaussian test upper bound}, which runs in $O(n)$ time and achieves strong detection whenever $\mu K/\sqrt n\to\infty$.

\subsubsection{Perfect matchings}\label{appendix: perfect matchings}
We consider the case where $\mathcal C$ consists of all perfect matchings in a complete bipartite graph $K_{m,m}$.
We have $n=m^2$, $N=m!$, and $K=m$.
As pointed out in \citet{addario2010combinatorial}, the overlap size $|S\cap S'|$ between two randomly chosen perfect matchings $S$ and $S'$ from $\mathcal C$ follows the distribution of the number of fixed points in a random permutation of size $m$; let us denote this random variable by $Z$, i.e., $Z:=\sum^m_{i=1} \mathds 1\{\pi(i)=i\}$, where $\pi$ is a uniformly random permutation of $[m]$.
Then, we have
$$
\Pr(Z=r)=\frac{\binom m r D_{m-r}}{m!}, \quad r=0,1,...,m,
$$
where $D_{m-r}$ is the number of derangements of size $m-r$, with $D_0=1$.
Hence,
\begin{align*}
	\ee_{\bfY\sim\bbQ}[L_\nu(\bfY)^2] &= \ee_{S,S'\sim \nu}[\exp(\mu^2 |S\cap S'|)]\\
	&= \ee_Z[\exp(\mu^2 Z)] \\
	&= \sum_{r=0}^m \exp(\mu^2 r) \frac{\binom m r D_{m-r}}{m!} \\
	&= \sum_{r=0}^m \exp(\mu^2 r) \frac{1}{r!} \sum_{j=0}^{m-r} \frac{(-1)^j}{j!} \\
	&= \sum_{\substack{r,j\geq0\\r+j\leq m}} \frac{\exp(\mu^2 r)(-1)^j}{r! j!}\\
	&= \sum_{k=0}^m \frac{1}{k!} \sum_{r=0}^k \binom k r (\exp(\mu^2))^r(-1)^{k-r}\quad \text{where}~k:=r+j\\
	&= \sum_{k=0}^m \frac{1}{k!}(\exp(\mu^2)-1)^k \\
	&\leq \exp(\exp(\mu^2)-1).
\end{align*}
Therefore, if 
$$\mu \leq \sqrt{\log(1+\log C)},
$$
then $\ee_{\bfY\sim\bbQ}[L_\nu(\bfY)^2]$ is bounded by $C$.
Thus, every bounded sequence $\mu=O(1)$ precludes strong detection, whereas the average test in \Cref{proposition: gaussian test upper bound} achieves strong detection whenever $\mu\to\infty$, since $K/\sqrt n=1$.

\begingroup
\makeatletter
\let\originalNote\@todonotes@drawMarginNote
\def\@todonotes@drawMarginNote{%
    \raisebox{3cm}[\height][\depth]{\originalNote}%
}
\makeatother
\endgroup

\paragraph{Negative association of matching indicators.}
\citet[Section~4.3]{addario2010combinatorial} establish the same statistical lower bound for perfect matchings by invoking Proposition~3.4, which requires negative association of the edge-inclusion indicators along a fixed reference matching.
As the calculation below shows, these indicators are not negatively associated when $m\geq2$.
The exact overlap calculation above nevertheless yields the stated lower bound without this assumption.
For completeness, we recall the definition of negative association and give a simple counterexample.

\begin{definition}[Negative association \citep{addario2010combinatorial}]
	A collection of random variables $Y_1,...,Y_n$ is said to be negatively associated if for any pair of disjoint subsets $A_1,A_2 \subset [n]$ and any two coordinatewise nondecreasing functions $f:\mathbb R^{|A_1|}\rarr \mathbb R$ and $g:\mathbb R^{|A_2|}\rarr \mathbb R$,
	$$
	\ee[f(Y_i:i\in A_1)g(Y_j:j\in A_2)] \leq \ee[f(Y_i:i\in A_1)]\ee[g(Y_j:j\in A_2)].
	$$
\end{definition}

For $m\geq2$, fix the reference matching $S'=\{(1,1),(2,2),...,(m,m)\}$ by symmetry, and let $S$ be uniformly distributed over $\mathcal C$.
For $(i,j)\in S'$, define the edge-inclusion indicator $I_{i,j}:=\mathds 1\{(i,j)\in S\}$.
Consider $I_{1,1}$ and $I_{2,2}$, and take $f(x)=g(x)=x$.
Then, we have

$$\ee[I_{1,1}I_{2,2}] = \Pr((1,1)\in S, (2,2)\in S) = \frac{(m-2)!}{m!} = \frac{1}{m(m-1)},$$
and
$$\ee[I_{1,1}]\ee[I_{2,2}] = \Pr((1,1)\in S)\Pr((2,2)\in S) = \frac{(m-1)!}{m!}\cdot \frac{(m-1)!}{m!} = \frac{1}{m^2}.$$
Since $1/(m(m-1)) > 1/m^2$ for all $m\geq 2,$ this shows that $I_{1,1}$ and $I_{2,2}$ are not negatively associated.

	\subsubsection{Stars}\label{appendix: star}
	Consider a network of $m$ nodes, where a \textit{star} is a subgraph of the complete graph consisting of all $K=m-1$ edges incident to a fixed node.
	Consider the set $\mathcal C$ of all stars in the complete graph.
	For $m\geq3$, we have $n=\binom m2$ and $N=m$.
	In this case, the overlap size $|S\cap S'|$ has the distribution
	$$
	|S\cap S'|=\begin{cases}
	m-1, &\text{with probability}~ \frac{1}{m},\quad \text{the centers coincide}\\
	1, &\text{with probability}~ \frac{m-1}{m},\quad \text{otherwise}.
	\end{cases}
	$$
	Therefore, we have
	\begin{align*}
		\ee_{\bfY\sim\bbQ}[L_\nu(\bfY)^2] &= \ee_{S,S'\sim \nu}[\exp(\mu^2 |S\cap S'|)] \\
		&= \frac{1}{m}\exp(\mu^2 (m-1)) + \frac{m-1}{m}\exp(\mu^2)\\
		&= \frac{1}{m}\big(\exp(\mu^2 (m-1))-1\big) + \frac{m-1}{m}\big(\exp(\mu^2)-1\big) + 1\\
		&\leq \frac{2}{m}\big(\exp(\mu^2 (m-1))-1\big)+1,
	\end{align*}
	where the inequality follows from $e^{(m-1)\mu^2}-1\geq (m-1)(e^{\mu^2}-1)$ by the convexity of $f(x)=e^{x}$.
	Thus, if
	$$\mu \leq \sqrt{\frac{1}{m-1}\log\bigg(1+\frac{m(C-1)}{2}\bigg)},
	$$
	then $\ee_{\bfY\sim\bbQ}[L_\nu(\bfY)^2]$ is bounded by $C$.
	Note that the scan test achieves the matching upper bound; all star sums can be accumulated in one pass over the $n$ observations, so this scan takes $O(n)$ time.

	\subsubsection{Spanning trees}\label{appendix: spanning tree}
	Consider a network of $m$ nodes, where a \textit{spanning tree} is a subgraph of the complete graph which is a tree and connects all $m$ nodes.
	We have $n=\binom m2$, $N=m^{m-2}$, and $K=m-1$.
	Note that while \citet{addario2010combinatorial} provide two different ways to obtain a statistical lower bound in this case, we provide an alternative approach focusing on combinatorial properties of spanning trees.
	A key step is to interpret each realization of $S\cap S'$ such that $|S\cap S'|=r$ as a forest with $m$ nodes and $r$ edges, including isolated vertices.
	For example, when $m=6$ and $r=3$, one possible realization of $S\cap S'$ is $\{(1,2),(2,3),(4,5)\}$, which can be interpreted as a forest $F$ with $c=3$ components whose vertex sets are $\big\{\{1,2,3\},\{4,5\},\{6\}\big\}$.
	We identify each such forest $F$ with its edge set; its vertex set is always $[m]$.
	Write $F\subset A$ for the event that all edges of $F$ belong to the edge set $A$; for example, $F\subset S\cap S'$ means that all edges in the forest $F$ are included in both $S$ and $S'$.
	Write $r(F)$ and $c(F)$ for the numbers of edges and components of $F$, respectively; then $c(F)=m-r(F)$.
	We introduce a useful lemma to count the number of spanning trees containing a given forest:
	\begin{lemma}[{\citealt[Theorem 6.1]{moon1970counting}}]\label{lemma: counting spanning trees containing a given forest}
		Let $F$ be a forest on vertex set $[m]$ with $c(F)$ components.
		If we denote the number of nodes in each component as $m_1,...,m_{c(F)}$, then the number of spanning trees containing $F$ is given by
		$$
		m^{c(F)-2}\prod_{i=1}^{c(F)} m_i.
		$$
	\end{lemma}
	Therefore, given a forest $F$, the probability that independent random spanning trees $S$ and $S'$ both contain $F$ is given by
	\begin{align*}
		\Pr(F \subset S\cap S') &= \Pr(F \subset S)\Pr(F \subset S') \\
		&= \bigg(\frac{m^{c(F)-2}\prod_{i=1}^{c(F)} m_i}{m^{m-2}}\bigg)^2 \\
		&\stackrel{(a)}= \bigg(\frac{\prod_{i=1}^{c(F)} m_i}{m^{r(F)}}\bigg)^2\\
		&\stackrel{(b)}\leq \bigg(\frac{2^{\sum_i (m_i-1)}}{m^{r(F)}}\bigg)^2\\
		&= \bigg(\frac{4}{m^2}\bigg)^{r(F)},
	\end{align*}
	where (a) follows from the relationship $c(F)=m-r(F),$ and (b) follows from the fact that $m_i \leq 2^{m_i-1}$ for all $m_i\geq 1$. 
	Now, denoting $Z:=|S\cap S'|$, observe that, for integers $0\leq R\leq m-1$, the $R$-th falling factorial of $Z$ can be expressed as
	$$
	(Z)_R = Z(Z-1)...(Z-R+1) = \sum_{I_R \in [m]_{2,R}} \mathds 1\{(i_1,j_1),...,(i_R,j_R) \in S\cap S'\},
	$$
	where we denote an ordered $R$-tuple of edges as $I_R:=\big((i_1,j_1),...,(i_R,j_R)\big)$, and the collection of ordered $R$-tuples of distinct edges as $[m]_{2,R}:=\big\{ I_R: (i_a,j_a) \in [m]_2, (i_a,j_a)\neq (i_b,j_b)~\text{for all}~ 1\leq a\neq b\leq R\big\}$.
	We use the conventions $(Z)_0=1$ and $[m]_{2,0}=\{()\}$.
	Then, the $R$-th factorial moment can be expressed as
	\begin{align*}
		\ee[(Z)_R] &= \sum_{I_R \in [m]_{2,R}} \Pr\big((i_1,j_1),...,(i_R,j_R) \in S\cap S'\big).
	\end{align*}
	Now, observe that if the edges $(i_1,j_1),...,(i_R,j_R)$ contain a cycle, then there is no spanning tree containing all those edges, and thus $\Pr\big((i_1,j_1),...,(i_R,j_R) \in S\cap S'\big)=0$.
	Therefore, the only nonzero contributions to the above summation come from those $I_R$ that form a forest with $R$ edges.
	Also, note that we can bound the size of the collection as $\big|[m]_{2,R}\big|\leq n(n-1)\cdots(n-R+1)\leq n^R$, where $n=\binom m2$.
	\begin{align*}
		\ee[(Z)_R] &= \sum_{I_R \in [m]_{2,R}} \Pr\big((i_1,j_1),...,(i_R,j_R) \in S\cap S'\big) \\
		&\leq n^R \bigg(\frac{4}{m^2}\bigg)^R \\
		&= \bigg(\frac{2(m-1)}{m}\bigg)^R\\
		&\leq 2^R.
	\end{align*}
	For $t\geq0$, we have
	$$
	\exp(tZ)=(1+(e^t-1))^Z = \sum_{R=0}^Z \frac{(e^t-1)^R}{R!} (Z)_R.
	$$
	Hence, we rewrite the moment generating function of $Z$ using the factorial moments as
	\begin{align*}
		\ee[\exp(tZ)]= \sum_{R=0}^{m-1} \frac{(e^t-1)^R}{R!} \ee[(Z)_R] \leq \sum_{R=0}^{m-1} \frac{(e^t-1)^R}{R!} 2^R \leq \sum_{R=0}^{\infty} \frac{\big(2(e^t-1)\big)^R}{R!}= \exp(2(e^t-1)).
	\end{align*}
	Therefore, we have
	\begin{align*}
		\ee_{\bfY\sim\bbQ}[L_\nu(\bfY)^2] = \ee_{S,S'\sim \nu}[\exp(\mu^2 |S\cap S'|)]\leq \exp(2(e^{\mu^2}-1)).
	\end{align*}
	Thus, if
	$$
	\mu \leq \sqrt{\log\bigg(1+\frac{\log C}{2}\bigg)},
	$$
	then $\ee_{\bfY\sim\bbQ}[L_\nu(\bfY)^2]$ is bounded by $C$.
	The matching signal scale is achieved by the average test, which strongly detects whenever $\mu\to\infty$, since $K/\sqrt n=(m-1)/\sqrt{\binom m2}\to\sqrt2$.

\subsection{Further examples}
As in the preceding subsection, let $\nu$ be uniform over $\mathcal C$, let $S,S'$ be independent draws from $\nu$, and fix a second-moment bound $C>1$ independent of the ambient dimension.
\subsubsection{\texorpdfstring{$k$-vertex stars}{k-vertex stars}}\label{appendix: k-size star}
	We consider a variant of the star graph in \Cref{appendix: star}, consisting of $k$ vertices and hence $K=k-1$ edges, where $3\leq k\leq m-1$.
	In this case, we have $n=\binom m2$ and $N= \binom m k \cdot k$.
	Note that the overlap size $|S\cap S'|$ between two randomly chosen $k$-vertex stars $S$ and $S'$ from $\mathcal C$ has the following conditional distributions:
	$$
	|S\cap S'|=\begin{cases}
	A\sim\mathrm{Hypergeometric}(m-1,k-1,k-1), & \text{if the centers coincide}\\
	B\sim\mathrm{Bernoulli}(\frac{k-1}{m-1}), & \text{if the center of $S'$ is a leaf of $S$}\\
	0, & \text{otherwise}.
	\end{cases}
	$$
	These three cases have probabilities $1/m$, $(k-1)/m$, and $(m-k)/m$, respectively.
	Therefore, denoting $p = \frac{k-1}{m-1}$, we have
	\allowdisplaybreaks
	\begin{align*}
		\ee_{\bfY\sim\bbQ}[L_\nu(\bfY)^2] &= \ee_{S,S'\sim \nu}[\exp(\mu^2 |S\cap S'|)]\\
	&=\frac{1}{m} \ee[\exp({\mu^2A})]+\frac{k-1}{m} \ee[\exp({\mu^2B})]+\frac{m-k}{m} \cdot 1\\
	&=\frac{1}{m} \ee[\exp({\mu^2A})]+\frac{k-1}{m}\big[p\{\exp({\mu^2})-1\}+1\big]+ \frac{m-k}{m}\\
	&=\frac{1}{m} \ee[\exp({\mu^2A})]+\frac{m-1}{m}+ \frac{(k-1)^2}{m(m-1)}\big\{\exp({\mu^2})-1\big\}\\
	&\stackrel{(a)}{\leq} \frac{1}{m}\big[p\{\exp(\mu^2)-1\}+1\big]^{k-1}+\frac{m-1}{m}+ \frac{(k-1)^2}{m(m-1)}\big\{\exp({\mu^2})-1\big\}\\
	&= 1+\frac{1}{m}\Big(\big[p\{\exp(\mu^2)-1\}+1\big]^{k-1}-1\Big)+ \frac{(k-1)^2}{m(m-1)}\big\{\exp({\mu^2})-1\big\}\\
	&= 1+\frac{1}{m}\Big(\exp\Big\{(k-1)\log\big[p\{\exp(\mu^2)-1\}+1\big]\Big\}-1\Big)+ \frac{(k-1)^2}{m(m-1)}\big\{\exp({\mu^2})-1\big\}\\
	&\stackrel{(b)}{\leq} 1+\frac{1}{m}\Big(\exp\big[(k-1)p\{\exp(\mu^2)-1\}\big]-1\Big)+ \frac{(k-1)^2}{m(m-1)}\big\{\exp({\mu^2})-1\big\}\\
	&=1+\frac{1}{m}\bigg(\exp\bigg[\frac{(k-1)^2}{m-1}\{\exp(\mu^2)-1\}\bigg]-1\bigg)+ \frac{(k-1)^2}{m(m-1)}\big\{\exp({\mu^2})-1\big\}\\
	&=1+\frac{1}{m}\big(e^z-1+z\big),\quad\bigg(z:=\frac{(k-1)^2}{m-1}\big\{\exp({\mu^2})-1\big\}\bigg),
	\end{align*}
	where $(a)$ uses the convex-order comparison of a hypergeometric random variable with a binomial random variable \citep[Theorem~4]{hoeffding1963probability}, and $(b)$ uses $\log(1+x)\leq x$.
	Then, since $z\leq e^z-1$, we have
	\begin{align*}
	\ee_{\bfY\sim\bbQ}[L_\nu(\bfY)^2]&\leq1+\frac{1}{m}\big(e^z-1+z\big)\\*
	&\leq 1+\frac{2}{m}(e^z-1)\\*
	&= 1+\frac{2}{m}\bigg(\exp\bigg[\frac{(k-1)^2}{m-1}\big\{\exp({\mu^2})-1\big\}\bigg]-1\bigg).
	\end{align*}
	Therefore, $\ee_{\bfY\sim\bbQ}[L_\nu(\bfY)^2]$ is bounded by $C$ if
	$$
	\mu \leq \sqrt{\log\bigg\{1 + \frac{m-1}{(k-1)^2}\log\bigg(1+\frac{m(C-1)}{2}\bigg)\bigg\}}.
	$$
	We further simplify this sufficient condition by considering two regimes.
	$(i)$ Suppose that $k^2=O(m^{1-\delta})$ for some fixed $\delta\in(0,0.1)$.
	Since
	$$
	\log\bigg(1+\frac{m-1}{(k-1)^2}\log\bigg(1+\frac{m(C-1)}{2}\bigg)\bigg) \geq \delta\log m+O(1),
	$$
	the quantity on the left is at least $\frac{\delta}{2}\log m$ for all sufficiently large $m$.
	Hence, $\mu\leq\sqrt{\frac{\delta}{2}\log m}$ is sufficient to ensure that the second moment remains bounded by $C$.
	This lower bound matches, up to a constant factor, the upper bound achieved by the outlier test, which runs in $O(n)$ time.
	$(ii)$ Suppose that $k^2=\omega(m\log m)$.
	In this regime,
	$$
	\frac{m-1}{(k-1)^2}\log\bigg(1+\frac{m(C-1)}{2}\bigg)=o(1).
	$$
	Applying $\log(1+x)\geq x/2$ for $0\leq x\leq 1$, we obtain
	$$
	\log\bigg(1+\frac{m-1}{(k-1)^2}\log\bigg(1+\frac{m(C-1)}{2}\bigg)\bigg) \geq \frac{(m-1)\log m}{4(k-1)^2}
	$$
	for all sufficiently large $m$. Consequently, it suffices to require
	$\mu\leq\frac{\sqrt{(m-1)\log m}}{2(k-1)}$.
	The resulting lower bound matches, up to a constant factor, the upper bound attained by the degree test in \Cref{remark: gaussian degree test}, which can be implemented in $O(n)$ time.

	\subsubsection{\texorpdfstring{Paths and cycles of length $k$}{Paths and cycles of length k}}\label{appendix: k-length path}
	We first consider the path case.
	For integers $1\leq k\leq m-1$, possibly depending on $m$, we consider the case where $\mathcal C$ consists of all paths of length $k$ in a complete graph with $m$ nodes.
	We have $n=\binom m2$ and $N=m(m-1)\cdots (m-k)/2$, with $K=k$.
	We take an approach similar to that in \Cref{appendix: spanning tree}.
	Note that each realization of $S\cap S'$ is a set of disjoint paths with a total of $r\leq k$ edges.
	We denote the corresponding forest on $[m]$ by $F$, and let $F_p$ be the subforest obtained by discarding all isolated vertices.
	Also, let $v(F_p)$ and $c(F_p)$ be the numbers of nodes and components in $F_p$, respectively.
	For example, when $m=5$ and $k=3$, one possible realization of $S\cap S'$ is $\{(1,2),(2,3)\}$, which can be interpreted as a forest with $c(F)=3$ components whose vertex sets are $\big\{\{1,2,3\},\{4\},\{5\}\big\}$; $F_p$ is the path $1$--$2$--$3$, with $c(F_p)=1$ and $v(F_p)=3$.
	Then, note that we have $r(F)=r(F_p)=v(F_p)-c(F_p)$, and we can count the number of paths of length $k$ containing a given $F_p$ with $v(F_p)\leq k+1$ as
	\begin{align*}
		\binom{m-v(F_p)}{k+1-v(F_p)} \cdot 2^{c(F_p)} (k+1-r(F_p))!  /2,
	\end{align*}
	where the first term accounts for the number of ways to choose the remaining nodes, the second term accounts for the number of directions for each component, the last term accounts for the number of ways to arrange all the nodes and components, and the division by 2 accounts for the fact that a path can be traversed in two directions.
	If $v(F_p)>k+1$, no path of length $k$ can contain $F_p$.
	Therefore, given a feasible $F_p$, the probability that a uniformly random path $S$ of length $k$ contains $F_p$ is bounded by
	\allowdisplaybreaks
	\begin{align*}
		\Pr(F_p \subset S) &= \frac{\binom{m-v(F_p)}{k+1-v(F_p)} \cdot 2^{c(F_p)} (k+1-r(F_p))!  /2}{m(m-1)\cdots (m-k)/2} \\*
		&=\frac{(m-v(F_p))!}{(k+1-v(F_p))! (m-k-1)!} \cdot\frac{(m-k-1)!}{m!} \cdot 2^{c(F_p)} (k+1-r(F_p))! \\*
		&= \frac{2^{c(F_p)}}{(m)_{v(F_p)}}\cdot\frac{(k+1-r(F_p))!}{(k+1-v(F_p))!} \\*
		&\leq \frac{2^{c(F_p)}}{(m)_{v(F_p)}}\cdot (k+1)^{c(F_p)}.
	\end{align*}
	Set $Z:=|S\cap S'|$.
	For $1\leq R\leq k$, let $c_R:=\min\{R,k+1-R\}$; a feasible forest with $R$ edges has at most $c_R$ nontrivial components.
	In the sums below, $F_p$ ranges over edge sets of disjoint nontrivial paths with $v(F_p)\leq k+1$, and $v=R+c$ whenever $r(F_p)=R$ and $c(F_p)=c$.
	Now, considering the $R$-th falling factorial of $Z$, we have
	\begin{align*}
		\ee[(Z)_R] &= \sum_{I_R \in [m]_{2,R}} \Pr\big((i_1,j_1),...,(i_R,j_R) \in S\cap S'\big) \\
		&= \sum_{I_R \in [m]_{2,R}} \Pr\big((i_1,j_1),...,(i_R,j_R) \in S\big)^2\\
		&= R! \sum_{r(F_p)=R}\Pr(F_p \subset S)^2\\
		&\leq R! \sum_{c=1}^{c_R} \sum_{\substack{r(F_p)=R\\c(F_p)=c}}\Pr(F_p \subset S)^2\\
		&\leq R! \sum_{c=1}^{c_R} \bigg(\frac{2^{c}}{(m)_{v}}\cdot (k+1)^{c}\bigg)^2 \cdot \big|\big\{F_p: r(F_p)=R,~ c(F_p)=c\big\}\big|.
	\end{align*}
	We count the size of the collection $\big\{F_p: r(F_p)=R,~ c(F_p)=c\big\}$ for $1\leq c\leq c_R$ as follows.
	We first choose and arrange $v=R+c$ nodes from $m$ nodes, which can be done in $(m)_v$ ways.
	Then, we partition those $v$ nodes into $c$ components containing at least one edge, which can be done in $\binom{R-1}{c-1}$ ways.
	Also, we divide by $2^c$ to account for the fact that each component can be traversed in two directions, and thus we only count one of them.
	Last, we divide by $c!$ to account for the indistinguishability of components.
	Therefore, we have
	\begin{align*}
		\big|\big\{F_p: r(F_p)=R,~ c(F_p)=c\big\}\big| &= \frac{(m)_v }{2^c c!}\binom{R-1}{c-1}.
	\end{align*}
	Then, observe that
	\begin{align*}
		\ee[(Z)_R] &\leq R! \sum_{c=1}^{c_R} \frac{2^{2c}}{(m)_{v}^{2}}\cdot (k+1)^{2c} \cdot \frac{(m)_v }{2^c c!}\binom{R-1}{c-1}\\
		&\leq R! \sum_{c=1}^{c_R} \frac{2^c(k+1)^{2c}}{(m)_vc!}\binom{R-1}{c-1}\\
		&\stackrel{(\dagger)}\leq R! \sum_{c=1}^{c_R} \frac{2^c(k+1)^{2c}}{(m/e)^vc!}\binom{R-1}{c-1}\\
		&\leq R! \sum_{c=1}^R \frac{2^ce^v(k+1)^{2c}}{m^vc!}\binom{R-1}{c-1},
	\end{align*}
	where $(\dagger)$ uses $(m)_v\geq(m!)^{v/m}\geq(m/e)^v$ for $v\leq k+1\leq m$.
	The last line extends the nonnegative sum to all $1\leq c\leq R$, with $v=R+c$ in each summand.
	Now, for $t\geq0$ with $e^t-1<m/e$, we bound the moment generating function of $Z$ as follows:
	\allowdisplaybreaks
	\begin{align*}
		\ee[\exp(tZ)] &= \sum_{R=0}^k \frac{(e^t-1)^R}{R!} \ee[(Z)_R] \\*
		&\leq 1+\sum_{R=1}^k (e^t-1)^R \sum_{c=1}^R \frac{2^ce^v(k+1)^{2c}}{m^vc!}\binom{R-1}{c-1}\\
		&= 1+\sum_{R=1}^k \bigg(\frac{e(e^t-1)}{m}\bigg)^R \sum_{c=1}^R \bigg(\frac{2e(k+1)^{2}}{m}\bigg)^c\frac{1}{c!}\binom{R-1}{c-1}\\
		&= 1+\sum_{c=1}^k \bigg(\underbrace{\frac{2e(k+1)^{2}}{m}}_{:=x}\bigg)^c\frac{1}{c!}\sum_{R=c}^k \binom{R-1}{c-1}\bigg(\underbrace{\frac{e(e^t-1)}{m}}_{:=y}\bigg)^R\\
		&\leq 1+\sum_{c=1}^k x^c \frac{1}{c!}\sum_{R=c}^\infty \binom{R-1}{c-1} y^R\\
		&=1+\sum_{c=1}^k x^c \frac{1}{c!}\sum_{R'=0}^\infty \binom{R'+c-1}{c-1} y^{R'+c}\quad \text{where}~R':=R-c\\
		&=1+\sum_{c=1}^k x^c \frac{y^c}{c!}\sum_{R'=0}^\infty \binom{R'+c-1}{c-1} y^{R'}\\
		&\stackrel{(\ddagger)}=1+\sum_{c=1}^k x^c \frac{1}{c!}\cdot \frac{y^c}{(1-y)^c} \\
		&\leq \exp\bigg(\frac{xy}{1-y}\bigg)\\
		&= \exp\bigg(\frac{2e^2(k+1)^2(e^{t}-1)}{m\{m-e(e^t-1)\}}\bigg),
	\end{align*}
	The infinite-series bound and the negative-binomial identity in $(\ddagger)$ require $0\leq y<1$.
	For $t\geq0$, this is equivalent to $e^t-1<m/e$, and hence to $e^{\mu^2}-1<m/e$ when $t=\mu^2$.
	We verify below that the final sufficient condition guarantees this restriction.
	Therefore, we have
		$$
		\ee_{\bfY\sim\bbQ}[L_\nu(\bfY)^2] = \ee_{S,S'\sim \nu}[\exp(\mu^2 |S\cap S'|)]\leq \exp\bigg(\frac{2e^2(k+1)^2(e^{\mu^2}-1)}{m\{m-e(e^{\mu^2}-1)\}}\bigg)=\exp\bigg(\frac{2e^2(k+1)^2z}{m(m-ez)}\bigg),
		$$
	where we denote $z:=e^{\mu^2}-1$.
	Note that, for $0\leq z<m/e$, the right-hand side is bounded by $C$ if and only if
	\begin{align}\label{eq:kpath condition}
	z\leq\frac{m^2\log C}{2e^2(k+1)^2+em\log C},
	\end{align}
	which is equivalent to
	$$
	\mu\leq\sqrt{\log\bigg(1+\frac{m^2\log C}{2e^2(k+1)^2+em\log C}\bigg)}.
	$$
	Observe that the condition $e^{\mu^2}-1<m/e$ is automatically satisfied when \Cref{eq:kpath condition} holds.

	We now examine this sufficient condition in two regimes.
	$(i)$ Assume that $k=O(m^{1-\delta})$ for some fixed $\delta\in(0,0.1)$.
	Since
	$$
	2e^2(k+1)^2+em\log C=O(m^{2-2\delta}),
	$$
	it follows that
	$$
	\log\bigg(1+\frac{m^2\log C}{2e^2(k+1)^2+em\log C}\bigg)\geq 2\delta\log m+O(1).
	$$
	Hence, for sufficiently large $m$, the condition
	$\mu\leq\sqrt{\delta\log m}$
	guarantees that the second moment is bounded by $C$.
	A matching upper bound up to a constant factor is obtained by the outlier test in \Cref{proposition: gaussian test upper bound}, which runs in $O(n)$ time.
	$(ii)$ Assume that $k=\Theta(m)$.
	Since $k\leq m-1$, we obtain
	$$
	\frac{m^2\log C}{2e^2(k+1)^2+em\log C}\geq \frac{\log C}{2e^2+e\log C/m}\geq \frac{\log C}{4e^2}
	$$
	for all sufficiently large $m$. Therefore, it is sufficient to have
	$$
	\mu\leq\sqrt{\log\bigg(1+\frac{\log C}{4e^2}\bigg)}.
	$$
	This establishes a constant-order statistical lower bound, which matches the signal scale achieved by the average test, since the test achieves strong detection whenever $\mu=\omega(m/k)=\omega(1)$.

	We now consider cycles of length $3\leq k\leq m$.
	The cycle case can be handled using a grouping argument similar to those employed in the proofs of our main results, such as the multiple planted cliques problem (see \Cref{proof: multiple cliques problem}).
	The main idea is to cover a cycle by two paths of length $\lceil k/2\rceil$ and analyze their overlaps instead of directly studying the overlap between cycles.
	Given a uniformly random cycle $S$, choose a starting edge and an orientation uniformly at random, and denote the successive edges in the chosen direction by $S_1,\ldots,S_k$.
	Let $q:=\lceil k/2\rceil$, as in \Cref{eq: multiple cliques grouping}, and set $\mathcal S_1:=\{S_1,\ldots,S_q\}$ and $\mathcal S_2:=\{S_{k-q+1},\ldots,S_k\}$.
	The two paths are edge-disjoint when $k$ is even and share one edge when $k$ is odd.
	Define $\mathcal S'_1$ and $\mathcal S'_2$ analogously for an independent cycle $S'$, and observe that
	$$
	S=\mathcal S_1\cup\mathcal S_2,
	\qquad
	S'=\mathcal S'_1\cup\mathcal S'_2.
	$$
	Following the same approach as in \Cref{eq: multiple cliques decomposition} and the subsequent analysis therein, and applying H\"older's inequality to the four pairwise overlaps, we obtain
	\begin{equation*}
	\begin{aligned}
	\ee[\exp(\mu^2 |S\cap S'|)]\leq\ee\big[\exp\big(4\mu^2 \big|\mathcal S_1\cap \mathcal S'_1\big|\big)\big].
	\end{aligned}
	\end{equation*}
	Note that each of $\mathcal S_1,\mathcal S_2,\mathcal S'_1$, and $\mathcal S'_2$ is marginally distributed as a uniformly random path of length $q$, and the blocks from $S$ are independent of those from $S'$.
	Then, applying the second moment bound established for the path case with path length $q$ and signal strength $2\mu$, we obtain a sufficient condition on $\mu$ that differs only by a constant factor.
	Finally, depending on the regime, the matching upper bounds are achieved by the same tests as in the path case.

\end{document}